\documentclass[onecolumn, 12 pt, doublespace, fullpage, a4paper]{report}

\usepackage{amssymb}
\usepackage{amsthm}
\usepackage[cmex10]{amsmath}
\usepackage{adjustbox} 
\usepackage{academicons}
\usepackage{breakcites} 
\usepackage{color}
\usepackage{epsfig}
\usepackage{epstopdf}
\usepackage{etoolbox}
\usepackage{enumitem}
\usepackage{float}
\usepackage{fancyhdr}
\usepackage{graphicx}
\usepackage{graphics}
\usepackage{latexsym,amsfonts}
\usepackage{longtable}
\usepackage{listings}
\usepackage{lscape} 
\usepackage{lipsum}
\usepackage{multirow}             
\usepackage{pdfpages}
\usepackage[figuresright]{rotating} 
\usepackage{sectsty}
\usepackage{setspace}
\usepackage{subfigure}
\usepackage{textcomp}
\usepackage{url} 
\usepackage{wrapfig} 
\usepackage{wasysym} 
\usepackage{xcolor}

\usepackage{natbib}
\usepackage{bibentry}
\usepackage{hyperref}
\definecolor{linkcolor}{rgb}{0.45,0,0.05} 
\hypersetup{ %
  colorlinks=true,
  linkcolor=linkcolor,
  citecolor=linkcolor,
  filecolor=linkcolor,
  urlcolor=linkcolor,
}
\newcommand{\linkcolor}{\color{linkcolor}}

\def\2{$^2$}			 
\def\3{$^3$}			 
\def\-2{$^{-2}$}		 
\def\-3{$^{-3}$}		 
\def\-1{$^{-1}$}		 

\chapterfont{\fontsize{14}{15}\selectfont}   
\sectionfont{\fontsize{14}{15}\selectfont}
\subsectionfont{\fontsize{14}{15}\selectfont}
\subsubsectionfont{\fontsize{14}{15}\selectfont}

\fancyheadoffset[L]{0.01mm}
\renewcommand{\headrulewidth}{0pt}

\makeatletter
\patchcmd{\@makechapterhead}{50\p@}{20pt}{}{}
\patchcmd{\@makeschapterhead}{50\p@}{20pt}{}{}
\makeatother

\fancypagestyle{plain}{
\fancyhf{} 
\fancyhead[C]{\thepage} 
\renewcommand{\headrulewidth}{0pt}
\renewcommand{\footrulewidth}{0pt}}

   \usepackage{ifthen,xkeyval,xfor,amsgen}
   \usepackage[acronym,toc, nogroupskip]{glossaries}
   \newglossary[slg]{symbols}{syi}{sbl}{List of Symbols}
 
   \makeglossaries

\newglossaryentry{symb:Pi}{
name=$\pi$, type=symbols,
description=A mathematical constant whose value is the ratio of any circle's circumference to its diameter,
sort=symbolpi
}

\newglossaryentry{symb:Phi}{
name=$\varphi$, type=symbols,
description=An angle,
sort=symbolphi
}

\newglossaryentry{symb:Lambda}{
name=$\lambda$, type=symbols,
description=Lambda indicates usually an eigenvalue in linear algebra,
sort=symbollambda
}

\newacronym{toc}{ToC}{Table of Contents}
\newacronym{los}{LoS}{List of Symbols}
\newacronym{loa}{LoA}{List of Abbreviations}
\newacronym{phd}{PhD}{Doctoral}
\newacronym{MS}{MS}{Masters}
\newacronym{M$}{MS}{Microsoft}
\newacronym{CD}{CD}{Compact Disc}
\newacronym{kaust}{KAUST}{King Abdullah University of Science and Technology}

\newacronym{AD}{AD}{Active Directory\protect\glsadd{glos:AD}}

\newglossaryentry{glos:AD}{
name=Active Directory,
description={Active Directory is the directory service for Windows based networks, that allows central organization and administration of any network resource. It allows a single-sign-on concept independent from network topologies or network protocols. As a prerequisite you need a Windows Server acting as Domain Controller. This computer stores all necessary data, e.g.~usernames and corresponding passwords}
}

\newglossaryentry{glos:RespF}{name=response file, description={A file 
that allows unattended software installation}}

\usepackage[lmargin=40mm, rmargin=25mm, vmargin=25mm, headsep=2.5mm]{geometry}
\newcommand{\mathsym}[1]{{}}
\newcommand{\unicode}[1]{{}}
\renewcommand{\thechapter}{\arabic{chapter}}
\renewcommand\bibname{\centering BIBLIOGRAPHY}

\usepackage{graphicx}
\usepackage{apptools}
\usepackage[flushleft]{threeparttable}
\usepackage{array,booktabs,makecell}
\usepackage{multirow}
\usepackage{nicefrac}
\usepackage{amsthm}
\usepackage{mathtools}
\usepackage{amsmath,amsfonts,bm,amssymb}
\usepackage{dsfont}
\usepackage{wrapfig}
\usepackage{caption}
\usepackage{siunitx}
\usepackage{thm-restate}
\usepackage{nccmath}
\usepackage{empheq}
\usepackage{bbm}
\usepackage{tabularx}

\usepackage{algorithmic}
\usepackage{algorithm}
\usepackage{filecontents}

\usepackage[capitalize,noabbrev]{cleveref}

\usepackage{color}
\usepackage{colortbl}
\definecolor{bgcolor}{rgb}{0.76,0.88,0.50}
\definecolor{bgcolor0}{rgb}{0.93,0.99,1}
\definecolor{bgcolor1}{rgb}{0.8,1,1}
\definecolor{bgcolor2}{rgb}{0.8,1,0.8}
\definecolor{bgcolor3}{rgb}{0.50,0.90,0.50}
\usepackage{tcolorbox}
\usepackage{pifont}

\definecolor{mydarkblue}{rgb}{0,0.16,0.54}
\definecolor{mydarkgreen}{RGB}{39,130,67}
\definecolor{mydarkorange}{RGB}{236,147,14}
\definecolor{mydarkred}{RGB}{192,47,25}
\definecolor{ruby}{RGB}{155,17,30}
\definecolor{chili}{RGB}{191,0,0}
\definecolor{sangria}{RGB}{146,0,10}
\definecolor{burgundy}{RGB}{128,0,32}
\definecolor{darkred}{RGB}{132,0,0} 
\definecolor{cherry}{RGB}{192,0,0} 
\definecolor{myaccent}{rgb}{0,0.45,0.45}

\newcommand{\highlightcolor}{\color{myaccent}}

\usepackage[textsize=tiny]{todonotes}

\usepackage{xspace}
\usepackage{relsize}
\newcommand{\algname}[1]{{\sf\relscale{0.95}#1}\xspace}

\newcommand{\norm}[1]{\left\| #1 \right\|}
\newcommand{\sqnorm}[1]{\left\| #1 \right\|^2}
\newcommand{\infnorm}[1]{\left\| #1 \right\|_{\infty}}
\newcommand{\flr}[1]{\left\lfloor #1\right\rfloor} 
\newcommand{\ceil}[1]{\left\lceil #1\right\rceil} 
\newcommand{\inp}[2]{\left\langle#1,#2\right\rangle} 
\newcommand{\abs}[1]{\left| #1 \right|}

\newcommand{\R}{\mathbb{R}} 
\newcommand{\N}{\mathbb{N}} 
\newcommand{\Z}{\mathbb{Z}} 
\newcommand{\E}[1]{\mathbb{E}\left[#1\right]}

\newcommand{\Exp}[1]{{\mathbb{E}}\left[#1\right]}
\newcommand{\ExpSub}[2]{{\mathbb{E}}_{#1}\left[#2\right]}
\newcommand{\ExpCond}[2]{{\mathbb{E}}\left[\left.#1\ \right\vert\ #2\right]}

\newcommand{\hmu}{\hat{\mu}}

\newcommand{\conf}{\mathrm{conf}}

\newcommand{\cA}{\mathcal{A}}

\newcommand{\cC}{\mathcal{C}}
\newcommand{\cD}{\mathcal{D}}
\newcommand{\cE}{\mathcal{E}}

\newcommand{\cN}{\mathcal{N}}
\newcommand{\cO}{\mathcal{O}}

\newcommand{\cR}{\mathcal{R}}

\newcommand{\bx}{\bm{x}}

\newcommand{\eqdef}{\coloneqq} 

\newcommand{\minimize}{\mathop{\mathrm{minimize}}}

\newcommand{\tauavg}{\tau_{\mathrm{avg}}}

\makeatletter
\newcommand{\vast}{\bBigg@{4}}

\DeclareMathOperator*{\argmin}{arg\,min}
\DeclareMathOperator*{\argmax}{arg\,max}

\def\<{\left\langle}
\def\>{\right\rangle}
\def\[{\left[}
\def\]{\right]}
\def\({\left(}
\def\){\right)}

\definecolor{myblue}{rgb}{0.3,0.25,0.2}
\definecolor{myorange}{rgb}{0.3,0.25,0.2}
\newcommand{\checkmarkgreen}{\textbf{\color{myblue}\ding{52}}}
\newcommand{\crossmarkred}{{\color{myorange}\ding{56}}}

\newcommand{\iasgd}{{\sf\relscale{0.95}IA\textsuperscript{2}SGD}\xspace}
\newcommand{\iasgdtitle}{IA\textsuperscript{2}SGD}
\newcommand{\ringleader}{\algname{Ringleader~ASGD}}
\newcommand{\ringleadertitle}{Ringleader~ASGD}
\newcommand{\malenia}{\algname{Malenia~SGD}}
\newcommand{\maleniatitle}{Malenia~SGD}
\newcommand{\herosgd}{\algname{Hero~SGD}}
\newcommand{\herosgdtitle}{Hero~SGD}
\newcommand{\minibatch}{\algname{Minibatch~SGD}}

\newcommand{\naiveminibatch}{\algname{Naive~Minibatch~SGD}}
\newcommand{\naiveminibatchtitle}{Naive~Minibatch~SGD}
\newcommand{\asyncsgdtitle}{Asynchronous~SGD}
\newcommand{\asyncsgd}{\algname{\asyncsgdtitle}}
\newcommand{\naiveasyncsgdtitle}{Naive~Asynchronous~SGD}
\newcommand{\naiveasyncsgd}{\algname{\naiveasyncsgdtitle}}

\newcommand{\saga}{\algname{SAGA}}

\newcommand{\ringmaster}{\algname{Ringmaster~ASGD}}
\newcommand{\ringmastertitle}{Ringmaster~ASGD}
\newcommand{\rennalatitle}{Rennala~SGD}
\newcommand{\rennala}{\algname{Rennala~SGD}}
\newcommand{\sgd}{\algname{SGD}}
\newcommand{\ata}{\algname{ATA}}
\newcommand{\hogwild}{\algname{HOGWILD!}}

\usepackage{thmtools}
\usepackage{thm-restate}
\usepackage{mdframed}

\theoremstyle{plain}
\newtheorem{theorem}{Theorem}[section]

\newtheorem{proposition}[theorem]{Proposition}

\newtheorem{remark}[theorem]{Remark}

\theoremstyle{definition} 

\newmdenv[
  font=\normalfont\bfseries,
  linecolor=black,
  linewidth=0.8pt,
  topline=false,
  bottomline=false,
  leftline=false,
  rightline=false,
  backgroundcolor=gray!13,
  skipabove=2pt,
  skipbelow=2pt,
  innertopmargin=-5pt,
  innerbottommargin=4pt,
  innerrightmargin=4pt,
  innerleftmargin=4pt,
]{myboxed}

\newmdtheoremenv[
  font=\normalfont\bfseries,
  linecolor=black,
  linewidth=0.8pt,
  topline=false,
  bottomline=false,
  leftline=false,
  rightline=false,
  backgroundcolor=gray!13,
  skipabove=2pt,
  skipbelow=2pt,
  innertopmargin=-5pt,
  innerbottommargin=4pt,
  innerrightmargin=4pt,
  innerleftmargin=4pt,
]{boxedtheorem}[theorem]{Theorem}

\newmdtheoremenv[
  font=\normalfont\bfseries,
  linecolor=black,
  linewidth=0.8pt,
  topline=false,
  bottomline=false,
  leftline=false,
  rightline=false,
  backgroundcolor=gray!13,
  skipabove=2pt,
  skipbelow=2pt,
  innertopmargin=-5pt,
  innerbottommargin=4pt,
  innerrightmargin=4pt,
  innerleftmargin=4pt,
]{boxedlemma}[theorem]{Lemma}

\newmdtheoremenv[
  font=\normalfont\bfseries,
  linecolor=black,
  linewidth=0.8pt,
  topline=false,
  bottomline=false,
  leftline=false,
  rightline=false,
  backgroundcolor=gray!13,
  skipabove=2pt,
  skipbelow=2pt,
  innertopmargin=-5pt,
  innerbottommargin=4pt,
  innerrightmargin=4pt,
  innerleftmargin=4pt,
]{boxeddefinition}[theorem]{Definition}

\newmdtheoremenv[
  font=\normalfont\bfseries,
  linecolor=black,
  linewidth=0.8pt,
  topline=false,
  bottomline=false,
  leftline=false,
  rightline=false,
  backgroundcolor=gray!13,
  skipabove=2pt,
  skipbelow=2pt,
  innertopmargin=-5pt,
  innerbottommargin=4pt,
  innerrightmargin=4pt,
  innerleftmargin=4pt,
]{boxedassumption}[theorem]{Assumption}

\newtheorem{innercustomthm}{Theorem}
\newenvironment{restate-theorem}[1]
  {\renewcommand\theinnercustomthm{#1}\innercustomthm}
  {\endinnercustomthm}

\newtheorem{innercustomlemma}{Lemma}
\newenvironment{restate-lemma}[1]
  {\renewcommand\theinnercustomlemma{#1}\innercustomlemma}
  {\endinnercustomlemma}

\newtheorem{innercustomproposition}{Proposition}
\newenvironment{restate-proposition}[1]
  {\renewcommand\theinnercustomproposition{#1}\innercustomproposition}
  {\endinnercustomproposition}

\newenvironment{restate-boxedtheorem}[1]
  {\renewcommand\theinnercustomthm{#1}\begin{myboxed}\begin{innercustomthm}}
  {\end{innercustomthm}\end{myboxed}}

\newenvironment{restate-boxedlemma}[1]
  {\renewcommand\theinnercustomlemma{#1}\begin{myboxed}\begin{innercustomlemma}}
  {\end{innercustomlemma}\end{myboxed}}

\newcommand*{\sketchproofname}{Sketch of Proof}

\numberwithin{algorithm}{section}
\numberwithin{figure}{section}

\begin{document}
\nobibliography*


\thispagestyle{empty}
\addvspace{5mm}  


\begin{center}
\begin{doublespace}
{\textbf{First Provably Optimal Asynchronous SGD \\ for Homogeneous and Heterogeneous Data}}
\end{doublespace}

\vspace{10mm}
{Dissertation by}\\
{Artavazd Maranjyan} 

\vspace{30mm}

{ In Partial Fulfillment of the Requirements}\\[12pt]
{ For the Degree of}\\[12pt]
{Doctor of Philosophy} \vfill
{King Abdullah University of Science and Technology }\\
{Thuwal, Kingdom of Saudi Arabia}
\vfill

\begin{onehalfspace}
{\copyright December 4, 2025}\\
Artavazd Maranjyan\\               
All rights reserved\\

\small \href{https://artomaranjyan.github.io}{https://artomaranjyan.github.io}\\

\end{onehalfspace}

\end{center}
\newpage


\doublespacing

\chaptertitlefont{\fontsize{14}{15}\selectfont\centering}

\begin{center}

\end{center}

\begin{center}
{{\bf\fontsize{14pt}{14.5pt}\selectfont \uppercase{ABSTRACT}}}
\end{center}


\addcontentsline{toc}{chapter}{Abstract}

\begin{center}
	\begin{doublespace}
{\fontsize{14pt}{14.5pt}\selectfont {First Provably Optimal Asynchronous SGD \\ for Homogeneous and Heterogeneous Data}}\\
		{\fontsize{14pt}{14.5pt}\selectfont {Artavazd Maranjyan}}\\
	\end{doublespace}
\end{center}

\singlespace 


Artificial intelligence has achieved remarkable success across domains such as language modeling, vision, and autonomous systems.
These breakthroughs are driven by increasingly large neural networks trained on massive datasets using thousands of GPUs or TPUs.
Such training runs can occupy entire data centers for weeks to months, consuming vast computational and energy resources.

While advances in hardware and data availability have made this scaling possible, the optimization algorithms for training have not evolved as quickly.
Most large-scale training still relies on \emph{synchronous} methods, where all workers must finish their tasks before the next iteration begins.
As the number of devices grows, inefficiencies caused by synchronization may grow, too: faster workers remain idle while waiting for slower ones, wasting both compute and energy.
In practice, it is nearly impossible for all workers to run at identical speeds\textemdash hardware faults and network delays inevitably create \emph{computational heterogeneity}.

At first sight, it may seem that there is an easy fix\textemdash removing synchronization\textemdash which allows all workers to operate continuously.
However, asynchrony introduces \emph{staleness}\textemdash some results are produced using outdated versions of the model\textemdash making algorithms harder to analyze, especially when delays stem from system-level variability rather than the algorithm itself.
Despite extensive research, the \emph{time complexity} of asynchronous methods is poorly understood.

This dissertation addresses this gap.
We develop a rigorous framework for asynchronous first-order stochastic optimization, isolating the core challenge these methods target: \emph{heterogeneous worker speeds}.
Within this framework, we study stochastic gradient descent (\sgd) and show that, with proper design, asynchronous \sgd can be \emph{provably optimal} in terms of time complexity, matching optimality results recently achieved by synchronous variants of \sgd only.

The first contribution, \ringmaster, achieves optimal time complexity in the \emph{homogeneous data} setting by selectively discarding stale updates.
The second, \ringleader, extends this result to the \emph{heterogeneous data} regime\textemdash typical of federated learning\textemdash using a structured gradient-table mechanism that coordinates model updates.
Finally, \ata improves resource efficiency by learning workers' compute-time distributions and allocating tasks adaptively, achieving near-optimal wall-clock time with far less computation.

Together, these results establish asynchronous optimization as a theoretically sound and practically efficient foundation for parallel and distributed learning---showing that coordination without synchronization is not only possible, but that such a strategy enjoys optimal time complexity in theory, while outperforming competing synchronous methods in practice.



\begin{center}

\end{center}
\begin{center}

{\bf\fontsize{14pt}{14.5pt}\selectfont \uppercase{Acknowledgements}}\\\vspace{1cm}
\end{center}

\addcontentsline{toc}{chapter}{Acknowledgements} 
I would like to begin by expressing my deepest gratitude to my advisor, Professor Peter Richt\'{a}rik.
His continuous support, encouragement, and insightful guidance made this dissertation possible.
Every discussion with him was a source of motivation and clarity, pushing me forward throughout this challenging journey.
From him, I learned how to think about complex ideas from first principles\textemdash to seek precision, simplicity, and mathematical elegance in research.
He taught me that doing impactful research often begins by asking simple yet foundational questions.
His creativity and depth of thought have profoundly influenced the way I approach problems.

I am also deeply grateful to my coauthors\textemdash Abdurakhmon Sadiev, Alexander Tyurin, El Mehdi Saad, Francesco Orabona, Laurent Condat, Mher Safaryan, and Omar Shaikh Omar\textemdash with whom I have been fortunate to collaborate closely during my PhD. 
Special thanks to Professor Francesco Orabona\textemdash his \emph{Online Learning} course at KAUST remains one of my favorites, and I have greatly enjoyed many insightful discussions with him.

Beyond my direct coauthors, I have been fortunate to work alongside many brilliant colleagues in our research group at KAUST and to share countless technical and non-technical discussions that have shaped my research interests and made my time here truly enjoyable:
Ahmad Rammal,
Alina Abdikarimova,
Alexander Gaponov,
Avetik Karagulyan,
Egor Shulgin,
Elnur Gasanov,
Grigory Malinovsky,
Hanmin Li,
Hussein Rammal,
Igor Klimczak,
Igor Sokolov,
Ivan Ilin,
Kai Yi,
Kaja Gruntkowska,
Konstantin Burlachenko,
Majied Ammar Mahran,
Samuel Horváth,
Sarit Khirirat,
Slavomír Hanzely,
Yassine Maziane,
Yury Demidovich, 
and Zhirayr Tovmasyan. 

I would also like to thank Professor Yi-Shuai Niu for hosting me during my research visit to Beijing, China.
That visit gave me the privilege of meeting inspiring professors and exploring several leading universities in Beijing.

I am especially thankful to Hrant Khachatrian and Mher Safaryan, without whose help and guidance I would not have known about KAUST\textemdash and perhaps would never have started my PhD journey here.
My sincere appreciation also goes to KAUST for providing an exceptional research environment and continuous institutional support.

Last but not least, I am deeply grateful to my family and friends for their constant love and encouragement throughout this journey.
Special thanks to my mother, whose endless support and selfless sacrifices made all of this possible.

\singlespacing

\begin{onehalfspacing}

\tableofcontents

\printglossary[type=\acronymtype,style=long3col, title=\centerline{LIST OF ABBREVIATIONS}, toctitle=List of Abbreviations, nonumberlist=true] 

\printglossary[type=symbols,style=long3col, title=\centerline{LIST OF SYMBOLS}, toctitle=List of Symbols, nonumberlist=true]

\phantomsection
\addcontentsline{toc}{chapter}{\listfigurename} 
\renewcommand*\listfigurename{LIST OF FIGURES}
\listoffigures

\cleardoublepage
\phantomsection
\addcontentsline{toc}{chapter}{\listtablename}  
\renewcommand*\listtablename{LIST OF TABLES}
\listoftables

\end{onehalfspacing}


\chaptertitlefont{\fontsize{14}{15}\selectfont}  

\chapter{Introduction}\label{chapter:intro}
\thispagestyle{empty}
Artificial Intelligence (AI), once primarily a research topic, now shapes applications that billions use daily.
The most visible examples today are conversational AI systems such as ChatGPT and other chatbots \citep{GPT3,GPT4,Gemini}.
Beyond chat, AI also powers recommendation systems \citep{covington2016deep,zhang2019deep}, autonomous driving \citep{bojarski2016end,chib2023recent}, image recognition \citep{krizhevsky2012imagenet,mauricio2023comparing}, and medical diagnostics \citep{esteva2017dermatologist,rajpoot2024integrated}.
What unites these applications is the reliance on large-scale deep learning models trained on vast amounts of data \citep{kaplan2020scaling}.

These models are so computationally demanding that they cannot be trained or even deployed efficiently on a single machine \citep{dean2012large, li2020pytorch}.
They require distributed training on a cluster of specialized hardware, often involving thousands of interconnected GPUs or TPUs \citep{jouppi2017datacenter, narayanan2021efficient}.
Large-scale data centers\textemdash massive facilities that interconnect this hardware\textemdash have therefore become essential infrastructure, built by companies such as Google, Meta, and Amazon.
State-of-the-art AI training runs today often consume weeks or months of computation on these data centers \citep{narayanan2021efficient, grattafiori2024llama}.

As both models and datasets continue to grow, so do the computational and energy demands of AI training.
Data centers are already among the most energy-intensive building types, and their demand is increasing faster than that of most other sectors, posing a major challenge for sustainable energy systems \citep{strubell2020energy,iea2025energy}.
Recent estimates suggest that by 2030, data center energy consumption could reach levels comparable to Japan's current annual electricity usage \citep{iea2025energy}.

Given the growing computational and energy demands of AI training, improving the efficiency of training large models has become a critical and highly active research area \citep{shoeybi2019megatron, rajbhandari2020zero, narayanan2021efficient, hu2022lora, wan2024efficient}.
Despite rapid advances in hardware and software, current large-scale training pipelines still have substantial room for improvement.
For example, the training of LLaMA-3 achieved only about 38-43\% Model FLOPs Utilization (a standard metric for measuring training efficiency \citep{chowdhery2023palm}) \citep{grattafiori2024llama}.
This means that for much of the time, expensive GPUs were idle, yet still consuming power, rather than performing useful computation.
This underutilization is like reserving an entire fleet of airplanes for a week, but flying only a few of them each day while the rest sit idling on the runway with engines running.
The cost of wasted resources in such a scenario is enormous\textemdash and in the case of large AI training runs, this inefficiency translates into billions of dollars in infrastructure expenses and energy consumption \citep{de2023growing, cottier2024rising, odonnell2025energy, elsworth2025measuring}.

The underutilization of resources in large-scale AI training comes from several sources, including the limits of current hardware and the communication overhead of transferring data between machines.
A particularly important source of inefficiency, however, lies on the algorithmic side.
Most existing training algorithms were not originally designed with today's massive computing clusters in mind, and they fail to make \emph{optimal} use of such large-scale systems.
In particular, they rely on synchronization across machines, an approach that becomes increasingly costly as the system grows.
%
\section{Synchronous Training: Meta-Algorithm}
Historically, large-scale AI training has relied on synchronous iterative algorithms, where all workers must complete their assigned tasks before the system can advance to the next iteration.
This strategy has remained the dominant approach for training large AI models \citep{goyal2017accurate, li2020pytorch, douillard2024diloco}.

To illustrate the principle behind these methods, we abstract away the details of the underlying optimization problem and present a general \emph{meta-algorithm} (\Cref{intro:alg:general_synchronous}) that captures the structure common to many synchronous approaches.
A precise definition of “synchronous algorithm” is deferred to \Cref{intro:sec:sync_vs_async}; here, the meta-algorithm is intentionally informal and serves to highlight common behavior and its bottlenecks.

\begin{algorithm}[H]
    \caption{Synchronous Meta-Algorithm (Server Perspective)}
    \label{intro:alg:general_synchronous}
    \begin{algorithmic}[1]
        \FOR{iteration $k = 0,1,2,\dots$}
            \STATE Broadcast the current model $x^k$ to all workers
                \label{meta:line:broadcast}
            \STATE Instruct all workers to compute their local results using $x^k$
                \highlightcolor
            \STATE Wait until all workers have finished, and their results have been received
                \label{meta:line:wait}
                \color{black}
            \STATE Aggregate the collected results
                \label{meta:line:aggregate}
            \STATE Update the model to obtain $x^{k+1}$
                \label{meta:line:update}
        \ENDFOR
    \end{algorithmic}
\end{algorithm}

The defining feature of these methods is the \emph{synchronization barrier}: the server must \emph{wait} for all workers to finish (Line~\ref{meta:line:wait}) before it can aggregate their results (Line~\ref{meta:line:aggregate}), compute the new model (Line~\ref{meta:line:update}), and broadcast it again to the workers (Line~\ref{meta:line:broadcast}).
We highlight Line~\ref{meta:line:wait} because it is the bottleneck: faster workers cannot proceed with new computations and remain idle until the slowest worker finishes.
This idleness is the central drawback of synchronous methods and a major source of underutilization in large-scale data centers.
\subsection{Challenges of Synchronization in Practice}
\label{subsec:sync_challenges}
One might hope to solve the synchronization problem by building data centers with identical hardware, ensuring all workers run at the same speed.
In practice, however, this strategy faces several fundamental obstacles.
\subsubsection{Stragglers Persist Even with Identical Hardware}
%
Even when all machines in a cluster are identical, \emph{stragglers}\textemdash workers that are slow due to compute or communication delays\textemdash still arise because of unpredictable issues such as hardware faults (e.g., overheating) and network problems \citep{dean2013tail,ananthanarayanan2013effectivestraggler,maranjyan2025mindflayer}.
A single malfunctioning component can cause one worker to lag significantly behind others, forcing the entire cluster to wait.
In other words, identical hardware specifications do not guarantee identical performance in practice.
\subsubsection{Hardware Upgrades Become Prohibitively Expensive}
%
Enforcing hardware homogeneity creates a second problem: it makes upgrading the data center extremely difficult.
When newer, more powerful processors become available, adding them to an existing cluster introduces heterogeneity.
Since faster machines must still wait for slower ones in synchronous training, much of the performance benefit of the upgrade is lost.
This forces operators into an expensive dilemma: either replace all machines simultaneously (which is extremely costly) or delay upgrades until a complete replacement is feasible.
Such constraints make it difficult to gradually expand capacity or adopt new technology as it becomes available.

\subsubsection{Federated and Edge Settings}
%
Finally, there are important scenarios where ensuring identical hardware is simply impossible.
In federated learning, training occurs across everyday devices such as smartphones, sensors, or other personal electronics \citep{konevcny2016federatedlearning, konevcny2016federatedoptimization, mcmahan2017communication}.
These devices vary dramatically in their computational power, battery life, internet connectivity, and availability.
Participants may join or leave the training process unpredictably, and connection quality can fluctuate significantly.
In such environments, waiting for the slowest device to complete each round is not merely inefficient\textemdash it may be impossible, as some devices might disconnect or fail to complete their work within a reasonable timeframe.

\section{Asynchronous Training: Meta-Algorithm}
Given the fundamental limitations of synchronization discussed above, researchers have proposed an algorithmic solution: simply remove the synchronization requirement altogether \citep{tsitsiklis1986distributed,recht2011hogwild,agarwal2011distributed}.
If the synchronization barrier is the source of inefficiency, then eliminating it should allow workers to operate continuously without idle time.

To illustrate this approach, we again present a general meta-algorithm (\Cref{intro:algo:general_asynchronous}) that captures the structure common to many asynchronous methods.
As before, we defer a precise definition of ``asynchronous algorithm'' to \Cref{intro:sec:sync_vs_async}.

\begin{algorithm}[H]
    \caption{Asynchronous Meta-Algorithm (Server Perspective)}
    \label{intro:algo:general_asynchronous}
    \begin{algorithmic}[1]
        \STATE Initialize model and send to all workers
        \WHILE{not terminated}
            \STATE Receive result from some worker $i$
                \highlightcolor
            \STATE No waiting -- proceed immediately
                \label{intro:alg:asynch:line:no_waiting}
                \color{black}
            \STATE Optionally update the model (based on some criterion)
                \label{intro:alg:asynch:line:optional_update}
            \STATE Send current model to worker $i$
                \highlightcolor
            \STATE Worker $i$ starts computing again
                \label{intro:alg:asynch:line:no_waiting_worker}
        \ENDWHILE
    \end{algorithmic}
\end{algorithm}

The key distinction between synchronous and asynchronous methods appears in Line~\ref{intro:alg:asynch:line:no_waiting}, 
where the server proceeds as soon as it receives a result from any worker, without waiting for the others.
Similarly, in Line~\ref{intro:alg:asynch:line:no_waiting_worker}, each worker continues computing immediately\textemdash without waiting for other workers.
This design eliminates the idle time that plagues synchronous approaches, allowing all workers to operate independently and continuously using their current model versions.

However, this independence comes at a cost.
When a worker begins its computation and later sends back its result, the model may have already been updated multiple times by other workers in the meantime.
This means the worker's result is based on an outdated, or \emph{stale}, version of the model.
Such staleness can degrade optimization progress, potentially slowing convergence or leading to suboptimal solutions.

Moreover, Line~\ref{intro:alg:asynch:line:optional_update} reveals that asynchronous methods offer additional flexibility\textemdash the server need not update the model after every received result.
These observations raise fundamental design questions: how should one design an asynchronous method to maximize its computational benefits while mitigating the effects of staleness?
\section{Dissertation Focus}
While asynchronous methods eliminate the key inefficiency of synchronous methods\textemdash keeping workers idle while waiting for slow workers\textemdash they introduce a new challenge: staleness in model updates.
This creates a fundamental paradox.
Staleness is not merely a side effect of asynchrony\textemdash it is the very mechanism that enables continuous worker utilization.
Workers can continue computing precisely because they do not wait for the model to be fully up-to-date.
Yet this same staleness means that updates are computed using outdated information, which can interfere with optimization progress and potentially harm convergence.
If workers operate on severely outdated models, their updates may no longer be beneficial and could actively degrade the training process.

Predicting and managing staleness is challenging on multiple fronts.
It is difficult to anticipate how much staleness a given system will experience, how that staleness will affect convergence, and which design choices lead to algorithms that work well in practice while maintaining theoretical guarantees.
Moreover, staleness makes the analysis considerably more difficult: the classical convergence analysis of synchronous methods no longer applies, and new theoretical machinery must be developed to understand asynchronous behavior.

\emph{The central goal of this dissertation is to develop a rigorous understanding of asynchronous optimization and to design provably efficient asynchronous algorithms for large-scale machine learning.}

\subsection{Scope: Data Parallelism}
Modern distributed training systems employ multiple forms of parallelism to handle large models and datasets.
These include \emph{data parallelism}, where workers process different subsets of data; \emph{model parallelism}, where different parts of the model are distributed across workers; and \emph{pipeline parallelism}, where different layers of the model are assigned to different stages of a pipeline \citep{narayanan2021efficient, grattafiori2024llama}.
In practice, state-of-the-art training systems often combine several or even all of these strategies simultaneously.

This dissertation focuses exclusively on \emph{data parallelism}, the most fundamental and widely applicable form of parallelism.
In data parallelism, all workers maintain a copy of the complete model but process different batches of training data.
We restrict our attention to this setting for two reasons.
First, data parallelism is the simplest form of parallelism and serves as a natural starting point for developing theoretical understanding.
Second, despite its apparent simplicity, even the asynchronous data-parallel setting remains poorly understood\textemdash fundamental questions about convergence, optimal algorithm design, and the effects of staleness remain open.
\section{Problem Formulation}
Having outlined the general focus of this dissertation, we now provide a concrete formulation of the problem we study.
We begin by presenting the optimization problem that lies at the core of distributed machine learning.
Then, since our primary interest is in the convergence time of algorithms, we formalize which operations take time and how these quantities are modeled in our analysis.
\subsection{Optimization Problem}
We consider a distributed learning setup with $n$ clients, where each client~$i$ has access to its own local data distribution~$\mathcal{D}_i$.
The goal is to learn a global model, parameterized by a vector $x \in \mathbb{R}^d$, where $d$ denotes the dimensionality of the model parameters, by solving the optimization problem
\begin{equation}\label{intro:eq:problem}
    \minimize\limits_{x \in \R^d}
        \left\{ f(x) \coloneqq \frac{1}{n}\sum_{i=1}^{n}f_i(x) \right\},
\end{equation}
where $f_i \colon \R^d \to \R$ denotes the local objective function of client $i$, and the global objective function $f$ is defined as their average.

The local objective functions can generally be written as
\begin{equation}\label{intro:eq:f_i}
    f_i(x) \coloneqq \ExpSub{\xi_i \sim \mathcal{D}_i}{f_i(x; \xi_i)}~,
\end{equation}
where $f_i(x; \xi_i)$ is the loss function of the model parameterized by $x$ evaluated on a data sample $\xi_i$.
Thus, the local objective is the expectation of this loss over the entire local data distribution $\mathcal{D}_i$.
This formulation corresponds to minimizing the population risk.

When the local data distribution is finite, i.e., $|\cD_i| < \infty$, the local objective \eqref{intro:eq:f_i} reduces to the average loss over the local training data points,
\begin{equation*}
    f_i(x) = \frac{1}{\abs{\cD_i}} \sum_{j=1}^{\abs{\cD_i}} f_i(x; \xi_{i,j}) ~.
\end{equation*}
\subsubsection{Homogeneous vs. Heterogeneous Data}
\label{intro:sec:homo_hetero}
In general, the local objective functions $f_i$ can differ significantly across clients depending on their local data distributions~$\mathcal{D}_i$.
This heterogeneity poses a fundamental challenge for solving the distributed optimization problem and is a key characteristic that distinguishes different practical settings.
\paragraph{Heterogeneous setting.}
In federated learning scenarios, data heterogeneity is unavoidable.
For example, when clients are personal devices such as smartphones, each device contains user-specific data that reflects the habits, preferences, and environments of a particular individual.
Such data cannot be redistributed or centralized without violating privacy constraints.
As a result, the local data distributions~$\mathcal{D}_i$ may differ drastically across clients\textemdash ranging from text typed in different languages to photos taken under different conditions\textemdash making the optimization problem inherently more challenging.
We refer to this as the \emph{heterogeneous data setting}.
\paragraph{Homogeneous setting.}
In contrast, data center environments offer significantly more freedom in organizing and distributing data.
In such settings, data can be freely shuffled, replicated, or centralized without privacy concerns.
It is therefore common practice to ensure that all workers have access to data drawn from the same underlying distribution, either by sharing a common dataset or by partitioning the data so that the resulting subsets are approximately statistically identical across workers.
Formally, this corresponds to assuming $\mathcal{D}_i = \mathcal{D}$ for all clients~$i$.

In addition, we typically employ the same loss function across all clients, i.e., $f_i(x; \xi_i) \equiv f(x; \xi)$. 
Under these assumptions, the optimization problem \eqref{intro:eq:problem} simplifies to
\begin{equation}\label{intro:eq:homogeneous_problem}
    \minimize\limits_{x \in \mathbb{R}^d} 
    \left\{f(x) \coloneqq \ExpSub{\xi \sim \cD}{f(x;\xi)} \right\}.
\end{equation}
We refer to this as the \emph{homogeneous data setting}.

Both settings are of practical importance and present distinct algorithmic challenges.
The homogeneous setting, while more structured, still requires the design of efficient distributed algorithms.
The heterogeneous setting, being more general, demands methods that are robust to statistical diversity across clients.
Importantly, these settings need to be studied separately: analyzing only the general heterogeneous case does not necessarily yield optimal solutions for the homogeneous case.
In the homogeneous setting, one can exploit the fact that all clients share the same data distribution, which enables sharper guarantees and often simpler algorithms.
In this dissertation, we develop and analyze algorithms tailored to both settings.
\subsection{Assumptions}
To make the problem setup concrete, we introduce several standard assumptions on the class of functions we study.
\begin{boxedassumption}[Smoothness]\label{intro:ass:lipschitz_constant}
    The function $f$ is differentiable, and its gradient is Lipschitz continuous.
    That is, there exists a constant $L>0$ such that, for all $x, y \in \mathbb{R}^d$, 
    \begin{equation*}
        \norm{\nabla f(x) - \nabla f(y)} \le L\norm{x - y}.
    \end{equation*}
\end{boxedassumption}
This assumption, often referred to as \emph{$L$-smoothness}, controls how quickly the gradient can change.
%
\begin{boxedassumption}[Lower boundedness]\label{intro:ass:lower_bound}
    The objective function is bounded below: there exists $f^* > -\infty$ such that $f(x) \geq f^*$ for all $x \in \mathbb{R}^d$.
\end{boxedassumption}
This assumption ensures that the optimization problem is well-posed: without 
it, function values could decrease without bound, making convergence 
guarantees meaningless.

We define $\Delta \coloneqq f(x^0) - f^*$, where $x^0$ is the initial point of the optimization algorithm.
The quantity $\Delta$ represents the initial suboptimality and serves as a measure of the inherent difficulty of the problem.
\subsubsection{Convergence Criterion}
Under these assumptions, the standard goal in nonconvex optimization is to find an $\varepsilon$-stationary point, that is, a (possibly random) vector $x$ satisfying
\begin{equation}\label{intro:eq:stationarity}
    \E{\sqnorm{\nabla f(x)}} \leq \varepsilon ~.
\end{equation}
This criterion is natural in the nonconvex setting: while finding global minima is generally intractable, first-order stationary points are attainable and often correspond to practically useful solutions.
\subsubsection{Stochastic First-Order Methods}
We focus on stochastic first-order methods, which rely solely on gradient information.
The earliest and most fundamental example is stochastic gradient descent (\sgd), introduced by \citet{robbins1951stochastic} and later popularized in machine learning by \citet{bottou2018optimization}.
Building on this foundation, adaptive methods such as \algname{Adam} \citep{kingma2014adam} and its variant \algname{AdamW} \citep{loshchilov2017decoupled} have become standard in training large-scale deep learning models nowadays.
In our analysis, we begin with \sgd, both for its simplicity and because it serves as the theoretical basis for many of these modern variants.
To study its behavior in the distributed setting, we impose additional assumptions on the stochastic gradients computed by the workers.
\begin{boxedassumption}
    \label{intro:ass:stochastic_variance_bounded}
    For each $i \in \{1,\ldots,n\}$ and every $\xi_i$, the function $f_i(x;\xi_i)$ is differentiable with respect to its first argument $x$.
    Moreover, the stochastic gradients are unbiased and have bounded variance $\sigma^2 \geq 0$, that is,
    \begin{gather*}
         \ExpSub{\xi_i \sim \cD_i}{\nabla f_i(x;\xi_i)} = \nabla f_i(x), 
            \quad \forall x \in \R^d, \;\; \forall i \in \{1,\ldots,n\}~,\\
         \ExpSub{\xi_i \sim \cD_i}{\sqnorm{\nabla f_i(x;\xi_i) - \nabla f_i(x)}} \leq \sigma^2,
            \quad \forall x \in \R^d, \;\; \forall i \in \{1,\ldots,n\}~.
    \end{gather*}
\end{boxedassumption}
The first condition ensures that stochastic gradients are unbiased estimators of the true gradient, while the second condition bounds their variance.
The parameter $\sigma$ quantifies the level of stochastic noise: when $\sigma = 0$, the gradients are exact.
\subsection{Worker Heterogeneity Model}
We are interested in the \emph{time complexity} of solving problem~\eqref{intro:eq:problem}, that is, the wall-clock time required to reach a desired level of stationarity \eqref{intro:eq:stationarity}.
To analyze time complexity in distributed systems, we model the computational characteristics of individual workers.
In this section, we formalize the notion of worker computational heterogeneity through their computation times.
\subsubsection{Fixed Computation Times}

We begin with the case where workers operate at constant computation speeds.
This setting provides an intuitive foundation for understanding the dynamics of asynchronous distributed optimization and isolates the effects of heterogeneity from other sources of variability.

Following the \textit{fixed computation model} of \citet{mishchenko2022asynchronous}, we assume:
\begin{equation}\label{intro:eq:fixed_time}
    \begin{minipage}{0.52\linewidth}
    Each worker~$i$ requires $\tau_i$~seconds\footnotemark\ to compute one~stochastic~gradient~$\nabla f_i(x;\xi_i)$.
    
    \vspace{4pt}
    
    Without loss of generality, we assume
    
    \centering
    $0 < \tau_1 \le \tau_2 \le \cdots \le \tau_n$~.
    \end{minipage}
\end{equation}
\footnotetext{Throughout this dissertation, we use ``seconds'' as a generic unit of time. The analysis applies to any consistent time unit.}
Differences in computation times $\{\tau_i\}_{i=1}^n$ may arise from diverse sources: heterogeneous hardware (different GPU generations), variation in local dataset sizes, or competition for shared computational resources.
This model captures the central challenge of asynchronous optimization: designing algorithms that efficiently exploit updates from workers that progress at different speeds.
\subsubsection{Time-Varying Computation Times}
In practice, computation times are rarely constant.
Even in data centers, workers may slow down due to shared resource usage, background processes, overheating, or occasional hardware failures.
This variability is even more pronounced in federated learning, where clients can face unstable network connections, limited battery life, or interruptions from user activity.

To capture these effects, we also consider models with time-varying computation times \citep{tyurin2024tighttimecomplexitiesparallel}.
A formal treatment of this more general setting is deferred to later chapters (see, for example, \Cref{sec:universal_computation_model}).
Crucially, we assume that variations in computation times are \emph{non-adversarial}: they may depend on system conditions or external factors, but not on the optimization algorithm itself.
This assumption reflects realistic scenarios where delays are caused by environment and hardware limitations rather than by deliberate interference.

\subsubsection{Stochastic Computation Times}
While the time-varying model introduced above can capture virtually any computation behavior, it provides no meaningful theoretical guarantees\textemdash its performance depends entirely on how the computation times evolve over time.
To obtain more interpretable results, we also consider a stochastic setting, where each worker's computation time is modeled as a random variable drawn from an underlying probability distribution.
This formulation allows us to derive guarantees in expectation, explicitly characterizing how the distributional parameters (e.g., mean, variance) affect the algorithm's behavior.
We adopt this stochastic assumption in parts of this dissertation; the precise formulation is given in \Cref{ata:modeling_assumptions}.
\subsubsection{Zero Communication and Auxiliary Costs}
For clarity of exposition, we assume that the only operation that consumes time is the computation of stochastic gradients.  
All other operations\textemdash including worker-to-server and server-to-worker communication, server-side model updates, and memory reads or writes\textemdash are assumed to take zero seconds.  
Of course, this assumption is not realistic, as communication delays and server overhead are non-negligible in real systems.
\subsubsection{Memory Locking and Inconsistent Reads}
A line of work initiated by \citet{recht2011hogwild} studies \emph{lock-free} updates on a shared parameter vector, where multiple workers read and write at the same time.
In that setting, while a worker is reading the server's parameter vector, other workers can change yet-unread coordinates; the result can be a local copy of the model that mixes values taken before and after those updates.

We do not consider such effects in our analysis, since we assume that communication, reads, and writes are instantaneous (zero-time) and \emph{atomic}.
Consequently, each read returns a single, well-defined model vector and each write completes indivisibly.
Under these assumptions, such inconsistencies cannot arise, and locking is unnecessary.

\subsubsection{Rationale Behind These Simplifications}
The assumptions introduced above intentionally eliminate all sources of time other than the workers' computation times.
We adopt this simplified setting because the primary challenge that asynchronous methods aim to address is \emph{heterogeneity in worker speeds}.
Even under such idealized conditions, the behavior of asynchronous algorithms is not yet fully understood: in many cases, they provide no advantage over their synchronous counterparts.
Establishing a rigorous understanding in this controlled environment is therefore a necessary first step before incorporating additional system-level complexities.
\section{Worker and Server Dynamics}\label{intro:sec:dynamics}
Before describing specific algorithms, let us outline the basic set of actions available to the workers and to the server.
All algorithms we construct later will be built from these basic operations.
\subsection{Worker Actions}
Each worker can perform only two types of actions at any given time:
\begin{itemize}
    \item \textbf{Idle:}
        The worker simply does nothing, waiting for instructions from the server.
    \item \textbf{Compute and send gradient:}
        The worker samples a new iid data point, computes a stochastic gradient at the model it was assigned, and immediately sends this gradient to the server.
        These stochastic gradients are unbiased, and their variance is bounded as assumed in \Cref{intro:ass:stochastic_variance_bounded}.
\end{itemize}
These are the only worker actions we consider.
One could imagine allowing workers to store gradients locally before sending them, but in our setting this makes no difference, since we assume communication is instantaneous.

A related approach is to perform local model updates, as commonly done in federated learning \citep{mcmahan2017communication}.
While this dissertation does not consider local steps, such extensions are possible \citep{tyurin2025birch,fradin2025local}.
\subsection{Server Actions}
The server has a richer set of possible actions, since it must coordinate the overall training process.
When a gradient arrives from a worker, the server can:
\begin{itemize}
    \item \textbf{Store the gradient:}
        Keep it in memory for potential use in a future update.
    \item \textbf{Discard the gradient:}
        For example, if it decides the gradient is outdated or of poor quality.
\end{itemize}
In addition to reacting to incoming gradients, the server can independently perform the following actions at any time:
\begin{itemize}
    \item \textbf{Update the model:}
        Update the model using stored gradients.
        This includes both immediate updates upon gradient arrival and delayed updates using multiple stored gradients.
    \item \textbf{Discard stored gradients:}
        Remove some of the previously stored gradients, for example if they have already been used sufficiently or are no longer needed.
    \item \textbf{Request gradient:}
        Send the current model to one or more workers, instructing them to compute gradients.
        In practice, the server may sometimes ask a worker to recompute a gradient at the same model it was already using, in which case resending the model is unnecessary.
        However, since we assume communication is instantaneous, we simply model requests as always sending the model.
    \item \textbf{Terminate one or more workers' computations:}
        If the server anticipates that the results of certain ongoing computations will be discarded, it may cancel them early to avoid wasted effort.
        The affected workers then simply return to the idle state until they receive new instructions.
\end{itemize}
Although these actions could, in principle, occur at any time, in practice an algorithm will specify precisely when they take place (e.g., upon receiving a gradient, after an update, etc.).
Thus, the server's behavior is structured rather than arbitrary.

The methods that operate within this framework will be referred to as \emph{distributed} or \emph{parallel} methods.
Both synchronous and asynchronous algorithms can be expressed in terms of these actions.
In the next section, we formalize how the distinction between synchronous and asynchronous methods arises within this framework.
\section{Synchronous vs. Asynchronous Algorithms: \\ Formal Definitions}\label{intro:sec:sync_vs_async}
We now formally define what we mean by \emph{synchronous} and \emph{asynchronous} algorithms.
\subsection{Synchronous Algorithms}
We call an algorithm \emph{synchronous} if, whenever the server performs a model update from stored stochastic gradients, it uses only gradients computed at the \emph{current} server model.
Equivalently, all gradients in an update correspond to the same model state.
Consequently, the update mirrors a standard \sgd step with an unbiased gradient estimator.

Formally, a synchronous update can be written as
\begin{equation}\label{intro:eq:synch_update}
    x^{k+1} = x^k - \gamma^k \sum_{i=1}^n w_i^k \sum_{j=1}^{b_i^k} \nabla f_i \(x^k;\xi_i^{k,j}\) ,
\end{equation}
where $\gamma^k>0$ is a stepsize, $w_i^k\ge 0$ are aggregation weights, and $b_i^k\in\{0,1,2,\ldots\}$ is the number of gradients contributed by worker $i$ in round $k$ (by convention, if $b_i^k=0$, then we take the sum as $0$).
All synchronous algorithms considered in this dissertation fit this template via different choices of $(\gamma^k, w_i^k, b_i^k)$\textemdash including \herosgd, \naiveminibatch, and \rennala (see \Cref{intro:sec:algos}).

Because updates must use only gradients computed at the \emph{same} model state, synchronous methods have inherent limitations.
When worker speeds differ, they either (i) force faster workers to wait until the slowest finishes, or (ii) discard late computations that miss the synchronization point.
Of course, this is not an issue if worker speeds are roughly the same, but this is hard to guarantee, as already discussed in \Cref{subsec:sync_challenges}.
\subsection{Asynchronous Algorithms}
By contrast, we call a method \emph{asynchronous} if its updates may use gradients evaluated at \emph{stale} model copies\textemdash that is, earlier versions of the model than the one currently held by the server.
While this need not occur in every iteration, asynchronous methods must generally account for such delays.

For simplicity, we assume that in each round, all gradients contributed by a given worker are computed using the \emph{same} model iterate.
This assumption covers all asynchronous methods analyzed in this dissertation, and relaxing it to allow mixtures of stale points per worker per round would be straightforward but unnecessary for our purposes.

Let $\delta_i^k \in \{0,1,2,\ldots\}$ denote the \emph{delay} at iteration $k$: the number of server updates that have occurred since worker $i$ received its model iterate.
With this notation, an asynchronous update takes the form
$$ 
    x^{k+1} = x^k - \gamma^k \sum_{i=1}^n w_i^k \sum_{j=1}^{b_i^k} \nabla f_i\(x^{k-\delta_i^k};\xi_i^{k-\delta_i^k, j}\),
$$
where $\gamma^k>0$ is the stepsize, $w_i^k\ge 0$ are aggregation weights, and $b_i^k\in\{0,1,2,\ldots\}$ is the number of gradients contributed by worker $i$ in round $k$.

Note that when all $\delta_i^k=0$, we recover the synchronous update \eqref{intro:eq:synch_update}; some algorithms occasionally realize such synchronous steps.
However, the goal is to provide guarantees for any admissible delay pattern.
An illustrative instance of an asynchronous method is \naiveasyncsgd, presented in \Cref{intro:sec:algos}.
%
\section{Algorithms for the Homogeneous Setting with Fixed Computation Times}
\label{intro:sec:algos}
Now that we have formalized the problem setting, we can begin exploring algorithmic solutions.
Recall that we have $n$ machines, each capable of computing unbiased stochastic gradients with bounded variance.
We assume that each machine's computation time is fixed, as specified by the fixed computation model in \eqref{intro:eq:fixed_time}, and that the server knows these times.

Our goal is to design algorithms that efficiently solve the optimization problem in this distributed setting.
We begin by focusing on the \emph{homogeneous data setting}, where we aim to solve problem \eqref{intro:eq:homogeneous_problem}.

If we were not concerned with time complexity, we could simply apply standard \sgd and achieve optimal iteration complexity \citep{ghadimi2013stochastic,arjevani2022lower}.
Consequently, in the single-machine setting, plain \sgd is optimal.
This observation motivates a natural starting point: begin with \sgd and investigate whether it can be adapted to achieve optimal time complexity in the distributed setting.
\subsection{\herosgdtitle}
\begin{figure*}[h]
    \centering
    \includegraphics[width=0.66\textwidth]{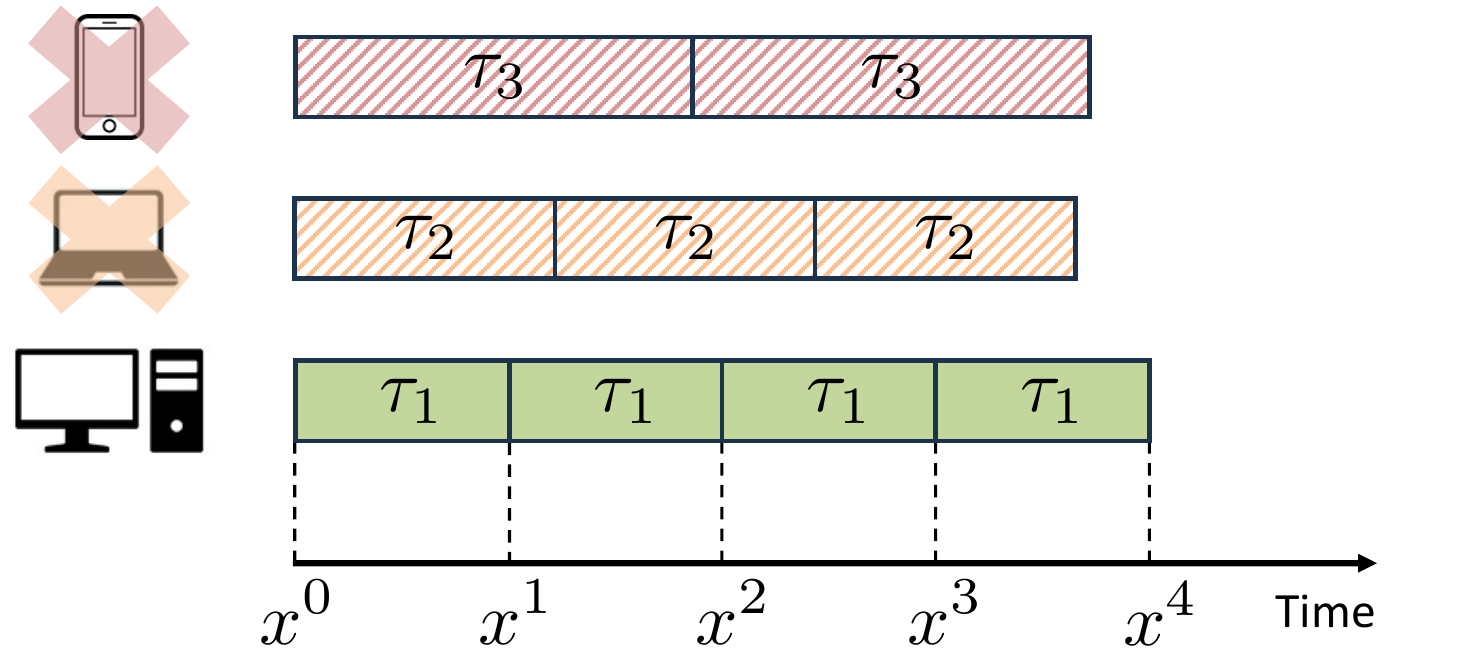}
    \caption{
        Illustration of \herosgd (\Cref{intro:alg:herosgd}) with 3 workers.
        Only the fastest worker (bottom) performs updates, while the others remain idle.
        Their idle computation capacity could, in principle, be utilized to compute additional gradients.
        }
    \label{intro:fig:herosgd}
\end{figure*}

The simplest approach to applying \sgd in our setting is to ignore the fact that we have $n$ machines available and use only a single one.
In that case, it is natural to select the fastest machine\textemdash the one with computation time $\tau_1$\textemdash and run standard \sgd on it.
We refer to this baseline as \herosgd,\footnote{The name is taken from the work by \citet{fradin2025local}.} and the corresponding algorithm is given in \Cref{intro:alg:herosgd}.
We also illustrate its behavior using a simple example with three workers in \Cref{intro:fig:herosgd}.

\begin{algorithm}[h]
    \caption{\herosgdtitle}
    \label{intro:alg:herosgd}
    \begin{algorithmic}[1]
        \STATE \textbf{Input:} initial point $x^0 \in \R^d$, stepsizes $\gamma^k > 0$
        \FOR{$k = 0, \dots, K - 1$}
            \STATE Request the fastest worker to compute a stochastic gradient $\nabla f\big(x^k; \xi^k\big)$
            \STATE Wait until the gradient arrives
            \STATE Update the model:
            $$
                x^{k+1} = x^k - \gamma^k \nabla f\(x^k; \xi^k\)
            $$
        \ENDFOR
    \end{algorithmic}
\end{algorithm}

This algorithm requires knowing which worker is the fastest.
If this information is unavailable, one could adopt a \emph{greedy} strategy: request all workers to compute a gradient and use the first result that arrives, discarding the rest.
Alternatively, one may attempt to learn the computation-time distributions of the workers and allocate tasks accordingly; we discuss such adaptive strategies later in \Cref{intro:sec:ata}.
In this section, however, we assume that computation times are fixed and known, so such considerations are unnecessary.

Since each iteration of this method takes $\tau_1$ seconds, and the iteration complexity of \sgd is known to be \citep{ghadimi2013stochastic}
$$
    K = \cO\!\left(\frac{L\Delta}{\varepsilon} + \frac{\sigma^2 L\Delta}{\varepsilon^2}\right),
$$
the total time complexity becomes
$$
    T = \tau_1 K~.
$$
While simple, this approach completely fails to exploit the computational power of the remaining $n-1$ machines, which remain idle throughout the training process.
\subsection{\naiveminibatchtitle}
\begin{figure*}[h]
    \centering
    \includegraphics[width=0.66\textwidth]{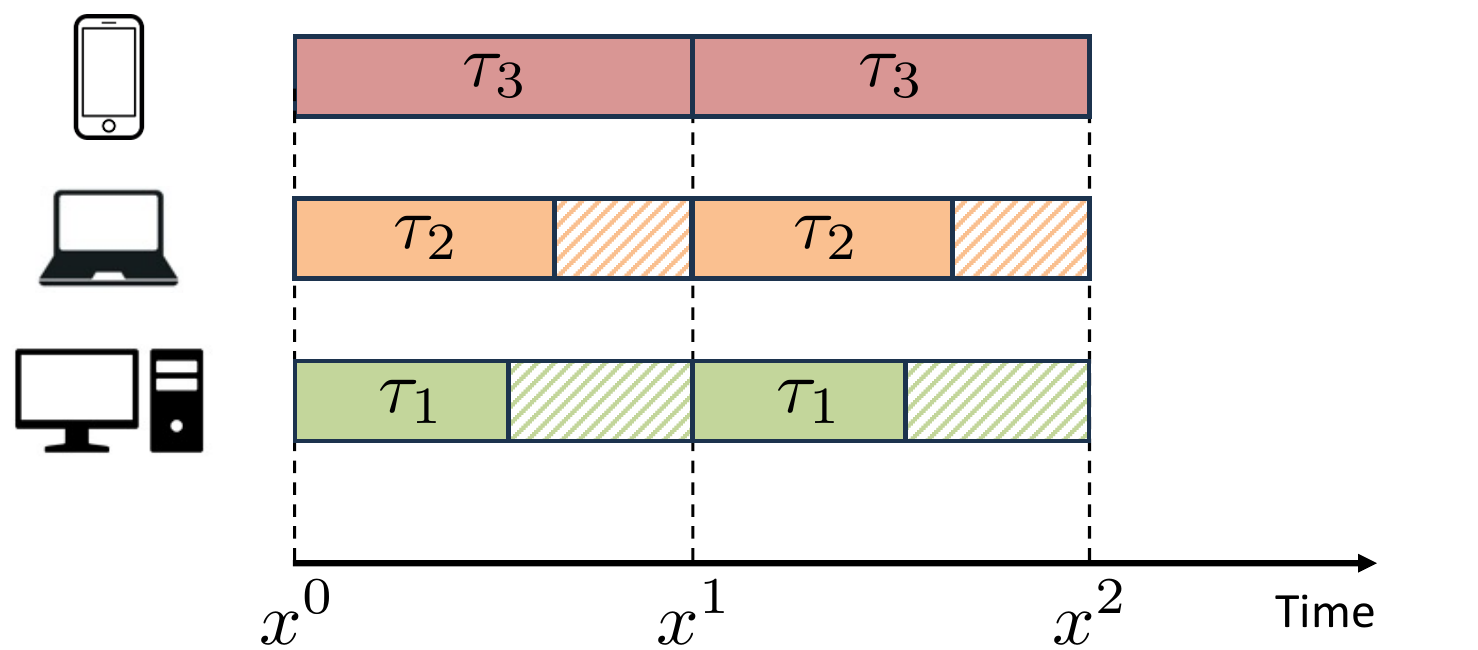}
    \caption{
        Illustration of \naiveminibatch (\Cref{intro:alg:naive_minibatch}) with 3 workers.
        The shaded regions indicate idle periods when faster devices must wait for the slowest one---wasting compute and doing nothing useful.
        }
    \label{intro:fig:naive_minibatch}
\end{figure*}

A natural next step is to utilize all $n$ machines.
The most straightforward way to do this is to have each machine compute one stochastic gradient per iteration and then update the model using their average.
This gives the \naiveminibatch algorithm, formalized in \Cref{intro:alg:naive_minibatch}.
Similar to the previous case, we illustrate the behavior of \naiveminibatch using a three-worker example in \Cref{intro:fig:naive_minibatch}.

\begin{algorithm}[h]
    \caption{\naiveminibatchtitle}
    \label{intro:alg:naive_minibatch}
    \begin{algorithmic}[1]
        \STATE \textbf{Input:} initial point $x^0 \in \R^d$, stepsizes $\gamma^k > 0$
        \FOR{$k = 0, \dots, K - 1$}
            \STATE Broadcast the current model $x^k$ to all workers
            \STATE Each worker $i$ computes a stochastic gradient $\nabla f\big(x^k; \xi_i^k\big)$
            \STATE Wait until all $n$ gradients have been received
            \STATE Aggregate the gradients and update the model:
            $$
                x^{k+1} = x^k - \gamma^k \frac{1}{n}\sum_{i=1}^n\nabla f\(x^k;\xi^k_i\)
            $$
        \ENDFOR
    \end{algorithmic}
\end{algorithm}

The iteration complexity of this method is similar to that of \herosgd.
Let $K(n)$ denote the number of iterations required by \naiveminibatch to reach an $\varepsilon$--stationary point; following \citet{ghadimi2013stochastic}, we have
\begin{equation}\label{intro:eq:naive_minibatch_iteration_complexity}
    K(n) = \cO\(\frac{L\Delta}{\varepsilon} + \frac{\sigma^2L\Delta}{n\varepsilon^2}\).
\end{equation}
The key difference is that the variance term is now divided by $n$, reflecting the variance reduction achieved by averaging $n$ independent gradient estimates.

However, each iteration requires $\tau_n$ seconds\textemdash the time needed for the slowest worker to complete its computation.
As a result, all faster workers remain idle until the slowest machine finishes before the next iteration can begin.
Therefore, the total time complexity is
$$
    \tau_n K(n)~.
$$
This synchronization overhead can be substantial in heterogeneous environments: when worker speeds differ greatly, the slowest device dictates the overall pace of training.  
In extreme cases, when $\tau_n$ is significantly larger than the others, the advantage of using multiple workers is almost entirely lost.
\subsection{\naiveasyncsgdtitle}
\begin{figure*}[h]
    \centering
    \includegraphics[width=0.66\textwidth]{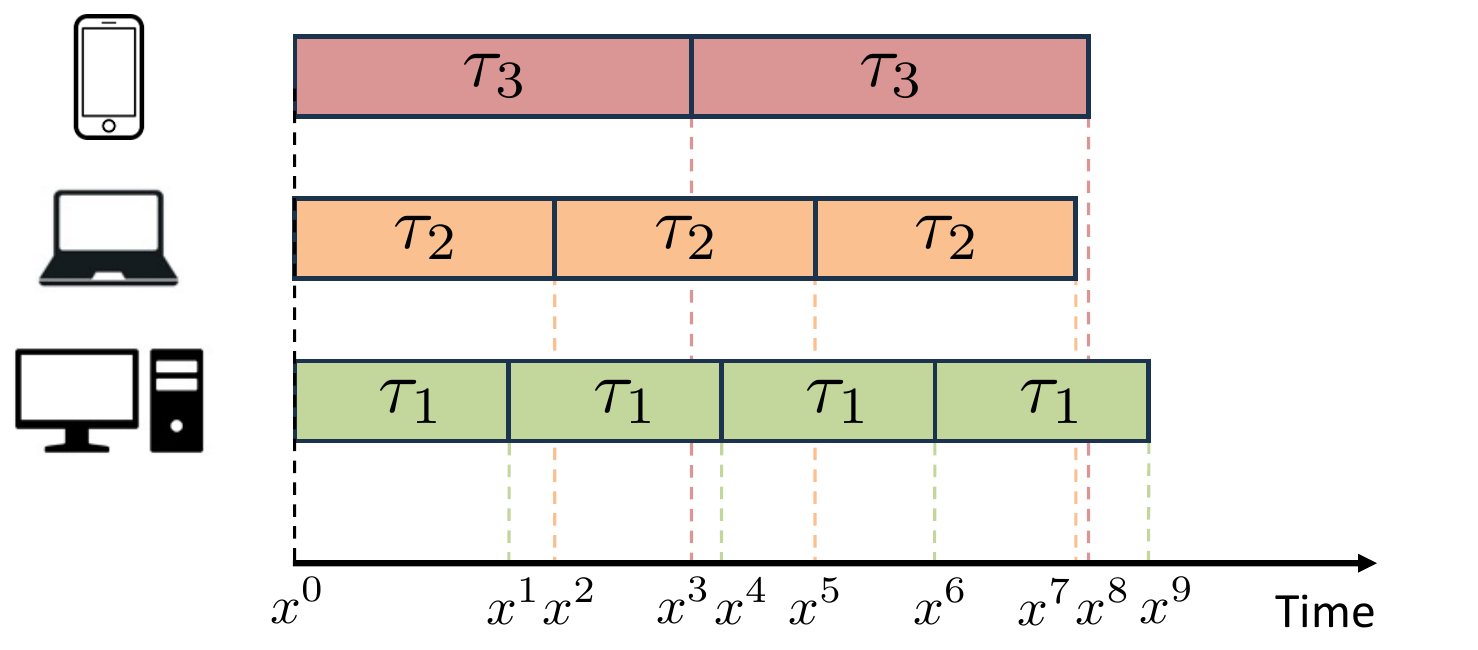}
    \caption{
        Illustration of \naiveasyncsgd (\Cref{intro:alg:asgd}) with 3 workers.
        Note the absence of shaded regions: all workers remain active, and every computed gradient contributes to a model update.
        }
    \label{intro:fig:asgd}
\end{figure*}

To address the inefficiency of idle workers in synchronous methods, researchers proposed the idea of \emph{asynchronous} \sgd \citep{recht2011hogwild,agarwal2011distributed}.
The key idea is simple: allow all workers to operate continuously and independently, without waiting for others.
Each worker computes stochastic gradients on its own copy of the model, and the server performs an update immediately whenever any worker finishes its computation.
We refer to this simple and naive algorithm, which updates the model as soon as any gradient arrives, as \naiveasyncsgd, and its pseudocode is given in \Cref{intro:alg:asgd}.
As before, we provide a three-worker illustration of \naiveasyncsgd in \Cref{intro:fig:asgd}.

\begin{algorithm}[h]
    \caption{\naiveasyncsgd}\label{intro:alg:asgd}
    \begin{algorithmic}[1]
        \STATE \textbf{Input:} initial point $x^0 \in \R^d$, stepsizes $\gamma^k > 0$
        \STATE Workers start computing stochastic gradients at $x^0$
        \FOR{$k = 0, \dots, K - 1$}
            \STATE A gradient $\nabla f\big(x^{k-\delta^k}; \xi^{k-\delta^k}_{i^k}\big)$ arrives from some worker $i^k$
            \STATE Update the model: $x^{k+1} = x^{k} - \gamma^k \nabla f\big(x^{k-\delta^k}; \xi^{k-\delta^k}_{i^k}\big)$
            \STATE Worker $i^k$ begins calculating $\nabla f\big(x^{k+1}; \xi^{k+1}_{i^k}\big)$
        \ENDFOR
    \end{algorithmic}
\end{algorithm}

In \Cref{intro:alg:asgd}, the term $\delta^k$ denotes the \emph{delay} of the gradient received at iteration $k$.
Specifically, $\delta^k$ counts how many model updates have occurred on the server since worker $i^k$ started computing its gradient.
Delays arise naturally because workers operate asynchronously\textemdash while one worker is computing, others may complete their gradients and perform model updates on the server.

While this algorithm ensures that all workers remain busy at all times, eliminating idle computation, it introduces a new challenge: \emph{staleness}.
Gradients computed on outdated models (large $\delta^k$) may misguide the updates, causing the optimization to progress in the wrong direction or converge more slowly.

Over the past decade, numerous variants of asynchronous \sgd have been proposed to mitigate this effect, primarily by improving how convergence bounds depend on the maximum delay \citep{cohen2021asynchronous, koloskova2022sharper, mishchenko2022asynchronous}.
Despite this progress, theoretical guarantees often remain pessimistic when delays are large.

A different line of research\textemdash initiated by \citet{tyurin2024optimal}\textemdash took an alternative path.
Rather than attempting to mitigate the delay issues inherent to asynchronous methods, this work turned attention back to the synchronous setting, aiming to squeeze as much performance as possible from it.
They proposed a method called \rennala.
\subsection{\rennalatitle}
\begin{figure*}[h]
    \centering
    \includegraphics[width=0.66\textwidth]{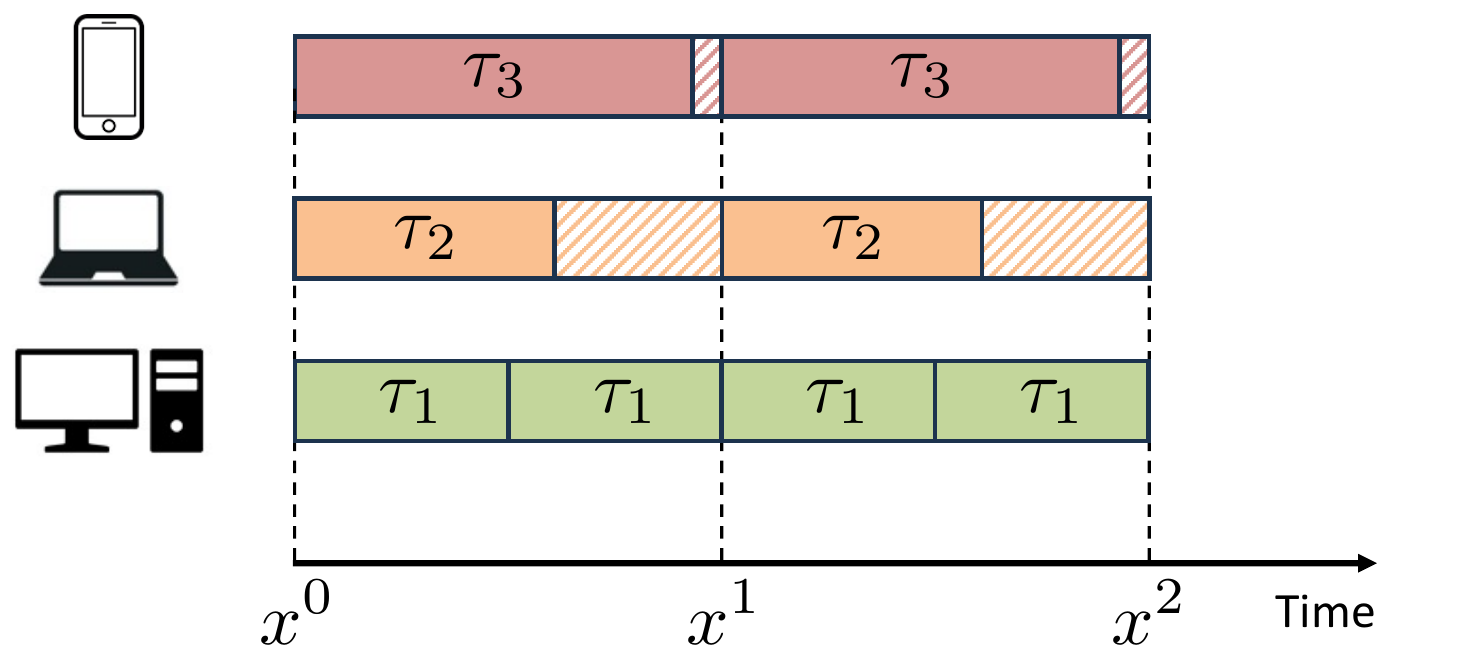}
    \caption{
        Illustration of \rennala (\Cref{intro:alg:rennala}) with 3 workers.
        The shaded regions indicate idle periods or gradient computations that are eventually discarded during synchronization.
        }
    \label{intro:fig:rennala}
\end{figure*}

The \rennala algorithm addresses the inefficiency of \naiveminibatch via a key insight: since it does not matter which worker computes a particular stochastic gradient, there is no need to collect exactly one gradient per worker.
Instead, the server continuously collects gradients as they arrive and performs updates as soon as ``enough'' have been gathered.
The algorithm is summarized in \Cref{intro:alg:rennala}.
Following the same setup, we show the behavior of \rennala with three workers in \Cref{intro:fig:rennala}.

\begin{algorithm}[h]
    \caption{\rennala \citep{tyurin2024optimal}}\label{intro:alg:rennala}
    \begin{algorithmic}[1]
        \STATE \textbf{Input:} initial point $x^0 \in \R^d$, stepsizes $\gamma^k > 0$, batch size $B\in \{1,2,\ldots\}$
        \FOR{$k = 0, \dots, K - 1$}
            \STATE Broadcast $x^k$ to all workers to compute stochastic gradients
            \label{intro:rennala:line:broadcast}
            \STATE Initialize $g^k = 0$ and $b = 0$
            \WHILE{$b < B$}
                \STATE A gradient $\nabla f\big(x^k; \xi_{i^{k,b}}^{k,b}\big)$ arrives from worker $i^{k,b}$
                \STATE $g^k = g^k + \nabla f\big(x^k; \xi_{i^{k,b}}^{k,b}\big)$
                \STATE Worker $i^{k,b}$ immediately begins computing a new gradient at $x^k$
                \STATE $b = b + 1$
            \ENDWHILE
            \STATE Update the model: 
            $$
                x^{k+1} = x^{k} - \gamma^k \frac{g^k}{B}
            $$
        \ENDFOR
    \end{algorithmic}
\end{algorithm}

The parameter $B$ denotes the number of gradients that the server aims to collect before performing an update.
These gradients are gathered asynchronously---as soon as they arrive---until $B$ in total have been received.
However, because \rennala is ultimately a \emph{synchronous} method, workers must synchronize once the server updates the model.
Specifically, as shown in Line~\ref{intro:rennala:line:broadcast}, when the server broadcasts the new model, all workers must \emph{discard} any ongoing gradient computations and restart from the updated point.
Alternatively, if the server knows the computation times of all workers, it can schedule the tasks so that the required $B$ gradients are collected in the minimum possible time; in that case, some workers may remain \emph{idle} once their assigned gradients are completed, waiting for the next synchronization.

The iteration complexity is analogous to the \naiveminibatch case \eqref{intro:eq:naive_minibatch_iteration_complexity}, since this is fundamentally minibatch \sgd with asynchronous batch collection:
$$
    K(B) = \cO\(\frac{L\Delta}{\varepsilon} + \frac{\sigma^2L\Delta}{B\varepsilon^2}\).
$$
Because $B$ gradients are averaged, the variance term now scales with $B$ instead of $n$.

The time per iteration, however, needs more reasoning.
As shown by \citet{tyurin2024optimal} (see also the alternative proof in \Cref{ringmaster:proof:time_R}), an upper bound on the time needed to collect $B$ gradients is given by
$$
    t(B) \eqdef 2 \min\limits_{m \in \{1,\ldots,n\}} \left(\frac{1}{m} \sum\limits_{i=1}^m \frac{1}{\tau_i}\right)^{-1} \left( 1 + \frac{B}{m} \right) . 
$$
The total time complexity therefore becomes
$$
    t(B)K(B)~.
$$
Remarkably, \citet{tyurin2024optimal} proved that by choosing
$$
    B=\max\left\{1,\ceil{\frac{\sigma^2}{\varepsilon}}\right\}~,
$$
\rennala attains the \emph{optimal time complexity}.
They established this result by deriving lower bounds for all first-order stochastic parallel methods applied to smooth nonconvex objectives\textemdash the exact setting considered in this dissertation\textemdash and showing that \rennala's complexity matches these bounds.

One important observation is that \rennala exhibits behavior similar to an asynchronous algorithm in that it reduces idle time by allowing workers to proceed without waiting for stragglers.
It almost fits the \emph{Asynchronous Meta-Algorithm} framework in \Cref{intro:algo:general_asynchronous}.
Since the update in \Cref{intro:algo:general_asynchronous} is optional, the server can simply accumulate gradients\textemdash just as \rennala does\textemdash and perform an update once $B$ gradients have been collected.
The only difference is that \Cref{intro:algo:general_asynchronous} cannot broadcast the updated model to all workers simultaneously.
However, instead of broadcasting to all workers at once, the server can send the updated model to each worker individually when they submit their next gradient.
In fact, this is how \rennala was originally introduced by \citet{tyurin2024optimal}.

Nevertheless, we do not classify \rennala as an asynchronous algorithm due to its update rule, as described in \Cref{intro:sec:sync_vs_async}.
\section{The Beginning of This Work}
Beyond establishing that \rennala achieves optimal time complexity, \citet{tyurin2024optimal} also showed---perhaps surprisingly---that no existing asynchronous \sgd algorithm attains this optimality.
These results raise a fundamental question:
\begin{center}
    \emph{Are asynchronous algorithms fundamentally flawed?}
\end{center}
If the optimal solution lies within the class of synchronous algorithms, should the community abandon asynchronous \sgd and focus exclusively on synchronous methods?
Was the long-standing enthusiasm for asynchrony misplaced?

This question forms the central motivation of this dissertation.
As we demonstrate in the chapters that follow, the answer is \emph{no}.
It is indeed possible to design asynchronous \sgd\ algorithms that achieve optimal time complexity.

Moreover, these asynchronous methods not only attain theoretically optimal time complexity but also outperform synchronous algorithms in practice.
As discussed in \Cref{intro:sec:sync_vs_async}, synchronous algorithms\textemdash even the optimal ones, such as \rennala---suffer from some computational waste: some devices remain idle, or their work is discarded during synchronization.
Asynchronous methods can potentially eliminate this inefficiency, which may lead to additional performance gains in real-world settings.

In this dissertation, we restore the theoretical and practical significance of asynchrony by introducing new asynchronous algorithms that are both provably optimal and empirically efficient.
\section{Contributions and Structure}

\begin{table*}
    \caption{
        Each chapter of this dissertation introduces a main algorithmic contribution, listed in the second column together with its reference.
        The third column indicates whether the method achieves optimal time complexity, i.e., whether it attains the theoretical lower bound on convergence time.
        The fourth column specifies the data regime in which each method operates\textemdash homogeneous or heterogeneous.
        The fifth column describes the worker computation time model used in the analysis, while the sixth column indicates whether the method is ``resource efficient''.
        The first two methods focus solely on minimizing wall-clock time and therefore utilize all available computational resources, whereas the third method, \ata, improves computational efficiency while losing only constant factors in the wall-clock time.
    }
    \label{table:chapters}
    \centering
    \begin{threeparttable}
      \resizebox{\textwidth}{!}{
        \begin{tabular}{cccccc}
            \toprule
                \bf Chapter 
                & \bf \makecell{Algorithm \\ and Reference}
                & \bf \makecell{Optimal \\ time \\ complexity}
                & \bf \makecell{Data \\ Heterogeneity}
                & \bf \makecell{Worker \\ Computation \\ Dynamics}
                & \bf \makecell{Resource \\ Efficient}
                \\
            \midrule
                \ref{chapter:Ringmaster}
                & \makecell{\ringmaster \\ \citep{maranjyan2025ringmaster}}
                & \checkmarkgreen 
                & \crossmarkred 
                & \makecell{Fixed and \\ Arbitrarily Changing}
                & \crossmarkred \\
            \midrule
                \ref{chapter:Ringleader}
                & \makecell{\ringleader \\ \citep{maranjyan2025ringleader}}
                & \centering \checkmarkgreen 
                & \checkmarkgreen 
                & \makecell{Fixed and \\ Arbitrarily Changing}
                & \crossmarkred \\
            \midrule
                \ref{chapter:ATA}
                & \makecell{\ata\textsuperscript{{\linkcolor (a)}} \\ \citep{maranjyan2025ata}}
                & \centering \checkmarkgreen\textsuperscript{{\linkcolor (b)}} 
                & \crossmarkred 
                & Stochastic
                & \checkmarkgreen \\
            \bottomrule
        \end{tabular}%
    }

    \resizebox{\textwidth}{!}{
        \begin{minipage}{\textwidth}
            \begin{tablenotes}[para,flushleft]
                \footnotesize
                \item[{\linkcolor (a)}]
                    \ata (Adaptive Task Allocation) is a meta-algorithm that can be applied on top of other methods\textemdash such as \ringmaster, \rennala, or \herosgd---to improve resource efficiency.
                \\
                \item[{\linkcolor (b)}]
                    Here, optimality is defined in a different sense: we compare the time complexity to that of a fixed competitor that knows the distribution of worker times, and show that \ata achieves performance within a constant multiplicative factor of the best fixed allocation in expectation.
            \end{tablenotes}
        \end{minipage}
    }

    \end{threeparttable}
\end{table*} 

This dissertation focuses on the design of asynchronous \sgd algorithms that achieve \emph{theoretically optimal time complexity}.
We develop optimal algorithms for both the homogeneous and heterogeneous data settings.
For the homogeneous case, we propose \ringmaster (\Cref{chapter:Ringmaster}); 
for the heterogeneous case, we introduce \ringleader (\Cref{chapter:Ringleader}).
The works of \citet{maranjyan2025ringmaster,maranjyan2025ringleader} introduced these algorithms, which, to the best of our knowledge, are the first \emph{theoretically optimal asynchronous \sgd methods}.
The key contributions of these chapters are summarized and discussed in more detail in \Cref{intro:sec:ringmaster} and \Cref{intro:sec:ringleader}, respectively.

Having established optimal methods in both settings, we next focus on improving their \emph{resource efficiency}.
While the optimal methods minimize the total \emph{time to convergence}, they do so by utilizing all available computational resources.
We show, however, that it is possible to slightly relax time optimality while achieving significant savings in resource usage.
To address this trade-off, we develop the \emph{Adaptive Task Allocation} (\ata) algorithm, presented in \Cref{chapter:ATA}, which introduces strategies that make \ringmaster and related methods in the homogeneous data regime more \emph{resource efficient}.
Further discussion of these ideas is provided in \Cref{intro:sec:ata}.

Finally, \Cref{table:chapters} summarizes the main contributions and distinctions between the chapters.
As summarized in Table \ref{table:chapters}, the dissertation progresses from optimal asynchronous methods (\ringmaster, \ringleader) toward resource-aware extensions (\ata).
The columns emphasize the distinct aspects of each method\textemdash optimality, heterogeneity, computation-time models, and efficiency.
In the following subsections, we provide a detailed overview of each chapter, clarifying the meaning of the table's columns and elaborating on specific contributions.
\subsection{\Cref{chapter:Ringmaster} -- Ringmaster ASGD: The First Asynchronous SGD with Optimal Time Complexity}\label{intro:sec:ringmaster}
In this chapter, we focus on the easier case\textemdash the homogeneous data setting.
As discussed earlier in \Cref{intro:sec:homo_hetero}, although this setting may seem simpler, it is by no means less important.
It is not merely a simplifying theoretical assumption; rather, it represents a distinct and practically relevant regime that deserves separate study.
Indeed, this setting closely reflects how most large-scale data centers operate today.

Our goal in this chapter is to design an asynchronous \sgd method that achieves \emph{optimal time complexity}\textemdash that is, whose upper bound matches the theoretical lower bound established by \citet{tyurin2024optimal}.
In other words, we aim to construct an algorithm that is as fast as theoretically possible.
Previously, such optimality had only been achieved by a \emph{synchronous} method---\rennala~\citep{tyurin2024optimal}.

We find that it is possible to achieve such optimality in an asynchronous setting by making a simple yet crucial modification to the \naiveasyncsgd algorithm.
Several recent works have attempted to improve \naiveasyncsgd, for example by introducing delay-adaptive stepsizes \citep{mishchenko2022asynchronous, koloskova2022sharper}.
While these methods outperform the synchronous \minibatch baseline, they remain suboptimal in general.
Our redesign takes a different approach: it leverages one of the possible server actions introduced in \Cref{intro:sec:dynamics}\textemdash the ability to \emph{discard} certain incoming stochastic gradients.

The key idea is the following.
The main issue in asynchronous methods is \emph{staleness}\textemdash some stochastic gradients are computed using outdated model parameters, which can harm convergence if the model is updated using them.
The vanilla \naiveasyncsgd algorithm (\Cref{intro:alg:asgd}) greedily applies all received gradients, regardless of their delay.
Instead, we propose to \emph{discard} gradients that are excessively stale, i.e., those whose delay exceeds a certain threshold.
By selecting this threshold carefully, we can provably attain optimal time complexity.

We call the resulting method \ringmaster, inspired by the metaphor of a circus ringmaster who brings order to chaos\textemdash taming the ``wild'' behavior of asynchronous workers and delayed gradients that characterizes \hogwild~\citep{recht2011hogwild}.
The algorithm is formalized as \Cref{intro:alg:ringmasternew}.
\begin{algorithm}[H]
    \caption{\ringmaster}
    \label{intro:alg:ringmasternew}
    \begin{algorithmic}[1]
        \STATE \textbf{Input:} initial point $x^0 \in \R^d$, stepsizes $\gamma^k > 0,$ delay threshold $R\ge 1$
        \STATE Set $k = 0$
        \STATE Workers start computing stochastic gradients at $x^0$
        \WHILE{not terminated}
            \STATE Gradient $\nabla f\big(x^{k-\delta^k}; \xi^{k-\delta^k}_{i}\big)$ arrives from worker $i$
            \IF{$\delta^k < R$}
                \STATE Update the model: 
                    $x^{k+1} = x^{k} - \gamma^k \nabla f\big(x^{k-\delta^k}; \xi^{k-\delta^k}_{i}\big)$
                \STATE Worker $i$ immediately begins calculating
                    $\nabla f\big(x^{k+1}; \xi^{k+1}_{i}\big)$
                \STATE Increment the iteration counter: $k = k + 1$
            \ELSE
                \STATE Discard the outdated gradient
                    $\nabla f\big(x^{k-\delta^k}; \xi^{k-\delta^k}_{i}\big)$
                \STATE Worker $i$ begins calculating
                    $\nabla f\big(x^{k}; \xi^{k}_{i}\big)$
            \ENDIF
        \ENDWHILE
    \end{algorithmic}
\end{algorithm}
Note that one can implement this more efficiently\textemdash instead of waiting for outdated gradients to arrive and then discarding them, the server can \emph{terminate} (another possible action described in \Cref{intro:sec:dynamics}) the computation of workers that are operating on excessively stale models.

Even though \ringmaster attains the same theoretical time complexity as the synchronous algorithm \rennala~(up to constant factors), it is often faster in practice due to two main factors:
%
%
%
\begin{itemize}
    \item \textbf{No synchronization waste.}  
    In \rennala, after each model update, all devices must synchronize with the server, discarding any partial computations that were still in progress.
    \ringmaster avoids this waste by synchronizing only those workers whose gradients have become stale.

    \item \textbf{More frequent updates.}  
    Whereas \rennala performs a single update after collecting $B$ gradients, \ringmaster updates the model whenever a fresh gradient arrives.
    This leads to faster progress, particularly in \emph{sparse models}, where each worker's gradient affects only a small subset of parameters with limited overlap across workers.
\end{itemize}
We validate this improvement through toy simulations presented in \Cref{ringmaster:sec:experiments}.

This chapter is based on the work of \citet{maranjyan2025ringmaster}:
\begin{quote}
    \bibentry{maranjyan2025ringmaster}
\end{quote}
\subsection{\Cref{chapter:Ringleader} -- Ringleader ASGD: The First Asynchronous SGD with Optimal Time Complexity \\ under Data Heterogeneity}\label{intro:sec:ringleader}

In this chapter, we turn from the homogeneous setting to the more general and more challenging \emph{heterogeneous data} case.
This setting is particularly important in applications such as \emph{federated learning}, where data is distributed across many clients (e.g., mobile devices or organizations), and each local dataset can follow a different distribution.
Moreover, these datasets are typically \emph{private} and cannot be shared or centralized due to privacy constraints.
Heterogeneity is therefore the rule rather than the exception in \emph{federated learning}.

Let us now discuss why this setting is more difficult than the homogeneous one.
Recall that our objective is to minimize the average of $n$ local functions:
\begin{equation*}
    \minimize_{x \in \mathbb{R}^d}
    \left\{ f(x) \coloneqq \frac{1}{n}\sum_{i=1}^{n} f_i(x) \right\}.
\end{equation*}
Suppose we apply the \naiveasyncsgd algorithm (\Cref{intro:alg:asgd}) in this setting.
The server updates the model whenever a gradient arrives from a worker.
However, in the heterogeneous case, these gradients correspond to different local functions $f_i$, and therefore may not accurately represent the global objective $f$.
As a result, naive asynchronous updates can introduce bias and lead to unstable training.

To address this, prior work on asynchronous optimization under data heterogeneity has typically used \emph{similarity assumptions} between the local objectives $f_i$ and the global objective $f$ (see, e.g.,~\citep{mishchenko2022asynchronous,koloskova2022sharper,nguyen2022federated,islamov2024asgrad}).
In practice, however, such assumptions rarely hold\textemdash especially in federated learning, where data of different users or devices can vary significantly (for example, in personalized text or image data).

The only known method in this heterogeneous data setting that achieves optimal time complexity \emph{without any similarity assumptions} is the synchronous method \malenia~\citep{tyurin2024optimal}\textemdash much like the situation in the homogeneous case.
Thus, the goal of this chapter is to design an \emph{asynchronous} method that achieves optimal time complexity \emph{without requiring any similarity assumptions}.
The hope is that such an \emph{asynchronous} method, while matching the optimal time complexity of the synchronous method \malenia, will perform better in practice by avoiding the computational waste caused by synchronization.

We accomplish this by leveraging another possible server action described in \Cref{intro:sec:dynamics}\textemdash the ability to \emph{store} received gradients rather than applying them immediately.
Specifically, we maintain a \emph{gradient table}, whose entries correspond to each worker's most recent gradients.
At each model update, the server averages all stored gradients, thereby incorporating information from all workers simultaneously.
This strategy eliminates the need for similarity assumptions, since the update uses a balanced view of the entire system.

This idea builds on the concept of maintaining and aggregating stored gradients, as in the \emph{Incremental Aggregated Gradient} (\algname{IAG})~\citep{blatt2007convergent,gurbuzbalaban2017convergence,vanli2018global} and \emph{Stochastic Averaged Gradient} (\algname{SAG})~\citep{roux2012stochastic,schmidt2017minimizing} methods.
An asynchronous extension of \algname{IAG}, called \iasgd, was recently proposed by \citet{wang2025incremental}, but their algorithm does not attain optimal time complexity.

While the gradient table idea removes the need for similarity assumptions, achieving optimal time complexity requires additional care.
As in the homogeneous setting, the main challenge lies in controlling \emph{staleness}.
Here, however, staleness arises because some entries in the gradient table may remain outdated if certain workers have not sent stochastic gradients for a long time.
If many iterations occur without refreshing these entries, the corresponding gradients become stale and degrade performance.

To address this issue, we modify the update mechanism to avoid performing a model update upon every received stochastic gradient.
Instead, we organize the \emph{updates} into \emph{rounds}.
In each round, the server performs $n$ updates\textemdash one per worker\textemdash and sends each worker the corresponding new model.
This ensures that, in the following round, the gradients computed by the workers are based on recent models.
Before the next round begins, the server clears the gradient table and starts collecting new gradients computed at these updated models.
This design ensures that the delays remain bounded, resulting in an optimal method.  
We call the resulting algorithm \ringleader---a natural successor to \ringmaster.

A nice property of our algorithm is that it never discards any gradients: all received gradients are stored in the table and eventually used for model updates.
In contrast, the optimal synchronous method \malenia~\citep{tyurin2024optimal} discards some workers' computations due to synchronization.
Therefore, although \ringleader achieves the same optimal time complexity (up to constants), it can be faster in practice\textemdash a result we demonstrate in the experimental section of the chapter (\Cref{sec:experiments}).

This chapter is based on the work of \citet{maranjyan2025ringleader}:
\begin{quote}
    \bibentry{maranjyan2025ringleader}
\end{quote}
\subsection{\Cref{chapter:ATA} -- ATA: Adaptive Task Allocation \\ for Efficient Resource Management \\ in Distributed Machine Learning}\label{intro:sec:ata}
In this chapter, we begin with the observation that the optimal parallel methods in the homogeneous data setting\textemdash \rennala~\citep{tyurin2024optimal} and \ringmaster---may result in significant resource waste.
The goal of these methods is purely to minimize training time, and to that end, they utilize all available resources at all times.
However, this often leads to redundant computations: a large portion of the work performed by slower machines may never be used.
We show that by slightly relaxing the objective of minimizing wall-clock training time, one can achieve a substantial reduction in resource waste.
The saved resources can then be allocated to other jobs or tasks.
\subsubsection{Where the Waste Comes From}
To illustrate this, consider \rennala, which provides a convenient example.
In each iteration, the server must collect exactly $B$ gradients before performing a model update (see \Cref{intro:alg:rennala}).
All $n$ workers continuously compute gradients until $B$ gradients have been received.
Once $B$ gradients arrive, the server updates the model and discards any ongoing computations.

For instance, consider a scenario with $n=1{,}010$ workers and $B=10$.
In this case, at least $1{,}000$ workers' computations will be discarded\textemdash since the server synchronizes after the first $10$ gradients arrive.
Ideally, we would avoid assigning tasks to those slower workers whose results will eventually be discarded.
However, this would require knowing in advance how long each worker will take to compute its gradient.
If these computation times were fixed and known, we could simply select the optimal subset of workers to minimize total computation time.
We instead consider a more challenging setup, in which the computation times are random variables.
Let us now formalize this idea.
\subsubsection{Problem Formulation}
We consider a system of $n$ workers, each capable of repeatedly performing the same task (e.g., computing a stochastic gradient).
In each round, the algorithm has a budget of $B$ tasks to distribute among the workers.
Allocating one unit of budget to a worker corresponds to assigning it one task.

We define the set of all feasible allocations as
$$
	\cA \eqdef \left\{ a \in \Z_+^n : \sum_{i=1}^n a_i = B \right\},
$$
where $a = (a_1, \dots, a_n)$ specifies how the total budget of $B$ tasks is distributed among the $n$ workers\textemdash that is, the $i$-th worker receives $a_i$ tasks\textemdash and $\mathbb{Z}_+$ denotes the set of nonnegative integers.
\subsubsection{Fixed Computation Times}
Under the fixed-time model described in \eqref{intro:eq:fixed_time}, where each worker $i$ requires $\tau_i$ seconds per task,
the time to receive $B$ completed gradients under allocation $a \in \mathcal{A}$ is
\begin{equation*}
	C(a) \eqdef \max_{i\in \{1,\ldots,n\}} \ a_i \tau_i ~.
\end{equation*}
Hence, the optimal allocation can be found by solving the combinatorial problem
$$
    a^* \in \argmin_{a \in\cA} C(a) ~.
$$
This allocation achieves the same wall-clock time as the greedy asynchronous collection strategy used in \Cref{intro:alg:rennala},
but with the minimum possible total computation, avoiding any wasted work.
\subsubsection{Random Computation Times}
In practice, however, worker computation times are not fixed and may change during the training.
We model this by assuming that the computation time of worker $i$ follows an unknown distribution $\nu_i$.
In round $k$, the time to complete the $u$-th task of worker $i$ is denoted by $X_i^{k,u}$, with $X_i^{k,u}\sim \nu_i$ iid across $u$.
If worker $i$ performs $a_i^k$ tasks in round $k$, the total computation time for that worker is
$$
	\sum_{u=1}^{a_i^k} X_i^{k,u}
$$
with the convention that this sum equals $0$ when $a_i^k = 0$.

Thus, the time to collect $B$ completed gradients under allocation $a^k$ is
$$
    C(a^k)
	    := \max_{i\in \{1,\ldots,n\}} \  \sum_{u=1}^{a_i^k} X_i^{k,u}~.
$$
Over $K$ rounds, the total expected computation time is
$$
	\cC_K := \E{ \sum_{k=1}^{K} C(a^k) } = \sum_{k=1}^{K} \E{C(a^k)}~.
$$
If the distributions ${\nu_i}$ were known, the optimal allocation would be obtained by solving
$$
    a^* \in \argmin_{a \in\cA} \ \E{C(a)} ~,
$$
and this same allocation would then be used in every round.
We refer to this strategy as the Optimal Fixed Task Allocation (OFTA).
\subsubsection{Goal: Adaptive Strategy without Knowing the Distributions}
In practice, the distributions ${\nu_i}$ are unknown.
Our goal is to design an adaptive allocation strategy whose total computation time $\mathcal{C}_K$ satisfies
\begin{equation}\label{intro:ata:eq:ata_obj}
    \cC_K \leq c \cdot \cC^*_K + \cE_K~,
\end{equation}
where $c \geq 1$ is a constant close to $1$, and $\cE_K$ is a negligible term compared to $\cC^*_K$ when $K\to \infty$.
This guarantee ensures that, asymptotically, our adaptive method performs nearly as well as the optimal strategy that knows all $\nu_i$.
\subsubsection{The ATA Algorithm}
We achieve this with our algorithm, Adaptive Task Allocation (\ata), which formulates the problem as a non-linear stochastic multi-armed bandit (MAB) problem \citep{LattimoreS18}.
In the standard MAB setting, an agent repeatedly chooses among several “arms” (here, workers), each associated with an unknown random outcome.
The agent must balance two objectives:
\begin{itemize}
    \item \textbf{Exploration:} learning the distribution of each arm.
    \item \textbf{Exploitation:} favoring the arms currently believed to yield better outcomes (faster workers, in our case).
\end{itemize}
In the MAB literature, this balance is often achieved using confidence bounds as estimators of each arm's mean reward (or cost).
In particular, algorithms based on Lower Confidence Bounds (LCBs) select arms using pessimistic estimates, ensuring sufficient exploration while minimizing regret \citep{auer2002finite}.

Following this principle, \ata maintains LCBs for each worker's expected computation time, updating them from observed samples.
At each iteration, given the LCB vector $s^k = (s_1^k, \ldots, s_n^k)$ and total budget $B$, the algorithm selects the next allocation by solving
$$
    a^k \in \argmin_{a \in \mathcal{A}} \ \max_{i \in \{1,\ldots,n\}} \  a_i^k s_i^k ~.
$$
This rule balances exploration and exploitation in a principled way and yields the performance guarantee stated in \eqref{intro:ata:eq:ata_obj}.

\subsubsection{The Asynchronous Case}
Since \ringmaster is the asynchronous version of \rennala, the \ata algorithm can also be applied on top of it.
In this setting, the allocation computed by \ata determines how many gradients each worker contributes, but instead of waiting to collect all gradients before updating the model, the server performs a model update immediately upon receiving each gradient.
At the same time, the server ensures that delays remain bounded by sending fresh models to any workers whose local copies become outdated during the process.

This chapter is based on the work of \citet{maranjyan2025ata}:
\begin{quote}
    \bibentry{maranjyan2025ata}
\end{quote}
\section{Publications Not Included in This Dissertation}
To maintain a coherent narrative and focus, several papers that I wrote during my PhD were not included in this dissertation.
For completeness, I also list earlier publications from before the PhD period.
All such works are summarized in \Cref{intro:table:excluded_papers}.
\begin{table}[thb]
    \caption{Publications not included in this dissertation.}
    \label{intro:table:excluded_papers}
    \centering
    \resizebox{\textwidth}{!}{
        \begin{tabular}{cc}
            \toprule
                \textbf{Reference and Title} & \textbf{During PhD} \\
            \midrule
                \makecell{\citet{maranjyan2025mindflayer} \\
                MindFlayer SGD: Efficient Parallel SGD in the Presence \\ of Heterogeneous and Random Worker Compute Times}
                & \checkmarkgreen \\
            \midrule
                \makecell{\citet{condat2025locodl} \\
                LoCoDL: Communication-Efficient Distributed Learning \\ with Local Training and Compression}
                & \checkmarkgreen \\
            \midrule
                \makecell{\citet{maranjyan2024differentially} \\
                Differentially Private Random Block Coordinate Descent}
                & \checkmarkgreen \\
            \midrule
                \makecell{\citet{maranjyan2025gradskip} \\
                GradSkip: Communication-Accelerated Local Gradient Methods \\ with Better Computational Complexity}
                & \crossmarkred \\
            \midrule
                \makecell{\citet{grigoryan2023menshov} \\
                Menshov-Type Theorem for Divergence Sets \\ of Sequences of Localized Operators}
                & \crossmarkred \\
            \midrule
                \makecell{\citet{grigoryan2021divergence} \\
                On the divergence of Fourier series in the general Haar system}
                & \crossmarkred \\
            \midrule
                \makecell{\citet{grigoryan2021unconditional} \\
                On the unconditional convergence of Faber-Schauder series in $L^1$}
                & \crossmarkred \\
            \bottomrule
        \end{tabular}
    }
\end{table}

Among these, let me briefly highlight one relevant work.
In \citet{maranjyan2025mindflayer}, we studied \rennala~\citep{tyurin2024optimal}—the optimal synchronous \sgd algorithm in the homogeneous data setting—under the assumption that each worker's computation time is a random variable.
We showed that \rennala can become significantly slower when the distributions of the workers' computation times have heavy tails.
To address this issue, we proposed a new algorithm, \algname{MindFlayer SGD}, which mitigates the effect of extreme delays by imposing a threshold on the maximum time allowed for each worker to compute its gradient.
Our analysis demonstrated that \algname{MindFlayer SGD} consistently outperforms \rennala and can achieve arbitrarily large speedups as the skewness of these distributions increases.
\section{Common Notations}
This section summarizes the common mathematical notations used throughout the dissertation.
\subsubsection{Basic Sets}
We denote by $\R$ the set of real numbers and by $\R_{+} \coloneqq [0, \infty)$ the set of nonnegative reals.
The set of natural numbers is denoted by $\N \coloneqq \{1, 2, \dots\}$, and the set of nonnegative integers by $\Z_{+} \coloneqq \{0, 1, 2, \dots\}$.
For a positive integer $n$, we use the shorthand notation $[n] \coloneqq \{1, 2, \dots, n\}$.
%
\subsubsection{Vectors and Norms}
Let $d \in \N$ denote the dimension of the vectors.
For vectors $a, b \in \R^d$, the standard inner product is defined as
$$
    \inp{a}{b} \coloneqq \sum_{i=1}^{d} a_i b_i~,
$$
where $a_i$ and $b_i$ are coordinates
and the corresponding Euclidean norm as
$$
    \norm{a} \coloneqq \sqrt{\inp{a}{a}}~.
$$
We also use $\norm{\cdot}_1$ and $\norm{\cdot}_\infty$ to denote the $\ell_1$ and $\ell_\infty$ norms, respectively:
$$
    \norm{a}_1 \coloneqq \sum_{i=1}^{d} \abs{a_i} ~,
    \qquad
    \norm{a}_\infty \coloneqq \max_{1 \le i \le d} \abs{a_i} ~.
$$
The symbol $\odot$ denotes the element-wise (Hadamard) product of two vectors.
For any scalar $x \in \R$, we write $(x)_{+} \coloneqq \max\{x, 0\}$.

For integers $a, b \in \Z$, we write $a \bmod b$ to denote the remainder when $a$ is divided by $b$, that is,
$$
    a \bmod b \coloneqq a - b \left\lfloor \frac{a}{b} \right\rfloor ~.
$$
\subsubsection{Probability and Expectation}
Expectation with respect to a random variable is denoted by $\E{\cdot}$.
For an event $\cE$, we write $\neg \cE$ for its complement (i.e., the event that $\cE$ does not occur).
We denote by $\mathds{1}(\cdot)$ the indicator function, which takes the value $1$ if the condition inside the parentheses is true and $0$ otherwise.

\subsubsection{Asymptotic Notation}
For functions $\phi, \psi : \mathcal{X} \to \R$, we use the standard asymptotic symbols:
\begin{itemize}[leftmargin=2em]
    \item $\phi = \cO(\psi)$ means that there exists a constant $C > 0$ such that 
    $$
        \phi(x) \le C\,\psi(x), \quad \forall x \in \mathcal{X}.
    $$
    \item $\phi = \Omega(\psi)$ means that there exists a constant $C > 0$ such that 
    $$
        \phi(x) \ge C\,\psi(x), \quad \forall x \in \mathcal{X}.
    $$
    \item $\phi = \Theta(\psi)$ means that both $\phi = \cO(\psi)$ and $\phi = \Omega(\psi)$ hold, 
    i.e., $\phi$ and $\psi$ are of the same order up to constant factors.
\end{itemize}


\chapter{Ringmaster ASGD: The First Asynchronous SGD \\ with Optimal Time Complexity}
\label{chapter:Ringmaster}
\thispagestyle{empty}
This chapter is based on the work of \citet{maranjyan2025ringmaster}:
\begin{quote}
    \bibentry{maranjyan2025ringmaster}
\end{quote}
\section{Introduction}
\label{ringmaster:sec:introduction}
We consider stochastic nonconvex optimization problems of the form
\begin{equation*}
    \minimize\limits_{x \in \R^d} \left\{f(x) \eqdef \ExpSub{\xi \sim \cD}{f(x;\xi)} \right\},
\end{equation*}
where $f \,:\, \R^d \times \mathbb{S}_{\xi} \to \R,$ $\R^d$ is a linear space, and $\mathbb{S}_{\xi}$ is a sample space.
In machine learning, $f(x;\xi)$ denotes the loss of a model parameterized by $x$ on a data sample $\xi$, and $\cD$ denotes the distribution of the training dataset.
In nonconvex optimization, our goal is to find an $\varepsilon$--stationary point, i.e., a (random) vector $x \in \R^d$ such that $\mathbb{E}[\|\nabla f(x)\|^2] \leq \varepsilon$.

We consider a setup involving $n$ workers (e.g., CPUs, GPUs, servers), each with access to the same distribution $\cD$.
Each worker is capable of computing independent, unbiased stochastic gradients with bounded variance (Assumption \ref{ringmaster:ass:stochastic_variance_bounded}).
We consider a setup with asynchronous, heterogeneous, and varying computation speeds.
We aim to account for all potential scenarios, such as random outages, varying computational performance over time, and the presence of slow or straggling workers \citep{dean2013tail}.

This setup is common in both data center environments \citep{dean2012large} and federated learning \citep{konevcny2016federated,mcmahan2016federated,kairouz2021advances} for distributed training. 
Although parallelism facilitates rapid convergence, variations in worker speeds make effective coordination more challenging.

Asynchronous Stochastic Gradient Descent is a popular approach for parallelization in such distributed settings.
We begin with a simple variant that keeps all workers continuously active and performs an update whenever a gradient arrives; we refer to this method as \naiveasyncsgd.
The method is outlined in \Cref{ringmaster:alg:asgd}.

\begin{algorithm}[H]
    \caption{\naiveasyncsgd}
    \label{ringmaster:alg:asgd}
    \begin{algorithmic}[1]
        \STATE \textbf{Input:} point $x^0 \in \R^d$, stepsizes $\gamma^k \geq 0$
        \STATE Workers start computing stochastic gradients at $x^0$
        \FOR{$k = 0, \dots, K - 1$}
            \STATE A gradient $\nabla f\big(x^{k-\delta^k}; \xi^{k-\delta^k}_{i^k}\big)$ arrives from worker $i^k$
            \STATE Update: $x^{k+1} = x^{k} - \gamma^k \nabla f\big(x^{k-\delta^k}; \xi^{k-\delta^k}_{i^k}\big)$
            \STATE Worker $i^k$ begins calculating $\nabla f\big(x^{k+1}; \xi^{k+1}_{i^k}\big)$
        \ENDFOR
    \end{algorithmic}
\end{algorithm}

This is a greedy and asynchronous method.
Once a worker finishes computation of the stochastic gradient, it immediately sends the gradient to the server, which updates the current iterate without waiting for other workers.
Notice that, unlike vanilla \sgd, the update is performed using the stochastic gradient calculated at the point $x^{k-\delta^k}$, where the index $k-\delta^k$ corresponds to the iteration when the worker started computing the gradient, which can be significantly outdated.
The sequence $\{\delta^k\}$ is a sequence of delays of asynchronous \sgd, where $\delta^k \geq 0$ is defined as the difference between the iteration when worker $i^k$ started computing the gradient and iteration $k$, when it was applied.

Asynchronous \sgd methods have a long history, originating in 1986 \citep{tsitsiklis1986distributed} and regaining prominence with the seminal work of \citet{recht2011hogwild}.
The core idea behind asynchronous \sgd is simple: to achieve fast convergence, all available resources are utilized by keeping all workers busy at all times.
This principle has been validated in numerous studies, showing that asynchronous \sgd can outperform naive synchronous \sgd methods \citep{feyzmahdavian2016asynchronous, dutta2018slow, nguyen2018sgd, arjevani2020tight, cohen2021asynchronous, mishchenko2022asynchronous, koloskova2022sharper, islamov2024asgrad, feyzmahdavian2023asynchronous}.

\begin{table}
    \caption{
        The time complexities of asynchronous stochastic gradient methods, which preform the step $x^{k+1} = x^{k} - \gamma^k \nabla f(x^{k-\delta^k}; \xi^{k-\delta^k}_{i^k}),$ to get an $\varepsilon$--stationary point in the nonconvex setting. 
        In this table, we consider the \emph{fixed computation model} from \Cref{ringmaster:sec:prel}. 
        Abbr.: $\sigma^2$ is defined as ${\rm \mathbb{E}}_{\xi}[\|\nabla f(x;\xi) - \nabla f(x)\|^2] \leq \sigma^2$ for all $x \in \R^d,$ $L$ is the smoothness constant of $f$, $\Delta \eqdef f(x^0) - f^*$, $\tau_i \in [0, \infty]$ is the time bound to compute a single stochastic gradient by worker $i^k$.
    }
    \label{ringmaster:table:complexities}
    \centering 
    \resizebox{\textwidth}{!}{
    \begin{threeparttable}  
        \begin{tabular}[t]{cccc}
        \toprule
            \bf Algorithm
            & \bf Worst-Case Time Complexity
            & \bf Optimal
            & \bf \makecell{Adaptive to \\ Varying Speeds} \\
        \midrule
            \makecell{\algname{Delay-Adaptive Asynchronous SGD} \\ \citep{koloskova2022sharper} \\ \citep{mishchenko2022asynchronous}} 
            & $\left(\frac{1}{n} \sum\limits_{i=1}^n \frac{1}{\tau_i}\right)^{-1}\left(\frac{L \Delta}{\varepsilon} + \frac{\sigma^2 L \Delta}{n \varepsilon^2}\right)$ 
            & \crossmarkred & \checkmarkgreen \\
        \midrule
            \makecell{\algname{Naive Optimal ASGD} \textbf{(new)} \\ 
            (\Cref{ringmaster:alg:ringmaster_short}; \\ \Cref{ringmaster:thm:ringmaster_short})}
            & $\min \limits_{m \in [n]} \left(\frac{1}{m}\sum\limits_{i=1}^m \frac{1}{\tau_i}\right)^{-1}\left(\frac{L \Delta}{\varepsilon} + \frac{\sigma^2 L \Delta}{m \varepsilon^2}\right)$ 
            & \checkmarkgreen
            & \crossmarkred \\
        \midrule
            \makecell{\ringmaster \textbf{(new)} \\ (Algorithms~\ref{ringmaster:alg:ringmasternew} or \ref{ringmaster:alg:ringmasternewstops}; \\ \Cref{ringmaster:thm:optimal_ringmaster})} 
            & $\min \limits_{m \in [n]} \left(\frac{1}{m}\sum\limits_{i=1}^m \frac{1}{\tau_i}\right)^{-1}\left(\frac{L \Delta}{\varepsilon} + \frac{\sigma^2 L \Delta}{m \varepsilon^2}\right)$ 
            & \checkmarkgreen & \checkmarkgreen\\
        \midrule
        \midrule
            \makecell{Lower Bound \\ \citep{tyurin2024optimal}} 
            & $\min \limits_{m \in [n]} \left(\frac{1}{m}\sum\limits_{i=1}^m \frac{1}{\tau_i}\right)^{-1}\left(\frac{L \Delta}{\varepsilon} + \frac{\sigma^2 L \Delta}{m \varepsilon^2}\right)$ 
            & \textbf{---} & \textbf{---} \\
        \bottomrule
        \end{tabular}
    \end{threeparttable}
    }
\end{table}

\subsection{Assumptions}\label{ringmaster:sec:assumption}
In this paper, we consider the standard assumptions from the nonconvex world.
\begin{boxedassumption}[Smoothness]\label{ringmaster:ass:lipschitz_constant}
    The function $f$ is differentiable, and its gradient is Lipschitz continuous.
    That is, there exists a constant $L>0$ such that, for all $x, y \in \R^d$, 
    \begin{equation*}
        \norm{\nabla f(x) - \nabla f(y)} \le L\norm{x - y}.
    \end{equation*}
 \end{boxedassumption}
This assumption, often referred to as \emph{$L$-smoothness}, controls how quickly the gradient can change.
\begin{boxedassumption}[Lower boundedness]\label{ringmaster:ass:lower_bound}
    The objective function is bounded below: there exists $f^* > -\infty$ such that $f(x) \geq f^*$ for all $x \in \R^d$.
\end{boxedassumption}
This assumption ensures that the optimization problem is well-posed: without it, function values could decrease without bound, making convergence guarantees meaningless.

We define $\Delta \coloneqq f(x^0) - f^*$, where $x^0$ is the initial point of the optimization algorithm.
The quantity $\Delta$ represents the initial suboptimality and serves as a measure of the inherent difficulty of the problem.
\begin{boxedassumption}
    \label{ringmaster:ass:stochastic_variance_bounded}
    For every $\xi$ the function $f(x;\xi)$ is differentiable with respect to its first argument $x$.
    Moreover, the stochastic gradients are unbiased and have bounded variance $\sigma^2 \geq 0$, that is,
    \begin{gather*}
         \ExpSub{\xi}{\nabla f(x;\xi)} = \nabla f(x), \quad \forall x \in \R^d~, \\
         \ExpSub{\xi}{\sqnorm{\nabla f(x;\xi) - \nabla f(x)}} \leq \sigma^2, \quad \forall x \in \R^d~.
    \end{gather*}
\end{boxedassumption}
\subsection{Related Work}
Despite the variety of asynchronous \sgd algorithms proposed over the years, a fundamental question remained unresolved: 
\textit{What is the optimal strategy for parallelization in this setting?} 

When we have one worker, the optimal number of stochastic gradients required to find an $\varepsilon$--stationary point is 
$$
    \Theta\left(\frac{L \Delta}{\varepsilon} + \frac{\sigma^2 L \Delta}{\varepsilon^2}\right),
$$ 
achieved by the vanilla \sgd method \citep{ghadimi2013stochastic,arjevani2022lower}.
In the parallel setting with many workers, several approaches have been proposed to obtain oracle lower bounds \citep{scaman2017optimal,woodworth2018graph, arjevani2020tight,lu2021optimal}.
Recent work by \citet{tyurin2024optimal, tyurin2024tighttimecomplexitiesparallel} addressed the question by establishing lower bounds for the \emph{time complexity} of asynchronous methods under the \emph{fixed computation model} and the \emph{universal computation model}.
Surprisingly, they demonstrated that none of the existing asynchronous \sgd methods are optimal.
Moreover, they introduced a minimax optimal method, \rennala, which achieves the theoretical lower bound for time complexity.

\begin{algorithm}[t]
    \caption{\rennala \citep{tyurin2024optimal}}
    \label{ringmaster:alg:rennala}
    \begin{algorithmic}[1]
        \STATE \textbf{Input:} point $x^0 \in \R^d$, stepsize $\gamma > 0$, batch size $B\in \N$
        \STATE Workers start computing stochastic gradients at $x^0$
        \FOR{$k = 0, \dots, K - 1$}
            \STATE Initialize $g^k = 0$ and $b = 0$
            \WHILE{$b < B$}
                \STATE A gradient $\nabla f\big(x^{k - \delta^{k,b}}; \xi^{k,b}\big)$ arrives from worker $i^{k,b}$
                \highlightcolor
                \IF{$\delta^{k,b} = 0$}
                    \color{black}
                    \STATE $g^k = g^k + \nabla f\big(x^{k - \delta^{k,b}}; \xi^{k,b}\big)$
                \ENDIF
                \color{black}
                \highlightcolor
                \STATE Worker $i^{k,b}$ immediately begins computing a new gradient at $x^k$
                \color{black}
                \STATE $b=b+1$
            \ENDWHILE
            \STATE Update the model:
                $$
                    x^{k+1} = x^{k} - \gamma \frac{g^k}{B}
                $$
        \ENDFOR
    \end{algorithmic}
\end{algorithm}

\rennala is \textit{semi-asynchronous} and can be viewed as \minibatch (which takes {\em synchronous} iteration/model updates) combined with an {\em asynchronous} minibatch collection mechanism.
Let us now explain how \rennala works, in the notation of \Cref{ringmaster:alg:rennala}, which facilitates comparison with further methods described in this work.
Due to the condition $\delta^{k,b} = 0$, which ignores all stochastic gradients calculated at the previous points, \rennala performs the step 
$$
    x^{k+1} = x^{k} - \gamma \frac{1}{B} \sum_{j=1}^{B} \nabla f\(x^{k}; \xi^{k,j}\) ~,
$$
where $\xi^{k,1}, \dots,\xi^{k,b}$ are independent samples from $\cD$ collected asynchronously across all workers.
Note that the workers compute the stochastic gradients at the {\em same} point $x^k$, with worker $i$ computing $B_i \geq 0$ gradients such that $\sum_{i=1}^n B_i = B$.

This approach has at least two fundamental drawbacks:
\begin{itemize}
    \item Once the fastest worker completes the calculation of the first stochastic gradient, $\nabla f\big(x^k; \xi^{k,1}\big)$, it begins computing another stochastic gradient at the same point $x^k$, even though it already possesses additional information from $\nabla f\big(x^k; \xi^{k,1}\big)$.
    \rennala does not update iterate $x^k$ immediately.
    \item Once \rennala finishes the inner while loop, it will ignore all stochastic gradients that were being calculated before the loop ended, even if a worker started the calculation just a moment before.
    In contrast, \naiveasyncsgd avoids these issues by fully utilizing all currently available information when asking a worker to calculate the next stochastic gradient and not ignoring any stochastic gradients.
\end{itemize}
These revelations raise an intriguing question: 
\textit{Is asynchronous parallelization fundamentally flawed?}
If the optimal solution lies in synchronous approaches, should the community abandon asynchronous \sgd and redirect its focus to developing synchronous methods?
Perhaps the widespread enthusiasm for asynchronous \sgd methods was misplaced.

Alternatively, could there be a yet-to-be-discovered variant of asynchronous \sgd that achieves optimal time complexity?
In this work, we answer this question affirmatively.
We reestablish the prominence of asynchronous \sgd by proposing a novel asynchronous optimization method that attains optimal time complexity.
\subsection{Contributions}
\label{ringmaster:sec:contributions}
Our contributions are summarized as follows:
\begin{itemize}
    \item We introduce a novel asynchronous stochastic gradient descent method, \ringmaster, described in \Cref{ringmaster:alg:ringmasternew} and \Cref{ringmaster:alg:ringmasternewstops}.
    This is the first asynchronous method to achieve optimal time complexity under arbitrary heterogeneous worker compute times (see \Cref{ringmaster:table:complexities}).
    Specifically, in Theorems~\ref{ringmaster:thm:optimal_ringmaster} and \ref{ringmaster:thm:optimal_ringmaster_dynamic}, we establish time complexities that match the lower bounds developed by \citet{tyurin2024optimal,tyurin2024tighttimecomplexitiesparallel}.
    \item Our work begins with another new  optimal method, \algname{Naive Optimal ASGD} (\Cref{ringmaster:alg:ringmaster_short}).
    We demonstrate that \algname{Naive Optimal ASGD} achieves optimality under the fixed computation model.
    However, we find that it is overly simplistic and lacks robustness in scenarios where worker computation times are chaotic and dynamic.
    To address this limitation, we designed \ringmaster, which combines the strengths of \algname{Naive Optimal ASGD}, previous non-optimal versions of \naiveasyncsgd methods \citep{cohen2021asynchronous,koloskova2022sharper,mishchenko2022asynchronous}, and the semi-synchronous \rennala \citep{tyurin2024optimal}.
    \item All our claims are supported by rigorous theoretical analysis showing the optimality of the method under virtually any computation scenario,  including unpredictable downtimes, fluctuations in computational performance over time, delays caused by slow or straggling workers, and challenges in maintaining synchronization across distributed systems (see Sections~\ref{ringmaster:sec:theory} and \ref{ringmaster:sec:dyn}).
    Using numerical experiments, we demonstrate that \ringmaster outperforms existing methods (see Section~\ref{ringmaster:sec:experiments}).
\end{itemize}
\section{Preliminaries and Naive Method}
\label{ringmaster:sec:prel}
To compare methods, we consider the \emph{fixed computation model} \citep{mishchenko2022asynchronous}.
In this model, it is assumed that
\begin{equation}\label{ringmaster:eq:worker-time}
    \begin{minipage}{0.52\linewidth}
    Each worker~$i$ requires $\tau_i$~seconds to compute one stochastic gradient $\nabla f(x;\xi)$.
    
    \vspace{4pt}
    
    Without loss of generality, we assume
    
    \centering
    $0 < \tau_1 \le \tau_2 \le \cdots \le \tau_n$~.
    \end{minipage}
\end{equation}
However, in \Cref{ringmaster:sec:dyn}, we will discuss how one can easily generalize our result to arbitrary computational dynamics, e.g., when the computation times are not bounded by the fixed values $\{\tau_i\}$, and can change in arbitrary/chaotic manner in time.
Under the fixed computation model, \citet{tyurin2024optimal} proved that the optimal time complexity lower bound is
\begin{align}
    \label{ringmaster:eq:rennala}
    T_{\textnormal{R}} 
    \eqdef \Theta\left(\min\limits_{m \in [n]} \left(\frac{1}{m} \sum\limits_{i=1}^{m} \frac{1}{\tau_{i}}\right)^{-1} \left(\frac{L \Delta}{\varepsilon} + \frac{\sigma^2 L \Delta}{m \varepsilon^2}\right)\right)
\end{align}
seconds achieved by \rennala (\Cref{ringmaster:alg:rennala}).
However, the best analysis of \naiveasyncsgd \citep{koloskova2022sharper,mishchenko2022asynchronous} with appropriate stepsizes achieves the time complexity (see Sec.~L of \citet{tyurin2024optimal})
\begin{align}
    \label{ringmaster:eq:async}
     T_{\textnormal{A}} \eqdef \Theta\left(\left(\frac{1}{n} \sum\limits_{i=1}^{n} \frac{1}{\tau_{i}}\right)^{-1} \left(\frac{L \Delta}{\varepsilon} + \frac{\sigma^2 L \Delta}{n \varepsilon^2}\right)\right).
\end{align}
Note that $T_{\textnormal{R}} \leq T_{\textnormal{A}}$; this is because $\min_{m \in [n]} g(m) \leq g(n)$ for any function $g \,:\, \N \to \R$.
Moreover, $T_{\textnormal{R}}$ can {\em arbitrarily} smaller.
To illustrate the difference, consider an example with $\tau_i = \sqrt{i}$ for all $i \in [n]$.
Then, 
$$
    T_{\textnormal{R}} = \Theta\left(\max\left\{\frac{\sigma L \Delta}{\varepsilon^{3/2}}, \frac{L \Delta \sigma^2}{\sqrt{n} \varepsilon^2}\right\}\right)
$$
and 
$$
    T_{\textnormal{A}} = \Theta\left(\max\left\{\frac{\sqrt{n} L \Delta}{\varepsilon}, \frac{L \Delta \sigma^2}{\sqrt{n} \varepsilon^2}\right\}\right),
$$
see the derivations in Section~\ref{ringmaster:sec:deriv}.
If $n$ is large, as is often encountered in modern large-scale training scenarios, $T_{\textnormal{A}}$ can be arbitrarily larger than $T_{\textnormal{R}}$.
Thus, the best-known variants of asynchronous \sgd are not robust to the scenarios when the number of workers is large and computation times are heterogeneous/chaotic.

\subsection{A Naive Optimal Asynchronous SGD}
We now introduce our first simple and effective strategy to improve the time complexity $T_{\textnormal{A}}$.
Specifically, we hypothesize that selecting a {\em subset} of workers at the beginning of the optimization process, instead of utilizing all available workers, can lead to a more efficient and stable approach.
As we shall show, this adjustment not only simplifies the computational dynamics, but also proves sufficient to achieve the optimal time complexity.

The idea is to select the fastest $[m] \eqdef \{1, 2, \dots, m\}$ workers, thereby ignoring the slow ones and eliminating delayed gradient updates. 
We demonstrate that the optimal algorithm involves finding the ideal number of workers $m$ and running delay-adaptive version of \naiveasyncsgd \citep{koloskova2022sharper, mishchenko2022asynchronous} on those workers.
The method is formalized in \Cref{ringmaster:alg:ringmaster_short}.
\begin{algorithm}[H]
    \caption{\algname{Naive Optimal ASGD}}
    \label{ringmaster:alg:ringmaster_short}
    \begin{algorithmic}[1]
        \STATE Find 
        $$
            m_{\star} \in \argmin\limits_{m \in [n]} \left(\frac{1}{m} \sum\limits_{i=1}^m \frac{1}{\tau_i}\right)^{-1} \left( 1 + \frac{\sigma^2}{m \varepsilon} \right)
        $$
        \STATE Run \naiveasyncsgd (\Cref{ringmaster:alg:asgd}) on $[m_{\star}]$ workers
    \end{algorithmic}
\end{algorithm}
The choice of $m_\star$ in \Cref{ringmaster:alg:ringmaster_short} effectively selects the fastest $m_\star$ workers only.
Note that it is possible for $m_\star$ to be equal to $n$, meaning that all workers participate, which occurs when all workers are nearly equally fast.
In this case, the harmonic mean in Line 1 of \Cref{ringmaster:alg:ringmaster_short} remains unchanged if all $\tau_i$s are equal, but the right-hand side decreases.
However, if some workers experience large delays, the harmonic mean in Line 1 increases as $m$ grows, introducing a trade-off between the two factors.
Conversely, if most workers are very slow, it may be optimal to have as few as one worker participating.

Next, we can easily prove that our algorithm, \algname{Naive Optimal ASGD} (\Cref{ringmaster:alg:ringmaster_short}), is optimal in terms of time complexity.
\begin{boxedtheorem}
    \label{ringmaster:thm:ringmaster_short}
    Consider the \emph{fixed computation model} \eqref{ringmaster:eq:worker-time}.
    Let Assumptions \ref{ringmaster:ass:lipschitz_constant}, \ref{ringmaster:ass:lower_bound}, and \ref{ringmaster:ass:stochastic_variance_bounded} hold.
    Then \algname{Naive Optimal ASGD} (\Cref{ringmaster:alg:ringmaster_short}) with $m_{\star}$ workers achieves the optimal time complexity \eqref{ringmaster:eq:rennala}.
\end{boxedtheorem}
\begin{proof}
    The proof is straightforward.
    Indeed, the time complexity \eqref{ringmaster:eq:async} of \Cref{ringmaster:alg:asgd} with $m_*$ workers is
        \begin{align*}
            \Theta\left(\left(\frac{1}{m_*} \sum\limits_{i=1}^{m_*} \frac{1}{\tau_{i}}\right)^{-1} \left(\frac{L \Delta}{\varepsilon} + \frac{\sigma^2 L \Delta}{m_* \varepsilon^2}\right)\right),
        \end{align*}
        which equals to \eqref{ringmaster:eq:rennala} due to the definition of $m_*$ in \Cref{ringmaster:alg:ringmaster_short}.
\end{proof}
To the best of our knowledge, this is the first variant of asynchronous \sgd that provides guarantees for achieving the optimal time complexity.

\subsection{Why Is \Cref{ringmaster:alg:ringmaster_short} Referred to as ``Naive''?}

Note that determining the optimal $m_{\star}$ requires the knowledge of the computation times $\tau_1, \dots, \tau_n$.
If the workers' computation times were indeed static in time, this would not be an issue, as these times could be obtained by querying a single gradient from each worker before the algorithm is run. However, in real systems, computation times are rarely static, and can vary with from iteration to iteration \citep{dean2013tail,chen2016revisiting, dutta2018slow, maranjyan2025mindflayer}, or even become infinite at times, indicating down-time.

Naively selecting the fastest $m_*$ workers at the start of the method and keeping this selection unchanged may therefore lead to significant issues in practice. The computational environment may exhibit adversarial behavior, where worker speeds change over time.
For instance, initially, the first worker may be the fastest, while the last worker is the slowest.
In such cases, \algname{Naive Optimal ASGD} would exclude the slowest worker.
However, as time progresses, their performance may reverse, causing the initially selected $m_*$ workers, including the first worker, to become the slowest.
This exposes a critical limitation of the strategy: it lacks robustness to time-varying worker speeds.

\section{Ringmaster ASGD}
We are now ready to present our new versions of asynchronous \sgd, called \ringmaster (\Cref{ringmaster:alg:ringmasternew} and \Cref{ringmaster:alg:ringmasternewstops}), which guarantee the \emph{optimal time complexity} without knowing the computation times a priori. 
Both methods are equivalent, up to  a minor detail that we shall discuss later. 
Let us first focus on \Cref{ringmaster:alg:ringmasternew}.

\begin{algorithm}[H]
    \caption{\ringmaster (without calculation stops)}
    \label{ringmaster:alg:ringmasternew}
    \begin{algorithmic}[1]
        \STATE \textbf{Input:} point $x^0 \in \R^d$, stepsize $\gamma > 0,$ delay threshold $R \in \N$
        \STATE Set $k = 0$
        \STATE Workers start computing stochastic gradients at $x^0$
        \WHILE{not terminated}
            \STATE Gradient $\nabla f\big(x^{k-\delta^k}; \xi^{k-\delta^k}_{i}\big)$ arrives from worker $i$
            \highlightcolor
            \IF{$\delta^k < R$}
                \color{black}
                \STATE Update the model:
                    $x^{k+1} = x^{k} - \gamma \nabla f\big(x^{k-\delta^k}; \xi^{k-\delta^k}_{i}\big)$
                \STATE Worker $i$ begins calculating 
                    $\nabla f\big(x^{k+1}; \xi^{k+1}_{i}\big)$
                \STATE Update the iteration number $k = k + 1$
            \ELSE
                \STATE Ignore the outdated gradient 
                    $\nabla f\big(x^{k-\delta^k}; \xi^{k-\delta^k}_{i}\big)$
                \STATE Worker $i$ begins calculating 
                    $\nabla f\big(x^{k}; \xi^{k}_{i}\big)$
            \ENDIF
        \ENDWHILE
    \end{algorithmic}
\end{algorithm}
\begin{algorithm}[H]
    \caption{\ringmaster (with calculation stops)}
    \label{ringmaster:alg:ringmasternewstops}
    \begin{algorithmic}[1]
        \STATE \textbf{Input:} point $x^0 \in \R^d$, stepsize $\gamma > 0$, delay threshold $R \in \N$
        \STATE Set $k = 0$
        \STATE Workers start computing stochastic gradients at $x^0$
        \WHILE{not terminated}
            \STATE {\highlightcolor Stop calculating stochastic gradients with delays $\geq R$, and start computing new ones at $x^k$ instead}
            \STATE Gradient $\nabla f\big(x^{k-\delta^k}; \xi^{k-\delta^k}_{i}\big)$ arrives from worker $i$
            \STATE Update the model:
                $x^{k+1} = x^{k} - \gamma \nabla f\big(x^{k-\delta^k}; \xi^{k-\delta^k}_{i}\big)$
            \STATE Worker $i$ begins calculating 
                $\nabla f\big(x^{k+1}; \xi^{k+1}_{i}\big)$
            \STATE Update the iteration number $k = k + 1$
        \ENDWHILE
    \end{algorithmic}
    (The core and essential modifications to Alg.~\ref{ringmaster:alg:asgd} are {\highlightcolor highlighted}.
    Alternatively, Alg.~\ref{ringmaster:alg:ringmasternew} is Alg.~\ref{ringmaster:alg:asgd} with a specific choice of adaptive stepsizes defined by \eqref{ringmaster:eq:adapative_step_size})
\end{algorithm}
\subsection{Description}
Let us compare \Cref{ringmaster:alg:ringmasternew} and \Cref{ringmaster:alg:asgd}:
the only difference, highlighted, lies in the fact that \Cref{ringmaster:alg:ringmasternew} disregards ``very old stochastic gradients'', which are gradients computed at points with significant delay.
Specifically, \Cref{ringmaster:alg:ringmasternew} receives $\nabla f\big(x^{k-\delta^k}; \xi^{k-\delta^k}_{i}\big)$, compares $\delta^k$ to our {\em delay threshold} hyperparameter $R$.
If $\delta^k \geq R$, then \Cref{ringmaster:alg:ringmasternew} completely ignores this very outdated stochastic gradient and requests the worker to compute a new stochastic gradient at the most relevant point $x^k$.

Notice that \ringmaster (\Cref{ringmaster:alg:ringmasternewstops}) is \Cref{ringmaster:alg:asgd} with the following adaptive stepsize rule:
\begin{equation}
\begin{aligned}
    \label{ringmaster:eq:adapative_step_size}
    &\gamma^k = 
    \begin{cases}
        \gamma, & \text{ if } \quad \bar{\delta}^k_i < R,\\
        0, &\text{ if }  \quad \bar{\delta}^k_i \geq R,
    \end{cases} 
 \\
    &\bar{\delta}^{k+1}_j = 
    \begin{cases}
     0, & \text{ if } \quad j = i, \\
        \bar{\delta}^{k}_j + 1,  & \text{ if }\quad  j\neq i \quad \text{and} \quad \bar{\delta}^k_i < R,\\
        \bar{\delta}^{k}_j, &  \text{ if } \quad j\neq i \quad \text{and} \quad \bar{\delta}^k_i \geq R,
    \end{cases}
\end{aligned}
\end{equation}
where $i$ is the  index of the worker whose stochastic gradient is applied at iteration $k$, and where we initialize the {\em virtual sequence of delays} $\{\bar{\delta}^{k}_j\}_{k,j}$ by setting  $\bar{\delta}^0_j = 0$ for all $j \in [n]$. 
\subsection{Delay Threshold} 
Returning to \Cref{ringmaster:alg:ringmasternew}, note that we have a parameter $R$ called the \emph{delay threshold}. 
When $R = 1$, the algorithm reduces to the classical \sgd method, i.e.,
$$
    x^{k+1} = x^k - \gamma \nabla f\(x^k; \xi^k_i\)~,
$$
since $\delta^k = 0$ for all $k \geq 0$.
In this case, the algorithm becomes highly conservative, ignoring all stochastic gradients computed at earlier points $x^{k-1}, \dots, x^0$. 
Conversely, if $R = \infty$, the method incorporates stochastic gradients with arbitrarily large delays, and becomes classical \naiveasyncsgd.
Intuitively, there should be a balance\textemdash a ``proper'' value of $R$ that would (i) prevent the method from being overly conservative, while (ii) ensuring stability by making sure that only informative stochastic gradients are used to update the model.
We formalize these intuitions by proposing an optimal $R$ in \Cref{ringmaster:sec:theory}.
Interestingly, the value of $R$ does {\em not} depend on the computation times.

\subsection{Why Do We Ignore the Old Stochastic Gradients?}
The primary reason is that doing so enables the development of the first optimal asynchronous \sgd method that achieves the lower bounds (see Sections~\ref{ringmaster:sec:theory} and \ref{ringmaster:sec:dyn}).
Ignoring old gradients allows us to establish tighter convergence guarantees.
Intuitively, old gradients not only fail to provide additional useful information about the function $f$, but they can also negatively impact the algorithm's performance.
Therefore, disregarding them in the optimization process is essential for achieving our goal of developing an optimal asynchronous \sgd method.

\subsection{Comparison to Rennala SGD}
Unlike \rennala, which combines a synchronous \minibatch update with an asynchronous minibatch collection strategy, \ringmaster is fully asynchronous, which, as we show in our experiments, provided further practical advantages:
\begin{itemize}
    \item \ringmaster updates the model  immediately  upon receiving a new and relevant (i.e., not too outdated) stochastic gradient. 
    This immediate update strategy is particularly advantageous for sparse models in practice, where different gradients may update disjoint parts of the model only, facilitating faster processing. 
    This concept aligns with the ideas of \citet{recht2011hogwild}, where lock-free asynchronous  updates were shown to effectively leverage sparsity for improved performance.
    \item \ringmaster ensures that stochastic gradients computed by fast workers are never ignored.
    In contrast, \rennala may discard stochastic gradients, even if they were recently initiated and would be computed quickly (see discussion in \Cref{ringmaster:sec:introduction}).
\end{itemize}
At the same time, \ringmaster adheres to the same principles that make \rennala optimal.
The core philosophy of \rennala is to prioritize fast workers by employing asynchronous batch collection, effectively ignoring slow workers.
This is achieved by carefully choosing the batch size $B$ in \Cref{ringmaster:alg:rennala}: large enough to allow fast workers to complete their calculations, but small enough to disregard slow workers and excessively delayed stochastic gradients.
Similarly, \ringmaster implements this concept using the delay threshold $R$.

In summary, while both \rennala and \ringmaster are optimal from a theoretical perspective, the complete absence of synchronization in the latter method intuitively makes it more appealing in practice, and allows it to collect further gains which are not captured in theory. As we shall see, this intuition is supported by our experiments.

\subsection{Comparison to Previous Asynchronous SGD Variants}
\citet{koloskova2022sharper} and \citet{mishchenko2022asynchronous} provided the previous state-of-the-art analysis of \naiveasyncsgd\ using delay-adaptive stepsizes.
However, as demonstrated in \Cref{ringmaster:sec:prel}, their analysis does not guarantee optimality.
This is due to at least two factors:
\begin{itemize}
    \item their approaches do {\em not} discard old stochastic gradients, instead attempting to utilize all gradients, even those with very large delays;
    \item although they select stepsizes $\gamma^k$ that decrease as the delays increase, this adjustment may not be sufficient to ensure optimal performance, and the choice of their stepsize may be suboptimal compared to \eqref{ringmaster:eq:adapative_step_size}.
\end{itemize}
\subsection{Stopping the Irrelevant Computations}
If stopping computations is feasible, we can further enhance \Cref{ringmaster:alg:ringmasternew} by introducing \cref{ringmaster:alg:ringmasternewstops}. 
Instead of waiting for workers to complete the calculation of outdated stochastic gradients with delays larger than $R$ which would not be used anyway, we propose to {\em stop/terminate} these computations immediately and reassign the workers to the most relevant point $x^k$. 
This adjustment provides workers with an opportunity to catch up, as they may become faster and perform computations more efficiently at the updated point.

\section{Theoretical Analysis}
\label{ringmaster:sec:theory}
We are ready to present the theoretical analysis of \ringmaster. 
We start with an \emph{iteration complexity} bound:
\begin{boxedtheorem}[Proof in \Cref{ringmaster:proof:ringmaster_iteration}]
    \label{ringmaster:thm:ringmaster_iteration}
    Under Assumptions \ref{ringmaster:ass:lipschitz_constant}, \ref{ringmaster:ass:lower_bound}, and \ref{ringmaster:ass:stochastic_variance_bounded}, let the stepsize in \ringmaster (\Cref{ringmaster:alg:ringmasternew} or \Cref{ringmaster:alg:ringmasternewstops}) be
    $$
        \gamma = \min \left\{ \frac{1}{2RL}, \frac{\varepsilon}{4L\sigma^2} \right\}.
    $$
    Then
    $$
       \frac{1}{K+1}\sum \limits_{k=0}^{K} \Exp{\sqnorm{\nabla f \(x^k\)}} \le \varepsilon~,
    $$
as long as
    \begin{align}
        \label{ringmaster:eq:HPMNVlgxoiRgUbRtj}
       K \geq \frac{8 R L \Delta}{\epsilon} + \frac{16 \sigma^2 L \Delta}{\epsilon^2}~,
    \end{align}
    where $R \in \{1, 2, \dots, \}$ is an arbitrary delay threshold.
\end{boxedtheorem}
The classical analysis of \naiveasyncsgd achieves the same convergence rate, with $R$ defined as $R \equiv \max_{k \in [K]} \delta^{k}$ \citep{stich2020error, arjevani2020tight}.
This outcome is expected, as setting $R = \max_{k \in [K]} \delta^{k}$ in \ringmaster makes it equivalent to classical \naiveminibatch, since no gradients are ignored.
Furthermore, the analyses by \citet{cohen2021asynchronous, koloskova2022sharper, mishchenko2022asynchronous} yield the same rate with $R = n$.
However, it is important to note that $R$ is a free parameter in \ringmaster that can be chosen arbitrarily.
While setting $R = \max_{k \in [K]} \delta^{k}$ or $R = n$ effectively recovers the earlier \emph{iteration complexities}, we show that there exists a {\em different} choice of $R$ leading to optimal time complexities (see Theorems~\ref{ringmaster:thm:optimal_ringmaster} and \ref{ringmaster:thm:optimal_ringmaster_dynamic}).

It is important to recognize that the \emph{iteration complexity} \eqref{ringmaster:eq:HPMNVlgxoiRgUbRtj} does {\em not} capture the actual  ``runtime'' performance of the algorithm.
To select an optimal value for $R$, we must shift our focus from \emph{iteration complexity} to \emph{time complexity}, which measures the algorithm's runtime.
We achieve the best practical and effective choice by optimizing the \emph{time complexity} over $R$.
In order to find the \emph{time complexity}, we need the following lemma.
\begin{boxedlemma}[Proof in Appendix~\ref{ringmaster:proof:time_R}]
    \label{ringmaster:lem:time_R}
    Let the workers' computation times satisfy the \emph{fixed computation model} \eqref{ringmaster:eq:worker-time}.
    Let $R$ be the delay threshold of \Cref{ringmaster:alg:ringmasternew} or \Cref{ringmaster:alg:ringmasternewstops}.
    The time required to complete any $R$ consecutive iterate updates of \Cref{ringmaster:alg:ringmasternew} or \Cref{ringmaster:alg:ringmasternewstops} is at most
    \begin{equation}\label{ringmaster:eq:time_R}
       t(R) 
       \eqdef 2 \min\limits_{m \in [n]} \left(\frac{1}{m} \sum\limits_{i=1}^m \frac{1}{\tau_i}\right)^{-1} \left( 1 + \frac{R}{m} \right).
    \end{equation}
\end{boxedlemma}
Combining \Cref{ringmaster:thm:ringmaster_iteration} and \Cref{ringmaster:lem:time_R}, we provide our main result.
\begin{boxedtheorem}[Optimality of \ringmaster]
    \label{ringmaster:thm:optimal_ringmaster}
    Let Assumptions \ref{ringmaster:ass:lipschitz_constant}, \ref{ringmaster:ass:lower_bound}, and \ref{ringmaster:ass:stochastic_variance_bounded} hold.
    Let the stepsize in \ringmaster (\Cref{ringmaster:alg:ringmasternew} or \Cref{ringmaster:alg:ringmasternewstops}) be
    $$
        \gamma = \min \left\{ \frac{1}{2RL}, \frac{\varepsilon}{4L\sigma^2} \right\}.
    $$
    Then, under the \emph{fixed computation model} \eqref{ringmaster:eq:worker-time}, \ringmaster achieves the optimal time complexity
    \begin{align}
        \label{ringmaster:eq:trfZhLNgGXylL}
        \cO\left(\min\limits_{m \in [n]} \left(\frac{1}{m} \sum\limits_{i=1}^{m} \frac{1}{\tau_{i}}\right)^{-1} \left(\frac{L \Delta}{\varepsilon} + \frac{\sigma^2 L \Delta}{m \varepsilon^2}\right) \right)
    \end{align}
    with the delay threshold
    \begin{align}
        \label{ringmaster:eq:TZkAIqt}
        R = \max\left\{1,\left\lceil \frac{\sigma^2}{\varepsilon}\right\rceil\right\}.
    \end{align}
\end{boxedtheorem}
Note that the value of $R$ does not in any way depend on the computation times $\{\tau_1,\dots,\tau_n\}$.
\begin{proof}
    From \Cref{ringmaster:thm:ringmaster_iteration}, the iteration complexity of \ringmaster is
    \begin{align}
        \label{ringmaster:eq:nyQNjhjnpDOOjflbnS}
        K = \left\lceil\frac{8 R L \Delta}{\epsilon} + \frac{16 \sigma^2 L \Delta}{\epsilon^2}\right\rceil.
    \end{align}
    Using \Cref{ringmaster:lem:time_R}, we know that \ringmaster requires at most $t(R)$ seconds to finish any $R$ consecutive updates of the iterates.
    Therefore, the total time is at most
    \begin{align*}
        t(R) \times \left\lceil\frac{K}{R}\right\rceil.
    \end{align*}
    Without loss of generality, we assume\footnote{Otherwise, using $L$--smoothness, $\|\nabla f(x^0)\|^2 \leq 2 L \Delta \leq \varepsilon$, and the initial point is an $\varepsilon$--stationary point.
    See \Cref{ringmaster:sec:l_init}.}
    that $L \Delta > \nicefrac{\varepsilon}{2}$.
    Therefore, using \eqref{ringmaster:eq:nyQNjhjnpDOOjflbnS}, we get
    \begin{align}
        \label{ringmaster:eq:dfgyMy}
        t(R) \times \left\lceil\frac{K}{R}\right\rceil = \cO\left(t(R) \left(\frac{L \Delta}{\epsilon} + \frac{\sigma^2 L \Delta}{R \epsilon^2}\right)\right).
    \end{align}
    It is left to substitute our choice \eqref{ringmaster:eq:TZkAIqt} into \eqref{ringmaster:eq:dfgyMy} to get \eqref{ringmaster:eq:trfZhLNgGXylL}.
    We take \eqref{ringmaster:eq:TZkAIqt}, noticing that this choice of $R$ minimizes \eqref{ringmaster:eq:dfgyMy} up to a universal constant.
\end{proof}

To the best of our knowledge, this is the first result in the literature to establish the optimality of a fully asynchronous variant of \sgd:
the time complexity \eqref{ringmaster:eq:trfZhLNgGXylL} is optimal and aligns with the lower bound established by \citet{tyurin2024optimal}.

The derived time complexity \eqref{ringmaster:eq:trfZhLNgGXylL} has many nice and desirable properties.  
First, it is robust to slow workers: if $\tau_n \to \infty$, the expression equals 
$$
    \min\limits_{m \in [{\highlightcolor n - 1}]} \left(\frac{1}{m} \sum \limits_{i=1}^{m} \frac{1}{\tau_{i}}\right)^{-1} \left(\frac{L \Delta}{\varepsilon} + \frac{\sigma^2 L \Delta}{m \varepsilon^2}\right),
$$
effectively disregarding the slowest worker.
Next, assume that $m_*$ is the smallest index that minimizes \eqref{ringmaster:eq:trfZhLNgGXylL}.
In this case, \eqref{ringmaster:eq:trfZhLNgGXylL} simplifies further to 
$$
    \left(\frac{1}{m_*} \sum \limits_{i=1}^{m_*} \frac{1}{\tau_{i}}\right)^{-1} \left(\frac{L \Delta}{\varepsilon} + \frac{\sigma^2 L \Delta}{m_* \varepsilon^2}\right).
$$
This shows that the method operates effectively as if only the fastest $m_*$ workers participate in the optimization process, resembling the idea from \Cref{ringmaster:alg:ringmaster_short}.
What is important, however, the method determines $m_*$ adaptively and automatically.
\subsection{The Choice of Threshold}
Another notable property of the method and the time complexity result in \Cref{ringmaster:thm:optimal_ringmaster} is that the threshold $R$ is independent of the individual computation times $\tau_i$.
As a result, the method can be applied across heterogeneous distributed systems (where workers have varying speeds) while still achieving the optimal time complexity given in \eqref{ringmaster:eq:dfgyMy}, up to a constant factor.

If a tighter, constant-level expression for the optimal threshold is desired, $R$ can be computed explicitly.
In that case, it will depend on the values of $\tau_i$.
The optimal $R$ solves
\begin{equation*}
    \argmin_{R \geq 1} \ t(R) \left(1 + \frac{\sigma^2}{R \varepsilon} \right) ,
\end{equation*}
which follows from \eqref{ringmaster:eq:dfgyMy} after discarding constants independent of $R$. Here, $t(R)$ denotes the total time required for $R$ consecutive iterations, and is upper bounded by \eqref{ringmaster:eq:time_R}.

Substituting this bound yields the following expression for the optimal $R$:
\begin{equation*}
    R = \max\left\{ \sigma \sqrt{\frac{m^*}{\varepsilon}}, 1 \right\},
\end{equation*}
where
\begin{equation*}
    m^* = \argmin\limits_{m \in [n]} \left( \frac{1}{m} \sum\limits_{i=1}^m \frac{1}{\tau_i} \right)^{-1} \left(1 + 2 \sqrt{\frac{\sigma^2}{m \varepsilon}} + \frac{\sigma^2}{m \varepsilon} \right) .
\end{equation*}
As the expression shows, the optimal threshold $R$ depends on $m^*$, which in turn is determined by the $\tau_i$ values.
\subsection{Proof Techniques}
Several mathematical challenges had to be addressed to achieve the final result.
Lemmas~\ref{ringmaster:lem:time_R} and \ref{ringmaster:lem:time_R_dyn} are novel, as estimating the time complexity of the \emph{asynchronous} \ringmaster method requires distinct approaches compared to the \emph{semi-synchronous} \rennala method because \ringmaster is more chaotic and less predictable.
Compared to the works of \citet{koloskova2022sharper,mishchenko2022asynchronous}, the proof of \Cref{ringmaster:thm:ringmaster_iteration} is tighter and more refined, as we more carefully analyze the sum
$$
    \sum_{k=0}^K \E{\sqnorm{x^k - x^{k-\delta^k}}}
$$
in \Cref{ringmaster:lemma:residual}.
We believe that the simplicity of \ringmaster, combined with the new lemmas, the refined analysis, and the novel choice of $R$ in \Cref{ringmaster:thm:optimal_ringmaster}, represents a set of non-trivial advancements that enable us to achieve the optimal time complexity.

\section{Optimality Under Arbitrary Computation Dynamics}\label{ringmaster:sec:dyn}
In the previous sections, we presented the motivation, improvements, comparisons, and theoretical results within the framework of the \emph{fixed computation model} \eqref{ringmaster:eq:worker-time}.
We now prove \ringmaster is optimal under virtually any computation behavior of the workers.
One way to formalize these behaviors is to use the \emph{universal computation model} \citep{tyurin2024tighttimecomplexitiesparallel}.

For each worker $i \in [n]$, we associate a \emph{computation power} function 
$$
    p_i \,:\, \R_{+} \rightarrow \R_{+} ~.
$$
The number of stochastic gradients computed by worker $i$ in $[T^0, T^1]$ is given by the integral of its computation power $p_i$, followed by a floor operation:
\begin{align}
  \label{ringmaster:eq:rSIiSfVcmivSKsfzoSA}
  \textnormal{``\# of stoch. grad. in $[T^0, T^1]$''} = \flr{\int_{T^0}^{T^1} p_i(\tau) d \tau}.
\end{align}
The computational power $p_i$ characterizes the behavior of workers, accounting for potential disconnections due to hardware or network delays, variations in processing capacity over time, and fluctuations or trends in computation speeds.
The integral of $p_i$ over a given interval represents the \emph{computation work} performed.
If $p_i$ is small within $[T^0, T^1]$, the integral is correspondingly small, indicating that worker $i$ performs less computation.
Conversely, if $p_i$ is large over $[T^0, T^1]$, the worker is capable of performing more computation.

For our analysis, we only assume  that $p_i$ is non-negative and continuous almost everywhere\footnote{Thus, it can be discontinuous and ``jump'' on a countable set.}.
We use this assumption non-explicitly when applying the Riemann integral.
The computational power $p_i$ can even vary randomly, and all the results discussed hold conditional on the randomness of $\{p_i\}$.

Let us examine some examples.
If worker $i$ remains inactive for the first $t$ seconds and then becomes active, this corresponds to $p_i(\tau) = 0$ for all $\tau \leq t$ and $p_i(\tau) > 0$ for all $\tau > t$.
Furthermore, we allow $p_i$ to exhibit periodic or even chaotic behavior\footnote{
For example, $p_i(t)$ might behave discontinuously as follows: $p_i(t) = 0.5t + \sin(10t)$ for $t \leq 10,$ $p_i(t) = 0$ for $10 < t \leq 20,$ and $p_i(t) = \max\{80 - 0.5 t, 0\}$.}.
The \emph{universal computation model} reduces to the \emph{fixed computation model} when $p_i(t) = \nicefrac{1}{\tau_i}$ for all $t \geq 0$ and $i \in [n]$.
In this case, the number of stochastic gradients computed in the interval $[T^0, T^1]$, given by \eqref{ringmaster:eq:rSIiSfVcmivSKsfzoSA}, becomes $\flr{\nicefrac{(T^1 - T^0)}{\tau_i}}$, indicating that worker $i$ computes one stochastic gradient after $T^0 + \tau_i$ seconds, two stochastic gradients after $T^0 + 2 \tau_i$ seconds, and so forth.

Note that \Cref{ringmaster:thm:ringmaster_iteration} is valid under any computation model.
Next, we introduce an alternative to \Cref{ringmaster:lem:time_R} for the \emph{universal computation model}.

\begin{boxedlemma}[Proof in Appendix~\ref{ringmaster:sec:time_R_dyn}]
    \label{ringmaster:lem:time_R_dyn}
    Let the workers' computation times satisfy the \emph{universal computation model}.
    Let $R$ be the delay threshold of \Cref{ringmaster:alg:ringmasternew} or \Cref{ringmaster:alg:ringmasternewstops}.
    Assume that some iteration starts at time $T^0$.
    Starting from this iteration, the $R$ consecutive iterate updates of \Cref{ringmaster:alg:ringmasternew} or \Cref{ringmaster:alg:ringmasternewstops} will be performed before the time
    \begin{equation*} 
           T(R,T^0) \eqdef \min \left\{T \geq 0 : \sum \limits_{i=1}^{n} \flr{\frac{1}{4} \int_{T^0}^{T} p_i(\tau) d \tau} \geq R\right\}.
    \end{equation*}
\end{boxedlemma}
This result extends \Cref{ringmaster:lem:time_R}. 
Specifically, if $p_i(t) = \nicefrac{1}{\tau_i}$ for all $t \geq 0$ and $i \in [n]$, it can be shown that $T(R,T^0) - T^0 = \Theta\left(t(R)\right)$ for all $T^0 \geq 0$.

Combining \Cref{ringmaster:thm:ringmaster_iteration} and \Cref{ringmaster:lem:time_R_dyn}, we are ready to present our main result.

\begin{boxedtheorem}[Optimality of \ringmaster; Proof in Appendix~\ref{ringmaster:sec:proof_dyn}]
    \label{ringmaster:thm:optimal_ringmaster_dynamic}
    Let Assumptions \ref{ringmaster:ass:lipschitz_constant}, \ref{ringmaster:ass:lower_bound}, and \ref{ringmaster:ass:stochastic_variance_bounded} hold.
    Let the stepsize in \ringmaster (\Cref{ringmaster:alg:ringmasternew} or \Cref{ringmaster:alg:ringmasternewstops}) be
    $$
        \gamma = \min \left\{ \frac{1}{2RL}, \frac{\varepsilon}{4L\sigma^2} \right\},
    $$ 
    and delay threshold 
    $$
        R = \max\left\{1,\left\lceil \frac{\sigma^2}{\varepsilon} \right\rceil\right\}.
    $$
    Then, under the \emph{universal computation model}, \ringmaster finds an $\varepsilon$--stationary point after at most $T^{\bar{K}}$ seconds, where 
    $$
    \bar{K} \eqdef \left\lceil\frac{48 L \Delta}{\epsilon}\right\rceil
    $$
    and $T^{\bar{K}}$ is the $\bar{K}$-th element of the following recursively defined sequence:
    \begin{align*}
        T^{k} \eqdef \min \left\{T \geq 0 : \sum \limits_{i=1}^{n} \flr{\frac{1}{4} \int_{T^{k - 1}}^{T} p_i(\tau) d \tau} \geq R\right\}
    \end{align*}
    for all $k \geq 1$ and $T^{0} = 0$.
\end{boxedtheorem}

Admittedly, the obtained theorem is less explicit than \Cref{ringmaster:thm:optimal_ringmaster}.
However, this is the price to pay for the generality of the result in terms of the computation time dynamics it allows. 
Determining the time complexity $T^{K}$ requires computing $T^{1}, T^{2},$ and so on, sequentially.
However, in certain scenarios, $T^{K}$ can be derived explicitly.
For example, if $p_i(t) = \nicefrac{1}{\tau_i}$ for all $t \geq 0$ and $i \in [n]$, then $T^{K}$ is given by \eqref{ringmaster:eq:trfZhLNgGXylL}.
Furthermore, $T^{1}$ represents the number of seconds required to compute the first $R$ stochastic gradients, $T^{2}$ represents the time needed to compute the first $2 \times R$ stochastic gradients, and so on.
Naturally, to compute $T^{2},$ one must first determine $T^{1},$ which is reasonable given the sequential nature of the process.
 
The obtained result is optimal, aligns with the lower bound established by \citet{tyurin2024tighttimecomplexitiesparallel}, and cannot be improved by any asynchronous parallel method.
\section{Experiments} 
\label{ringmaster:sec:experiments}
Asynchronous \sgd has consistently demonstrated its effectiveness and practicality, achieving strong performance in various applications \citep{recht2011hogwild,lian2018asynchronous,mishchenko2022asynchronous}, along with numerous other studies supporting its utility.
Our main goal of this paper is to refine the method further and establish that it is not only practical but also \emph{theoretically optimal}.

We conduct a toy experiment with our proposed \ringmaster method and compare it against two baselines: the previous \naiveasyncsgd (referred to as \algname{Delay-Adaptive ASGD}, following~\citet{mishchenko2022asynchronous}) and \rennala. 
The optimization task is based on a convex quadratic function $f : \R^d \to \R$ defined as
$$
    f(x) = \frac{1}{2} x^\top A x - b^\top x~, \qquad \forall x \in \R^d ~,
$$
where
\begin{align*}
    A = \frac{1}{4}
    \begin{bmatrix}
    2 & -1 &  & 0 \\
    -1 & \ddots & \ddots &  \\
    & \ddots & \ddots & -1 \\
    0 & & -1 & 2 \\
    \end{bmatrix}
    \in \mathbb{R}^{d \times d}~,
    \quad
    & b = \frac{1}{4}
    \begin{bmatrix}
    -1 \\
    0 \\
    \vdots \\
    0 \\
    \end{bmatrix}
    \in \mathbb{R}^d~.
\end{align*}
We set $d = 1729$ and $n = 6174$.
Each of the $n$ workers computes unbiased stochastic gradients of the form  
$$
    \nabla f(x, \xi) = \nabla f(x) + \xi~,
$$
where $\xi \sim \mathcal{N}(0, 0.01^2)$.

\begin{figure}[ht]
    \begin{center}
    \centerline{\includegraphics[width=0.75\textwidth]{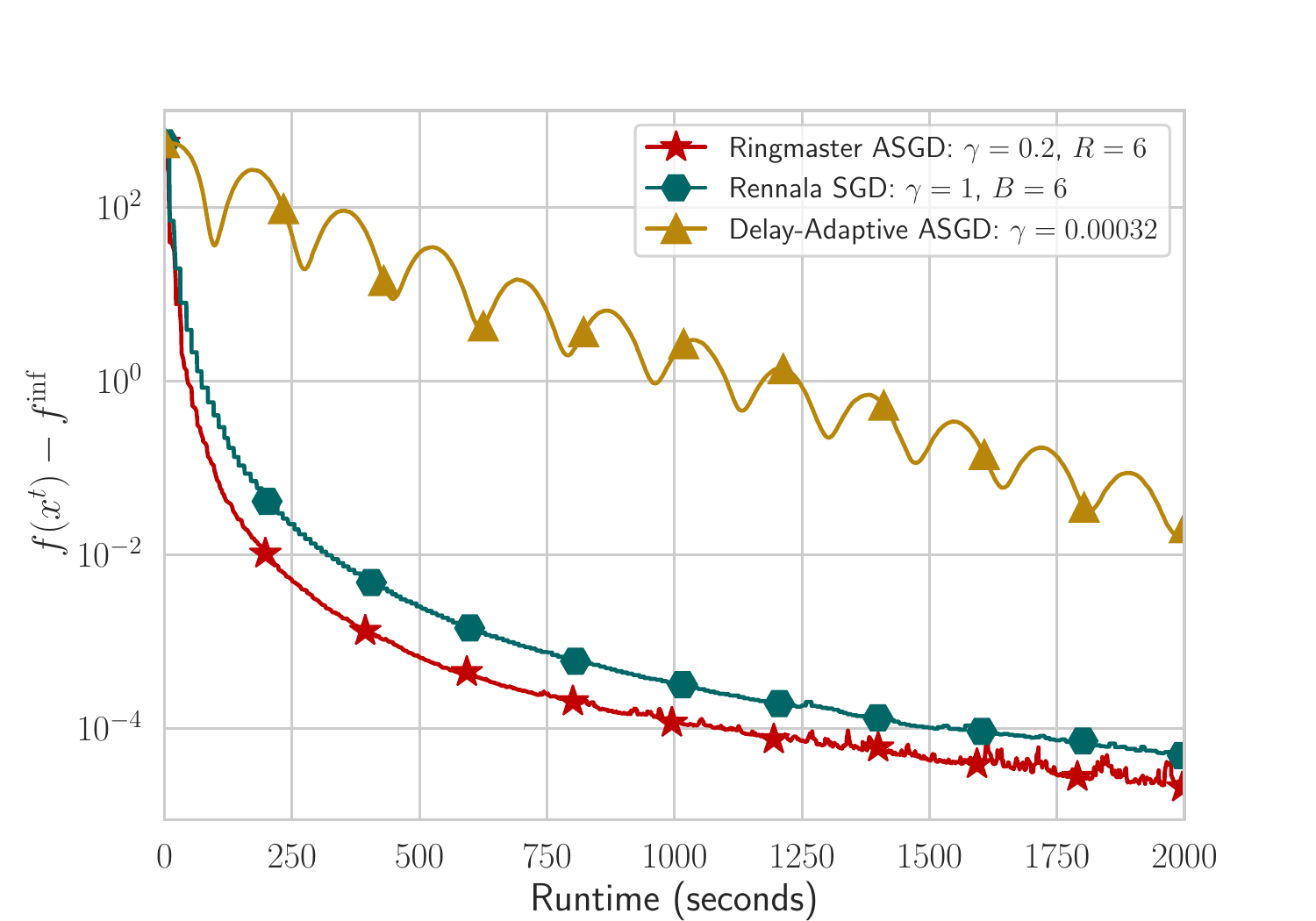}}
    \caption{Experiment with $n = 6174$ and $d = 1729$ showing the convergence of \ringmaster, \algname{Delay-Adaptive ASGD}, and \rennala.}
    \label{ringmaster:fig:1}
    \end{center}
\end{figure}

The experiments were implemented in Python.
The distributed environment was emulated on machines with Intel(R) Xeon(R) Gold 6248 CPU @ 2.50GHz.
The computation times for each worker are simulated as $\tau_i = i + \abs{\eta_i}$ for all $i \in [n]$, where $\eta_i \sim \cN\(0,i\)$.
We tuned the stepsize from the set $\left\{5^p : p \in [-5, 5]\right\}$.
Both the batch size for \rennala and the delay threshold for \ringmaster were tuned from the set $\left\{\lceil \nicefrac{n}{4^p} \rceil : p \in \Z_+\right\}$.
The experimental results are shown in \Cref{ringmaster:fig:1}.

The obtained result confirms that \ringmaster is indeed faster than \algname{Delay-Adaptive ASGD} and \rennala in the considered setting.
One can see the numerical experiments support that our theoretical results, and we significantly improve the convergence rate of the previous version of \naiveasyncsgd (\algname{Delay-Adaptive ASGD}).

\subsection{Neural Network Experiment}
\begin{figure}[ht]
    \begin{center}
    \centerline{\includegraphics[width=0.75\textwidth]{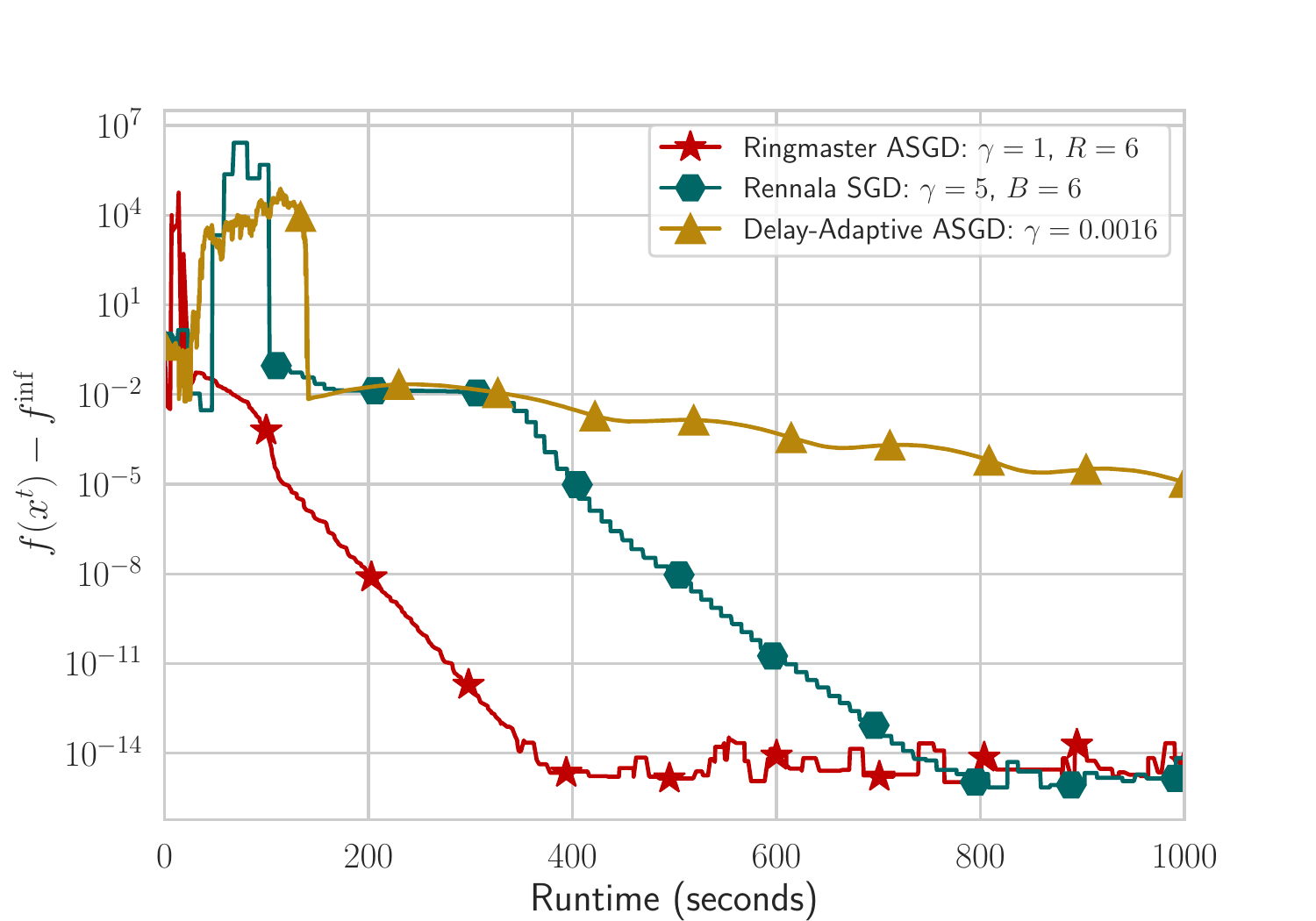}}
    \caption{
        We run an experiment on a small 2-layer neural network with ReLU activation on the MNIST dataset, showing that our method, \ringmaster, is more robust and outperforms \algname{Delay-Adaptive ASGD} and \rennala.
    }
    \label{ringmaster:fig:mnist}
    \end{center}
\end{figure}

To show that our method also works well for neural networks, we trained a small 20-layer neural network with ReLU activation on the MNIST dataset \citep{lecun1998gradient}.
We used the same number of workers as in the previous experiment ($n = 6174$) and kept the same time distributions.
The results are shown in \Cref{ringmaster:fig:mnist}.

\section{Additional Related Work}

In federated learning, communication is often more time-consuming than local computation, which motivates the use of local training strategies \citep{mcmahan2016federated,mishchenko2022proxskip,malinovsky2022variance,maranjyan2025gradskip}.
Similarly, in our setting, instead of repeatedly computing gradients at the same model point, workers can perform several local gradient descent steps and send the updated model to the server when requested.
Asynchronous variants of such schemes, including those with explicit time-complexity analyses, have been explored in prior works \citep{tyurin2025birch,fradin2025local}.
\section{Conclusion and Future Work}
\label{ringmaster:sec:sonclusion}

In this work, we developed the first asynchronous \sgd method, named \ringmaster, that achieves optimal time complexity.
By employing a carefully designed algorithmic approach, which can be interpreted as asynchronous \sgd with adaptive step sizes \eqref{ringmaster:eq:adapative_step_size}, we successfully reach this goal.
By selecting an appropriate delay threshold $R$ in \Cref{ringmaster:alg:ringmasternew}, the method attains the theoretical lower bounds established by \citet{tyurin2024optimal,tyurin2024tighttimecomplexitiesparallel}.

Future work can explore heterogeneous scenarios, where the data on each device comes from different distributions, which are particularly relevant for federated learning (FL) \citep{konevcny2016federated,mcmahan2016federated,kairouz2021advances}. 
In FL, communication costs also become crucial, making it worthwhile to consider similar extensions \citep{alistarh2017qsgd,tyurin2024shadowheart,tyurin2024optimalgraph}.
It would be also interesting to design a fully asynchronous analog to the method from \citep{tyurin2024freya}.
Additionally, as discussed by \citet{maranjyan2025mindflayer}, computation and communication times can be treated as random variables.
It would be valuable to investigate these cases further and derive closed-form expressions for various distributions.

\chapter{Ringleader ASGD: The First Asynchronous SGD \\ with Optimal Time Complexity under Data Heterogeneity}
\label{chapter:Ringleader}
\thispagestyle{empty}
This chapter is based on the work of \citet{maranjyan2025ringleader}:
\begin{quote}
    \bibentry{maranjyan2025ringleader}
\end{quote}
\section{Introduction}
Modern machine learning increasingly depends on large-scale distributed training across clusters with hundreds or even thousands of GPUs \citep{shoeybi2019megatron, GPT3, narayanan2021efficient}.
However, classical synchronous training methods struggle to scale in these settings, as device failures, network instabilities, and synchronization overheads introduce significant inefficiencies \citep{chen2016revisiting, grattafiori2024llama}.
These issues become even more pronounced in environments with heterogeneous computational power, such as Federated Learning (FL), where devices range from high-end data center GPUs to resource-constrained edge hardware \citep{konevcny2016federated, mcmahan2016federated, li2020federated, kairouz2021advances}.
Because synchronous methods are bottlenecked by the slowest participants, faster devices remain idle, leading to severe underutilization of computational resources when stragglers\textemdash nodes slowed down by computation or communication\textemdash lag significantly behind.

One way to reduce synchronization bottlenecks is to equip data centers with homogeneous GPUs.
However, this approach is prohibitively expensive and difficult to scale: upgrading to faster GPUs would require replacing all devices simultaneously, since heterogeneous hardware cannot be combined efficiently.
Even then, homogeneity does not eliminate synchronization issues, as hardware failures and device dropouts during training still cause stragglers and idle time.
Moreover, this solution applies only to controlled data center environments and is infeasible in FL, where edge devices are outside the server's control.

A more promising approach is to shift from hardware solutions to algorithmic ones by adopting asynchronous optimization methods.
These methods remove the need for synchronization, allowing fast workers to contribute updates without waiting for slower ones \citep{tsitsiklis1986distributed, recht2011hogwild, agarwal2011distributed, dean2012large, li2014communication}.
Despite their appeal, asynchronous methods are more difficult to analyze.
In particular, a meaningful analysis would require studying \textit{time to convergence}, rather than iteration complexity only.
While iteration complexity is the traditional metric in optimization, it does not necessarily reflect real-world training speed in parallel settings: a method that performs more iterations may finish faster in wall-clock time if those iterations can be computed without waiting for slow workers.
This distinction raises a fundamental question: \textit{among all parallel methods, which ones are provably fastest in theory?}
To make this question precise, we restrict our attention to smooth nonconvex problems and to stochastic first-order methods, encompassing algorithms with or without synchronization.
This will be the only setting considered in this paper.

\begin{table*}
    \caption{
        Comparison of time complexities for parallel first-order methods under the fixed computation time model, where each worker $i$ takes a fixed time $\tau_i$ to compute a stochastic gradient, with the times ordered so that $\tau_n$ is the largest \eqref{eq:fixed_time}.
        We denote by $\tauavg \eqdef \tfrac{1}{n}\sum_{i=1}^n \tau_i$ the average computation time across all workers.
        The table shows how the time complexity of different algorithms depends on key problem parameters: the initial function suboptimality $\Delta \eqdef f(x^0) - f^*$ (\Cref{ass:lower_bound}), the target stationarity $\varepsilon$, the variance bound of the stochastic gradients $\sigma^2$ (\Cref{ass:stochastic_variance_bounded}), and smoothness constants.
        Specifically, $L_f$ is the smoothness constant of $f$ (\Cref{eq:lipschitz_gradient}); $L_{\max} \eqdef \max_{i\in[n]} L_{f_i}$ with $L_{f_i}$ the smoothness constant of $f_i$; and $L$ is a constant associated with our new smoothness-type assumption (\Cref{ass:lipschitz_constant}).
        They satisfy $L_f \le L \le L_{\max}$ (\Cref{lemma:smoothness_relation}).
        All stated time complexities hide universal constant factors.\\
        Each column indicates whether a method satisfies the following desirable properties:
        \textbf{Optimal:} achieves the theoretical lower bound derived by \citet{tyurin2024optimal} for parallel first-order stochastic methods in heterogeneous data setting.
        \textbf{No sync.:} does not require synchronization and is therefore \textit{asynchronous}.
        \textbf{No idle workers:} all workers remain busy without waiting, so computational resources are fully utilized.
        \textbf{No discarded work:} no computation is wasted, and no worker is stopped mid-computation.
        Our new method, \ringleader, is the first asynchronous method to achieve optimal time complexity, while also ensuring full resource utilization (no idle workers) and no discarded computations / work.
    }
    \label{table:fixedtime}
    \centering
    \begin{threeparttable}
      \resizebox{\textwidth}{!}{
        \begin{tabular}{cccccc}
            \toprule
                \bf Algorithm
                & \makecell{\bf Time \\ \bf Complexity}
                & \bf Optimal
                & \makecell{\bf No \\ \bf Sync.}
                & \makecell{\bf No \\ \bf Idle \\ \bf Workers}
                & \makecell{\bf No \\ \bf Discarded \\ \bf Work} \\
            \midrule
                \makecell{\naiveminibatch \\ (\Cref{sec:minibatch})}
                & $\frac{L_f\Delta}{\varepsilon}\(\tau_n + \tau_n \frac{\sigma^2}{n\varepsilon}\)$
                & \crossmarkred
                & \crossmarkred
                & \crossmarkred
                & \checkmarkgreen \\
            \midrule
                \makecell{\iasgd \\ \citep{wang2025incremental} \\ (\Cref{sec:iasgd_time_complexity})}
                & $\frac{L_{\max}\Delta}{\varepsilon}\(\tau_n + \tau_n \frac{\sigma^2}{n\varepsilon}\)$ \textsuperscript{{\linkcolor (\dag)}}
                & \crossmarkred
                & \checkmarkgreen
                & \checkmarkgreen
                & \checkmarkgreen \\
            \midrule
                \makecell{\malenia \\ \citep{tyurin2024optimal}}
                & $\frac{L_f\Delta}{\varepsilon}\(\tau_n + \tauavg \frac{\sigma^2}{n\varepsilon}\)$
                & \checkmarkgreen
                & \crossmarkred
                & \checkmarkgreen
                & \crossmarkred \\
            \midrule
                \makecell{\ringleader \textbf{(new)} \\ (\Cref{algo:Ringleader}; \\ \Cref{thm:time_complexity})}
                & $\frac{L\Delta}{\varepsilon}\(\tau_n + \tauavg \frac{\sigma^2}{n\varepsilon}\)$
                & \checkmarkgreen \textsuperscript{{\linkcolor (\ddag)}}
                & \checkmarkgreen
                & \checkmarkgreen
                & \checkmarkgreen\\
            \bottomrule
        \end{tabular}%
    }

    \resizebox{\textwidth}{!}{
        \begin{minipage}{\textwidth}
            \begin{tablenotes}[para,flushleft]
                \footnotesize
                \item[{\linkcolor (\dag)}]
                    The analysis presented by \citet{wang2025incremental} is carried out under the assumption that each $f_i$ is smooth with the same smoothness constant, which is equivalent to requiring that each $f_i$ is $L_{f_i}$--smooth and then using the upper bound $L_{\max}$ for all $L_{f_i}$.
                    However, this strict assumption is not necessary: they could instead use our relaxed smoothness-type assumption (\Cref{ass:lipschitz_constant}), in which case the constant improves to $L$, and their analysis remains unchanged.
                \\
                \item[{\linkcolor (\ddag)}]
                    The time complexities of \ringleader and \malenia differ in the smoothness constant only.
                    Since \malenia is optimal, \ringleader is also optimal whenever $L$ exceeds $L_f$ by at most a universal constant factor, that is, $L = \cO(L_f)$.
            \end{tablenotes}
        \end{minipage}
    }

    \end{threeparttable}
\end{table*} 

Recently, \citet{tyurin2024optimal} studied this very regime, where they derived lower bounds.
They then proposed two algorithms: \rennala, designed for the \textit{homogeneous data setting}, where all workers draw samples from the same distribution, and \malenia, for the \textit{heterogeneous data setting}, where data distributions differ across workers.
They showed that both methods are optimal\textemdash achieving the lower bounds\textemdash and, perhaps surprisingly, both are synchronous (they periodically synchronize the workers).
The key idea in both is to fully utilize the available computational resources by keeping workers continuously busy: each worker computes independently, and synchronization occurs only after a sufficient number of gradient computations have been accumulated.

At first, the result of \citet{tyurin2024optimal} suggested a rather pessimistic outlook for asynchronous methods: despite their practical appeal, they showed that existing asynchronous methods are not optimal and that the method achieving the lower bound is \textit{synchronous}.
This created the view that optimality is inherently tied to synchronization.
However, this view was overturned by \citet{maranjyan2025ringmaster}, who, in the \textit{homogeneous data setting}, introduced \ringmaster---the first asynchronous \sgd method to achieve the same optimal time complexity as the synchronous \rennala.
Although both methods share the same theoretical guarantees, \ringmaster can be faster than \rennala in practice, since it avoids synchronization and benefits from more frequent updates.

Nevertheless, the work of \citet{maranjyan2025ringmaster} established optimality in the \textit{homogeneous data setting} only.
The question of whether some variant of a parallel method that does not rely on synchronization (i.e., is asynchronous) can also be optimal in the more general \textit{heterogeneous data setting} remained open.
In this work, we close this gap and answer the question affirmatively.

The heterogeneous data setting is both important and practically relevant.  
In FL, for instance, such heterogeneity arises naturally as participants hold distinct datasets \citep{zhao2018federated, li2020federated, tan2022towards}.
Yet this setting is significantly more challenging than the homogeneous one.
The standard philosophy of asynchronous \sgd---updating the model after every gradient computation\textemdash can be harmful here: fast workers contribute updates more frequently, causing the optimization process to become biased toward their local data.
To mitigate this, most existing asynchronous methods address this issue by assuming similarity across client distributions \citep{mishchenko2022asynchronous, koloskova2022sharper, nguyen2022federated, islamov2024asgrad}.
While this assumption simplifies the analysis, it is often unrealistic in practice, where clients may involve fundamentally different populations (e.g., hospitals with distinct demographics, mobile users in different countries, or financial institutions under varied regulations).

Recent work by \citet{wang2025incremental} took an important step toward removing these restrictive assumptions by proposing Incremental Aggregated Asynchronous \algname{SGD} (\iasgd), a method that provably converges without similarity assumptions.
However, their method achieves the same time complexity as standard \naiveminibatch (see the first two rows of \Cref{table:fixedtime})\textemdash the simplest synchronous \sgd baseline, which waits to collect one gradient from each worker before every update\textemdash thus failing to provide the computational advantages that motivate asynchronous approaches in the first place.

To the best of our knowledge, the only method proven to be optimal in the \textit{heterogeneous data setting} is the synchronous algorithm \malenia of \citet{tyurin2024optimal}, which, notably, does not rely on similarity assumptions.
However, synchronization is a major bottleneck in practice: although synchronous and asynchronous methods can share the same theoretical complexity, asynchronous methods are often faster in practice because they avoid costly synchronization and benefit from more frequent updates, as demonstrated in the homogeneous case by \citet{maranjyan2025ringmaster}.

This raises a fundamental question: \textit{Is it possible to design an asynchronous method that requires no similarity assumptions while still achieving optimal time complexity?}
In this paper, we answer this question affirmatively by introducing \ringleader, the first asynchronous \sgd method that achieves optimal time complexity\footnote{
    Throughout the paper, we refer to our method as optimal.
    Formally, this holds whenever the constant $L$\textemdash associated with our new smoothness-type assumption (\Cref{ass:lipschitz_constant})\textemdash is at most a constant factor larger than the smoothness constant $L_f$ used in the derived lower bounds \citep{tyurin2024optimal}.
    See \Cref{table:fixedtime} for details.
    }
in the \textit{heterogeneous data setting}.
Importantly, \ringleader attains this without relying on restrictive similarity assumptions.
\subsection{Contributions}
Our main contributions are the following:
\begin{itemize}
    \item \textbf{Optimal asynchronous SGD under data heterogeneity.}
    We prove that \ringleader (\Cref{algo:Ringleader}) is, to the best of our knowledge, the first asynchronous method in the heterogeneous setting under the fixed computation model \eqref{eq:fixed_time} that matches the lower bounds for parallel methods established by \citet{tyurin2024optimal} (\Cref{table:fixedtime}).
    Crucially, it does not rely on any similarity assumptions across clients' datasets.
    \item \textbf{Additional useful properties.}  
    Beyond achieving optimal time complexity, our method \ringleader satisfies two additional desirable properties: (i) all workers remain continuously active (\textit{no idle workers}), and (ii) every computed gradient is incorporated into the update (\textit{no discarded work}).
    These properties are crucial in practice, as they ensure maximum resource utilization: all workers contribute at all times, and their computations are never wasted.
    \Cref{table:fixedtime} compares \ringleader against benchmark algorithms with respect to these properties.
    \item \textbf{Parameter-free design.}
    In contrast to the optimal synchronous method \malenia \citep{tyurin2024optimal}, which requires prior knowledge of the gradient variance bound and target accuracy, our method operates in the fixed computation time model without such parameters (except for the stepsize, needed only to match the optimal theoretical rate).
    This makes it far more practical for real-world deployments, where these quantities are typically unknown or difficult to estimate.
    The same parameter-free improvement can also be extended to \malenia, as we discuss in \Cref{sec:malenia_param_free}.
    \item \textbf{Universal computation model.}
    In \Cref{sec:arbitrary_time}, we extend our analysis beyond the fixed computation time model to the general setting of arbitrarily varying computation times, accommodating virtually any computational behavior, including stochastic or adversarial patterns, while retaining optimality time complexity.
    \item \textbf{Empirical validation.}
    In \Cref{sec:experiments}, we evaluate \ringleader against benchmark methods on illustrative toy problems.
    The results validate our theoretical findings and demonstrate clear practical advantages over the baselines.
\end{itemize}
\section{Related Work}
Research on asynchronous stochastic gradient methods dates back to the seminal work of \citet{tsitsiklis1986distributed}, and gained renewed momentum with the introduction of the \algname{HOGWILD!} algorithm \citep{recht2011hogwild}.
\algname{HOGWILD!} is fundamentally an asynchronous coordinate descent method: updates are performed lock-free with inconsistent reads and writes, and its convergence guarantees rely on sparsity assumptions that are rarely realistic in modern large-scale machine learning.
Subsequent refinements of this paradigm include \citep{reddi2015variance, zhao2016fast, pan2016cyclades, mania2017perturbed, leblond2018improvedAsynchronous, nguyen2018sgd, zhou2018simple}, but these works remain tied to the coordinate descent setting with inconsistent memory accesses, and thus differ substantially from our focus.

Closer to our setting are works where updates are based on gradients that are applied consistently.
Early contributions, typically under the homogeneous-data assumption (all workers sample from the same distribution), include the work of \citet{agarwal2011distributed}, who studied convex objectives, as well as later extensions to the non-convex case such as the work of \citet{lian2015asynchronous} and \citet{dutta2018slow}, the latter analyzing exponentially distributed computation times.
Other relevant results in this line include \citep{feyzmahdavian2016asynchronous, zheng2017asynchronous, arjevani2020tight, feyzmahdavian2023asynchronous}, all of which assume fixed delays.
More recently, delay-adaptive methods have been proposed, aiming to improve performance by down-weighting very stale gradients \citep{cohen2021asynchronous, wu2022delayadaptive, koloskova2022sharper, mishchenko2022asynchronous}.

Particularly relevant to our work are asynchronous variants of \saga.
\citet{leblond2018improvedAsynchronous} developed a shared-parameter version in the spirit of \algname{HOGWILD!}, while \citet{glasgow2022asynchronous} studied a distributed setting that employs full gradients, in contrast to our stochastic-gradient perspective.

A large body of recent work investigates asynchronous methods in federated learning (FL), where clients hold data from heterogeneous distributions.
Notable contributions include \citep{aytekin2016analysis, xie2019asynchronous, mishchenko2022asynchronous, koloskova2022sharper, wang2022asyncfeded, wang2022asynchronous, glasgow2022asynchronous, fraboni2023general, zhang2023no, wang2024achieving, islamov2024asgrad, alahyane2025optimizing}.

More broadly, \citet{assran2020advances} provide a comprehensive survey of asynchronous optimization methods.

There is another line of work that began with \citet{tyurin2024optimal}, who established lower bounds for parallel methods and proposed optimal synchronous algorithms together with an asynchronous counterpart.
Several follow-up papers extended this semi-asynchronous framework to other settings \citep{tyurin2024freya,tyurin2024shadowheart, tyurin2024optimalgraph, maranjyan2025mindflayer}.

Finally, beyond asynchronous approaches, several synchronous methods address system heterogeneity by adapting local training to worker speeds.
The canonical method, \algname{FedAvg} \citep{mcmahan2017communication}, performs multiple local steps on each worker.
Variants have adapted the number of local steps to match workers' computation speeds \citep{li2020federated, maranjyan2025gradskip}, effectively balancing task assignments across heterogeneous systems.
More recently, \citet{maranjyan2025ata} proposed adapting the number of local steps dynamically, without requiring prior knowledge of worker speeds.
\section{Problem Setup}
We consider a distributed learning setting with $n$ workers, where each worker~$i$ possesses its own local data distribution~$\cD_i$.
Our goal is to solve the following distributed optimization problem:
\begin{equation}\label{eq:problem}
    \minimize\limits_{x \in \R^d}
    \left\{ f(x) \eqdef \frac{1}{n}\sum_{i=1}^{n}f_i(x) \right\},
    \quad\text{where}\quad
    f_i(x) \eqdef \ExpSub{\xi_i \sim \cD_i}{f_i(x; \xi_i)}.
\end{equation}
Here $f_i \colon \R^d \to \R$ denotes the local objective of worker $i$, defined as the expectation of the sample loss $f_i(x;\xi_i)$ over data points $\xi_i$ drawn from its local distribution $\cD_i$.
\subsection{Worker Heterogeneity Model}\label{sec:compute_times}
We first focus on the case where workers have constant computation speeds, as this setting is more intuitive and serves as a foundational model for understanding the dynamics of asynchronous distributed optimization.
The extension to arbitrary computation times is presented in \Cref{sec:arbitrary_time}.

Following the \textit{fixed computation model} \citep{mishchenko2022asynchronous}, we formalize:
\begin{equation}\label{eq:fixed_time}
    \begin{minipage}{0.52\linewidth}
    Each worker~$i$ requires $\tau_i$~seconds\footnotemark\ to compute one~stochastic~gradient~$\nabla f_i(x;\xi_i)$.
    
    \vspace{4pt}
    
    Without loss of generality, we assume
    
    \centering
    $0 < \tau_1 \le \tau_2 \le \cdots \le \tau_n$~.
    \end{minipage}
\end{equation}
\footnotetext{
    One could alternatively assume that each worker requires \textit{at most} $\tau_i$ seconds. 
    Under this formulation, all of our upper bounds would still hold; however, the lower bound would no longer be valid. 
    For this reason, we adopt the assumption that each worker requires exactly $\tau_i$ seconds.}
We assume communication is infinitely fast (taking $0$ seconds), both from workers to the server and from the server to workers\footnote{
    Alternatively, one could define $\tau_i$ as the time required for a worker to both compute a gradient and communicate it to the server, while keeping server-to-worker communication infinitely fast.
    Our upper bounds would still hold under this formulation, but the lower bounds would no longer apply, so we use the simpler model.}.
This is a modeling choice\textemdash arguably the simplest one\textemdash and has been the standard assumption in prior work \citet{mishchenko2022asynchronous, koloskova2022sharper, tyurin2024optimal, maranjyan2025ringmaster}, even if not always stated explicitly.
A related study by \citet{tyurin2024shadowheart} considers the case where communication is non-negligible and proposes techniques to address it.

Finally, we denote by $\tauavg \eqdef \tfrac{1}{n}\sum_{i=1}^n \tau_i$ the average computation time across all workers.
\subsection{Notations}
We denote the standard inner product in $\R^d$ by 
$$
    \inp{x}{y} = \sum_{i=1}^{d} x_i y_i~,
$$ 
and the corresponding Euclidean norm by $\norm{x} \eqdef \sqrt{\inp{x}{x}}$.
We use $[n] \eqdef \{1,2,\dots,n\}$ to denote the index set, and $\Exp{\cdot}$ for mathematical expectation.
For functions $\phi, \psi : \mathcal{Z} \to \R$, we write $\phi = \cO(\psi)$ if there exists a constant $C > 0$ such that $\phi(z) \leq C \psi(z)$ for all $z \in \mathcal{Z}$.
\subsection{Assumptions}
\label{sec:assumption}
We consider the following standard assumptions for the nonconvex setting.
\begin{boxedassumption}
    \label{ass:stochastic_variance_bounded}
    For each $i \in [n]$ and every $\xi_i$, the function $f_i(x;\xi_i)$ is differentiable with respect to its first argument $x$.
    Moreover, the stochastic gradients are unbiased and have bounded variance $\sigma^2 \geq 0$, that is,
    \begin{gather*}
         \ExpSub{\xi_i \sim \cD_i}{\nabla f_i(x;\xi_i)} = \nabla f_i(x), 
         \quad \forall x \in \R^d, \;\; \forall i \in [n]~,\\
         \ExpSub{\xi_i \sim \cD_i}{\sqnorm{\nabla f_i(x;\xi_i) - \nabla f_i(x)}} \leq \sigma^2,\quad \forall x \in \R^d, \;\; \forall i \in [n]~.
    \end{gather*}
\end{boxedassumption}
\begin{boxedassumption}\label{ass:lipschitz_constant}
    Each function $f_i$ is differentiable.
    There exists a constant $L>0$ such that, for all $x \in \R^d$ and $y_1,\dots,y_n \in \R^d$, 
    $$
        \sqnorm{\nabla f(x) - \frac{1}{n} \sum_{i=1}^n \nabla f_i(y_i)}
        \le \frac{L^2 }{n} \sum_{i=1}^n \sqnorm{x - y_i}.
    $$
\end{boxedassumption}
Recall the standard definition of smoothness
\begin{boxeddefinition}[Smoothness]\label{eq:lipschitz_gradient}
    A differentiable function $\phi \colon \R^d \to \R$ is called $L_\phi$--smooth if
    \begin{equation*}
        \norm{\nabla \phi(x) - \nabla \phi(y)} \le L_\phi \norm{x - y}, \quad \forall x,y \in   \R^d~.
    \end{equation*}
    By convention, $L_\phi$ denotes the smallest such constant.
\end{boxeddefinition}
Note that \Cref{ass:lipschitz_constant} is stronger than requiring $f$ itself to be $L_f$--smooth, yet weaker than all $f_i$ being $L_{f_i}$--smooth.
The following lemma establishes the relation among the constants $L_f$, $L$, and $L_{f_i}$.
\begin{boxedlemma}[Smoothness Bounds; Proof in \Cref{proof:smoothness}]\label{lemma:smoothness_relation}
    Suppose \Cref{ass:lipschitz_constant} holds with constant $L>0$.
    Then $f$ is $L_f$--smooth with $L_f \le L$.
    Moreover, if each $f_i$ is $L_{f_i}$--smooth, then \Cref{ass:lipschitz_constant} is satisfied, and we have
    $$
        L_f \;\le\; L \;\le\; \sqrt{\frac{1}{n}\sum_{i=1}^n L_{f_i}^2} \;\le\; \max_{i \in [n]} L_{f_i} =: L_{\max}~.
    $$
    Finally, if all $f_i$ are identical, i.e., $f_i = f$ for all $i \in [n]$, then $L = L_f$.
\end{boxedlemma}
To the best of our knowledge, prior work on asynchronous \sgd under data heterogeneity has always assumed smoothness of each $f_i$ \citep{koloskova2022sharper,nguyen2022federated,wang2025incremental}.
\begin{boxedassumption}\label{ass:lower_bound}
    There exists $f^* > -\infty$ such that $f(x) \geq f^*$ for all $x \in \R^d$.
    We define $\Delta \eqdef f(x^0) - f^*,$ where $x^0$ is the starting point of the optimization methods.
\end{boxedassumption}
Under these assumptions; our objective is to find an $\varepsilon$--stationary point\textemdash a (possibly random) vector $x$ satisfying $\mathbb{E}[\|\nabla f(x)\|^2] \leq \varepsilon$.
\section{Background and Motivation}
In this section, we review relevant existing methods in distributed optimization and discuss their limitations to motivate our algorithmic design.
We begin with \naiveminibatch as the canonical synchronous baseline, then consider \malenia \citep{tyurin2024optimal}, the first synchronous \sgd method to attain optimal time complexity.
We then turn to the challenges of asynchronous approaches, focusing on the recent \iasgd \citep{wang2025incremental} and its limitations.
Finally, we outline how these limitations can be addressed and introduce the core idea behind our algorithm.
\subsection{\naiveminibatchtitle}
\label{sec:minibatch}
\naiveminibatch provides the most straightforward approach to solving problem \eqref{eq:problem} in a distributed setting.
\paragraph{Algorithm Description.}
At each iteration $k$, the algorithm performs the following update:
$$
    x^{k+1} = x^k - \gamma \frac{1}{n} \sum_{i=1}^n \nabla f_i \(x^k; \xi_i^k\).
$$
The algorithm operates synchronously: at each iteration, the server waits to receive one stochastic gradient from each of the $n$ workers, all of which are evaluated at the current iterate $x^k$.
Once all gradients are collected, the server constructs an unbiased minibatch estimator of the full gradient by averaging these stochastic gradients and performs a standard gradient descent step.
\paragraph{Limitations.}
This synchronous approach has a significant computational bottleneck.
Each iteration requires waiting for the slowest worker to complete its gradient computation, resulting in the iteration time $\tau_n \eqdef \max_{i \in [n]} \tau_i$ \eqref{eq:fixed_time}.
Consequently, faster workers remain idle while waiting for stragglers, which leads to inefficient utilization of available computational resources.
\paragraph{Theoretical Performance.}
In terms of convergence rate, \naiveminibatch achieves the iteration complexity
$$
    \cO\( \frac{L_f\Delta}{\varepsilon}\(1 + \frac{\sigma^2}{n\varepsilon}\)\)
$$
to reach an $\varepsilon$--stationary point \citep{cotter2011better, goyal2017accurate, gower2019sgd}.
The corresponding time complexity becomes
$$
    \cO\(\frac{\tau_n L_f\Delta}{\varepsilon}\(1 + \frac{\sigma^2}{n\varepsilon}\)\) .
$$
This motivates the development of methods that can better utilize fast workers without waiting for stragglers.
\subsection{\maleniatitle}
\malenia \citep{tyurin2024optimal} resolves the straggler issue in \naiveminibatch by ensuring continuous worker utilization, rather than waiting for the slowest worker in each iteration.
However, it is fundamentally a \minibatch algorithm: in every iteration, it constructs an unbiased gradient estimator from multiple gradients and then performs a synchronous update.
The key distinction lies in how the batch is collected\textemdash \malenia gathers gradients asynchronously, allowing potentially multiple contributions from the same worker within a single iteration, unlike \naiveminibatch.
\paragraph{Algorithm Description.}
At each iteration $k$, all workers continuously compute stochastic gradients at the current point $x^k$.
The server accumulates these gradients until the stopping condition
\begin{equation}\label{eq:malenia_condition}
    \(\frac{1}{n} \sum_{i=1}^n \frac{1}{b_i^k} \)^{-1} \ge \max\left\{1,\frac{\sigma^2}{n\varepsilon}\right\}
\end{equation}
is satisfied, where $b_i^k$ denotes the number of gradients received from worker $i$.
Once this condition is met, the server performs the update
$$
    x^{k+1}
    = x^k - \gamma \bar g^k
    \eqdef x^k - \gamma \frac{1}{n} \sum_{i=1}^n \bar g_i^k~,
$$
where $\bar g_i^k$ is the average of the stochastic gradients received from worker $i$
$$
    \bar g_i^k = \frac{1}{b_i^k} \sum_{j=1}^{b_i^k} \nabla f_{i}\(x^{k}; \xi_{i}^{k,j}\).
$$
Because stochastic gradients are collected asynchronously, but the update is performed only after all required gradients are received, the algorithm can be regarded as semi-asynchronous\textemdash asynchronous in gradient collection but synchronous in parameter updates.
\paragraph{Stopping Condition Rationale.}
The left-hand side of condition \eqref{eq:malenia_condition} appears in the variance of the gradient estimator $\bar g^k$.
Ensuring that this quantity is sufficiently large allows the algorithm to control the variance at each step.
Moreover, condition \eqref{eq:malenia_condition} guarantees that every worker computes at least one gradient, which in turn yields a smaller variance than that of the estimator in \naiveminibatch.
\paragraph{Theoretical Performance.}
This asynchronous gradient collection strategy achieves the optimal time complexity
$$
    \cO\(\frac{L_f\Delta}{\varepsilon}\(\tau_n + \tauavg \frac{\sigma^2}{n\varepsilon}\)\),
$$
as shown by \citet{tyurin2024optimal}.
The key improvement over \naiveminibatch is that the variance term $\sigma^2$ is now multiplied by $\tau_{\text{avg}}$ instead of $\tau_n$.
Since $\tau_{\text{avg}} \ll \tau_n$ in computationally heterogeneous environments, \malenia can potentially provide substantial speedup in highly heterogeneous regimes.

The algorithm's benefits are most pronounced in the high-noise settings (where $\nicefrac{\sigma^2}{n \varepsilon}$ is large).
In low-noise scenarios or when $\sigma = 0$, \malenia offers no advantage over \naiveminibatch, since collecting multiple gradients per worker provides no benefit in terms of variance reduction.
\paragraph{Limitations.}
The main limitation of \malenia is its synchronous nature.
After collecting the necessary gradients, the server must broadcast the updated model to all workers simultaneously.
This all-at-once synchronization creates substantial communication overhead, which can become a major scalability bottleneck in large-scale or bandwidth-limited environments.

Moreover, the synchronization process forces the server to discard ongoing gradient computations from workers who are actively computing when the broadcast occurs.
Since workers operate continuously during the gradient collection phase, they must abandon their current computations upon receiving the new model, wasting valuable computational resources that could otherwise contribute to convergence.

Additionally, \malenia requires prior knowledge of the noise level $\sigma$ and the target accuracy $\varepsilon$ to evaluate the stopping condition \eqref{eq:malenia_condition}.
This dependence on problem-specific parameters, which are often unknown or difficult to estimate in practice, significantly limits its practical applicability.

These synchronization bottlenecks motivate the development of asynchronous optimal methods.
Beyond avoiding coordination overhead, asynchronous approaches enable immediate model updates upon gradient arrival, which can accelerate convergence through more frequent parameter updates.
This immediate processing is particularly beneficial for sparse models where gradients typically affect disjoint parameter subsets, allowing parallel updates without interference \citep{recht2011hogwild}.
\subsection{Toward Asynchronous Methods}
A naive approach to making optimization asynchronous would be to update the model immediately upon receiving any gradient, using
$$
    x^{k+1} = x^k - \gamma \nabla f_{i^k}\(x^{k-\delta^k};\xi_{i^k}^{k-\delta^k}\),
$$
where $i^k$ denotes the worker that sent the gradient at iteration $k$, and $\delta^k \ge 0$ is its delay, i.e., the number of server updates that occurred while the gradient was being computed.
Delays arise naturally in asynchronous execution: fast workers return gradients quickly and proceed with updated models, while slower workers compute on stale iterates; by the time they return their gradients, several server updates may already have occurred.
Consequently, $\delta^k$ is only determined when the server receives a gradient: at that point, it knows both the model iterate used by the worker and the current server iterate, and $\delta^k$ is simply the difference between the two.

\paragraph{Limitations.}
This approach suffers from a fundamental bias problem when workers have heterogeneous data distributions.
Faster workers send updates more frequently, causing their local objectives to dominate the optimization and pull the model toward their own minimizers.
Slower workers, when their updates finally arrive, push the model back toward different solutions.
As a result, the iterates may oscillate or stagnate, failing to converge to a stationary point.

This bias also makes theoretical analysis difficult.
Classical \algname{SGD}-style proofs rely on one-step progress toward minimizing the global function, but here each update direction reflects a different local objective.
Without additional data similarity assumptions \citep{mishchenko2022asynchronous,koloskova2022sharper,nguyen2022federated,islamov2024asgrad}, it becomes impossible to extend the analysis to the global function\textemdash yet such assumptions are rarely realistic when data can be arbitrarily heterogeneous across machines or organizations.

The root cause is that each update uses a gradient from one worker only.
A better strategy is to incorporate information from all workers, even if some gradients are computed at stale iterates.
This idea motivates methods such as Incremental Aggregated Gradient (\algname{IAG}) \citep{blatt2007convergent,gurbuzbalaban2017convergence,vanli2018global} and Stochastic Averaged Gradient (\algname{SAG}) \citep{roux2012stochastic,schmidt2017minimizing}, which maintain and aggregate gradients from all workers.





\subsection{\iasgdtitle}
As discussed above, the key insight is to construct a gradient estimator using information from all workers at each model update.
\iasgd \citep{wang2025incremental} achieves this by maintaining a gradient table on the server, similar to \algname{SAG} or \algname{IAG}, but with asynchronous table updates.
\paragraph{Algorithm Overview.} 
The server maintains a gradient table $\{g_i\}_{i=1}^n$ that stores the most recent gradient received from each of the $n$ workers.
The table is initialized by collecting one gradient from each worker at the starting point $x^0$.
The server then computes the first model update 
$$
    x^1 = x^0 - \gamma \bar g \eqdef x^0 - \gamma \frac{1}{n} \sum_{i=1}^n g_i
$$
and broadcasts $x^1$ to all workers.
Subsequently, the workers compute stochastic gradients in parallel, and their corresponding table entries are updated asynchronously as soon as the computations finish.

At each iteration $k$, the server performs the following steps:
\begin{enumerate}
    \item Receive gradient $\nabla f_{i^k}\Big(x^{k-\delta_{i^k}^k};\xi_{i^k}^{k-\delta_{i^k}^k}\Big)$ from worker $i^k$
    \item Update the gradient table entry: 
        $g_{i^k} = \nabla f_{i^k}\Big(x^{k-\delta_{i^k}^k};\xi_{i^k}^{k-\delta_{i^k}^k}\Big)$
    \item Perform model update: 
        $x^{k+1} = x^k - \gamma \bar g = x^k - \gamma \frac{1}{n} \sum_{i=1}^n g_i$
    \item Send the updated iterate $x^{k+1}$ to worker $i^k$ for its next computation
\end{enumerate}
The gradient estimator $\bar g$ combines the most recent gradient from each worker, ensuring that every update reflects information from the entire set of workers despite asynchronous execution.
In this way, the method retains the statistical benefits of using all workers' data while allowing them to operate independently, thereby avoiding the synchronization bottlenecks that limit scalability.

Note that, due to asynchrony, the stochastic gradients stored in the table generally correspond to different iterates of the model.
We therefore record each worker $i$'s delay $\delta_i^k$ to track the iterate at which its gradient was computed.

\paragraph{Theoretical Performance.}
The iteration complexity of this algorithm was established by \citet{wang2025incremental}, and by a straightforward conversion to the fixed computation time model (see \Cref{sec:iasgd_time_complexity}), we obtain the corresponding time complexity
$$
    \cO\(\frac{\tau_nL_{\max}\Delta}{\varepsilon}\(1 + \frac{\sigma^2}{n\varepsilon}\)\),
$$
which matches the complexity of \naiveminibatch.
This indicates that asynchronous execution alone does not provide computational advantages over synchronous approaches.
Thus, a fundamental challenge lies in how gradient delay affects convergence, which we address next.
\subsection{Motivation}
The primary limitation of asynchronous algorithms stems from gradient delay, which can significantly degrade convergence performance.
Large delays can cause the optimization steps to follow suboptimal trajectories, which disrupts convergence.

This delay problem occurs even in homogeneous settings where all functions $f_i$ are equal ($f_i \equiv f$ for all $i\in [n]$).
The state-of-the-art solution for this case, \algname{Ringmaster ASGD} \citep{maranjyan2025ringmaster}, achieves optimal time complexity by explicitly controlling delay to prevent it from becoming large.
\algname{Ringmaster ASGD} accomplishes this by discarding gradients that arrive with large delays.

Unfortunately, this gradient-discarding strategy fails in the heterogeneous setting of \iasgd.
The fundamental issue is that slow workers inevitably experience large delays due to their computational constraints.
If we ignore or discard their delayed gradients, the corresponding table entries remain outdated and may never be updated again if subsequent gradients also arrive late and are discarded.
This creates a persistent information bottleneck that degrades the quality of the gradient estimator and harms convergence.

This issue suggests we should prevent such situations from occurring by controlling the number of updates performed using fast workers.
The simplest approach would be to ignore some updates from fast workers, but this contradicts the core principle of asynchronous methods whereby all workers compute gradients continuously.

Instead, our approach \textit{buffers} the gradients produced by fast workers rather than applying them immediately, similar to the strategy in \malenia.
By buffering gradients and performing a model update only once a sufficient number have been collected, we control the delays induced by slow workers while keeping all workers continuously active.
This buffering mechanism provides an additional advantage: when multiple gradients computed at the same iterate are aggregated, they yield lower-variance estimates, thereby improving convergence.
\section{Ringleader ASGD}
We now formally introduce our algorithm, \ringleader.

\ringleader builds upon three key insights from existing methods.
First, inspired by \iasgd, we maintain a gradient table\footnote{
    It is not strictly necessary to maintain a table on the server.
    An alternative, as in \iasgd \citep{wang2025incremental}, is to store one vector on each worker along with an additional vector on the server.
    However, this can be problematic when workers have limited memory.
    For clarity and simplicity, we adopt the server-side table formulation in our description.}
to ensure that information from all workers is incorporated into every update, which eliminates the need for data similarity assumptions between workers.
Second, following \ringmaster, we recognize that controlling gradient delay is essential for efficient asynchronous optimization.
Third, drawing from \malenia, we use buffering of stochastic gradients\textemdash rather than discarding delayed ones\textemdash to control delays while preserving valuable computations, enabling continuous utilization of all resources.

An important property of the algorithm is that all workers remain continuously active, computing stochastic gradients.
As soon as a worker finishes computing a gradient, it immediately sends it to the server.
Recall that we assumed communication is instantaneous, i.e., takes zero time (\Cref{sec:compute_times}).
When the server receives a gradient, it either buffers it for later use or applies it immediately to perform a model update.
If the gradient is buffered, no further action is taken and the worker simply continues computing and sending new gradients.
If the server decides to perform an update, it updates the model and sends the updated model back to the worker that provided the gradient.
This server-to-worker communication is also assumed instantaneous, after which the worker resumes computing gradients at the new model point, ensuring that workers are never idle.

The algorithm proceeds in rounds.
In each round, exactly $n$ model updates are performed\textemdash one for each worker.
Specifically, when a worker sends a stochastic gradient, the server may apply an update and return the updated model to that worker, but it ensures that each worker receives an updated model at most once per round.
Repeating this procedure $n$ times ensures that every worker obtains exactly one fresh model per round, which in turn keeps delays bounded.
To avoid discarding the computations of fast workers, the server buffers their gradients and applies them only at the appropriate moment, thereby guaranteeing that each round consists of exactly $n$ updates.

Each round consists of two phases:
\begin{itemize}
    \item \textbf{Phase~1:} Buffer stochastic gradients in a table until at least one gradient from each worker is available.
    \item \textbf{Phase~2:} Perform exactly $n$ updates (one for each worker), then discard the old stochastic gradients from the table and return to Phase~1 to start collecting fresh ones.
\end{itemize}
The complete algorithm is formalized in \Cref{algo:Ringleader}.

\begin{algorithm}[H]\caption{\ringleader (server algorithm)}
    \label{algo:Ringleader}
    \begin{algorithmic}[1]
        \STATE \textbf{Input:} Stepsize $\gamma > 0$, initial point $x^0 \in \R^d$
        \STATE \textbf{Initialization:} Broadcast $x^0$ to all workers, which then start running \Cref{algo:Ringleader_worker} in parallel
        \STATE Set $k = 0$, $S = \emptyset$; initialize $G_i = 0$, $b_i = 0$ for all $i \in [n]$
        \WHILE{not terminated}
        \STATE {\highlightcolor \textbf{--- Phase~1: await stochastic gradients from all workers ---}}
        \WHILE{$S \ne [n]$}
            \STATE Receive stochastic gradient $g_j^k$ \big(at $x^{k-\delta_j^k}$\big) from some worker $j \in [n]$
                \label{line:receive_ph1}
            \STATE $G_j = G_j + g_j^k$;\;  $b_j = b_j + 1$;\; $S = S \cup \{j\}$
                \label{line:update_ph1}
        \ENDWHILE
        \STATE {\highlightcolor \textbf{--- Phase~2: perform exactly one update for every worker ---}}
        \STATE Update the model: $x^{k+1} = x^k - \gamma \frac{1}{n}\sum_{i=1}^n \nicefrac{G_i}{b_i}$
            \label{line:step_last_worker}
        \STATE Broadcast $x^{k+1}$ to worker $j$
            \hfill $\diamond$ Last worker to complete Phase~1
            \label{line:broadcast_last_worker}
        \STATE $k=k+1$;\; $S = S \setminus \{j\}$
            \label{line:counter_last_worker}
        \STATE $g_i^+ = 0$, $b_i^+ = 0$ for all $i \in [n]$;\; $S^+ = \emptyset$
            \hfill $\diamond$ Initialize temporary buffer
            \label{line:buffer_initialization}
        \WHILE{$S \ne \emptyset$}
            \STATE Receive stochastic gradient $g_j^k$ \big(at $x^{k-\delta_j^k}$\big) from some worker $j \in [n]$
            \IF{$j \in S$}
                \STATE $G_j = G_j + g_j^k$;\; $b_j = b_j + 1$
                    \label{line:update_table_ph2}
                \STATE Update the model: $x^{k+1} = x^k - \gamma \frac{1}{n}\sum_{i=1}^n \nicefrac{G_i}{b_i}$
                    \label{line:step_ph2}
                \STATE Broadcast $x^{k+1}$ to worker $j$
                    \label{line:broadcast_ph2}
                \STATE $k=k+1$;\; $S = S \setminus \{j\}$
                    \label{line:counter_ph2}
            \ELSE
                \STATE $G_j^+ = G_j^+ + g_j^k$;\; $b_j^+ = b_j^+ + 1$;\; $S^+ = S^+ \cup \{j\}$ 
                    \hfill $\diamond$ Buffer
                    \label{line:buffer}
            \ENDIF
        \ENDWHILE
        \STATE $G_i = G_i^+$;\; $b_i = b_i^+$ for all $i\in [n]$;\; $S = S^+$
            \label{line:transfer_table}
    \ENDWHILE
\end{algorithmic}
\end{algorithm}
\begin{algorithm}[H]\caption{Worker $i$'s subroutine}
    \label{algo:Ringleader_worker}
    \begin{algorithmic}[1]
        \STATE \textbf{Input:} Model $x$
        \WHILE{not terminated}
            \STATE Compute $g_i = \nabla f_i(x; \xi_i)$ using a freshly sampled data point $\xi_i \sim \cD_i$ 
            \STATE Send $g_i$ to the server 
        \ENDWHILE
    \end{algorithmic}
\end{algorithm}
\subsection{Detailed Description}
\paragraph{Initialization.}
The algorithm begins by broadcasting the initial point $x^0 \in \R^d$ to all workers, which then start executing the worker subroutine (\Cref{algo:Ringleader_worker}).
Each worker continuously computes stochastic gradients at its current point and sends them to the server until instructed to stop, at which point the server provides a new point to resume from.
This design ensures that workers remain fully utilized and never idle.

The server maintains a gradient table with entries $\{(G_i, b_i)\}_{i=1}^n$, all initialized to zero.
Here, $G_i$ accumulates gradients, while $b_i$ tracks the number of stochastic gradients received from worker $i$ to form proper minibatch estimators, with $b_i = 0$ for all $i \in [n]$ at the start.

Before Phase~1 begins, we also initialize the set $S = \emptyset$, which tracks which table entries contain at least one stochastic gradient.
Since no gradients have yet arrived, $S$ is empty.
\paragraph{Phase~1 --- Gradient Collection.}
In this phase, the server continuously receives stochastic gradients from the workers and stores them in the gradient table $\{(G_i, b_i)\}_{i=1}^n$.
We denote by $g_j^k$ the stochastic gradient sent by worker $j$ at iteration $k$, which is computed at point $x^{k-\delta_j^k}$ using an i.i.d. sample $\xi_j \sim D_j$.

We do not specify how the delays $\delta_j^k$ evolve, since this information is not needed to run the algorithm: whenever necessary, a delay can be obtained as the difference between the current model iterate and the iterate at which the stochastic gradient was computed.
The delays will only play a role in the analysis, not in the execution of the method.

Upon receiving $g_j^k$ from worker $j$ (Line~\ref{line:receive_ph1}), the server updates the corresponding table entry and the stochastic gradient counter as follows (Line~\ref{line:update_ph1})
$$
    G_j = G_j + g_j^k ~, \quad b_j = b_j + 1 ~, \quad S = S \cup \{j\} ~.
$$
This process continues until $S = [n]$, i.e., the server has collected at least one gradient from every worker.
No model updates are performed during this phase, and workers do not receive new points; hence, all stochastic gradients from a given worker are computed at the same point.
\paragraph{Phase~2 --- Sequential Updates.}
In this phase, the server performs exactly one model update for each worker $i$, for a total of $n$ updates.
Phase~2 starts with the last worker that completed Phase~1, i.e., the worker whose gradient made the table complete.
The server first computes an update by averaging the accumulated stochastic gradients in the table $\{(G_i, b_i)\}_{i=1}^n$ and taking a descent step with this estimate (Line~\ref{line:step_last_worker}).
The updated model is then sent to this worker (Line~\ref{line:broadcast_last_worker}), the worker is removed from the set $S$, and the iteration counter is incremented (Line~\ref{line:counter_last_worker}).

Next, the server must update the remaining $n-1$ workers.
These updates are performed sequentially as soon as each worker finishes its current computation.
During this waiting period, new gradients may arrive from workers not in $S$\textemdash e.g. for example, the last updated worker may send another stochastic gradient before the other workers complete their computation.
Since discarding these gradients would waste information, they are instead buffered for the next round.
\paragraph{Temporary Table Management.}
To achieve this, the server maintains a temporary table $\{(G_i^+, b_i^+)\}_{i=1}^n$, initialized to zero (Line~\ref{line:buffer_initialization}), together with a set $S^+$ that records which workers have contributed to the table.
Whenever a gradient arrives from a worker not in $S$, it is stored in the temporary table (Line~\ref{line:buffer}).

If instead the gradient comes from a worker $j \in S$\textemdash i.e., one of the workers whose model we still need to update\textemdash the server first updates the main table $\{(G_i, b_i)\}_{i=1}^n$ with this new stochastic gradient (Line~\ref{line:update_table_ph2}).
It then performs a model update by again averaging the accumulated stochastic gradients in the table (Line~\ref{line:step_ph2}), broadcasts the new model to worker $j$ (Line~\ref{line:broadcast_ph2}), increments the iteration counter, and removes $j$ from the set $S = S \setminus \{j\}$ (Line~\ref{line:counter_ph2}).
\paragraph{Preparing for the Next Round.}
Once all workers in $S$ have been updated and Phase~2 is complete ($S=\emptyset$), the server prepares for the next round by copying the contents of the temporary table $\{(G_i^+, b_i^+)\}_{i=1}^n$ into the main table $\{(G_i, b_i)\}_{i=1}^n$ (Line~\ref{line:transfer_table}).
The set $S$ is also reset to $S = S^+$, since these workers already contributed gradients at their updated models.
Entries in the main table corresponding to workers not in $S^+$ remain zero, as the temporary table was initialized with zeros at the start of Phase~2 (Line~\ref{line:buffer_initialization}).

The server can now proceed to the next round by returning to Phase~1 and beginning a new gradient collection phase.
\subsection{Delay Control Analysis}
The structure of \ringleader, with its two phases, is specifically designed to prevent the unbounded delays that arise in standard asynchronous methods.
To understand why this works, consider that in each round we perform exactly $n$ updates\textemdash one per worker\textemdash before moving to the next round.
This ensures that no worker can fall more than one full round behind the current iteration.
The precise bound on delays is given in the following lemma
%
\begin{boxedlemma}[Bounded Delay]\label{lemma:delay}
   In \ringleader (\Cref{algo:Ringleader}), the delays $\delta_i^k$ satisfy
   $$
       \delta_i^k \le 2n-2~,
   $$
   for any worker $i \in [n]$ and any iteration $k \geq 0$.
\end{boxedlemma}
\begin{proof}
    We prove this by analyzing the structure of \ringleader.
    The algorithm operates in rounds, where each round consists of Phase~1 (gradient collection) followed by Phase~2 (sequential updates).
    In each Phase~2, the server performs exactly $n$ updates, one for each worker.
    Phase~2 begins at iterations $0, n, 2n, 3n, \ldots$, i.e., at multiples of $n$.
    \paragraph{First round (iterations $0$ to $n-1$):}
    Initially, all workers compute gradients at the point $x^0$, so during iterations $0, 1, \ldots, n-1$, the server receives gradients computed at $x^0$.
    For any iteration $k$ in this range, the server processes stochastic gradients computed at point $x^{k-\delta_i^k}$, so $\delta_i^k = k \le n-1$.
    Thus, delays simply increment during this first Phase~2.

    At the end of this round, each worker $i$ has received a new model $x^j$ for some $j \in \{1, 2, \ldots, n\}$, and these update iterations are distinct across workers.
    \paragraph{Second round (iterations $n$ to $2n-1$):}
    At the start of the second Phase~2 (at iteration $n$), the gradient table contains gradients computed at points $x^{n-\delta_i^n}$ for worker $i$, where $\delta_i^n \in \{0, 1, \ldots, n-1\}$.
    These delay values are distinct across workers since each worker received its update at a different iteration in the previous round.
    
    During iterations $n$ to $2n-1$, these delays increase by $1$ at each iteration for the same reason as in the first Phase~2, giving $\delta_i^{2n-1} \in \{n-1, n, \ldots, 2n-2\}$ at the end of this round.
    At the same time, all workers receive new points to compute gradients from, so during the next Phase~2, the delays will again be distinct for all workers and in $\{0, 1, \ldots, n-1\}$.
    \paragraph{General pattern:}
    By induction, at the beginning of each round starting at iteration $cn$ (for integer $c \ge 1$), the delays $\delta_i^{cn}$ take distinct values in $\{0, 1, \ldots, n-1\}$.
    During each Phase~2, these delays increase by at most $n-1$, giving the bound
    $$
        \delta_i^k \le (n-1) + (n-1) = 2n-2~.
    $$
\end{proof}
\subsection{Comparison to \iasgdtitle}
Our method is a delay-controlled version of \iasgd.
We can recover \iasgd by removing Phase~1 (gradient collection) and Phase~2 (structured updates), and thus perform updates naively\textemdash immediately upon gradient arrival.
In contrast, our algorithm operates in structured rounds, performing exactly one update per worker in each round, which provides the crucial delay control that \iasgd lacks.

In \iasgd, delays for slow workers can grow unboundedly because the server continuously updates the model using gradients from fast workers, causing slow workers to fall increasingly behind.
Our method prevents this issue by buffering the gradients from fast workers rather than immediately applying these gradients, to ensure that all workers receive updated models within $n$ subsequent iterations.
\subsection{Comparison to \maleniatitle}
\malenia also operates as an algorithm with two phases.
In Phase~1, \malenia collects gradients using a similar method to our approach, but uses a different termination condition \eqref{eq:malenia_condition} that requires knowledge of the noise parameter $\sigma$ and the target stationarity level $\varepsilon$, making it impractical.
In Phase~2, \malenia performs a \textit{synchronous} update by averaging all collected gradients and broadcasting the new model to \textit{all} workers simultaneously before returning to Phase~1.
This synchronization forces \malenia to discard ongoing gradient computations from workers that are active during the broadcast.

In contrast, our method performs Phase~2 \textit{asynchronously}: we update workers sequentially as they complete their current gradient computations, which ensures that no computational work is wasted.

Regarding \malenia's termination condition \eqref{eq:malenia_condition}, in \Cref{sec:malenia_param_free} we demonstrate that this condition can be replaced with our simpler requirement of obtaining at least one gradient from every worker.
With this modification, \malenia remains optimal in the fixed-time regime \eqref{eq:fixed_time} while becoming parameter-free, which eliminates the need for prior knowledge of $\sigma$ and $\varepsilon$.
In this parameter-free variant, the only difference between \malenia and \ringleader lies in Phase~2: we perform updates asynchronously without discarding gradients, while \malenia operates synchronously.
\section{Theoretical Results}
Before presenting the theoretical results, we first write the algorithm in a compact form.
The gradients for each worker in the table are all computed at the same point; for worker $i$ at iteration $k$, the point is $x^{k-\delta_i^k}$.
The update rule can be written compactly as
$$
    x^{k+1} = x^k - \gamma \bar g^k~,
$$
where the gradient estimator $\bar g^k$ is defined by
$$
    \bar g^k
    \eqdef \frac{1}{n} \sum_{i=1}^n \bar g_i^k
    \eqdef \frac{1}{n} \sum_{i=1}^n \frac{1}{b_i^k} \sum_{j=1}^{b_i^k} g_i^{k,j} ~.
$$
Since multiple gradients may be received from the same worker, we denote by $g_i^{k,j}$ the $j$-th gradient from worker $i$ at iteration $k$.
Here the index $j$ corresponds to the i.i.d. sampled data point, and more concretely
$$
    g_i^{k,j} \eqdef \nabla f_i\(x^{k-\delta_i^k};\xi_i^{k-\delta_i^k, j}\)~.
$$
The quantity $b_i^k$ denotes the number of gradients from worker $i$ stored in the table at iteration $k$, i.e., the value of $b_i$ in Lines~\ref{line:step_last_worker} and \ref{line:step_ph2}.
Thus, the pair $(b_i^k, \delta_i^k)$ fully determines the method's behavior at iteration $k$.

Note that the sequence $\{b_i^k\}$ depends only on the computation times $\{\tau_i\}$ and the algorithm design (i.e., the stopping rule for collecting gradients).
Once these are fixed, all $b_i^k$ for every $i \in [n]$ and iteration $k$ are determined.
Crucially, the values of $b_i^k$ do not depend on the method's hyperparameters $\gamma$, $x^0$, or on the variance parameter $\sigma$ or the stationarity level $\varepsilon$.
\subsection{Iteration Complexity}
Our convergence analysis follows the structured approach employed by \citet{maranjyan2025ringmaster}, which decomposes the proof into two key components: a descent lemma that tracks the progress of the objective function and a residual estimation lemma that controls the accumulated delays in the system.

We begin by establishing notation for the harmonic means of the batch sizes across rounds:
$$
    B^k = \( \frac{1}{n} \sum_{i=1}^n \frac{1}{b_i^k} \)^{-1}, \quad \text{and} \quad B = \inf_{k \ge 0} B^k ~.
$$
Note that $B \geq 1$, since by the algorithm's design each $b_i^k \geq 1$.
A sharper bound on $B$ will be established later in \Cref{lemma:time_for_n_iter}.

The first lemma quantifies how much the objective function decreases at each iteration, accounting for both the standard gradient descent progress and the additional complexities introduced by asynchronous updates.
\begin{boxedlemma}[Descent Lemma; Proof in \Cref{proof:descent}]\label{lemma:descent}
    Under Assumptions~\ref{ass:stochastic_variance_bounded} and~\ref{ass:lipschitz_constant}, if the stepsize in \Cref{algo:Ringleader} satisfies $\gamma \le \nicefrac{1}{4L}$, then the following inequality holds
    \begin{align*}
        \E{f\(x^{k+1}\)}
        &\le \E{f\(x^{k}\)}
            - \frac{\gamma}{2} \E{ \sqnorm{\nabla f\(x^{k}\)} }
            - \frac{\gamma}{4} \E{\sqnorm{\frac{1}{n} \sum_{i=1}^n \nabla f_i\(x^{k-\delta_i^k}\)}} \\
            &\quad + \frac{\gamma L^2}{2n} \sum_{i=1}^n \E{\sqnorm{x^{k} - x^{k-\delta_{i}^k}}}
                    + \frac{3\gamma^2 L \sigma^2}{2B} \\
            &\quad + \gamma^2 L \sum_{\ell = k-(k \bmod n)}^{k-1} \E{\sqnorm{\frac{1}{n}\sum_{i=1}^n \nabla f_i\(x^{\ell-\delta_i^\ell}\)}}.
    \end{align*}
\end{boxedlemma}
This descent lemma shares a similar structure with its counterpart in the homogeneous setting analyzed by \citet{maranjyan2025ringmaster}, but with a crucial additional term.
The final summation term in the upper bound captures the effect of using stale gradients from the gradient table\textemdash a phenomenon we refer to as ``table delay''.
This term is absent in the homogeneous case because no gradient table is maintained.
Indeed, when $n=1$, our setting reduces to the homogeneous case, the gradient table becomes unnecessary, and this additional term vanishes, recovering the original descent lemma established by \citet{maranjyan2025ringmaster}.

Next, similar to the work of \citet{maranjyan2025ringmaster}, we derive a lemma to bound the term involving the difference between current and old points.
\begin{boxedlemma}[Residual Estimation; Proof in \Cref{proof:residual}]
    \label{lemma:residual}
    Under \Cref{ass:stochastic_variance_bounded}, the iterates of \ringleader (\Cref{algo:Ringleader}) with stepsize $\gamma \le \nicefrac{1}{4nL}$ satisfy the following bound
    \begin{equation*}
    \frac{1}{K} \sum_{k=0}^{K-1} \frac{1}{n} \sum_{i=1}^n \E{\sqnorm{x^{k} - x^{k-\delta_i^k}}}
        \le \frac{2\gamma n}{LK}\sum_{k=0}^{K-1} \E{\sqnorm{ \frac{1}{n}\sum_{j=1}^n \nabla f_j\(x^{k-\delta_j^k}\) }}
            + \frac{2\gamma \sigma^2}{LB}~.
    \end{equation*}
\end{boxedlemma}
Finally, we get the iteration complexity combining these two lemmas.
\begin{boxedtheorem}[Iteration Complexity]\label{theorem:convergence}
    Under Assumptions \ref{ass:stochastic_variance_bounded}, \ref{ass:lipschitz_constant}, and \ref{ass:lower_bound}, let the stepsize in \ringleader (\Cref{algo:Ringleader}) be
    $$
        \gamma = \min \left\{\frac{1}{8nL}, \frac{\varepsilon B}{10 L \sigma^2 } \right\}.
    $$
    Then,
    $$
        \frac{1}{K}\sum_{k=0}^{K-1} \E{ \sqnorm{\nabla f\(x^{k}\)} } 
        \le \varepsilon~,
    $$
    for
    $$
        K 
        \ge \frac{32 nL\Delta }{\varepsilon} + \frac{40 L\Delta \sigma^2}{B\varepsilon^2}
        = \cO\left(\frac{nL\Delta}{\varepsilon} \left(1 + \frac{\sigma^2}{Bn\varepsilon}\right)\right).
    $$
\end{boxedtheorem}
\begin{proof}    
    We start by averaging the inequality from \Cref{lemma:descent} over $K$ iterations and dividing both sides by $\nicefrac{\gamma}{2}$
    \begin{align*}
        \frac{1}{K}\sum_{k=0}^{K-1} & \E{ \sqnorm{\nabla f\(x^{k}\)} } 
            + \frac{1}{2K}\sum_{k=0}^{K-1} \E{\sqnorm{\frac{1}{n} \sum_{i=1}^n \nabla f_i\(x^{k-\delta_i^k}\)}} \\
        &\le \frac{2\Delta}{\gamma K}
            + \frac{3 \gamma L \sigma^2}{B} \\
            &\quad + \frac{L^2}{n} \frac{1}{K}\sum_{k=0}^{K-1} \sum_{i=1}^n \E{\sqnorm{x^{k} - x^{k-\delta_{i}^k}}} \\
            &\quad + 2 \gamma L \frac{1}{K}\sum_{k=0}^{K-1} \sum_{\ell = k- (k \bmod n)}^{k-1} \E{\sqnorm{\frac{1}{n}\sum_{i=1}^n \nabla f_i\(x^{\ell-\delta_i^\ell}\)}}.
    \end{align*}
    We now bound the third term on the right using \Cref{lemma:residual}
    \begin{align*}
        \frac{1}{K}\sum_{k=0}^{K-1} & \E{ \sqnorm{\nabla f\(x^{k}\)} } 
            + \frac{1}{2K}\sum_{k=0}^{K-1} \E{\sqnorm{\frac{1}{n} \sum_{i=1}^n \nabla f_i\(x^{k-\delta_i^k}\)}} \\
        &\le \frac{2\Delta}{\gamma K}
            + \frac{3 \gamma L \sigma^2}{B}
            + \frac{2\gamma L \sigma^2}{B} \\
            &\quad + 2\gamma L n \frac{1}{K}\sum_{k=0}^{K-1} \E{\sqnorm{ \frac{1}{n}\sum_{j=1}^n \nabla f_j\(x^{k-\delta_j^k}\) }} \\
            &\quad + 2 \gamma L \frac{1}{K}\sum_{k=0}^{K-1} \sum_{\ell = k- (k \bmod n)}^{k-1} \E{\sqnorm{\frac{1}{n}\sum_{i=1}^n \nabla f_i\(x^{\ell-\delta_i^\ell}\)}} \\
        &\le \frac{2\Delta}{\gamma K}
            + \frac{5 \gamma L \sigma^2}{B} \\
            &\quad + 2\gamma L n\frac{1}{K}\sum_{k=0}^{K-1} \E{\sqnorm{ \frac{1}{n}\sum_{j=1}^n \nabla f_j\(x^{k-\delta_j^k}\) }} \\
            &\quad + 2 \gamma Ln \frac{1}{K}\sum_{k=0}^{K-1} \E{\sqnorm{\frac{1}{n}\sum_{i=1}^n \nabla f_i\(x^{k-\delta_i^k}\)}} .
    \end{align*}
    Now, using the bound $\gamma \le \nicefrac{1}{8nL}$, we obtain
    \begin{equation*}
        \frac{1}{K}\sum_{k=0}^{K-1} \E{ \sqnorm{\nabla f\(x^{k}\)} }
        \le \frac{2\Delta}{\gamma K}
            + \frac{5 \gamma L \sigma^2}{B}~.
    \end{equation*}
    Finally, plugging in the stepsize and the expression for $K$ ensures the right-hand side is bounded by $\varepsilon$.
\end{proof}
For parallel and asynchronous methods, iteration complexity is less important than time complexity.
What truly matters is how quickly we can finish training.
We are willing to perform more iterations and extra computation if it means completing the process faster.
Having established the iteration complexity, we now turn to the time complexity.
\subsection{Time Complexity}
Since the algorithm operates in rounds with $n$ steps per round, and its iteration complexity is already known, it remains to determine the duration of each round.
We have the following lemma
\begin{boxedlemma}\label{lemma:time_for_n_iter}
    Each block of $n$ consecutive iterations (each round) of \Cref{algo:Ringleader} takes at most $2\tau_n$ seconds.
    Moreover, we have
    $$
        B \ge \frac{\tau_n}{2} \left( \frac{1}{n} \sum_{i=1}^n \tau_i \right)^{-1} = \frac{\tau_n}{2\tauavg}~.
    $$
\end{boxedlemma}
\begin{proof}
   The upper bound of $2\tau_n$ follows from the structure of the algorithm, which consists of two phases.
   In the first phase, the server waits until all workers complete at least one gradient computation, which takes at most $\tau_n$ seconds.
   In the second phase, the server applies the received gradients and waits for all ongoing computations to finish\textemdash which again takes at most $\tau_n$ seconds.
   Thus, the total time for $n$ iterations is bounded by $2\tau_n$.

   We now prove the second part of the lemma.
   Recall that
   $$
      B = \inf_{k \ge 0} B^k = \inf_{k \ge 0} \left( \frac{1}{n} \sum_{i=1}^n \frac{1}{b_i^k} \right)^{-1} \ge \left( \frac{1}{n} \sum_{i=1}^n \frac{1}{b_i} \right)^{-1},
   $$
   where we define
   $$
      b_i = \inf_{k \ge 0} b_i^k~.
   $$
   We are interested in the number of gradients stored in the table at iteration $k$.
   This count includes gradients computed during Phase~1 plus one additional gradient from Phase~2 (except for the worker that finished Phase~1 last).

   Since every worker needs to compute at least one gradient during Phase~1, the slowest worker will take $\tau_n$ seconds to complete single gradient computation.
   During this $\tau_n$-second interval, faster workers $i < n$ may still be finishing gradients from the previous round's Phase~2 before starting their Phase~1 computations for the current round.

   After completing any remaining Phase~2 work (which takes at most $\tau_i$ seconds), worker $i$ has at least $\tau_n - \tau_i$ seconds remaining to compute additional gradients for the current round's Phase~1.
   The number of gradients that worker $i$ can compute satisfies
   $$
      b_i \ge \max\left\{1, \left\lceil \frac{\tau_n - \tau_i}{\tau_i} \right\rceil \right\} \ge \max\left\{1, \frac{\tau_n}{\tau_i} - 1 \right\}.
   $$
   For workers $i$ where $\tau_n \geq 2\tau_i$, we have
   $$
      \frac{\tau_n}{\tau_i} - 1 \ge \frac{\tau_n}{2\tau_i}~,
   $$
   and hence
   $$
      b_i \ge \frac{\tau_n}{2\tau_i}~.
   $$
   Plugging this bound into the expression for $B$ gives the claimed result.
\end{proof}
Based on this lemma, we derive the time complexity guarantee of our algorithm
\begin{boxedtheorem}\label{thm:time_complexity}
    Let Assumptions \ref{ass:lipschitz_constant}, \ref{ass:lower_bound}, and \ref{ass:stochastic_variance_bounded} hold.
    Let the stepsize in \ringleader (\Cref{algo:Ringleader}) be
    $\gamma = \min \left\{\nicefrac{1}{8nL}, \nicefrac{\varepsilon B}{10 L \sigma^2 } \right\}$.
    Then, under the \emph{fixed computation model} \eqref{eq:fixed_time}, \ringleader achieves the optimal time complexity
    $$
        \cO\(\frac{L\Delta}{\varepsilon}\( \tau_n + \tauavg\frac{\sigma^2}{n\varepsilon} \) \) .
    $$
\end{boxedtheorem}
\begin{proof}
    We start with the iteration complexity from \Cref{theorem:convergence}
    $$
        K 
        \ge \frac{32 nL\Delta }{\varepsilon} + \frac{40 L\Delta \sigma^2}{B\varepsilon^2}
        = \cO\(\frac{nL\Delta}{\varepsilon} \(1 + \frac{\sigma^2}{B n\varepsilon}\) \).
    $$
    The time to do $n$ steps is at most $2\tau_n$ form \Cref{lemma:time_for_n_iter}.
    Then the time complexity is 
    $$
        2\tau_n \times \frac{K}{n} 
        = \cO\(\tau_n\frac{L\Delta}{\varepsilon} \( 1 + \frac{\sigma^2}{Bn\varepsilon} \) \) .
    $$
    It remains to put $B \ge \nicefrac{\tau_n}{2\tauavg}$ from \Cref{lemma:time_for_n_iter}.
\end{proof}
The obtained time complexity consists of two key terms that illuminate the algorithm's behavior.
The first term depends on the slowest device, which is fundamental since all devices must contribute to solving the problem.
The second term, however, involves $\tauavg$ rather than $\tau_n$ as in \naiveminibatch (see \Cref{table:fixedtime})\textemdash this substitution captures the core benefit of asynchronous execution.
Specifically, this advantage becomes pronounced when $\sigma$ is relatively large.
Intuitively, in high-noise regimes, collecting many gradients from workers is essential for convergence, and asynchronous methods can leverage faster workers more effectively.
Conversely, in low-noise settings, fewer gradient evaluations suffice for good performance, making \naiveminibatch already quite effective and rendering the additional complexity of asynchrony unnecessary.

This time complexity result matches the lower bound derived by \citet{tyurin2024optimal}, thus making \ringleader the first asynchronous algorithm to achieve optimality.
\section{Experiments}
\label{sec:experiments}
To validate our theoretical results we perform a toy simulation.

We consider image classification on MNIST \citep{lecun2010mnist} and on Fashion-MNIST \citep{xiao2017fashion} with standard normalization for both datasets.
To enable equal client sizes, we first trim the dataset so that the total number of examples is divisible by the number of clients $n=100$.
To obtain heterogeneous datasets across clients, we then partition the trimmed set using an \emph{equal-size Dirichlet} procedure with concentration parameter $\alpha=0.1$ \citep{yurochkin2019bayesian}.
For each client $j\in[n]$, we draw proportions $p_j \sim \mathrm{Dirichlet}(\alpha,\ldots,\alpha)$ over the classes and compute a rounded class-allocation vector whose entries sum exactly to $\nicefrac{N}{n}$, where $N$ is the trimmed dataset size.
This creates non-IID data where each client has a skewed distribution over classes (with $\alpha=0.1$, clients typically observe only 1-2 classes frequently).

When assigning samples, we take the requested number from each class pool for client $j$.
If a class pool does not have enough remaining examples to meet the requested amount, the client receives whatever is left from that class and the shortfall is \emph{topped up} using samples from other classes that still have available examples.

Our model is a two-layer MLP $\mathrm{Linear}(d,128)\!\to\!\mathrm{ReLU}\!\to\!\mathrm{Linear}(128,10)$ trained with mean cross-entropy.
Stochastic gradients at the clients use a minibatch of size $4$, while reported gradient norms are computed on the \emph{full} dataset.

We emulate heterogeneous compute by assigning each client $i$ a base delay and a random jitter:
$$
    \tau_i \;=\; i \;+\; |\eta_i|~, \qquad 
    \eta_i \sim \cN(0,i)~, 
    \quad \text{for all } i\in [n]~.
$$  
For each method we tune the stepsize $\gamma$ within a fixed wall-clock budget.
We sweep  
$$
    \gamma \in \{0.001,\,0.005,\,0.01,\,0.02,\,0.05,\,0.1,\,0.2,\,0.5,\,1.0,\,2.0\}~,
$$  
and then fix the best $\gamma$ per method for evaluation.

We report the full-batch gradient-norm squared $\| \nabla f(x^k) \|^2$, versus wall-clock time.
Each method is run $30$ times over different seeds.
We report the \emph{median} with \emph{interquartile range (IQR)}.
To reduce high-frequency noise, we apply a centered moving-average smoothing to the aggregated curves (post-aggregation), while keeping the initial point unchanged.

\begin{figure}[h]
    \centering
    \begin{minipage}{0.48\linewidth}
        \centering
        \includegraphics[width=\linewidth]{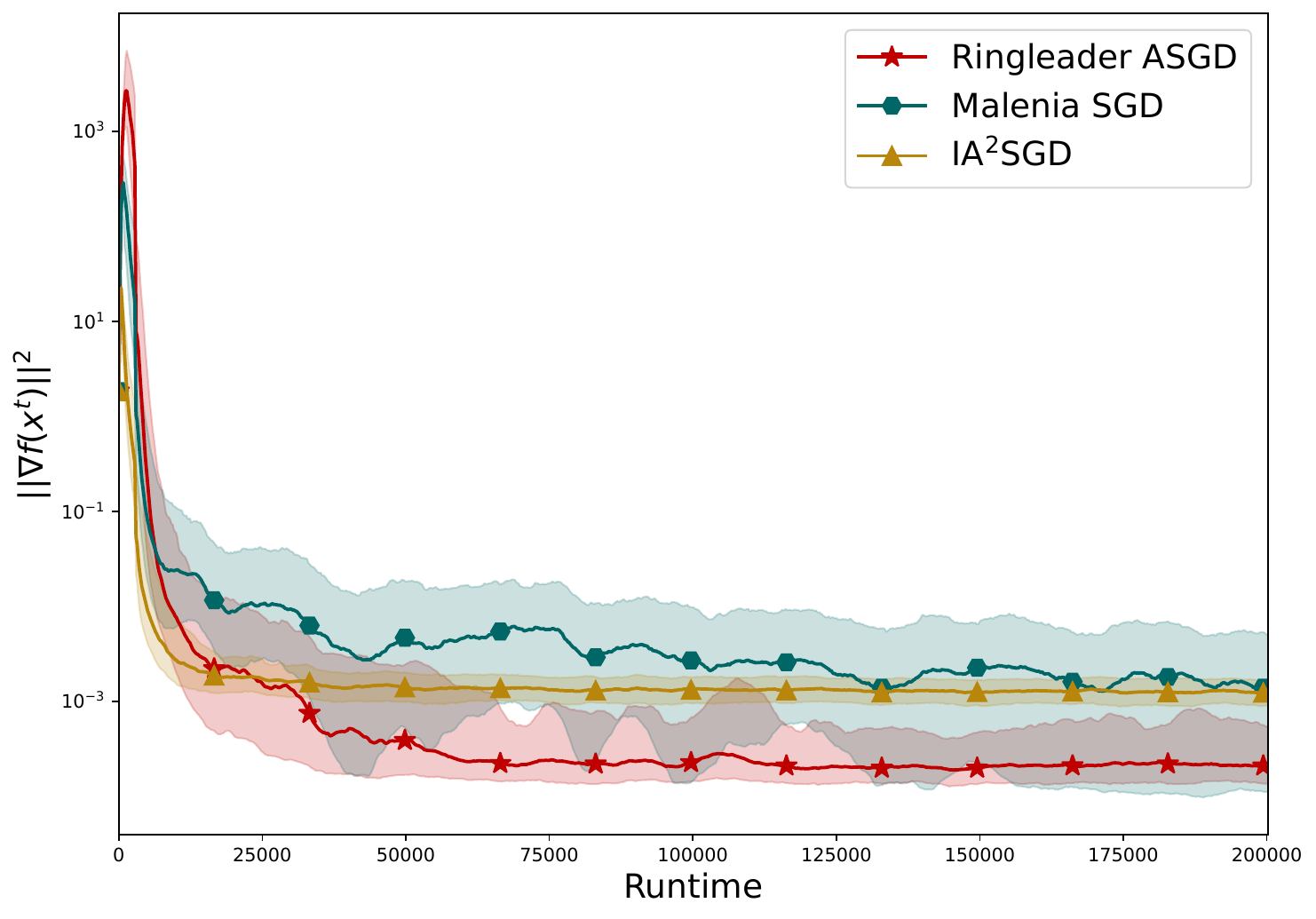}
        \caption*{(a) MNIST}
    \end{minipage}
    \hfill
    \begin{minipage}{0.48\linewidth}
        \centering
        \includegraphics[width=\linewidth]{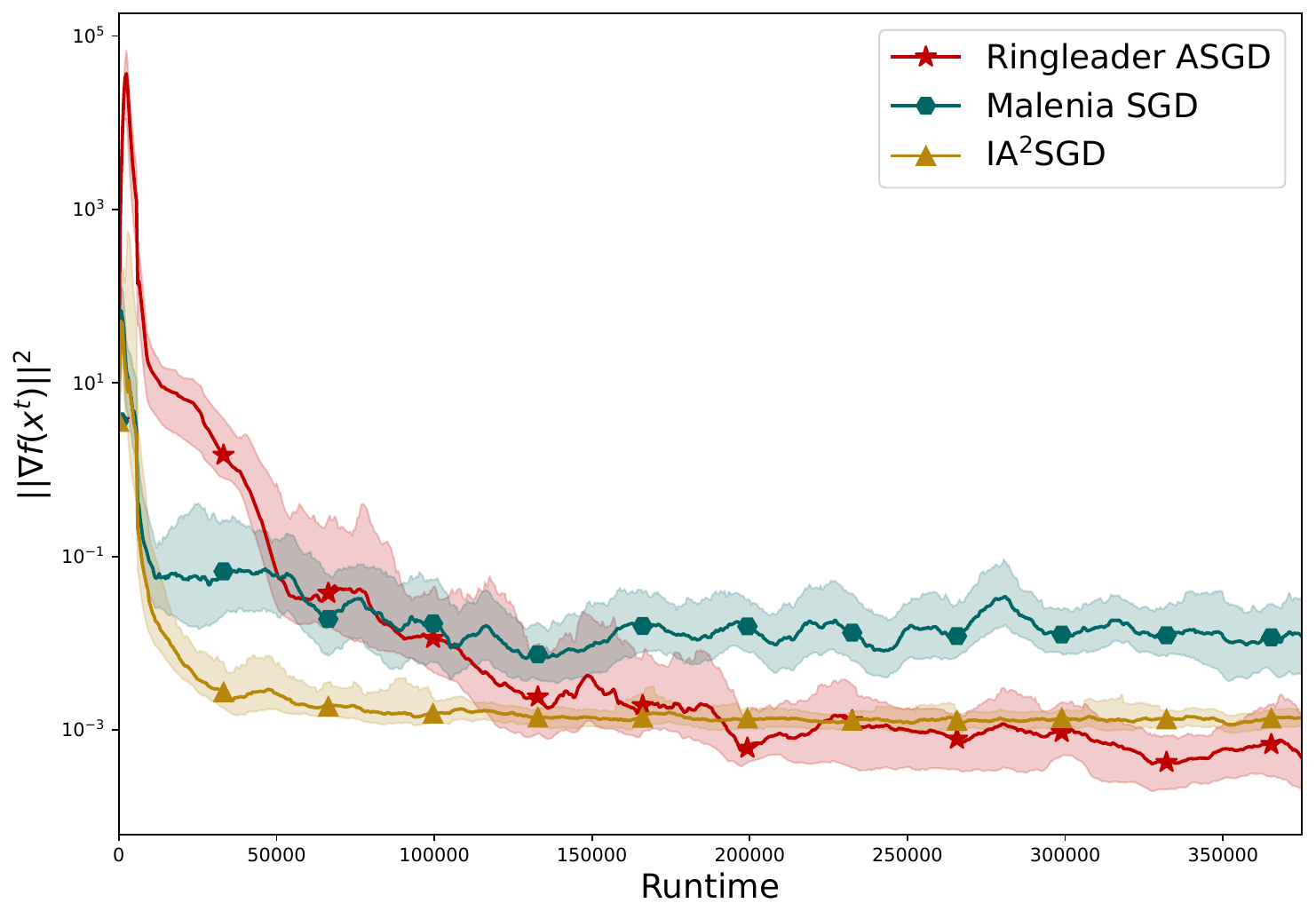}
        \caption*{(b) Fashion-MNIST}
    \end{minipage}
    \caption{
        Convergence comparison showing median squared gradient norm $\|\nabla f(x^k)\|^2$ (solid lines) with interquartile ranges (shaded regions) plotted against wall-clock time, averaged over 30 random seeds.
        \textbf{Setup:} A two-layer MLP with architecture $\mathrm{Linear}(d,128)$ $\rightarrow$ $\mathrm{ReLU}$ $\rightarrow$ $\mathrm{Linear}(128,10)$ is trained on (a) MNIST and (b) Fashion-MNIST datasets.
        \textbf{Client delays:} Heterogeneous delays are simulated as $\tau_i = i + |\eta_i|$, where $\eta_i \sim \mathcal{N}(0, i)$ for each client $i \in [n]$, with $n = 100$.
        \textbf{Results:} With optimally tuned stepsizes, \ringleader converges faster than both \malenia and \iasgd, even though \ringleader and \malenia share the same theoretical time-complexity guarantees.
    }
    \label{fig:median}
\end{figure}

\Cref{fig:median} shows the results.
We observe that \ringleader converges faster compared to \malenia and \iasgd.  
Although both \ringleader and \malenia have the same theoretical time complexity, in practice \ringleader performs better because Phase 2 involves $n$ asynchronous updates instead of a single synchronous one.
This design enables more optimization steps within the same wall-clock budget, which is especially advantageous when updates are sparse.

\section{Conclusion}
We have introduced \ringleader, the first asynchronous stochastic gradient method to achieve optimal time complexity under arbitrary data heterogeneity and arbitrarily heterogeneous computation times in distributed learning, without requiring similarity assumptions between workers' datasets.

Its core innovation is a two-phase structure within each round: the model is updated once per worker (for a total of $n$ updates), while a buffering mechanism manages gradient delays and preserves the efficiency of asynchronous execution.
By maintaining a gradient table and alternating between gradient collection and sequential updates, \ringleader prevents the unbounded delays common in naive asynchronous methods.
Every gradient received by the server is either used in the current round or stored for future use, ensuring no computation is wasted.

Our analysis shows that \ringleader matches the optimal time complexity bounds established by \citet{tyurin2024optimal}.
In contrast to the optimal but synchronous \malenia method, \ringleader is asynchronous and requires no prior knowledge of problem parameters in the algorithm design, making it practical for real-world deployments.

Finally, with a minor modification, \ringleader also achieves optimality in the more general setting of arbitrarily varying computation times (\Cref{sec:arbitrary_time}).

\chapter{
    ATA: Adaptive Task Allocation for Efficient\\
    Resource Management in Distributed Machine Learning
    }
\label{chapter:ATA}
\thispagestyle{empty}
This chapter is based on the work of \citet{maranjyan2025ata}:
\begin{quote}
    \bibentry{maranjyan2025ata}
\end{quote}
\section{Introduction}\label{ata:section:introduction}
In this work, we address a very general yet fundamental and important problem arising in various contexts and fields.
In particular, there are $n$ workers/nodes/devices collaborating to run some iterative algorithm which has the following structure:
\begin{itemize}
    \item 
        In order to perform a single iteration of the algorithm, a certain number ($B$) of {\em tasks} needs to be performed.
    \item 
        Each task can be computed by any worker, and the tasks are not temporally related.
        That is, they can be computed in any order, in parallel, and so on.
    \item 
        Whenever a worker is asked to perform a single task, the task will take a certain amount of time, modeled as a nonnegative random variable drawn from an unknown distribution specific to that worker.
        The stochastic assumption makes sense because in real systems computation times are not fixed and can vary with each iteration \citep{dean2013tail,chen2016revisiting, dutta2018slow, maranjyan2025mindflayer}.
    \item 
        Each worker can only work on a single task at a time.
        That is, a worker processes all tasks it has to perform sequentially.
        Different workers work in parallel.
\end{itemize}
%
%
\subsection{Greedy Task Allocation}
A natural goal in this setup is to make sure all tasks are completed as fast as possible (in expectation), which minimizes the (expected) time it takes for a single iteration of the algorithm to be performed provided that the task completion time is the dominant time factor of the iteration.
Provided we are willing to waste resources, there is a simple solution to this problem, a Greedy Task Allocation (\algname{GTA}) strategy, which follows this principle:
{\em Make sure all workers are always busy working on some task, and stop once $B$ tasks have been completed.}
In \algname{GTA}, we initially ask all $n$ workers to start working on a task, and as soon as some worker is done with a task, we ask it to start completing another task.
This process is repeated until $B$ tasks have been completed.
\subsection{Wastefulness of GTA}
While \algname{GTA} minimizes the completion time, it can be immensely wasteful in terms of the total worker utilization time needed to collect all $B$ tasks.
Indeed, consider the scenario with $n=1{,}010$ workers and $B=10$ tasks.
In this case, \algname{GTA} will lead to at least $n-B = 1{,}000$ unnecessary tasks being run in each iteration!
This is highly undesirable in situations where the workers are utilized across multiple other jobs besides running the iterative algorithm mentioned above.
\subsection{Goal of this work}
The goal of our work is to design new task allocation strategies, with rigorous theoretical support, that would attempt to minimize the expected completion time subject to the {\em constraint} that such wastefulness is completely eliminated.
That is, we ensure that no more than $B$ tasks are completed in each round.
\subsection{A Motivating Example: Optimal Parallel SGD}
A key inspiration for our work, and the prime example of the general task collection problem described above, relates to recent development in the area of parallel stochastic gradient descent (\algname{SGD}) methods.
Consider the problem of finding an approximate stationary point of the optimization problem
$$
    \minimize\limits_{x \in \R^d} \left\{f(x) \eqdef \ExpSub{\xi \sim {\cD}}{f_{\xi}(x)} \right\},
$$
where $f_{\xi}:\R^d\to \R$ are smooth nonconvex functions, and $f$ is assumed to be bounded from below.
We assume that
$$
    \ExpSub{\xi \sim {\cD}}{\sqnorm{f_{\xi}(x) - \nabla f(x)}} \leq \sigma^2
$$
for all $x\in \R^d$.

In a recent breakthrough, \citet{tyurin2024optimal} developed a parallel \algname{SGD} method, {\em optimal} in terms of a novel notion of complexity called {\em time complexity}, for solving the above problem with $n$ parallel workers, assuming that it takes $\tau_i>0$ seconds to worker $i$ to compute a stochastic gradient of $f$ (this corresponds to a task).
Their method, \rennala, corresponds to \algname{Minibatch SGD} of minibatch size $B$ (which depends on the target accuracy and $\sigma$ only), with the $B$ tasks (stochastic gradients) completed via \algname{GTA}.
While minimax optimal in terms of time complexity, the \algname{GTA} task allocation strategy employed within \rennala can be wasteful, as explained above.

Recently, \citet{maranjyan2025ringmaster} proposed \ringmaster, a fully asynchronous \algname{SGD} method, matching the optimal time complexity of \rennala and achieving optimality for arbitrary compute time patterns associated with the tasks (stochastic gradients), including random, as considered in our setup.
However, \ringmaster also employs a greedy task allocation strategy, leading to wastefulness.

Numerous other parallel/distributed methods involve the implementation of a task allocation strategy, including stochastic proximal point methods (task = evaluation of the stochastic prox operator), higher-order methods (task = evaluation of stochastic Hessian), and beyond.
So, by addressing the general task allocation problem, we aim to tame the inherent resource wastefulness of all these methods.
\subsection{Contributions}
In this work, we formalize the task allocation problem as a \emph{combinatorial online learning problem with partial feedback and non-linear losses}.
Then, we introduce \ata, a lower-confidence bound-based algorithm designed to solve the proposed allocation problem.
\ata is agnostic to workers' computation times, and our theoretical analysis demonstrates that the total computation time achieved by our methods remains within a small multiplicative factor of the optimal computation time (i.e., the one attainable with full knowledge of the workers' arm distributions).
Additionally, we present \algname{ATA-Empirical}, a variant of \ata that leverages a novel data-dependent concentration inequality and achieves better empirical results.
Finally, we validate our approach through numerical simulations.
\section{Related Work}
\label{ata:section:related}

Most of the literature on asynchronous methods focuses on demonstrating advantages over their synchronous counterparts.
For the simplest method, \algname{SGD}, this was only recently established by \citet{tyurin2024optimal}.
With this result in place, the community can now shift its focus to reducing the overhead of asynchrony.
Our work may be the first step in this direction.

In federated learning (FL) \citep{konevcny2016federated,mcmahan2016federated,kairouz2021advances}, several works account for system heterogeneity.
The most well-known FL method, \algname{FedAvg} \citep{mcmahan2017communication}, operates by performing multiple local steps on workers, where each step can be viewed as a task.
Some works adjust the number of local steps based on worker computation times \citep{li2020federated, maranjyan2025gradskip}, effectively adapting task assignments to worker speed.
However, these methods rely on prior knowledge of these times rather than learning them adaptively, as we do.

We reformulate our problem as an online bandit problem.
The literature on bandit algorithms is vast, and we refer the reader to the work by \citet{LattimoreS18} for an introduction to this subject.
Our algorithm is based on the approach of using Lower Confidence Bounds (LCBs) on the true means of the arms.
This idea, originally proposed by \citet{auer2002finite} for the classical Multi-Armed Bandit (MAB) setting, has since been widely adopted in the stochastic combinatorial bandits literature \citep{gai2012combinatorial, chen2013combinatorial, combes2015combinatorial, kveton2015tight}.
Using LCBs instead of the empirical estimates of the means allows an optimal trade-off between exploration and exploitation.

The ``greedy'' approach we employ, which involves selecting the action that minimizes the loss function based on lower confidence bounds instead of the unknown means, is a standard technique in the literature \citep{chen2013combinatorial, lin2015stochastic}.
However, note that our larger action space and the discontinuity of our loss function necessitates a more tailored analysis.
To the best of our knowledge, this is the first work addressing a non-continuous loss function in a stochastic combinatorial MAB-like framework.
To overcome this challenge, we exploit the specific structures of our loss function and action space to control the number of rounds where suboptimal actions are chosen.
Additionally, our procedure is computationally efficient.
\section{Problem Setup}\label{ata:sec:setup}
In this section, we formally describe the problem setup.
\subsection{Task Allocation Protocol}
We consider a system of $n$ workers, each responsible for repeatedly performing the same task (e.g., computing a gradient).
In each round, the allocation algorithm has a budget of $B$ units, which must be distributed among the $n$ workers.
Allocating one unit corresponds to one execution of the task.

Let $K$ denote the total number of rounds, which is assumed to be unknown to the learner.
For worker $i \in [n] := {1,2,\dots,n}$, let $X_i^{k,u}$ be the computation time required in round $k \in [K]$ to complete its $u$-th task.
Thus, the total computation time for worker $i$ to perform $a_i^k$ tasks in round $k$ is
$$
	\sum_{u=1}^{a_i^k} X_i^{k,u}
$$
with the convention that this sum equals $0$ when $a_i^k = 0$.

At each round $k$, the algorithm selects an allocation vector $a^k = (a_1^k, \ldots, a_n^k) \in \Z_{+}^n$ such that $\norm{a^k}_1 = \sum_{i=1}^n a_i^k = B$, based only on the information available up to round $k-1$.
The feedback consists of $a_i^k$ observed times for all the chosen workers.
We will denote the action set by 
$$
	\cA := \left\{ a \in \Z_{+}^n : \norm{a}_1 = B \right\}~.
$$
The objective of the allocation strategy in each round $k$ is to minimize the total computation time.
Formally, we define 
$$
    C : \cA \to \R_+~,
$$
as the computation time the optimizer must wait to receive $B$ completed tasks (e.g., gradients) under an allocation vector $a \in \cA$, given by
\begin{equation}
    \label{ata:eq:def_c}
	C(a^k)
	:= \max_{i\in \text{supp}(a^k)} \  \sum_{u=1}^{a_i^k} X_i^{k,u}~.
\end{equation}
\subsection{Modeling Assumptions}\label{ata:modeling_assumptions}
We assume that the computation time of each worker $i \in [n]$ are i.i.d. drawn from a random variable $X_i$ following a probability distribution $\nu_i$.
We denote by $\mu = (\mu_1, \dots, \mu_n)$ the vector of unknown means.
Hence, the random variables $X_i^{k,u}$ with $u \in \{1, \dots, a_i^k\}$ are $a_i^k$ i.i.d. samples drawn from $\nu_i$.

We assume that the distribution $\nu_i$ of the computation times to be sub-exponential random variables.
To quantify this assumption, we recall the definition of the sub-exponential norm, also known as the Orlicz norm, for a centered real-valued random variable $X$:
\begin{equation}\label{ata:eq:def_se}
	\norm{X}_{\psi_1} := \inf\left\{C > 0: \E{\exp\( \frac{\abs{X}}{C} \)} \le 2\right\}~.
\end{equation}
Hence, formally we make the following assumption.
\begin{boxedassumption}
    \label{ata:a:sube}
	For each $i \in [n]$, $X_i$ is a positive random variable with 
    $$
        \|X_i - \mu_i\|_{\psi_1} \le \alpha_i~, \quad \alpha_i \ge 0 ~.
    $$
\end{boxedassumption}
In the remainder of this paper we denote $\alpha \eqdef \max_{i \in [n]} \alpha_i$. 

The considered class encompasses several other well-known classes of distributions in the literature, such as support-bounded and sub-Gaussian distributions.
Moreover, it includes exponential distributions, which are frequently used in the literature to model waiting or computation times in queuing theory and resource allocation in large distributed systems \citep{gelenbe2010analysis, gross2011fundamentals,hadjis2016omnivore, mitliagkas2016asynchrony, dutta2018slow, nguyen2022federated}.

\subsection{Objective of the Allocation Algorithm}
The main objective of this work is to develop an online allocation strategy that achieves a small expected total computation time, defined as 
$$
	\cC_K := \sum_{k=1}^{K} \E{C(a^k)} ~.
$$
If the distributions of the arms were known in advance, the optimal allocation
$$
	a^* \in \argmin_{a \in \cA} \ \E{C(a)}
$$
would be selected to minimize the expected computation time per round, and the same allocation would be used consistently across all $K$ rounds.
The corresponding optimal total computation time is then
$$
	\cC_K^{*} := K \E{C(a^*)}.
$$
Our goal is to design a strategy that ensures the computation time $\cC_K$ remains within a small multiplicative factor of the optimal time $\cC^*_K$, up to an additional negligible term.
Specifically, we aim to satisfy
\begin{equation}
    \label{ata:eq:ata_obj}
	\cC_K \leq c \cdot \cC^*_K + \cE_K~,
\end{equation}
where $c \geq 1$ is a constant close to $1$, and $\cE_K$ is a negligible term compared to $\cC^*_K$ when $K\to \infty$.
This would assure us that in the limit we are a constant multiplicative factor away from the performance of the optimal allocation strategy that has full knowledge of the distributions of the computational times of the workers.

Finding a strategy solving the objective in \eqref{ata:eq:ata_obj} presents several technical challenges.
First, the action space $\cA$ is discrete, and the nonlinearity of the computation time functions $C(\cdot)$ prevents reducing our objective to a convex problem.
Second, the size of $\cA$ is combinatorial, growing on the order of $\binom{n + B - 1}{B}$, which necessitates exploiting the inherent problem structure to develop efficient strategies.
Third, because the workers' computation times are stochastic, any solution must account for uncertainty.
Finally, the online setting forces the learner to balance exploration and exploitation under a limited allocation budget of $B$ units per round and partial feedback\textemdash only the computation times of workers who receive allocations are observed.
This last point naturally suggests adopting a MAB approach.

In the next section, we show how to reduce this problem to a MAB problem and how to efficiently solve it.
\section{Adaptive Task Allocation}
\label{ata:section:ata}
Here, we first show how to reduce the problem in \eqref{ata:eq:ata_obj} to a \emph{non-linear} stochastic MAB problem.
Then, we propose an efficient algorithm for this formulation.
\subsection{Reduction to Multi-Armed Bandit and Proxy Loss}\label{ata:sec:reduction}
The stochastic MAB problem is a fundamental framework in sequential decision-making under uncertainty.
It involves a scenario where an agent must choose among a set of arms, each associated with an unknown reward distribution.
The agent aims to maximize cumulative reward (or equivalently minimize the cumulative loss) over time by balancing exploration (gathering information about the reward distributions) and exploitation (leveraging the best-known arm).
The challenge lies in the trade-off between exploring suboptimal arms to refine reward estimates and exploiting the arm with the highest observed reward, given the stochastic nature of the outcomes.
Using the terminology from bandit literature, here we will refer to each worker as an ``arm''.

However, differently from the standard MAB problem, we have a harder problem because $\mathbb{E}[C(a^k)]$ depends on the joint distribution of all the arms in the support of $a^k$, rather than on their expectations only.
This dependency potentially renders the task of relying on estimates of $\mathbb{E}[C(a)]$ for $a \in \cA$ computationally challenging due to the combinatorial nature of the set $\cA$.

To solve this issue, our first idea is to introduce a \emph{proxy loss} $\ell : \cA\times \R_+^n \to \R_+$, defined as
\begin{equation}
    \label{ata:eq:def_l}
	\ell(a,\mu) \eqdef  \max_{i \in [n]} \  a_{i} \mu_i~.
\end{equation}
Due to the convexity of $C(\cdot)$, the introduced proxy-loss underestimates the expected computation time.
However, in \Cref{ata:sec:proof_2} we prove that this quantity also upper bounds the expected computation time up to a constant that depends on the distribution of the arms.
In particular, for any $a \in \cA$, we show that
\begin{equation}
    \label{ata:eq:enc}
	\ell(a, \mu) \le \E{C(a)} \le (1+4\eta \ln(B)) \ell(a, \mu)~,
\end{equation}
where $\eta$ is defined as
\begin{equation}
    \label{ata:def:eta}
	\eta := \max_{i \in [n]} \frac{\alpha_i}{\mu_i}~.
\end{equation}
In words, $\eta$ provides an upper bound on the ratio between the standard deviation and the mean of the arms.
Note that in the literature, it is common to consider exponential, Erlang, or Gamma distributions, where the ratio $\eta$ is typically\footnote{For $\mathrm{Gamma}(\alpha, \lambda)$, $\sigma / \mu = 1 / \sqrt{\alpha}$, so the claim holds for $\alpha \geq 1$.} bounded by $1$.

The bound above will allow us to derive guarantees on the total computation time of an allocation strategy based on its guarantees for the proxy loss $ \ell(\cdot)$, up to a factor of the order $1 + 4\eta \ln(B)$.
We remark that in the special case where the arms' distributions are deterministic ($\eta = 0$) or the query budget is unitary ($B = 1$), the two targets $\mathbb{E}[C(a)]$ and $\ell$ exactly coincide.
\subsection{Comparison with the Combinatorial Bandits Setting}
Our setting is closely related to the Combinatorial Multi-Armed Bandits (CMAB) framework \citep{cesa2012combinatorial}, particularly due to the combinatorial nature of the action space and the semi-bandit feedback, where the learner observes outcomes from all chosen arms.
However, our formulation differs in two significant ways.
First, while CMAB typically involves selecting a subset of $n$ arms, resulting in an action space with a maximum size of $2^n$, our action space $\cA$ has a cardinality of $\binom{n + B - 1}{B}$.
The ratio between these two can be extremely large, potentially growing exponentially with $n$.
Second, although most works in this domain assume a linear loss function in the arms' means, some notable exceptions address non-linear reward functions \citep{chen2013combinatorial, lin2015stochastic, chen2016combinatorial, wang2018thompson}.
However, these approaches generally rely on assumptions such as smoothness, Lipschitz continuity, or higher-order differentiability of the reward function.
In contrast, our loss function $\ell(\cdot, \mu)$ is not continuous with respect to the arms' means.
Finally, motivated by the practical requirements of our setting, we place a strong emphasis on computational efficiency that rules out most of the approaches based on CMAB.

\subsection{Adaptive Task Allocation Algorithm}
Now, we introduce our Adaptive Task Allocation algorithm (\ata).
\ata does not require prior knowledge of the horizon $K$ and only relies on an upper bound $\alpha$ satisfying $\alpha \ge \max_{i \in [n]} \norm{X_i - \mu_i}_{\psi_1}$ for the Orlicz norms of the arm distributions. Recall that $\norm{X_i - \mu_i}_{\psi_1} \le 2 \norm{X_i}_{\psi_1}$, so an upper bound on $\norm{X_i}_{\psi_1}$ also provides one for $\|X_i - \mu_i\|_{\psi_1}$.
The core idea of the procedure is to allocate the workers based on \emph{lower confidence bound estimates} on the arm means $(\mu_i)_{i \in [n]}$, in order to balance exploration and exploitation.

For each arm $i \in [n]$ and round $k \in [K]$, let $K_i^k$ represent the number of samples collected from the distribution of arm $i$ up to round $k$.
At each round $k$, we compute an empirical mean, denoted by $\hat{\mu}_i^k$, using the $K_i^k$ samples obtained so far.
Based on these empirical means, we define the lower confidence bounds $s_i^k$ as
\begin{equation}
\label{ata:eq:s}
	s_i^k = \left(\hat{\mu}_i^k - \text{conf}(i,k) \right)_{+}~,
\end{equation}
where $(x)_{+} = \max\{x, 0\}$ and $\text{conf}(\cdot, \cdot)$ is defined as
\begin{equation*}
	\conf(i, k) =
	\begin{cases} 
		2\alpha\(\sqrt{\frac{\ln(2k^2)}{K_i^k}}+\frac{\ln(2k^2)}{K_i^k}\), & K_i^k \geq 1~, \\
		+\infty, &  K_i^k = 0~.
	\end{cases}
\end{equation*}
The term $\text{conf}(\cdot, \cdot)$ is derived from a known concentration inequality for sub-exponential variables with an Orlicz norm bounded by $\alpha$ (\Cref{ata:prop:concentration} in the Appendix).

Given the confidence bounds $s^k := (s_1^k, \dots, s_n^k)$, the learner selects the action $a^k \in \cA$ at round $k$ that minimizes the loss $\ell(\cdot, s^k)$, defined in \eqref{ata:eq:def_l}.
While nonconvex, we show in \Cref{ata:sec:RAS} that this optimization problem can be solved using a recursive routine, whose computational efficiency is
$$
    \cO\(n \ln\(\min\{B, n\}\) + \min\{B, n\}^2\) ~.
$$
\begin{algorithm}[t]
	\caption{\ata (Adaptive Task Allocation)}
    \label{ata:alg:ata}
	\begin{algorithmic}[1]
        \STATE \textbf{Input}: allocation budget $B$, $\alpha>0$
        \STATE \textbf{Initialize}: empirical means $\hmu_i^1 = 0$, usage counts $K_i^1 = 0$, and usage times $T_i^1 = 0$, for all $i \in [n]$
		\FOR{$k = 1,\ldots, K$}
            \STATE Compute LCBs $(s_i^k)$ for all $i \in [n]$ using \eqref{ata:eq:s}
            \STATE Find allocation:
            $$
                a^k \in  \argmin_{a\in \cA} \ell(a, s^k)
            $$
            \STATE Allocate $a_i^k$ tasks to each worker $i \in [n]$
            \STATE {\highlightcolor Update optimization parameters}
            \STATE For all $i$ such that $a_i^k \neq 0$, update:
            \begin{align*}
                K_i^{k+1} &= K_i^k + a_i^k \\
                T_i^{k+1} &= T_i^k + \sum_{j=1}^{a_i^k} X_i^{k,j} \\
                \hmu_i^{k+1} &= \frac{T_i^{k+1}}{K_i^{k+1}}
            \end{align*}
		\ENDFOR
	\end{algorithmic}
\end{algorithm}
As last step, the feedback obtained after applying the allocation $a^k$ is used to update the lower confidence bounds.
The complete pseudocode for \ata is provided in \Cref{ata:alg:ata}.
\begin{remark}
	Line 7 of the algorithm acts as a placeholder for the optimization method, where the optimization parameters are updated using the quantities computed by the workers (e.g., gradients in the case of \algname{SGD}).
    In this view, the allocation algorithm is independent of the specifics of the chosen optimization algorithm. Refer to \Cref{ata:section:other_methods} for further details.
\end{remark}
\subsection{Upper-Bound on the Total Computation Time}
We provide guarantees for \ata in the form of an upper bound on the expected total computation time required to perform $K$ iterations of the optimization procedure.
Recall that the proxy loss $\ell(\cdot, \mu)$ and the expected computation time are related through \eqref{ata:eq:enc}.
This relationship and \Cref{ata:thm:main} allow us to derive guarantees on the expected total computation time, denoted by 
$$
	\cC_K := \sum_{k=1}^{K} \E{C(a^k)} ~.
$$
Recall that the optimal expected total computation time in this framework is given by 
$$
	\cC_K^{*} := K \E{C(a^*)} ~,
$$
with 
$$
	a^* \in \argmin_{a \in \cA} \ \E{C(a)} ~.
$$
We now state the result that provides an upper bound on the total expected computation time achieved by \ata.
\begin{boxedtheorem}[Proof in \Cref{ata:sec:proof_2}]
    \label{ata:cor:main}
	Suppose \Cref{ata:a:sube} holds and let $\eta := \max_{i \in [n]} \nicefrac{\alpha_i}{\mu_i}$.
	Then, the total expected computation time after $K$ rounds, using the allocation prescribed by \ata with inputs $(B, \alpha)$ satisfies
	$$
		\cC_K \le \left(1+4\eta~\ln(B) \right)\cC_K^* + \cO(\ln K) ~.
	$$
\end{boxedtheorem}
\begin{remark}
	The $\cO(\cdot)$ term hides an instance dependent factor. We will give its full specifics in the regret upper bound of \Cref{ata:thm:main}.
\end{remark}
The bound in \Cref{ata:cor:main} shows that the total expected computation time of \ata remains within a multiplicative factor of $1 + 4\eta\ln(B)$ of the optimal computation time $\cC_K^*$, with an additional remainder term that scales logarithmically with $K$. 
Since $\cC_K = \Omega(K)$, this additive term is negligible compared to $\cC^*_K$.
In practical scenarios, where computation time follows common distributions such as exponential or Gamma, the factor $\eta$ is typically of order $1$, and $\ln(B)$ remains relatively small for the batch sizes commonly used in optimization algorithms like \algname{SGD}.

The reader might wonder if the more ambitious goal of deriving bounds with a multiplicative factor of exactly $1$ is achievable.
However, achieving this goal would require significantly more precise estimates of the expected computation time $\mathbb{E}[C(a)]$ for all $a \in \cA$.
Since $\mathbb{E}[C(a)]$ depends on the joint distribution of all workers in the support $a$, obtaining such precise estimates would come at the cost of computational efficiency in the allocation strategy.

We note that it is unsurprising that $\eta$ appears in the upper bound of \Cref{ata:cor:main}, since having a heavier-tailed distribution increases the gap between $\ell(a, \mu)$ and $\mathbb{E}[C(a)]$ through the convexity of $C(\cdot)$.
Instead, the factor $\ln(B)$ arises because $C(\cdot)$ is expressed as the maximum of up to $B$ random variables.
Moreover, in the edge cases where $\eta = 0$ (deterministic case) or $B=1$ (linear cost function), we guarantee that the expected computation time is at most an \emph{additive} factor away from the optimal one.
\section{Empirical Adaptive Task Allocation}\label{ata:sec:ata-em}
The \ata procedure is based on a lower confidence bound approach that relies on concentration inequalities.
These bounds play a key role in performance, as sharper concentration bounds lead to more accurate estimates and reduce exploration of suboptimal options.
Since workers' computation times follow sub-exponential distributions, their concentration behavior is determined by the Orlicz norm of the corresponding variables.
In \ata, the only prior knowledge available is an upper bound on the largest Orlicz norm among all arms.
When the Orlicz norms of the arms' distributions vary significantly, this uniform bound may result in loose confidence intervals and inefficient exploration. 

To address this issue, we introduce \algname{ATA-Empirical}, which better adapts to the distribution of each arm, particularly its Orlicz norm.
This adaptation is achieved through a novel data-dependent concentration inequality for sub-exponential variables.
Unlike \ata, which depends on the maximum Orlicz norm, \algname{ATA-Empirical} accounts for the individual Orlicz norms of all arms, denoted by $(\alpha_i)_{i\in [n]}$.
This improvement is reflected in the upper bounds on regret presented in \Cref{ata:section:theory}.
In practice, this leads to improved performance at least some settings, as shown in our simulations in \Cref{ata:section:experiments}.
However, this increased adaptivity comes with a trade-off since \algname{ATA-Empirical} requires an upper bound on the quantity $\eta = \max_i \nicefrac{\alpha_i}{\mu_i}$, rather than a bound on the largest Orlicz norm.
That said, for many distributions of interest, the ratios $\nicefrac{\alpha_i}{\mu_i}$ across different arms tend to be of the same order, whereas their Orlicz norms can vary significantly.
The \algname{ATA-Empirical} procedure differs from \ata only in the lower confidence bounds it uses.
These bounds are derived from the novel concentration inequality in Lemma~\ref{ata:lem:0} and are defined for arm $i \in [n]$ at round $k \in [K]$ as
\begin{equation}\label{ata:def:s}
	\hat{s}_i^k =
	\hat{\mu}_i^k \left[ 1-2\,\eta\left(\sqrt{\frac{\ln(2k^2)}{K_i^k}}+ \frac{\ln(2k^2)}{K_i^k} \right) \right]_{+},
\end{equation}
where $\eta = \max_{i \in [n]} \nicefrac{\alpha_i}{\mu_i}$.\\

The expected total computation time $\cC_K$ of \algname{ATA-Empirical} satisfies the same guarantee presented in \Cref{ata:thm:main}, but we obtain an improved multiplicative factor of the additive logarithmic term.
The precise expressions of these factors are provided in the next section, and they show that the guarantees of \algname{ATA-Empirical} adapt to the Orlicz norms $\norm{X_i}_{\psi_1}$ of each arm, while the guarantees of \ata depend on the maximum Orlicz norm $\max_i\norm{X_i}_{\psi_1}$.
\section{Theoretical Results}
\label{ata:section:theory}

In this section, we sketch the derivation of \Cref{ata:cor:main} for \ata and \algname{ATA-Empirical}, through a regret analysis on the proxy losses.
We define the expected cumulative regret of the proxy loss $\ell(\cdot, \mu)$ after $K$ rounds
\begin{equation}
    \label{ata:eq:proxy_loss}
	\cR_K := \sum \limits_{k=1}^{K} \E{\ell(a^k, \mu)} - K \cdot \ell(\bar{a}, \mu)~,
\end{equation}
where $\bar{a} \in \argmin_{a \in \cA} \ell(a, \mu)$ represents the optimal allocation over the workers.
If multiple optimal actions exist, we consider the one returned by the optimization sub-routine used in \ata (line 5 of Algorithm \ref{ata:alg:sgd-ata}).

We derive upper bounds on the expected cumulative regret $\cR_K$.
Based on these bounds, we provide the guarantees on the expected total computation time required to complete $K$ iterations of the optimization process.
\subsection{Guarantees for ATA}
For each worker $i \in [n]$, recall that $\bar{a}_i$ denote the prescribed allocation of the optimal action $\bar{a}$.
Define $k_i$ as the smallest integer satisfying
\begin{equation}
    \label{ata:def:ki}
	(\bar{a}_i + k_i) \mu_i > \ell(\bar{a}, \mu)~.
\end{equation}
From the definition above, it follows that if the learner plays an action $a^k$ at round $k$ such that $a_i^k \ge \bar{a}_i + k_i$, then $\ell(a^k, \mu) \ge \ell(\bar{a}, \mu)$.
Thus, $k_i$ can be interpreted as the smallest number of additional units allocated to worker $i$ that result in a suboptimal loss.
Moreover, for every worker $i \in [n]$, we have $k_{i} \in \{1, 2\}$ (see Lemma~\ref{ata:lem:1} in the Appendix).

The next result provides an upper bound on the expected regret of \ata.
\begin{boxedtheorem}[Proof in \Cref{ata:proof:thm:main}]
	\label{ata:thm:main}
	Suppose that Assumption~\ref{ata:a:sube} holds.
	Then, the expected regret of \ata with inputs $(B, \alpha)$ satisfies
	\begin{align*}
		\cR_K &\le 2n\max_{i \in [n]} \{B\mu_i -\ell(\bar{a}, \mu)\} \\
		& \qquad +c\cdot\sum \limits_{i=1}^{n} \frac{\alpha^2(\bar{a}_i+k_i)(B \mu_i - \ell(\bar{a}, \mu)) }{\left((\bar{a}_i+k_i)\mu_i - \ell(\bar{a}, \mu)\right)^2}\cdot \ln K ~,
	\end{align*}
	where $\alpha \ge \max_{i \in [n]} \norm{X_i-\mu_i}_{\psi_1}$, and $c$ is a numerical constant.
\end{boxedtheorem}
The first term in the regret upper bound is independent on the number of rounds $K$.
The second term, however, grows logarithmically with $K$, which aligns with the behavior observed in stochastic bandit problems in the literature.

In the case where $B=1$, our setting reduces to the problem of regret minimization for the standard multi-armed bandits.
Observe that in this case $\ell(\bar{a}, \mu) = \min_{i\in[n]}\mu_i$, $k_i=1$ for all $i \in [n]$.
Therefore, the guarantees of \Cref{ata:thm:main} recover the known optimal bound 
$$
    \cO \( \sum_i \frac{\ln(K)}{\Delta_i} \)
$$
of the standard MAB setting, where $\Delta_i := \mu_i - \min_j \mu_j$.
\paragraph{Proof sketch.}
In standard and combinatorial MAB problems, regret bounds are typically derived by controlling the number of rounds in which the learner selects suboptimal arms.
These bounds are often of the order $\nicefrac{\ln(K)}{\Delta^2}$, where $\Delta$ denotes the suboptimality gap and quantifies the exploration cost required to distinguish optimal actions from suboptimal ones.

In our setting, the problem is more complex since the learner must not only choose which arms to pull but also determine the allocation of resources across selected arms.
With this in mind, we develop the following key arguments leading to the bound in \Cref{ata:thm:main}.  

We define \textit{over-allocation} for worker $i$ at round $k$ as the event where $a_i^k \geq \bar{a}_i + k_i$.
By definition of $k_i$ (see \eqref{ata:def:ki}), this implies that $\ell(a^k, \mu) > \ell(\bar{a}, \mu)$.
We define a \textit{bad round} as a round where $\ell(a^k, \mu) > \ell(\bar{a}, \mu)$, and we say that a bad round is \textit{triggered by arm $i$} when $a_i^k \mu_i = \ell(a^k, \mu) > \ell(\bar{a}, \mu)$.
Then, the proof revolves around establishing an upper bound on the total number of bad rounds.

To derive this bound, we consider the number of samples required to verify that the mean computation time of worker $i$ under over-allocation exceeds the optimal waiting time $\ell(\bar{a}, \mu)$.
Specifically, we need to test whether the mean of the corresponding distribution, at least $(\bar{a}_i + k_i)\mu_i$, surpasses the threshold $\ell(\bar{a}, \mu)$.
This is equivalent to testing whether 
$$
	\left\{ \mu_i > \frac{\ell(\bar{a}, \mu)}{\bar{a}_i + k_i} \right\}.
$$
Using the concentration inequality applied in our analysis, the number of samples required for this test is of the order:  
\begin{equation}
    \label{ata:eq:n_rounds}
	\alpha_i^2 \left(\mu_i - \frac{\ell(\bar{a}, \mu)}{\bar{a}_i + k_i}\right)^{-2} = \frac{\alpha_i^2 (\bar{a}_i + k_i)^2}{\left((\bar{a}_i + k_i)\mu_i - \ell(\bar{a}, \mu)\right)^2}~.
\end{equation}
During rounds where worker $i$ is over-allocated, the learner collects at least $\bar{a}_i + k_i$ samples from the corresponding distribution.
Therefore, the total number of rounds required to accumulate enough samples to stop over-allocating worker $i$ can be upper-bounded by
$$
    \frac{\alpha_i^2 (\bar{a}_i + k_i)}{\left((\bar{a}_i + k_i)\mu_i - \ell(\bar{a}, \mu)\right)^2} ~.
$$
In the regret bound of \Cref{ata:thm:main}, the term $\alpha^2$ appears instead of $\alpha_i^2$ because the learner's prior knowledge is limited to an upper bound $\alpha \ge \max_i \norm{X_i-\mu_i}_{\psi_1}$ on the maximal Orlicz norm of the arm distributions.
Finally, considering that the worst-case excess loss incurred when over-allocating worker $i$ is $B\mu_i - \ell(\bar{a}, \mu)$, we obtain the stated bound.
\subsection{Guarantees for ATA-Empirical}
We present theoretical guarantees for \algname{ATA-Empirical} by providing an upper bound on the expected cumulative regret \eqref{ata:eq:proxy_loss}.
As discussed in \Cref{ata:section:ata}, \algname{ATA-Empirical} leverages lower confidence bounds derived from a novel data-dependent concentration inequality introduced below.
The proof of this result is detailed in \Cref{ata:sec:technical}.
\begin{boxedlemma}
    \label{ata:lem:0}
	Let $X_1, \dots, X_n$ be i.i.d.\ positive random variables with mean $\mu$, such that $\alpha = \norm{X_1-\mu}_ {\psi_1}<+\infty$. Let $\hat{X}_n$ denote the empirical mean. For $\delta \in (0,1)$, let
	$$
		C_{n,\delta} := 2 \sqrt{\frac{\log(\nicefrac{2}{\delta})}{n}}+2 \frac{\log(\nicefrac{2}{\delta})}{n}~,
	$$
	where $\eta = \nicefrac{\alpha}{\mu}$.
	Then, with probability at least $1-\delta$, we have
	$$
		\mu \ge \hat{X}_n \left(1- \eta C_{n,\delta}\right)_{+}~.
	$$
	Moreover, if $\eta C_{n,\delta} \le \nicefrac{1}{4}$, then, we have with probability at least $1-\delta$, we have
	$$
		\hat{X}_n \left(1- \eta C_{n,\delta}\right)_{+} \le \mu \le \hat{X}_n \left(1+\frac{4}{3}\eta C_{n,\delta}\right)~.
	$$
\end{boxedlemma}
Using the concentration inequality above, we construct the lower confidence bounds $\hat{s}_i^k$ as defined in \eqref{ata:def:s}.
The following theorem provides an upper bound on the regret of \algname{ATA-Empirical}.
\begin{boxedtheorem}[Proof in \Cref{ata:proof:thm:main2}]
	\label{ata:thm:main2}
	Suppose that \Cref{ata:a:sube} holds.
	Then, the expected regret of \algname{ATA-Empirical} with inputs $(B, \eta)$, satisfies
	\begin{align*}
		\cR_K &\le 2n\max_{i \in [n]} \{B\mu_i -\ell(\bar{a}, \mu)\}\\
		& +c \eta^2\cdot\sum \limits_{i=1}^{n} \frac{\mu_i^2(\bar{a}_i+k_i)(B \mu_i - \ell(\bar{a}, \mu)) }{\left((\bar{a}_i+k_i)\mu_i - \ell(\bar{a}, \mu)\right)^2}\cdot \ln K~,
	\end{align*}
	where $\eta \ge \max_{i \in [n]} \nicefrac{\alpha_i}{\mu_i}$ and $c$ is a numerical constant.
\end{boxedtheorem}
Comparing the bounds for \algname{ATA-Empirical} and \ata, we observe a key differences.
Unlike the bound in \Cref{ata:thm:main}, which incurs a squared maximal Orlicz norm penalty of $\alpha^2$ for all terms in the upper bound, \algname{ATA-Empirical} benefits from its adaptive nature, leading to a term-specific factor of $\eta^2 \mu_i^2$.
In the case where the arm distributions have a ratio $\nicefrac{\alpha_i}{\mu_i}$ of the same order (such as the exponential distributions), the bound of \Cref{ata:thm:main2} shows that \algname{ATA-Empirical} adapts to the quantities $\alpha_i$ as we have, in the last case, $\eta \mu_i = \alpha_i$.

\section{Experiments}
\label{ata:section:experiments}

\begin{figure*}[t]
    \centering
    { %
    \setlength{\tabcolsep}{-4pt} 
    \renewcommand{\arraystretch}{0} 
    \begin{tabular}{@{}cccc@{}}
        \includegraphics[width=0.27\textwidth]{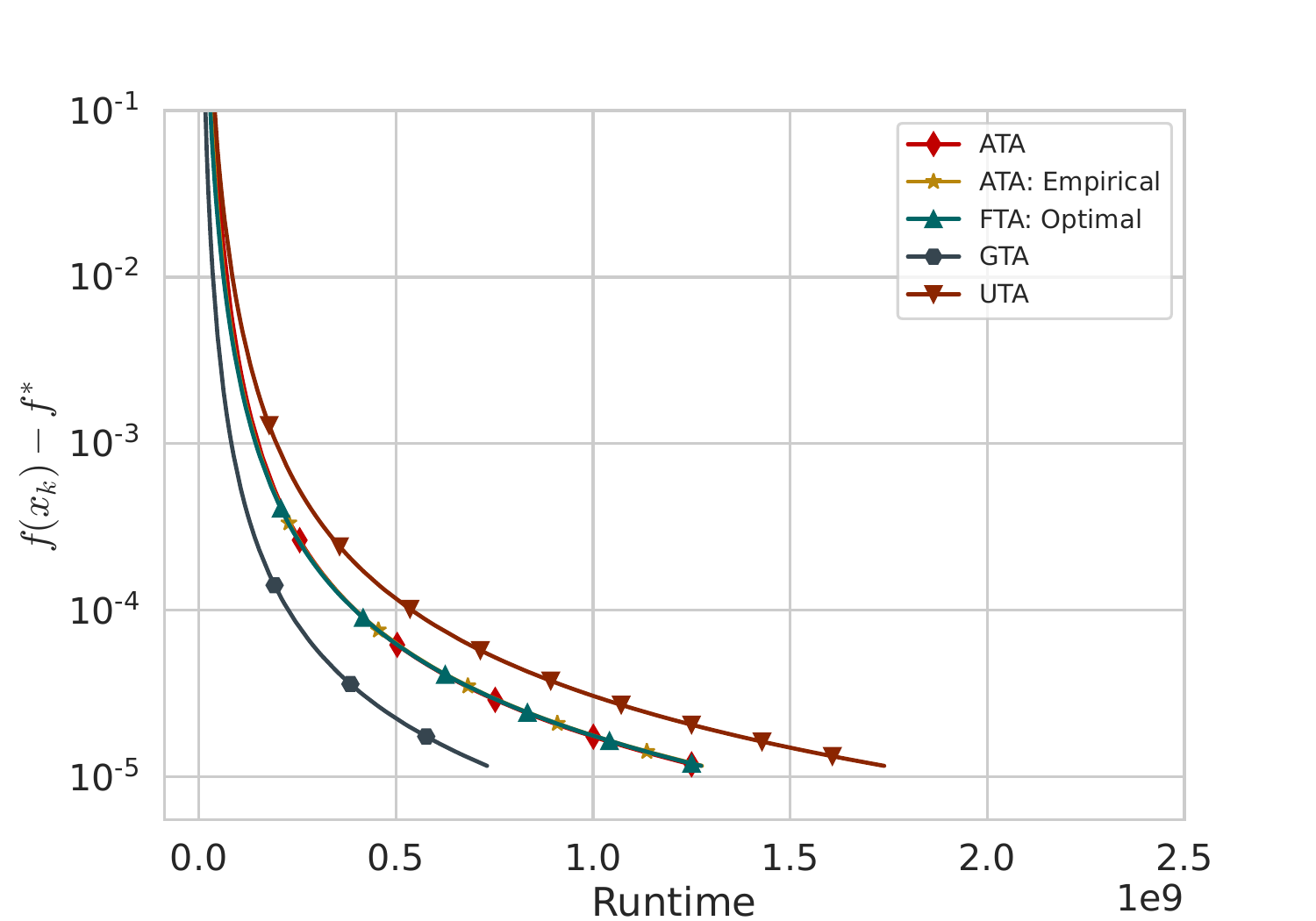} \hspace{-2mm} &
        \includegraphics[width=0.27\textwidth]{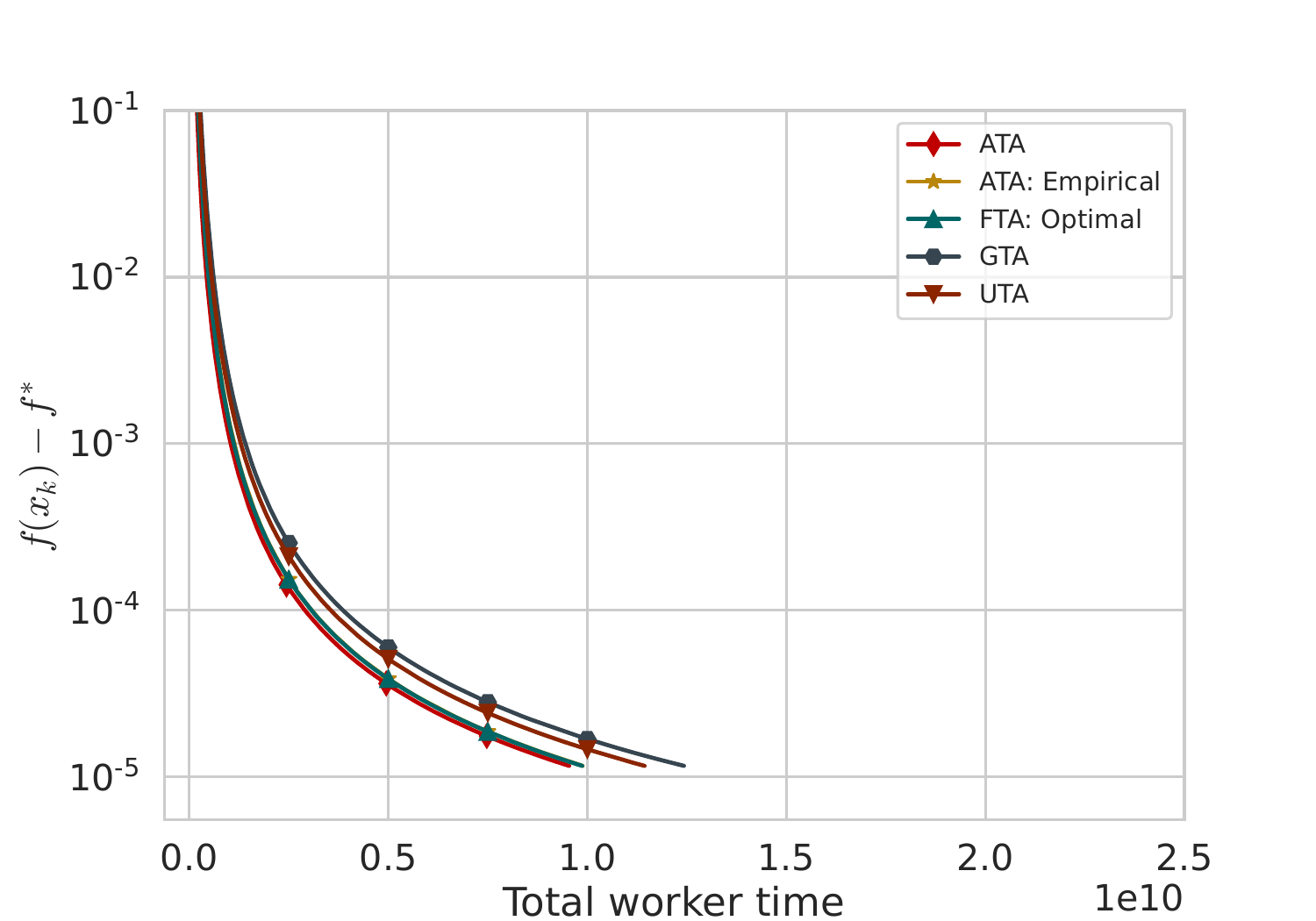} \hspace{-2mm} &
        \includegraphics[width=0.27\textwidth]{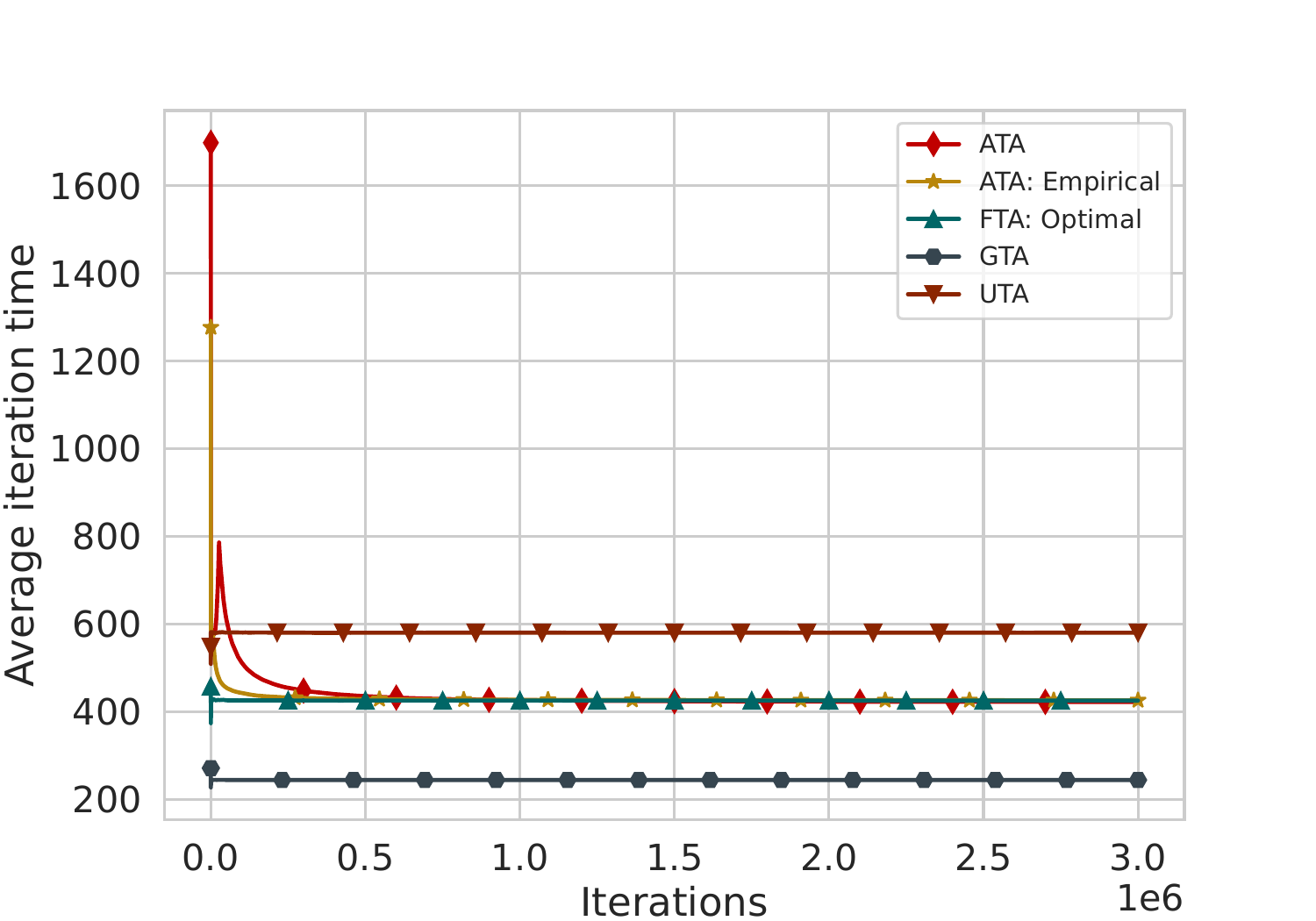} \hspace{-2mm} &
        \includegraphics[width=0.27\textwidth]{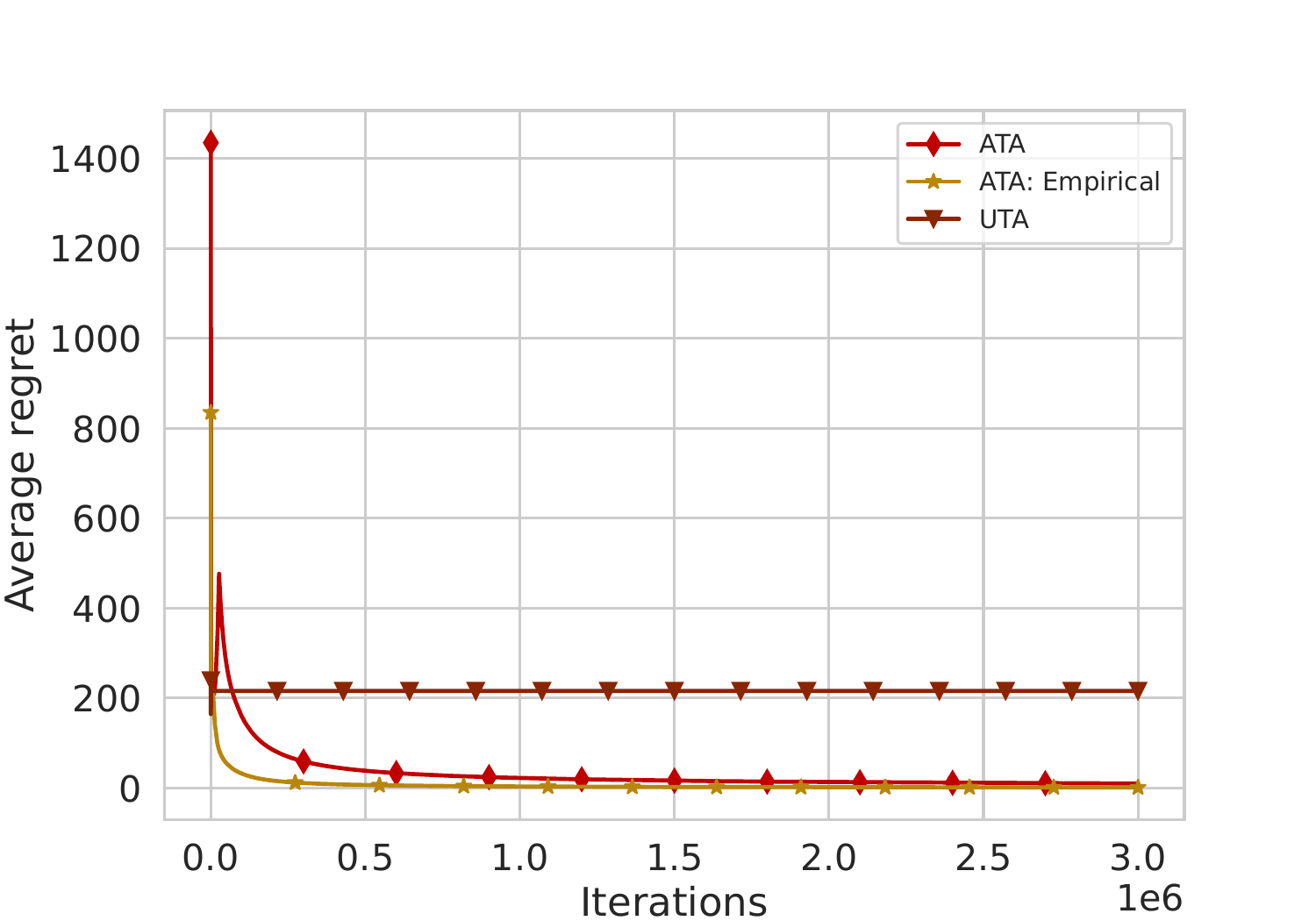} \\
        \includegraphics[width=0.27\textwidth]{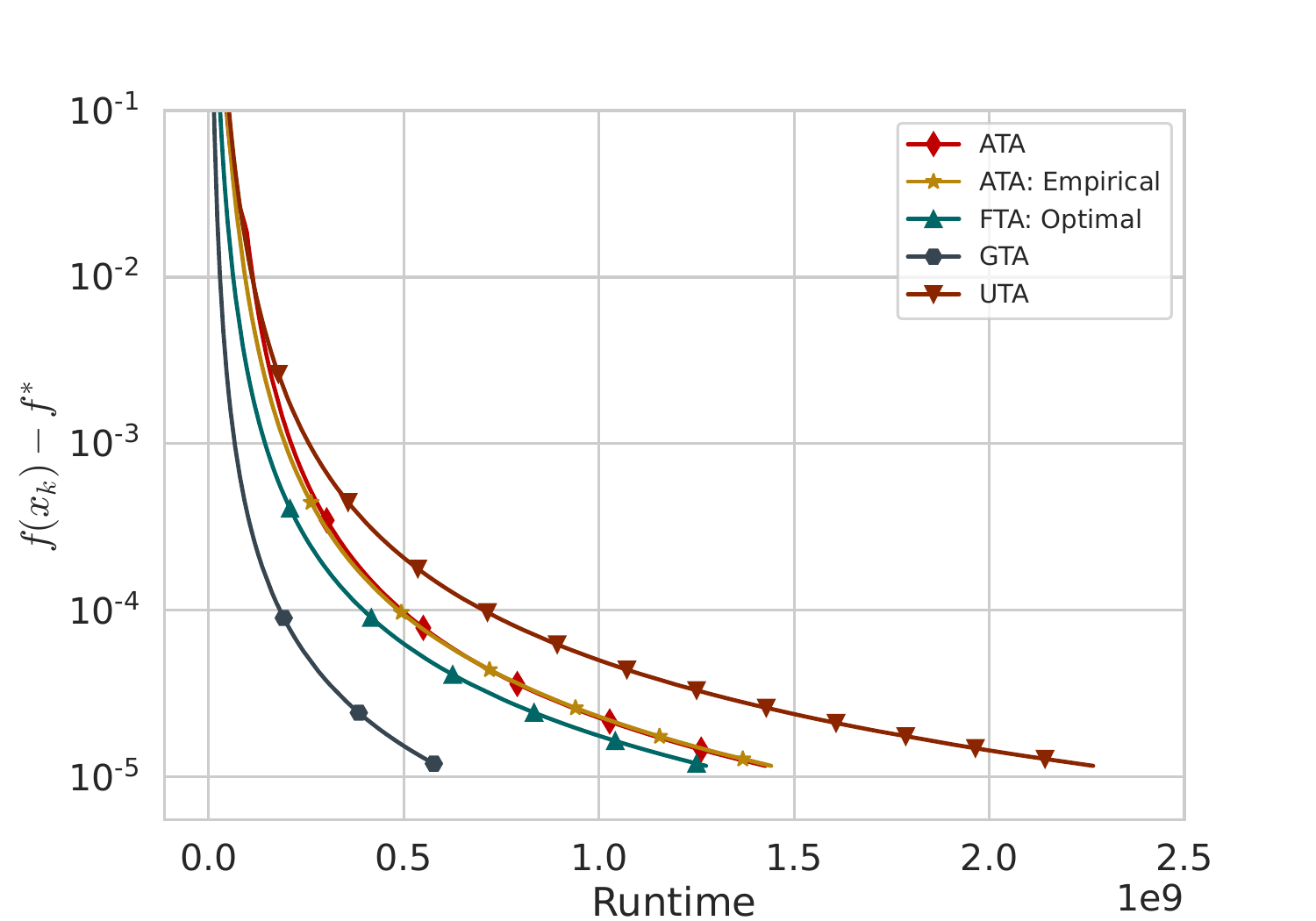} &
        \includegraphics[width=0.27\textwidth]{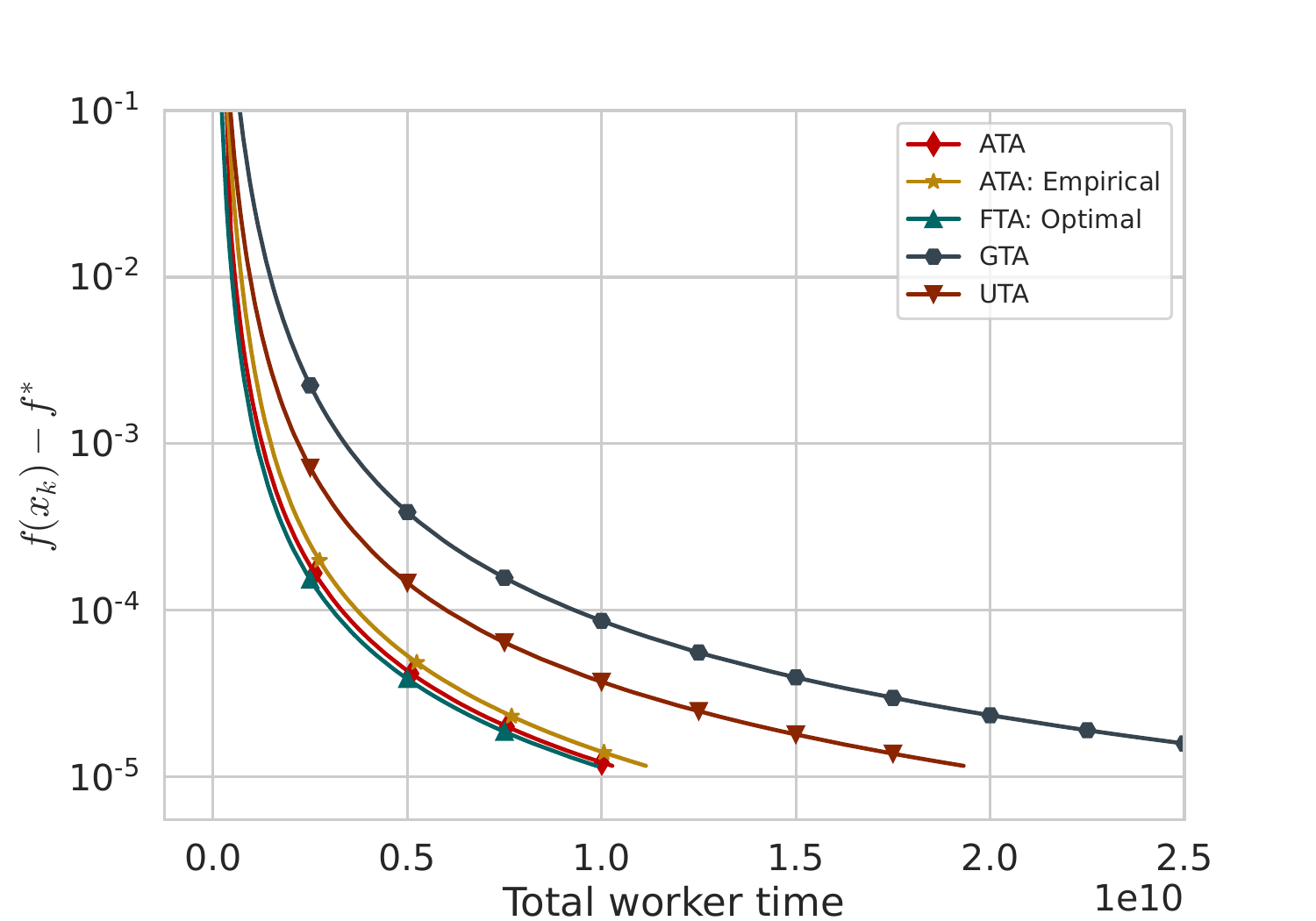} &
        \includegraphics[width=0.27\textwidth]{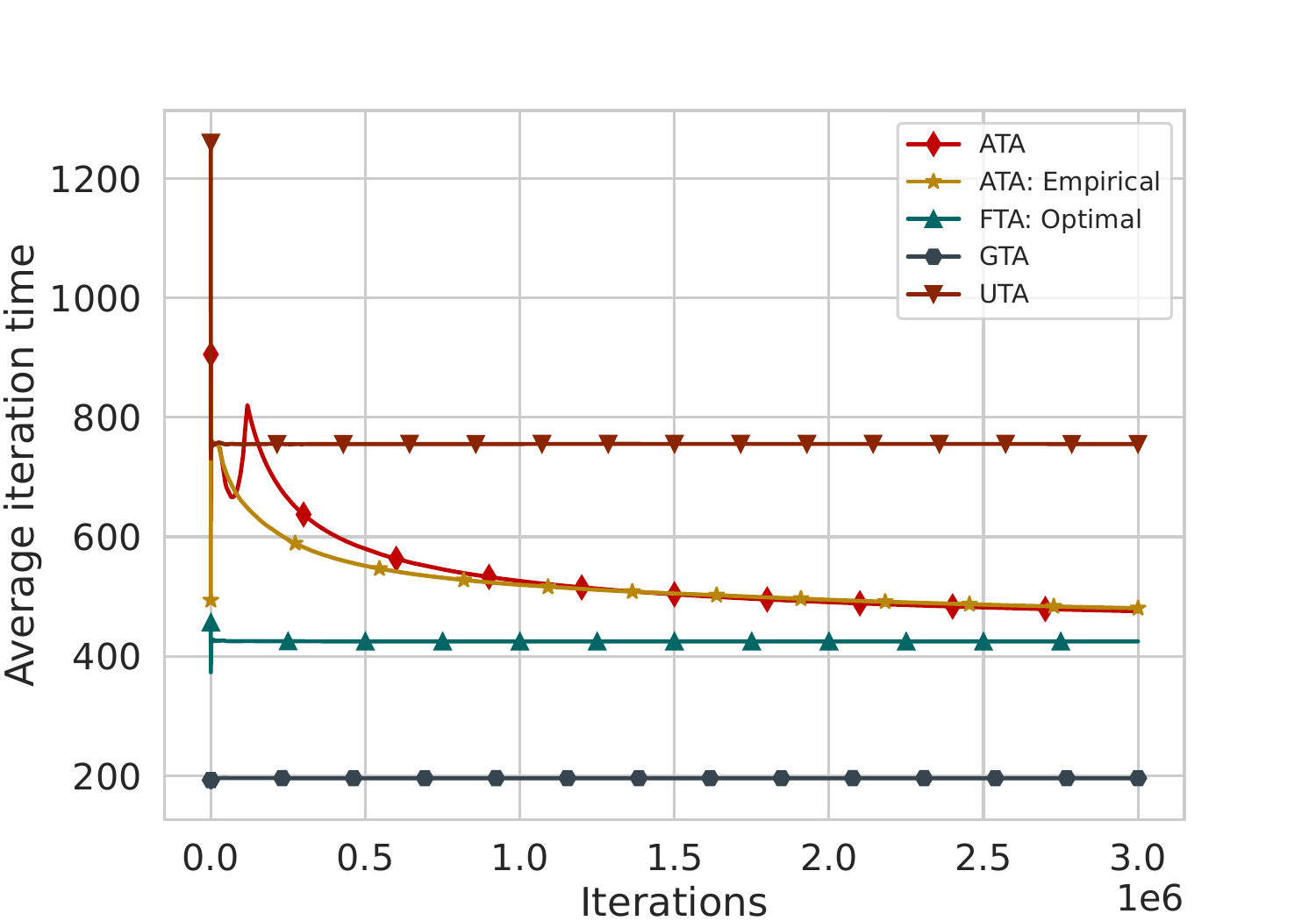} &
        \includegraphics[width=0.27\textwidth]{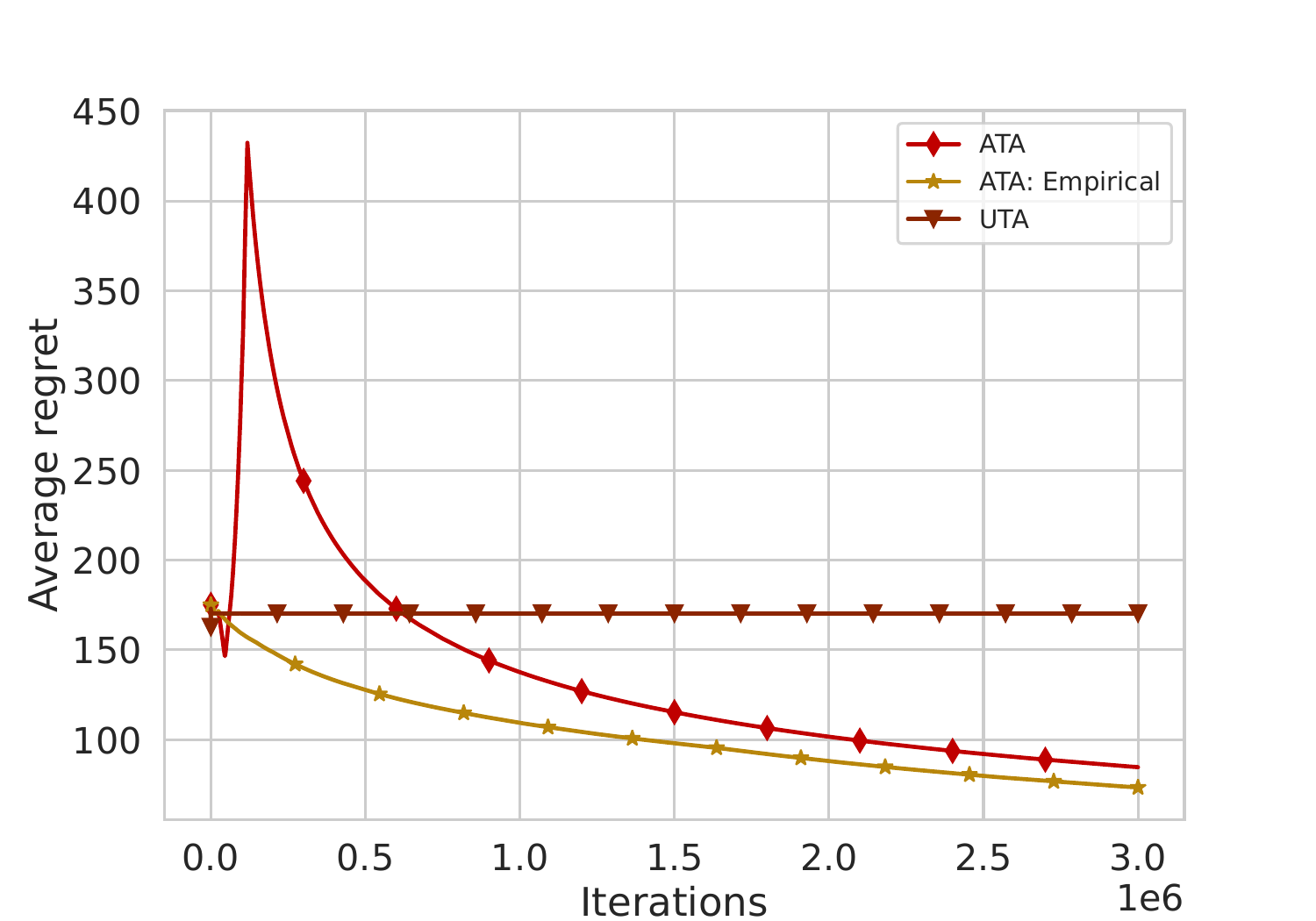} \\
        \includegraphics[width=0.27\textwidth]{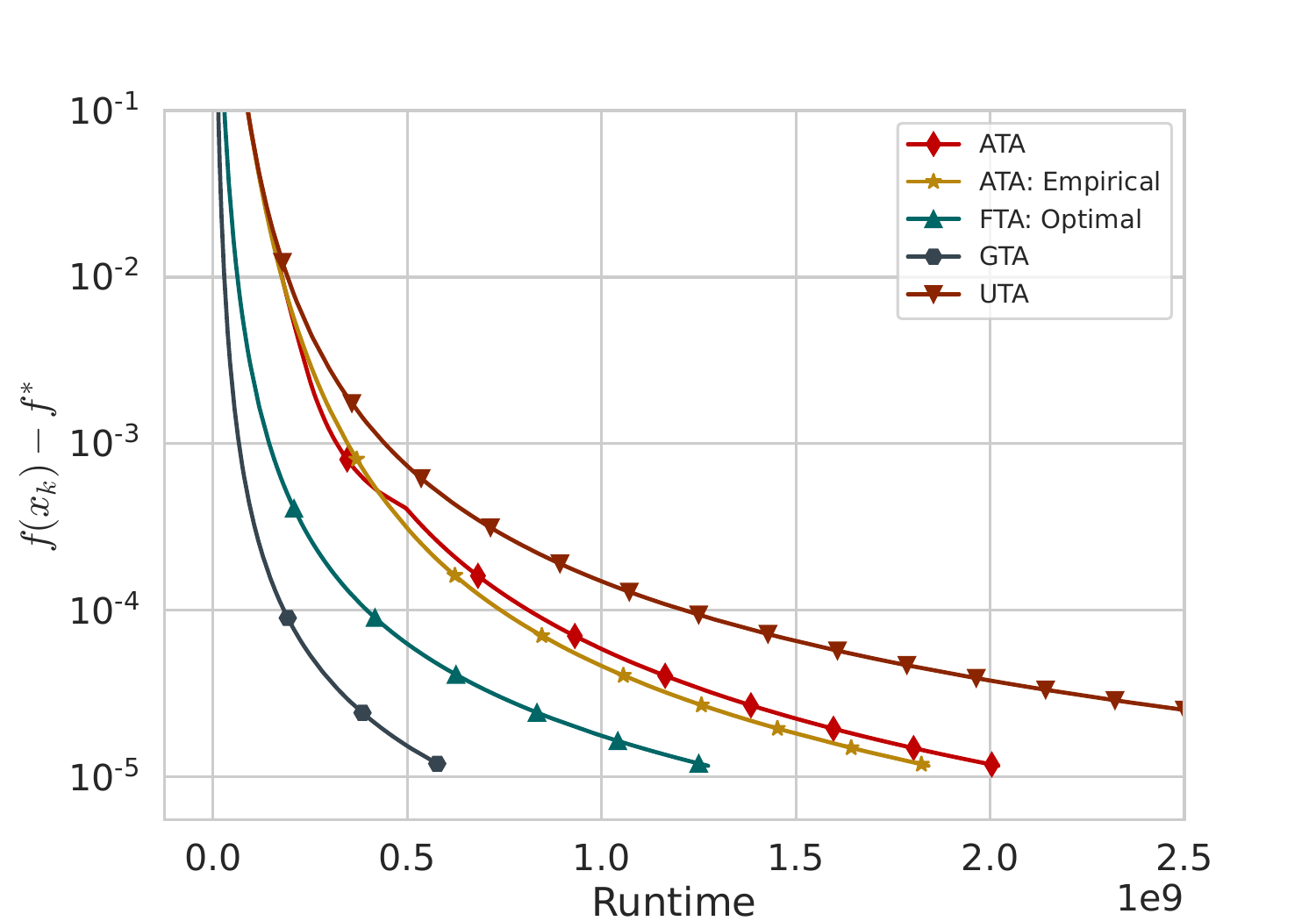} &
        \includegraphics[width=0.27\textwidth]{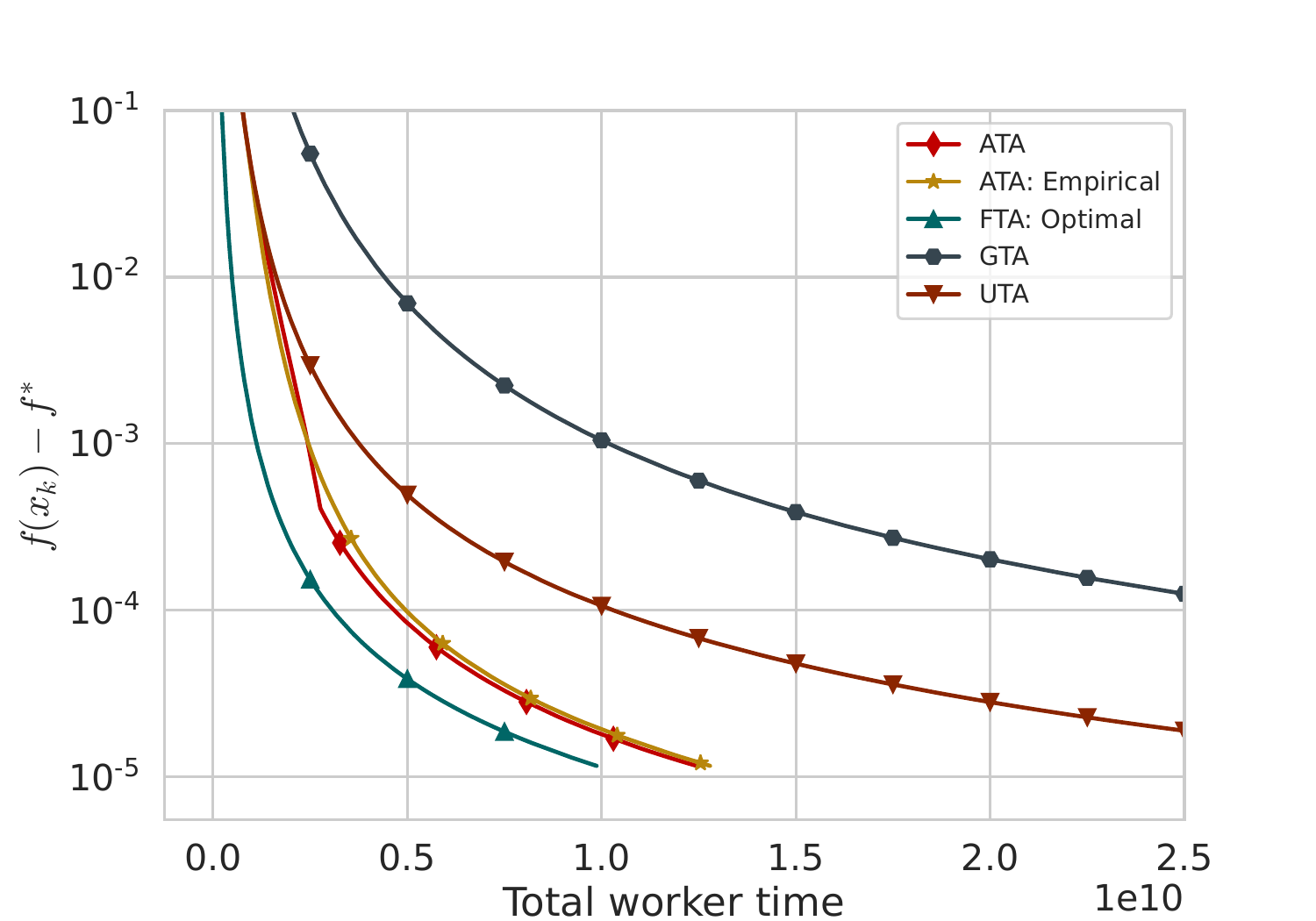} &
        \includegraphics[width=0.27\textwidth]{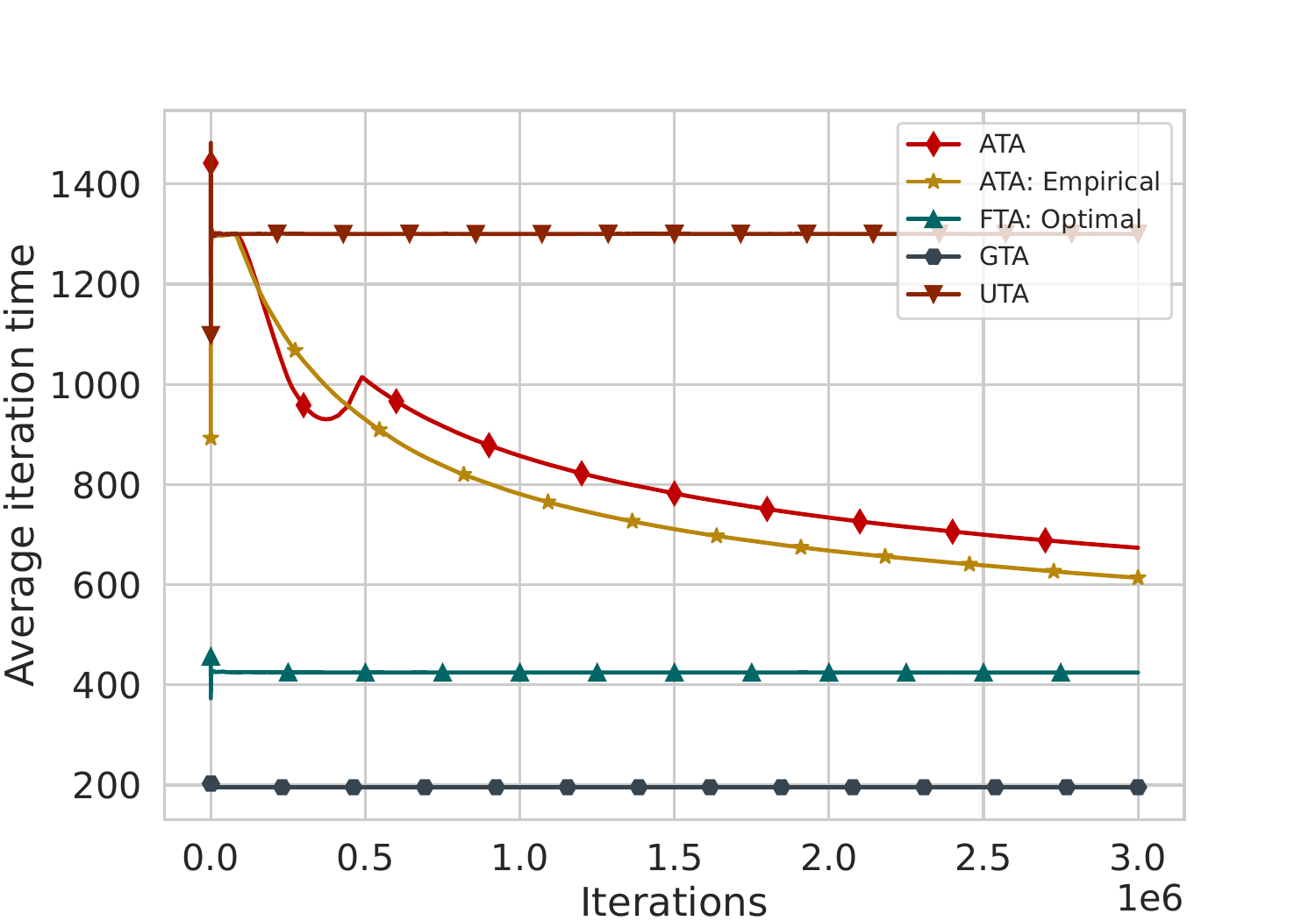} &
        \includegraphics[width=0.27\textwidth]{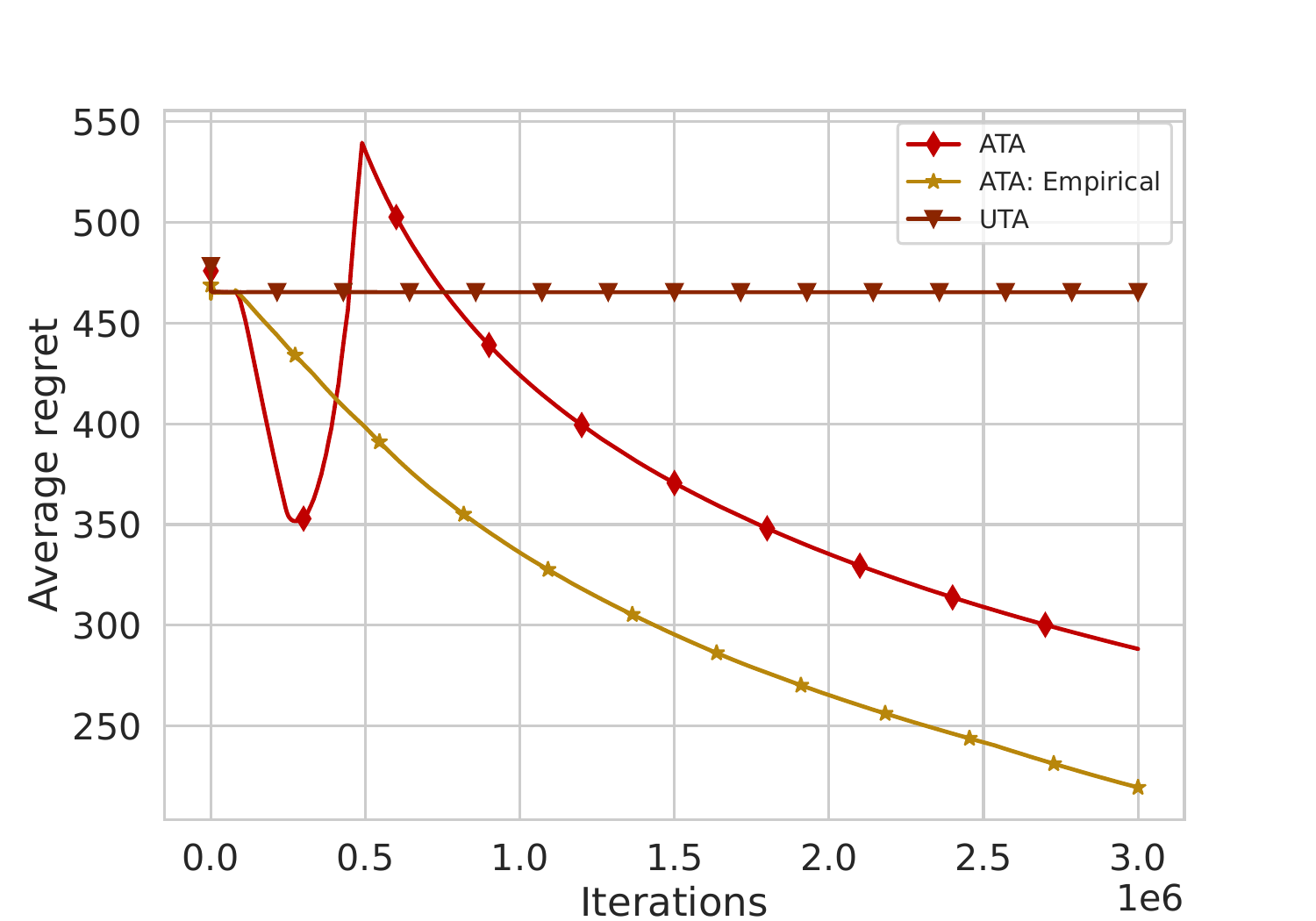} \\
        \includegraphics[width=0.27\textwidth]{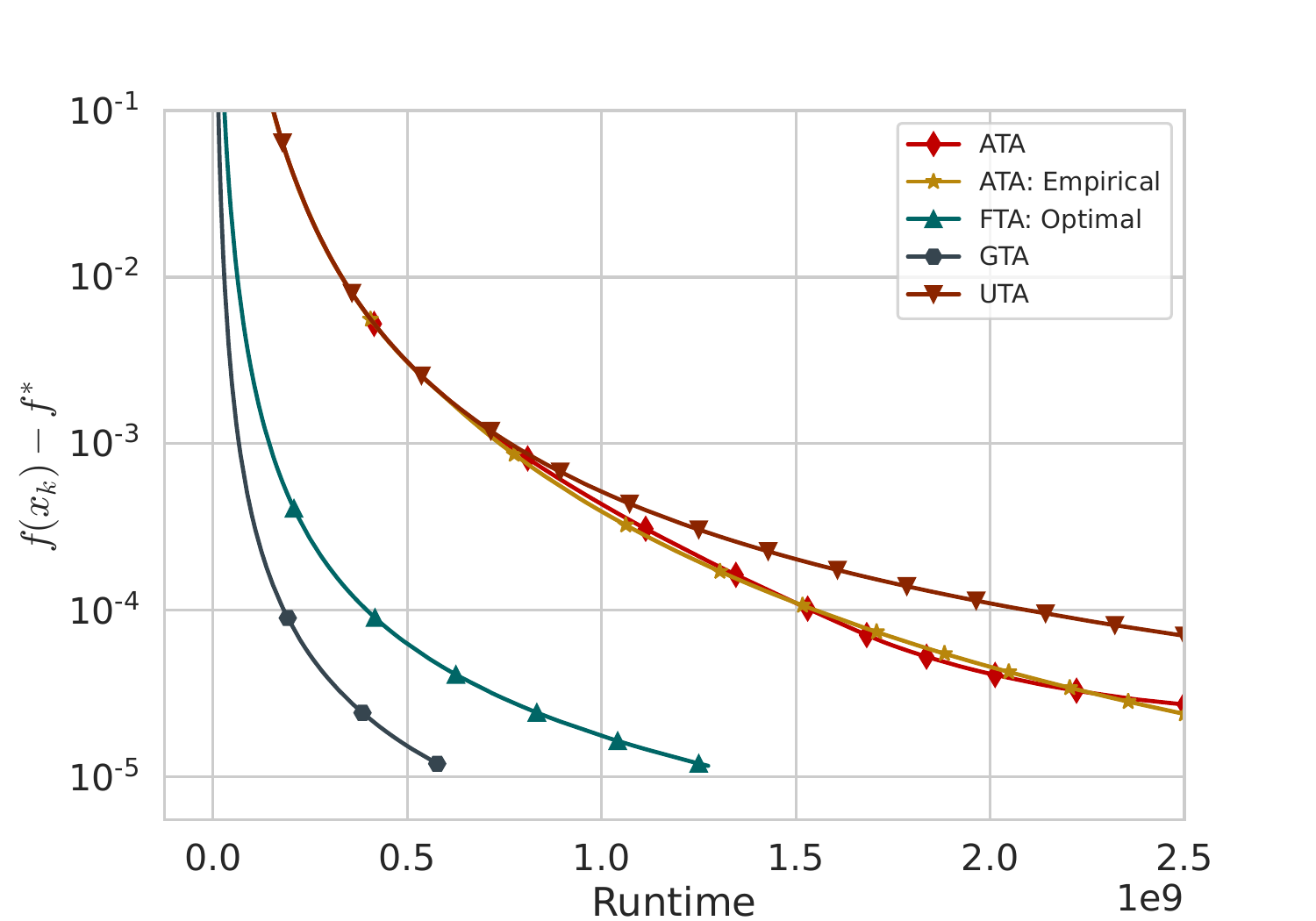} &
        \includegraphics[width=0.27\textwidth]{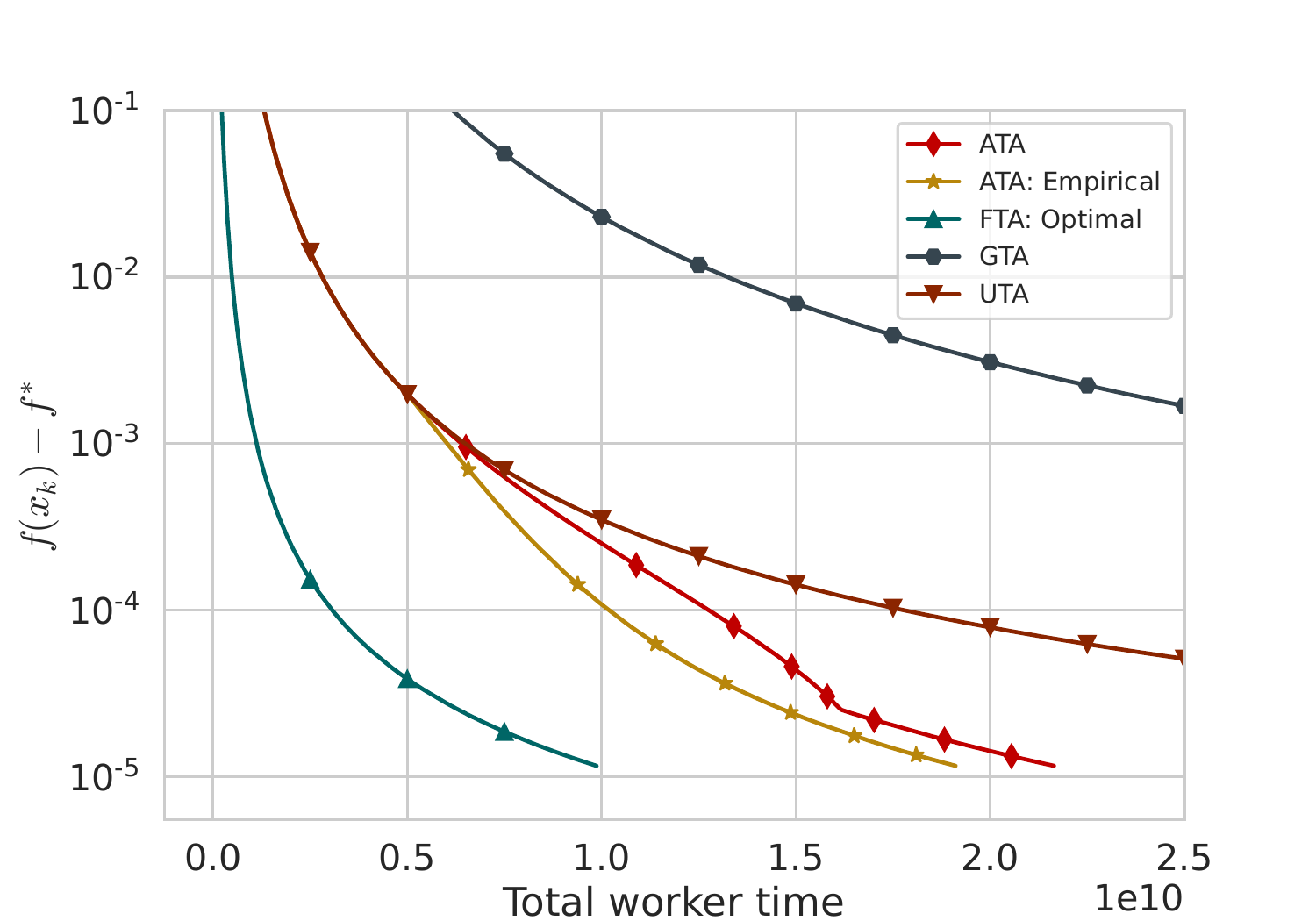} &
        \includegraphics[width=0.27\textwidth]{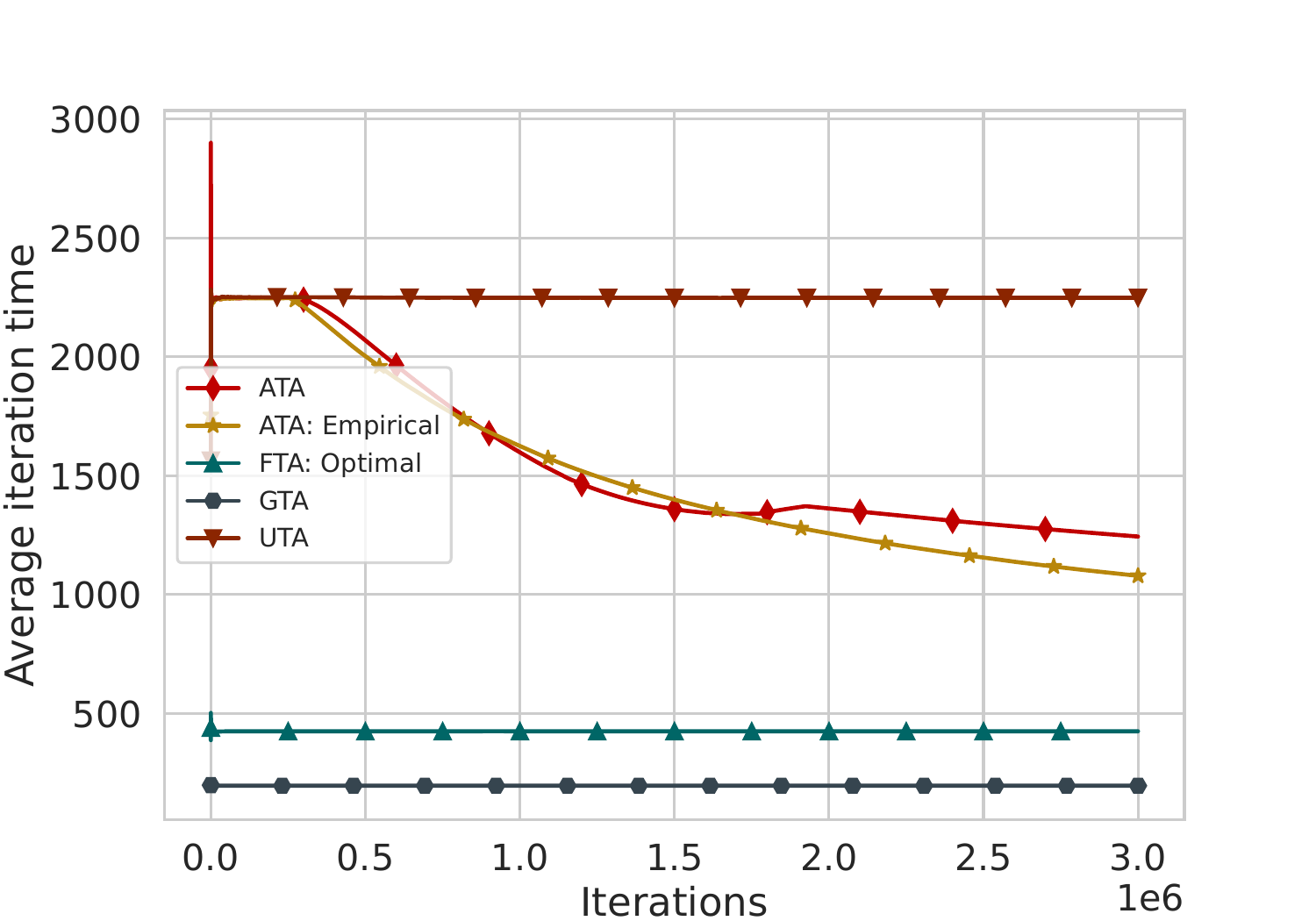} &
        \includegraphics[width=0.27\textwidth]{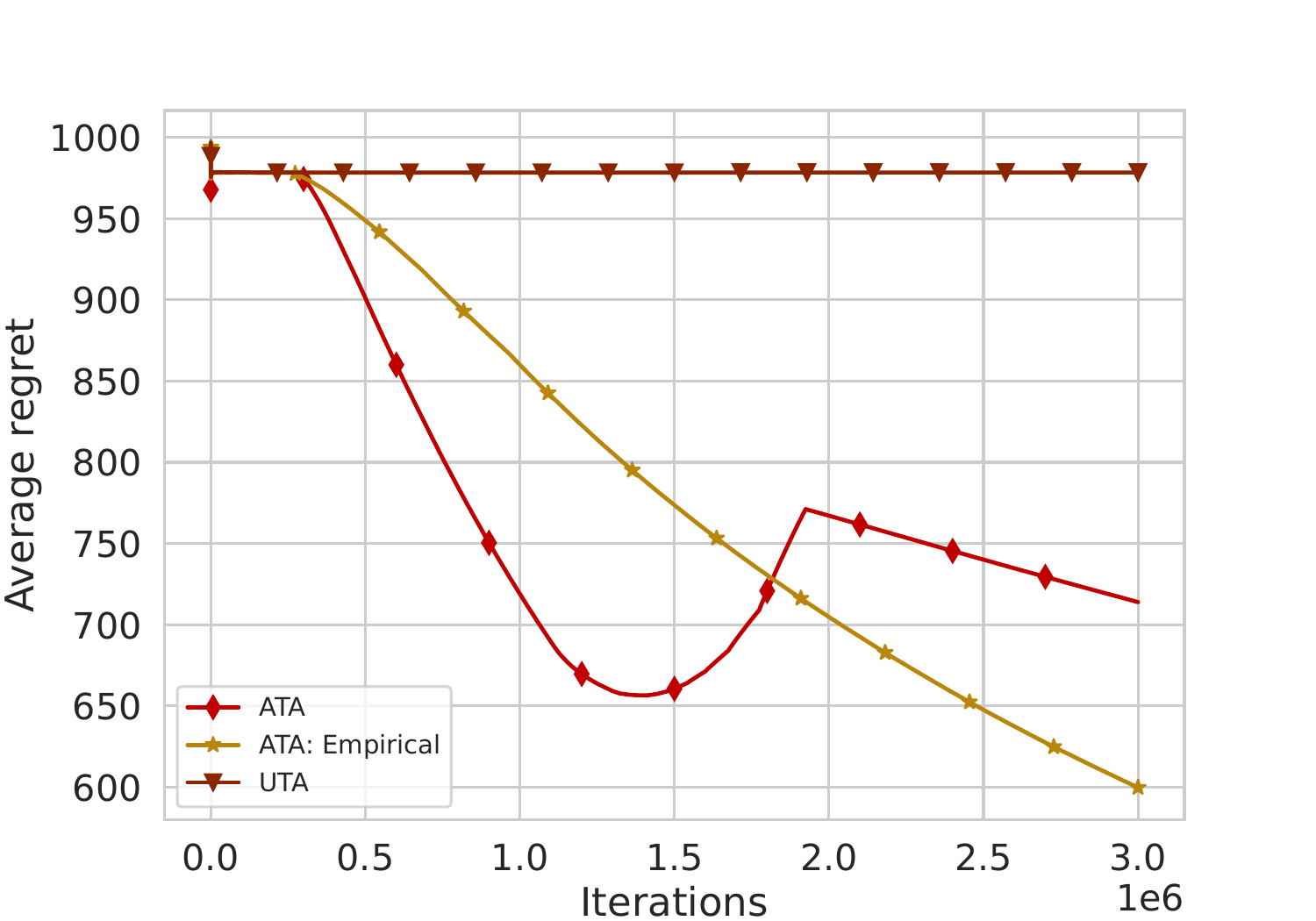} \\
    \end{tabular}
    } %
    \caption{
        Each row increases the number of workers by a factor of 3, starting from $17$, that is, $n = 17, 51, 153, 459$ from top to bottom.
        The first column shows runtime vs. suboptimality.  
        The second column also plots suboptimality, but against total worker time, i.e., $\sum_{i=1}^n T_i^k$ in \Cref{ata:alg:ata}.
        The third column presents the average iteration time, given by $\nicefrac{C_k}{k}$ over all iterations $k$.  
        The last column displays the averaged cumulative regret, as defined in \eqref{ata:eq:proxy_loss}.
    }
    \label{ata:fig:sqrt}
\end{figure*}

In this section, we validate our algorithms by simulating a scenario with $n$ workers, where we solve a simple problem using \algname{SGD}.
In each iteration, we collect $B=23$ gradients from the workers and perform a gradient descent step.

The objective function $f \,:\, \R^d \to \R$ is a convex quadratic defined as
$$
    f(x) = \frac{1}{2} x^\top A x - b^\top x~,
$$
where
\begin{align*}
    A &= \frac{1}{4}
    \begin{bmatrix}
    2 & -1 &  & 0 \\
    -1 & \ddots & \ddots &  \\
    & \ddots & \ddots & -1 \\
    0 & & -1 & 2 \\
    \end{bmatrix}
    \in \R^{d \times d} ~, \quad
    b = \frac{1}{4}
    \begin{bmatrix}
    -1 \\
    0 \\
    \vdots \\
    0 \\
    \end{bmatrix}
    \in \R^d~.
\end{align*}
We denote $f^*$ as the minimum value of the function $f$.
Each of the $n$ workers is able to calculate unbiased stochastic gradients $g(x)$ that satisfy
$$
    \E{\sqnorm{g(x) - \nabla f(x)}} \le 0.01^2 ~.
$$
This is achieved by adding Gaussian noise to the gradients of $f$.

The computation time for worker $i$ is modeled by the distribution
$$
    \nu_i = 29\sqrt{i} + \mathrm{Exp}\(29\sqrt{i}\),
$$
for all $i \in [n]$, where $\mathrm{Exp}(\beta)$ denotes the exponential distribution with scale parameter $\beta$.
The expected value of this distribution is $\mu_i = 2 \cdot 29 \sqrt{i}$.
Furthermore, the Orlicz norm satisfies the bound $\alpha_i \le 2\mu_i$.

We consider three benchmark algorithms.
\algname{GTA-SGD}, originally introduced as \rennala by \citet{tyurin2024optimal}.
Additionally, we include \algname{OFTA} (Optimal Fixed Task Allocation), which assumes the oracle knowledge of the mean computation times and uses the optimal allocation $\bar{a}$ in \eqref{ata:eq:proxy_loss} in each iteration, and \algname{UTA} (Uniform Task Allocation), which distributes $B$ tasks uniformly among the $n$ workers.
If $n>B$, then in \algname{UTA} we select $B$ workers at random, each one tasked to calculate one stochastic gradient.
Our algorithms aim to achieve a performance close to the one of \algname{OFTA}, without any prior knowledge of the true means.

For \ata, we set $\alpha = \alpha_n = 4 \cdot 29 \sqrt{n}$, while for \algname{ATA-Empirical} we use $\eta = 1$.
The results of our experiments are shown in \Cref{ata:fig:sqrt}.
As expected, \algname{GTA} is the fastest in terms of runtime (first column), but it performs poorly in terms of total worker time (second column).
This is because it uses all devices, most of which perform useless computations that are never used, leading to worse performance as the number of workers increases.
In fact, its performance can become arbitrarily worse.
On the other hand, \algname{OFTA} performs best in terms of total worker time.
Although it is slower in terms of runtime, the difference is by a constant factor that does not increase as $n$ grows.
This is because additional workers are less efficient and do not provide significant benefits for \algname{GTA}.

Turning our attention to our algorithms, both \ata and \algname{ATA-Empirical} initially behave like \algname{UTA}, as it is expected by the need to perform an initial exploration phase with uniform allocations.
However, after this phase, they begin to converge to the performance of \algname{OFTA}.

The last two columns contain plots that confirm our theoretical derivations.
The third plot validates \Cref{ata:cor:main}, showing that \ata and \algname{ATA-Empirical} converge to \algname{OFTA} up to a constant.
The final column shows the averaged cumulative regret, vanishing over time as predicted by Theorems~\ref{ata:thm:main}~and~\ref{ata:thm:main2}.

\begin{table}[h!]
	\caption{
        Ratios of total worker times and runtimes required to achieve $f(x) - f^{*} < 10^{-5}$.
        For total worker time, we divide the total worker time of \algname{GTA} by the corresponding total worker times of the other algorithms listed.
        For runtime, we do the opposite, dividing the runtime of the other algorithms by the runtime of \algname{GTA}, since \algname{GTA} is the fastest.
        To simplify the naming, we refer to \algname{ATA-Empirical} as \algname{ATA-E}.
        }
	\label{ata:table:sqrt}

	\begin{center}
	\begin{small}
	\begin{sc}
	\begin{tabular}{l|ccc|ccc}
        \toprule
        \multirow{2}{*}{$n$} & \multicolumn{3}{c|}{Tot. worker time ratio} & \multicolumn{3}{c}{Runtime ratio} \\
        \cmidrule(lr){2-4} \cmidrule(lr){5-7}
        & \ata & \algname{ATA-E} & \algname{OFTA} & \ata & \algname{ATA-E} & \algname{OFTA} \\
        \midrule
        $17$ & $1.3$ & $1.26$ & $1.26$ & $1.73$ & $1.75$ & $1.74$ \\
        $51$ & $2.91$ & $2.69$ & $3.03$ & $2.43$ & $2.45$ & $2.17$ \\
        $153$ & $7.22$ & $7.02$ & $9.1$ & $3.44$ & $3.14$ & $2.17$ \\
        $459$ & $12.45$ & $14.1$ & $27.3$ & $6.36$ & $5.51$ & $2.17$ \\
        \bottomrule
	\end{tabular}
	\end{sc}
	\end{small}
	\end{center}
\end{table}

In \Cref{ata:table:sqrt}, we compare the results numerically.
Both the total worker time ratio and runtime ratio increase as $n$ grows.
The total worker time ratio increases because \algname{GTA} becomes less efficient, using more resources than necessary.
The runtime ratio grows for \ata and \algname{ATA-Empirical} since a larger number of workers requires more exploration.
However, for \algname{OFTA} this ratio remains unchanged, as discussed earlier.

We remark that in these experiments we started all runs for \ata and \algname{ATA-Empirical} without prior knowledge of the computation time distribution.
However, in real systems, where these algorithms are used multiple times, prior estimates of computation times from previous runs could be available.
With this information, \ata and \algname{ATA-Empirical} would be much faster, as they would spend less time on exploration, approaching the performance of \algname{OFTA} in a faster way.
We validate this through experiments presented in \Cref{ata:sec:prior_knowledge}.

In \Cref{ata:sec:linear_noise}, we conducted similar experiments with a different time distribution, where the mean times vary linearly across the arms.
In \Cref{ata:sec:heterogeneous_distributions}, we examine scenarios with varying client time distributions.
Additionally, in \Cref{ata:sec:regret}, we analyze regret performance, confirming its logarithmic behavior as predicted by Theorems~\ref{ata:thm:main} and \ref{ata:thm:main2}.
Finally, in \Cref{ata:sec:real_dataset}, we trained a simple CNN on the CIFAR-100 dataset \citep{krizhevsky2009learning} using Adam \citep{kingma2014adam}.
\chapter{Concluding Remarks}\label{chapter:conclusion}
\thispagestyle{empty}
Asynchronous methods address the fundamental inefficiency of synchronous distributed training, where synchronization barriers cause faster workers to idle while waiting for slower ones.
However, prior to this work, asynchronous methods lacked rigorous theoretical foundations and had not demonstrated any theoretical advantages over synchronous methods.

To address this gap, we began with a clean theoretical setting that isolates the core problem asynchronous methods are designed to solve: heterogeneous computation times across workers.
We studied the data-parallel regime\textemdash where all workers maintain copies of the full model but process different data\textemdash using first-order stochastic methods for smooth nonconvex optimization.
This setting allowed us to abstract away factors orthogonal to asynchrony (such as communication delays and memory consistency) and focus on the fundamental question: can asynchronous methods be provably better?

In this setting, existing analyses often relied on restrictive assumptions, and crucially, no asynchronous algorithm was known to achieve optimal time complexity\textemdash a benchmark attained only by synchronous methods \citep{tyurin2024optimal}.

This dissertation establishes that asynchronous methods can indeed achieve optimal time complexity.
Through careful algorithm design and rigorous analysis, we demonstrated that staleness\textemdash the central challenge in asynchronous optimization\textemdash can be controlled through principled mechanisms such as selective gradient discarding and structured gradient buffering.
Our work provides the foundational theory showing that asynchronous methods can be both theoretically optimal and practically efficient in the fundamental setting of data-parallel stochastic optimization.
\section{Summary of Contributions}
Within this theoretical framework, we developed three main contributions:
\subsection{\ringmastertitle \ (\Cref{chapter:Ringmaster})}
We introduced the first asynchronous \sgd algorithm that achieves optimal time complexity in the homogeneous data setting.
\ringmaster (\Cref{ringmaster:alg:ringmasternewstops}) can be viewed as \asyncsgd (\Cref{ringmaster:alg:asgd}) with selective gradient discarding: it rejects gradients whose delay exceeds a threshold $R$, filtering out stale updates that would otherwise harm convergence.
With an appropriate choice of $R$, \ringmaster achieves optimal time complexity in both the fixed computation time \eqref{ringmaster:eq:worker-time} and arbitrarily varying computation time (\Cref{ringmaster:sec:dyn}) settings.
Moreover, our experiments demonstrate that \ringmaster can be faster than the optimal synchronous method \rennala \citep{tyurin2024optimal} in practice, due to more frequent model updates and the absence of synchronization overhead.

\subsection{\ringleadertitle \ (\Cref{chapter:Ringleader})}
We extended optimal asynchronous methods to the heterogeneous data regime, where local objectives differ across workers\textemdash a setting typical of federated learning \citep{konevcny2016federatedlearning, konevcny2016federatedoptimization, mcmahan2017communication}.

The key algorithmic innovation combines two mechanisms.
First, we maintain a \emph{gradient table} that stores gradients per worker, and use the average of all stored gradients for model updates.
This approach eliminates the need for restrictive similarity assumptions between local objectives, which rarely hold in federated learning due to diverse client data.
Second, we control staleness through \emph{structured rounds}: each round performs exactly $n$ updates (one per worker), ensuring that both the local models and the gradients in the table remain fresh.
This design keeps delays bounded without requiring explicit delay thresholds.

With these mechanisms, \ringleader achieves optimal time complexity in both the fixed and arbitrarily varying computation time settings.
Importantly, every gradient received by the server is either used immediately or stored for future updates\textemdash no computation is wasted, unlike synchronous methods that discard work at synchronization barriers.
This makes \ringleader the first asynchronous method proven optimal for heterogeneous data without similarity assumptions.
Moreover, our numerical experiments (\Cref{sec:experiments}) demonstrate that \ringleader can outperform the theoretically optimal synchronous method \malenia \citep{tyurin2024optimal} in practice, precisely because it avoids the computational waste inherent to synchronization.

\subsection{Adaptive Task Allocation (\Cref{chapter:ATA})}
This chapter addresses a different but equally important concern: \emph{computational efficiency}.
While the previous optimal methods minimize wall-clock time, they achieve this by utilizing all available computational resources.
In large-scale systems, such full utilization can lead to unnecessary and wasteful computation without proportional benefit.
This raises a natural question: \emph{Can we reduce computational waste without significantly increasing the total wall-clock time?}

To investigate this, we focus on the homogeneous data setting and build upon \ringmaster and \rennala.
We assume that worker computation times follow unknown probability distributions.
If these distributions were known, one could allocate the optimal number of tasks to each worker to minimize the expected wall-clock time.
We use this optimal fixed allocation as our competitor\textemdash the benchmark against which we measure performance.

Our goal is to design an algorithm that performs nearly as well as the best fixed competitor, but without prior knowledge of the worker-time distributions.
We achieve this with our proposed method, \ata.
\ata formulates task allocation as a multi-armed bandit problem, learning the workers' compute-time distributions online while balancing exploration and exploitation.
We prove that \ata matches the performance of the best fixed competitor up to a constant multiplicative factor in expectation, a standard guarantee in online learning known as constant regret.
\section{Future Directions}
Several important directions remain for future research.
We organize these into two categories: extensions specific to individual algorithms developed in this dissertation, and broader research directions that could apply to all the methods presented here or to asynchronous optimization more generally.
\subsection{Algorithm-Specific Extensions}
The following directions represent natural extensions of individual algorithms that would address their current limitations or enhance their specific capabilities.

\paragraph{Utilizing discarded gradients in \ringmastertitle.}
Currently, \ringmaster completely discards gradients whose delay exceeds the threshold $R$.
An interesting question is whether these stale gradients could be incorporated with reduced weights rather than being discarded entirely, similar to delay-adaptive asynchronous \sgd methods \citep{koloskova2022sharper,mishchenko2022asynchronous}.
Such an approach might improve wall-clock convergence time by making better use of available computation.
Moreover, in data center settings, consistently slow workers may never contribute if their gradients always exceed the delay threshold, meaning their data is never incorporated into training.
Developing mechanisms to include such workers\textemdash perhaps through delay-dependent weighting or adaptive thresholds\textemdash could improve both data utilization and model quality while maintaining near-optimal convergence guarantees.

\paragraph{Client sampling in \ringleadertitle.}
\ringleader requires gradients from \emph{all} workers before proceeding with updates, which may be impractical in federated learning where clients frequently become unavailable due to network issues, battery constraints, or user inactivity.
A natural extension would allow \ringleader to work with only a subset of available clients in each round, rather than waiting for all workers.
However, this modification likely requires reintroducing similarity assumptions between local objectives to maintain convergence guarantees.
A complementary research direction would be to establish lower bounds for asynchronous methods under these similarity assumptions.

\paragraph{Closing the optimality gap in \ringleadertitle.}
Our analysis shows that \ringleader attains existing lower bounds when the new smoothness constant (\Cref{ass:lipschitz_constant}) is within a constant factor of the standard smoothness parameter.
However, a gap remains: it is unclear whether this constant factor is fundamental or whether \ringleader's complexity can be further improved.
Future work could address this by either (i) relaxing the assumptions under which optimality is established or (ii) deriving tighter lower bounds tailored to asynchronous methods in the heterogeneous setting.
We believe the latter is more promising\textemdash new, sharper lower bounds for asynchronous algorithms may ultimately determine whether \ringleader is truly optimal or if additional improvements are still possible.

\paragraph{Adaptive task allocation for heterogeneous data.}
The \ata algorithm was developed for the homogeneous data setting, where it improves the computational efficiency of \ringmaster and \rennala \citep{tyurin2024optimal}.
A natural extension would be to adapt these ideas to the heterogeneous data regime, thereby enhancing the computational efficiency of \ringleader and \malenia \citep{tyurin2024optimal}.
Such an extension would need to account for the different criteria for stopping gradient collection, which in this case depend on the harmonic mean of the gradient counts.

\paragraph{Adaptive task allocation under adversarial delays.}
The existing \ata framework assumes worker computation times follow unknown but fixed distributions.
In federated learning, however, worker behavior can be highly unpredictable: devices may become unavailable, experience extreme delays, or exhibit non-stationary patterns that do not follow any stable distribution.
Developing allocation strategies robust to such adversarial or non-stochastic delays would significantly enhance practical applicability.
This direction connects to the adversarial online learning and bandit algorithm literature \citep{auer2002nonstochastic}.

\subsection{Broader Research Directions}
Beyond extensions of individual algorithms, several fundamental questions could enhance all the methods in this dissertation or advance asynchronous optimization more broadly.

\paragraph{Beyond SGD.}
Our analysis has focused exclusively on stochastic gradient descent, the most fundamental first-order optimization method.
However, modern deep learning frequently relies on more advanced optimizers such as \algname{Adam} \citep{kingma2014adam}, \algname{AdamW} \citep{loshchilov2017decoupled}, and preconditioned gradient methods like \algname{Shampoo} \citep{gupta2018shampoo}, as well as newer algorithms such as \algname{Muon} \citep{jordan2024muon}.
Extending our asynchronous framework to these methods while preserving theoretical guarantees would significantly broaden the practical scope of this work.

\paragraph{Extensions to other forms of parallelism.}
This dissertation focused exclusively on data parallelism.
However, modern large-scale training systems increasingly combine data parallelism with model and pipeline parallelism \citep{narayanan2021efficient,shoeybi2019megatron}.
Extending asynchronous optimization methods to such hybrid settings is a natural next step.
Of particular interest is developing time-complexity analyses for asynchronous pipeline parallelism, where different stages of a model are executed on distinct workers and subject to inter-stage delays.
Several recent works have explored asynchronous pipeline training in this context \citep{yang2019pipemare,ajanthan2025nesterov,kong2025clapping}.

\paragraph{Communication costs and network delays.}
Throughout this work, we assumed instantaneous communication\textemdash a simplification that allowed us to isolate the fundamental challenges of asynchrony.
In federated and geographically distributed systems, however, communication often dominates computation time.
To mitigate this, existing approaches employ either communication compression \citep{alistarh2017qsgd} or local training \citep{mcmahan2016federated,mishchenko2022proxskip,malinovsky2022variance,maranjyan2025gradskip} to reduce the frequency or volume of transmitted updates.
Combining asynchrony with these techniques\textemdash local training or communication compression\textemdash and analyzing the resulting time complexity remains an important open problem.
Initial progress has been made in both directions, with recent results on asynchronous compression \citep{tyurin2024freya} and asynchronous local training \citep{tyurin2025birch,fradin2025local}.

\paragraph{Lock-free shared memory algorithms.}
Throughout this work, we assumed atomic reads and writes with instantaneous communication.
An alternative model, studied in works such as \hogwild \citep{recht2011hogwild}, considers shared memory systems where multiple workers can read and write simultaneously without locks.
In this setting, a worker reading the parameter vector may observe inconsistent values\textemdash a mixture of the current server state and partial updates from other workers.
Recent works \citep{mania2017perturbed, leblond2018improvedAsynchronous, nguyen2018sgd} have provided improved analyses of \hogwild-style algorithms.
It would be interesting to investigate whether the techniques developed in this dissertation can be adapted to this lock-free setting and whether optimal time complexity can be achieved despite inconsistent reads.

\paragraph{Large-scale empirical validation.}
The experiments in this dissertation were primarily small-scale demonstrations designed to validate theoretical predictions.
A valuable direction for future work is large-scale empirical validation in realistic settings: training large language models or other modern architectures in actual data center or federated learning environments.
Such experiments would reveal practical bottlenecks not captured by our theoretical model, and would help bridge the gap between theory and deployment.

\section{Closing Remarks}

This dissertation establishes that asynchronous methods can achieve optimal time complexity, providing the first rigorous theoretical foundation for their use in distributed optimization.
While our results focus on the fundamental setting of data-parallel stochastic gradient descent applied to smooth nonconvex objectives, they mark an essential first step toward a comprehensive theory of asynchronous optimization.

Beyond the specific optimality results, this work introduces broadly applicable techniques: principled mechanisms for controlling staleness (e.g., selective gradient discarding and structured gradient buffering), a rigorous framework for analyzing convergence and time complexity under asynchrony, and strategies for balancing wall-clock time with computational cost.
Together, these ideas lay the groundwork for extending asynchronous methods to more complex and realistic environments.

The significance of this work extends beyond theory.
Asynchronous algorithms have historically seen limited practical adoption, hindered by weak theoretical guarantees and limited software support.
By demonstrating that they can be both theoretically optimal and practically efficient, this dissertation provides strong motivation for developing robust implementations and integrating these methods into large-scale training systems.
As AI systems continue to grow in scale and as training infrastructure becomes increasingly heterogeneous\textemdash spanning data centers, edge devices, and federated networks\textemdash the importance of efficient asynchronous optimization will only continue to rise.

The central message of this dissertation is one of possibility rather than limitation.
Synchronous algorithms offer predictability but waste resources through forced idleness.
Asynchronous algorithms eliminate this waste but introduce the challenge of managing stale information.
This dissertation has shown that the trade-off is not fundamental: through careful algorithm design, staleness can be controlled while maintaining both theoretical optimality and practical efficiency.
The future of large-scale machine learning lies not in forcing all workers to march in lockstep, but in enabling them to contribute asynchronously\textemdash coordinated through principled mechanisms that ensure convergence and efficiency.
As demonstrated here, such coordination is not only possible\textemdash it is optimal.


\begin{onehalfspacing}
	\renewcommand*\bibname{\centerline{REFERENCES}} 
    \phantomsection
	\addcontentsline{toc}{chapter}{References}
	\newcommand{\BIBdecl}{\setlength{\itemsep}{0pt}}
		\bibliography{bib}
        \bibliographystyle{iclr2026_conference}
\end{onehalfspacing}



\appendix
		\newpage
		\begingroup
			\begin{center}
			\vspace*{2\baselineskip}
			{ \textbf{{\large APPENDICES}}} 
            \phantomsection
			\addcontentsline{toc}{chapter}{Appendices} 
			\end{center}
            

\chapter{Appendix for \Cref{chapter:Ringmaster}}
\label{appendix_ringmaster}
\thispagestyle{empty}
\section{Proof of \Cref{ringmaster:lem:time_R}}
\label{ringmaster:proof:time_R}
To analyze the time complexity of \ringmaster, we first establish a bound on the time needed to perform a fixed number of updates.
For completeness, we restate the lemma below before providing its proof.
\begin{restate-boxedlemma}{\ref{ringmaster:lem:time_R}}
    Let workers' computation times satisfy the \emph{fixed computation model} \eqref{ringmaster:eq:worker-time}.
    Let $R$ be the delay threshold of \Cref{ringmaster:alg:ringmasternew} or \Cref{ringmaster:alg:ringmasternewstops}.
    The time required to complete any $R$ consecutive iterates updates of \Cref{ringmaster:alg:ringmasternew} or \Cref{ringmaster:alg:ringmasternewstops} is at most
    \begin{equation}
        \tag{\ref{ringmaster:eq:time_R}}
        t(R) \eqdef 2 \min\limits_{m \in [n]} \left(\frac{1}{m} \sum\limits_{i=1}^m \frac{1}{\tau_i}\right)^{-1} \left( 1 + \frac{R}{m} \right) .
    \end{equation}
\end{restate-boxedlemma}
\begin{proof}
    We now focus on \Cref{ringmaster:alg:ringmasternewstops}.
    Let us consider a simplified notation of $t(R)$:
    \begin{align*}
        t \eqdef t(R) = 2 \min_{m\in[n]} \left(\sum_{i=1}^m \frac{1}{\tau_i}\right)^{-1} (R + m) = 2 \min\limits_{m \in [n]} \left(\frac{1}{m} \sum\limits_{i=1}^m \frac{1}{\tau_i}\right)^{-1} \left( 1 + \frac{R}{m} \right).
    \end{align*}
    Let us fix any iteration $k$, and consider the consecutive iterations from $k$ to $k + R - 1$.
    By the design of \Cref{ringmaster:alg:ringmasternewstops}, any worker will be stopped at most one time, meaning that at most one stochastic gradient will be ignored from the worker.
    After that, the delays of the next stochastic gradients will not exceed $R - 1$ during these iterations.
    As soon as some worker finishes calculating a stochastic gradient, it immediately starts computing a new stochastic gradient.
    
    Using a proof by contradiction, assume that \Cref{ringmaster:alg:ringmasternewstops} will not be able to finish iteration $k + R - 1$ by the time $t$.
    At the same time, by the time $t$, all workers will calculate at least
    \begin{align}
        \label{ringmaster:eq:RQKxmnzvrIUDaroZE}
        \sum_{i=1}^{n} \max\left\{\flr{\frac{t}{\tau_i}} - 1, 0\right\}
    \end{align}
    stochastic gradients with delays less than $R$ because $\flr{\nicefrac{t}{\tau_i}}$ is the number of stochastic gradients that worker $i$ can calculate after $t$ seconds.
    We subtract $1$ to account for the fact that one stochastic gradient with a large delay can be ignored.
    Let us take an index
    \begin{align*}
        j^* = \argmin_{m\in[n]} \left(\sum_{i=1}^m \frac{1}{\tau_i}\right)^{-1} (R + m) .
    \end{align*}
    Since $\flr{x} \geq x - 1$ for all $x \geq 0$, we get
    \begin{align*}
        \sum_{i=1}^{n} \max\left\{\flr{\frac{t}{\tau_i}} - 1, 0\right\} 
        &\geq \sum_{i=1}^{j^*} \max\left\{\flr{\frac{t}{\tau_i}} - 1, 0\right\} \geq \sum_{i=1}^{j^*} \left(\flr{\frac{t}{\tau_i}} - 1\right) \\
        &\geq \sum_{i=1}^{j^*} \frac{t}{\tau_i} - 2 j^* \\
        &= 2 \left(\sum_{i=1}^{j^*} \frac{1}{\tau_i}\right) \left(\left(\sum_{i=1}^{j^*} \frac{1}{\tau_i}\right)^{-1} (R + j^*)\right) - 2 j^* \\
        &= 2 R + 2 j^* - 2 j^* \geq R.
    \end{align*}
    We can conclude that by the time \eqref{ringmaster:eq:time_R}, the algorithm will calculate $R$ stochastic gradients with delays less than $R$ and finish iteration $k + R - 1$, which contradicts the assumption. 
    
    The proof for \Cref{ringmaster:alg:ringmasternew} is essentially the same.
    In \Cref{ringmaster:alg:ringmasternew}, one stochastic gradient with a large delay for each worker will be ignored, and all other stochastic will be used in the optimization process.
\end{proof}

\section{Proof of \Cref{ringmaster:lem:time_R_dyn}}
\label{ringmaster:sec:time_R_dyn}
To handle the dynamic setting where worker speeds vary over time, we now prove an analogue of \Cref{ringmaster:lem:time_R} under the universal computation model.
\begin{restate-boxedlemma}{\ref{ringmaster:lem:time_R_dyn}}
    Let the workers' computation times satisfy the \emph{universal computation model}.
    Let $R$ be the delay threshold of \Cref{ringmaster:alg:ringmasternew} or \Cref{ringmaster:alg:ringmasternewstops}.
    Assume that some iteration starts at time $T^0$.
    Starting from this iteration, the $R$ consecutive iterate updates of \Cref{ringmaster:alg:ringmasternew} or \Cref{ringmaster:alg:ringmasternewstops} will be performed before the time
    \begin{equation*} 
           T(R,T^0) \eqdef \min \left\{T \geq 0 : \sum \limits_{i=1}^{n} \flr{\frac{1}{4} \int_{T^0}^{T} p_i(\tau) d \tau} \geq R\right\}.
    \end{equation*}
\end{restate-boxedlemma}
\begin{proof}
    Let 
    \begin{align*}
        T \eqdef T(R,T^0) = \min \left\{T \geq 0 : \sum_{i=1}^{n} \flr{\frac{1}{4} \int_{T^0}^{T} v_i(\tau) d \tau} \geq R\right\}.
    \end{align*}
    At the beginning, this proof follows the proof of \Cref{ringmaster:lem:time_R} up to \eqref{ringmaster:eq:RQKxmnzvrIUDaroZE}.
    Here we also use a proof by contradiction and assume that \Cref{ringmaster:alg:ringmasternewstops} will not be able to do $R$ consecutive iterative updates by the time $T$.

    Instead of \eqref{ringmaster:eq:RQKxmnzvrIUDaroZE}, all workers will calculate at least
    \begin{align*}
        N \eqdef \sum_{i=1}^{n} \max\left\{\flr{\int_{T^0}^{T} v_i(\tau) d \tau} - 1, 0\right\} 
    \end{align*}
    stochastic gradients by time $T$ because, due to \eqref{ringmaster:eq:rSIiSfVcmivSKsfzoSA}, $\flr{\int_{T^0}^{T} v_i(\tau) d \tau}$ is the number of stochastic gradients that worker $i$ will calculate in the interval $[T^0, T]$.
    We define $V_i(T) \eqdef \int_{T^0}^{T} v_i(\tau) d \tau$.
    Thus,
    \begin{align*}
        N = \sum_{i=1}^{n} \max\left\{\flr{V_i(T)} - 1, 0\right\} \textnormal{ and } T = \min \left\{T \geq 0 : \sum_{i=1}^{n} \flr{\frac{V_i(T)}{4}} \geq R\right\}.
    \end{align*}
    Let us additionally define
    \begin{align}
        S \eqdef \left\{i \in [n] \,:\, \frac{V_i(T)}{4} \geq 1\right\}.
        \label{ringmaster:eq:peCcI}
    \end{align} 
    Note that
    \begin{align}
        \label{ringmaster:eq:XwmTgavNaBJOKaRgFMT}
        \sum_{i=1}^{n} \flr{\frac{V_i(T)}{4}} = \sum_{i \in S} \flr{\frac{V_i(T)}{4}},
    \end{align}
    since $\flr{\nicefrac{V_i(T)}{4}} = 0$ for all $i \not\in S$.
    Using simple algebra, we get
    \begin{align*}
        N &= \sum_{i=1}^{n} \max\{\flr{V_i(T)} - 1, 0\} \geq \sum_{i \in S} \max\{\flr{V_i(T)} - 1, 0\} \\
        &\geq \sum_{i \in S} \flr{V_i(T)} - |S| \geq \sum_{i \in S} V_i(T) - 2 |S|~,
    \end{align*}
    where the last inequality due to $\flr{x} \geq x - 1$ for all $x \in \R$.
    Next, using \eqref{ringmaster:eq:peCcI}, we have
    \begin{align*}
        N &\geq \frac{1}{2} \sum_{i \in S} V_i(T) + \frac{1}{2} \sum_{i \in S} 4 - 2 |S| = \sum_{i \in S} \frac{V_i(T)}{2} \\
        &\geq \sum_{i \in S} \flr{\frac{V_i(T)}{4}} \overset{\eqref{ringmaster:eq:XwmTgavNaBJOKaRgFMT}}{=} \sum_{i=1}^{n} \flr{\frac{V_i(T)}{4}} \geq R~,
    \end{align*}
    where the last inequality due to the definition of $T$.
    As in the proof of \Cref{ringmaster:lem:time_R}, we can conclude that by the time $T = T(R,T^0)$, the algorithm will calculate $R$ stochastic gradients with delays less than $R$ and finish iteration $k + R - 1$, which contradicts the assumption.
\end{proof}

\section{Proof of \Cref{ringmaster:thm:ringmaster_iteration}}
\label{ringmaster:proof:ringmaster_iteration}
We restate below the result that bounds the number of iterations required by \ringmaster to reach the desired level of stationarity, and then provide its proof.
\begin{restate-boxedtheorem}{\ref{ringmaster:thm:ringmaster_iteration}}
    Under Assumptions \ref{ringmaster:ass:lipschitz_constant}, \ref{ringmaster:ass:lower_bound}, and \ref{ringmaster:ass:stochastic_variance_bounded}, let the stepsize in \ringmaster (\Cref{ringmaster:alg:ringmasternew} or \Cref{ringmaster:alg:ringmasternewstops}) be
    $$
        \gamma = \min \left\{ \frac{1}{2RL}, \frac{\varepsilon}{4L\sigma^2} \right\}.
    $$
    Then, the following holds
    $$
        \frac{1}{K+1}\sum_{k=0}^{K} \Exp{\sqnorm{\nabla f \(x^k\)}} \le \varepsilon~,
    $$
    for 
    \begin{align}
        \tag{\ref{ringmaster:eq:HPMNVlgxoiRgUbRtj}}
        K \geq \frac{8 R L \Delta}{\epsilon} + \frac{16 \sigma^2 L \Delta}{\epsilon^2}~,
    \end{align}
    where $R \in \{1, 2, \dots \}$ is an arbitrary delay threshold.
\end{restate-boxedtheorem}
We adopt the proof strategy of \citet{koloskova2022sharper} and rely on the following two lemmas.
\begin{boxedlemma}[Descent Lemma; Proof in Appendix~\ref{ringmaster:proof:descent}]
    \label{ringmaster:lemma:descent}
    Under Assumptions \ref{ringmaster:ass:lipschitz_constant} and \ref{ringmaster:ass:stochastic_variance_bounded}, if the stepsize in \ringmaster (\Cref{ringmaster:alg:ringmasternew} or \Cref{ringmaster:alg:ringmasternewstops}) satisfies $\gamma \le \nicefrac{1}{2L}$, the following inequality holds:
    \begin{align*}        
        \ExpSub{k+1}{f\(x^{k+1}\)} 
        \leq f\(x^{k}\) 
        &- \frac{\gamma}{2} \sqnorm{\nabla f\(x^{k}\)} 
        - \frac{\gamma}{4}\sqnorm{\nabla f\(x^{k-\delta^k}\)} \\
        &+ \frac{\gamma L^2}{2}\sqnorm{x^k - x^{k-\delta^k}} 
        + \frac {\gamma^2 L}{2} \sigma^2~,
    \end{align*}
    where $\ExpSub{k+1}{\cdot}$ represents the expectation conditioned on all randomness up to iteration $k$.
\end{boxedlemma}
\begin{boxedlemma}[Residual Estimation; Proof in Appendix~\ref{ringmaster:proof:residual}]
    \label{ringmaster:lemma:residual}
    Under Assumptions \ref{ringmaster:ass:lipschitz_constant} and \ref{ringmaster:ass:stochastic_variance_bounded}, the iterates of \ringmaster (\Cref{ringmaster:alg:ringmasternew} or \Cref{ringmaster:alg:ringmasternewstops}) with stepsize $\gamma \leq \nicefrac{1}{2RL}$ satisfy the following bound:
    $$
    \frac{1}{K+1} \sum_{k=0}^K \E{\sqnorm{x^k - x^{k-\delta^k}}} \leq \frac{1}{2 L^2(K+1)} \sum_{k=0}^K \E{\sqnorm{\nabla f\(x^{k-\delta^k}\)}} + \frac{\gamma}{L} \sigma^2.
    $$
\end{boxedlemma}
With these results, we are now ready to prove \Cref{ringmaster:thm:ringmaster_iteration}.
\begin{proof}[Proof of \Cref{ringmaster:thm:ringmaster_iteration}]
We begin by averaging over $K+1$ iterations and dividing by $\gamma$ in the inequality from \Cref{ringmaster:lemma:descent}:
\begin{align*}
    \frac{1}{K+1} \sum_{k=0}^K & \( \frac{1}{2} \E{ \sqnorm{\nabla f\(x^k\)}} + \frac{1}{4} \E{ \sqnorm{\nabla f(x^{k-\delta^k})}} \) 
    \leq \frac{\Delta}{\gamma(K+1)} + \frac{\gamma L}{2} \sigma^2 \\ 
    &\hspace{13em} + \frac{1}{K+1} \frac{L^2}{2} \sum_{k=0}^K \E{ \sqnorm{x^k - x^{k-\delta^k}}}.
\end{align*}
Next, applying \Cref{ringmaster:lemma:residual} to the last term, we have:
\begin{align*}
    \frac{1}{K+1} \sum_{k=0}^K & \( \frac{1}{2} \E{ \sqnorm{\nabla f\(x^k\)}} + \frac{1}{4} \E{ \sqnorm{\nabla f(x^{k-\delta^k})}} \) 
    \leq \frac{\Delta}{\gamma(K+1)} + \frac{\gamma L}{2} \sigma^2 \\ 
    &\hspace{9em} + \frac{1}{4(K+1)} \sum_{k=0}^K \E{\sqnorm{\nabla f\(x^{k-\delta^k}\)}} + \frac{\gamma L}{2}\sigma^2.
\end{align*}
Simplifying further, we obtain:
$$
    \frac{1}{K+1} \sum_{k=0}^K \E{ \sqnorm{\nabla f\(x^k\)} } \leq \frac{2\Delta}{\gamma(K+1)} + 2\gamma L \sigma^2.
$$
Now, we choose the stepsize $\gamma$ as
$$
    \gamma = \min \left\{ \frac{1}{2RL}, \frac{\varepsilon}{4L\sigma^2} \right\} \leq \frac{1}{2 R L} ~.
$$
With this choice of $\gamma$, it remains to choose
$$
    K \geq \frac{8 \Delta R L}{\epsilon} + \frac{16 \Delta L \sigma^2}{\epsilon^2}
$$ 
to ensure
$$
    \frac{1}{K+1} \sum_{k=0}^{K} \E{\sqnorm{ \nabla f\(x^k\) }} \leq \epsilon~.
$$
This completes the proof.
\end{proof}
\subsection{Proof of \Cref{ringmaster:lemma:descent}}
\label{ringmaster:proof:descent}
We restate below the descent lemma used in the convergence analysis of \ringmaster and then provide its proof.
\begin{restate-boxedlemma}{\ref{ringmaster:lemma:descent}}[Descent Lemma]
    Under Assumptions \ref{ringmaster:ass:lipschitz_constant} and \ref{ringmaster:ass:stochastic_variance_bounded}, if the stepsize in \ringmaster (\Cref{ringmaster:alg:ringmasternew} or \Cref{ringmaster:alg:ringmasternewstops}) satisfies $\gamma \le \nicefrac{1}{2L}$, the following inequality holds:
    \begin{align*}
        \ExpSub{k+1}{f\(x^{k+1}\)} &
        \leq f\(x^{k}\) 
        - \frac{\gamma}{2} \sqnorm{\nabla f\(x^{k}\)} 
        - \frac{\gamma}{4}\sqnorm{\nabla f\(x^{k-\delta^k}\)} \\
        & + \frac{\gamma L^2}{2}\sqnorm{x^k - x^{k-\delta^k}} 
        + \frac{\gamma^2 L}{2} \sigma^2,
    \end{align*}
    where $\ExpSub{k+1}{\cdot}$ represents the expectation conditioned on all randomness up to iteration $k$.
\end{restate-boxedlemma}


\begin{proof}
Assume that we get a stochastic gradient from the worker with index $i_k$ when calculating $x^{k+1}$.
Since the function $f$ is $L$-smooth (Assumption~\ref{ringmaster:ass:lipschitz_constant}), we have \citep{nesterov2018lectures}:
\begin{align*}
    \ExpSub{k+1}{f\(x^{k+1}\)} 
    &\le f\(x^{k}\) 
        - \gamma \underbrace{\ExpSub{k+1}{\< \nabla f\(x^{k}\), \nabla f\(x^{k-\delta^k}; \xi^{k-\delta^k}_{i_k}\) \>}}_{=: t_1} \\
        &+ \; \frac{L}{2} \gamma^2 \underbrace{\ExpSub{k+1}{\sqnorm{\nabla f\(x^{k-\delta^k}; \xi^{k-\delta^k}_{i_k}\)}}}_{=: t_2}.
\end{align*}
Using Assumption~\ref{ringmaster:ass:stochastic_variance_bounded}, we estimate the second term as
\begin{align*}
    t_1 
    &= \< \nabla f\(x^{k}\), \nabla f\(x^{k-\delta^k}\) \> \\
    &= \frac{1}{2} \[\sqnorm{\nabla f\(x^{k}\)} 
        + \sqnorm{\nabla f\(x^{k-\delta^k}\)} 
        - \sqnorm{\nabla f\(x^{k}\) - \nabla f\(x^{k-\delta^k}\)} \].
\end{align*}
Using the variance decomposition equality and Assumption~\ref{ringmaster:ass:stochastic_variance_bounded}, we get
\begin{eqnarray*}
t_2 
&=& \ExpSub{k+1}{\sqnorm{\nabla f\(x^{k-\delta^k}; \xi^{k-\delta^k}_{i_k}\) - \nabla f\(x^{k-\delta^k}\)}} 
    + \sqnorm{\nabla f\(x^{k-\delta^k}\)} \\
&\le& \sigma^2 
    + \sqnorm{\nabla f\(x^{k-\delta^k}\)} .
\end{eqnarray*}
Combining the results for $t_1$ and $t_2$, and using $L$--smoothness to bound \\
$\|\nabla f(x^{k}) - \nabla f(x^{k-\delta^k})\|^2$, we get
\begin{align*}
    \ExpSub{k+1}{f\(x^{k+1}\)} 
    &\leq f\(x^{k}\) 
    - \frac{\gamma}{2} \sqnorm{\nabla f\(x^{k}\)} 
    - \frac{\gamma}{2} (1 - \gamma L) \sqnorm{\nabla f\(x^{k-\delta^k}\)} \\
    &+ \frac{\gamma L^2}{2}\sqnorm{x^k - x^{k-\delta^k}} 
    + \frac{\gamma^2 L}{2} \sigma^2.
\end{align*}
Finally, applying the condition $\gamma \le \nicefrac{1}{2 L}$ completes the proof.
\end{proof}

\subsection{Proof of \Cref{ringmaster:lemma:residual}}
\label{ringmaster:proof:residual}
We restate below the lemma that bounds the expected residual between consecutive iterates of \ringmaster and then provide its proof.
\begin{restate-boxedlemma}{\ref{ringmaster:lemma:residual}}[Residual Estimation]
    Under Assumptions \ref{ringmaster:ass:lipschitz_constant} and \ref{ringmaster:ass:stochastic_variance_bounded}, the iterates of \ringmaster (\Cref{ringmaster:alg:ringmasternew} or \Cref{ringmaster:alg:ringmasternewstops}) with stepsize $\gamma \leq \nicefrac{1}{2RL}$ satisfy the following bound:
    $$
        \frac{1}{K+1} \sum_{k=0}^K \E{\sqnorm{x^k - x^{k-\delta^k}}} \leq \frac{1}{2 L^2(K+1)} \sum_{k=0}^K \E{\sqnorm{\nabla f\(x^{k-\delta^k}\)}} + \frac{\gamma}{L} \sigma^2.
    $$
\end{restate-boxedlemma}


\begin{proof}
Assume that we get a stochastic gradient from the worker with index $i_k$ when calculating $x^{k+1}$.
We begin by expanding the difference and applying the tower property, Assumption~\ref{ringmaster:ass:stochastic_variance_bounded}, Young's inequality, and Jensen's inequality:
\begin{align*}
    \E{ \sqnorm{x^k - x^{k-\delta^k}} } 
    & = \E{ \sqnorm{ \sum_{j=k-\delta^k}^{k-1} \gamma \nabla f\(x^{j-\delta_j}; \xi^{j-\delta_j}_{i_j}\) } } \\
    & \leq 2 \E{ \sqnorm{ \sum_{j=k-\delta^k}^{k-1} \gamma \nabla f\(x^{j-\delta_j}\) } } \\
        &\quad + 2 \E{ \sqnorm{ \sum_{j=k-\delta^k}^{k-1} \gamma \left(\nabla f\(x^{j-\delta_j}; \xi^{j-\delta_j}_{i_j}\) - \nabla f\(x^{j-\delta_j}\)\right) } }\\
    & \leq 2 \E{ \sqnorm{ \gamma \sum_{j=k-\delta^k}^{k-1}  \nabla f\(x^{j-\delta_j}\) } } + 2 \delta^k \gamma^2 \sigma^2\\
    & \leq 2 \delta^k \gamma^2 \sum_{j=k-\delta^k}^{k-1} \E{ \sqnorm{\nabla f\(x^{j-\delta_j}\)}} + 2 \delta^k \gamma^2 \sigma^2 ~.
\end{align*}
Using that $\gamma \leq \nicefrac{1}{2RL}$ and $\delta^k \leq R - 1$ (by the design), we obtain:
\begin{align*}
    \E{ \sqnorm{x^k - x^{k-\delta^k}}} 
    & \leq \frac{1}{2 L^2 R} \sum_{j=k-\delta^k}^{k-1} \E{ \sqnorm{\nabla f\(x^{j-\delta_j}\)}} + \frac{\gamma}{L} \sigma^2 ~.
\end{align*}
Next, summing over all iterations $k = 0, \dots, K$, we get:
\begin{align*}
    \sum_{k=0}^K \E{ \sqnorm{x^k - x^{k-\delta^k}}} 
    & \leq \frac{1}{2 L^2 R} \sum_{k=0}^K \sum_{j=k-\delta^k}^{k-1} \E{ \sqnorm{\nabla f\(x^{j-\delta_j}\)} } + \(K+1\) \frac{\gamma}{L} \sigma^2 ~.
\end{align*}
Observe that each squared norm $\| \nabla f (x^{j-\delta_j}) \|^2$ in the right-hand sums appears at most $R$ times due to the algorithms' design.
Specifically, $\delta^k \leq R - 1$ for all $k \geq 0$, ensuring no more than $R$ squared norms appear in the sums.
Therefore:
\begin{align*}
    \sum_{k=0}^K \E{ \sqnorm{x^k - x^{k-\delta^k}} }
    & \leq \frac{1}{2 L^2} \sum_{k=0}^K \E{ \sqnorm{\nabla f\(x^{k-\delta^k}\)} } + \(K+1\) \frac{\gamma}{L} \sigma^2.
\end{align*}
Finally, dividing the inequality by $K+1$ completes the proof.
\end{proof}
\section{Proof of \Cref{ringmaster:thm:optimal_ringmaster_dynamic}}
\label{ringmaster:sec:proof_dyn}
We restate below the theorem establishing the convergence of \ringmaster under the universal computation model and then provide its proof.
\begin{restate-boxedtheorem}{\ref{ringmaster:thm:optimal_ringmaster_dynamic}}
    Let Assumptions \ref{ringmaster:ass:lipschitz_constant}, \ref{ringmaster:ass:lower_bound}, and \ref{ringmaster:ass:stochastic_variance_bounded} hold.
    Let the stepsize in \ringmaster (\Cref{ringmaster:alg:ringmasternew} or \Cref{ringmaster:alg:ringmasternewstops}) be
    $$
        \gamma = \min \left\{ \frac{1}{2RL}, \frac{\varepsilon}{4L\sigma^2} \right\},
    $$ 
    and delay threshold 
    $$
        R = \max\left\{1,\left\lceil \frac{\sigma^2}{\varepsilon} \right\rceil\right\}.
    $$
    Then, under the \emph{universal computation model}, \ringmaster finds an $\varepsilon$--stationary point after at most $T^{\bar{K}}$ seconds, where 
    $$
    \bar{K} \eqdef \left\lceil\frac{48 L \Delta}{\epsilon}\right\rceil
    $$
    and $T^{\bar{K}}$ is the $\bar{K}$-th element of the following recursively defined sequence:
    \begin{align*}
        T^{k} \eqdef \min \left\{T \geq 0 : \sum \limits_{i=1}^{n} \flr{\frac{1}{4} \int_{T^{k - 1}}^{T} p_i(\tau) d \tau} \geq R\right\}
    \end{align*}
    for all $k \geq 1$ and $T^{0} = 0$.
\end{restate-boxedtheorem}
\begin{proof}
    From \Cref{ringmaster:thm:ringmaster_iteration}, the iteration complexity of \ringmaster is
    \begin{align}
    K = \left\lceil\frac{8 R L \Delta}{\epsilon} + \frac{16 \sigma^2 L \Delta}{\epsilon^2}\right\rceil.
    \end{align}
    Without loss of generality, we assume that
    $L \Delta > \nicefrac{\varepsilon}{2}$.\footnote{Otherwise, using $L$--smoothness, $\norm{\nabla f(x^0)}^2 \leq 2 L \Delta \leq \varepsilon,$ and the initial point is an $\varepsilon$--stationary point.}
    Thus, 
    $$
        \left\lceil\frac{8 R L \Delta}{\epsilon} + \frac{16 \sigma^2 L \Delta}{\epsilon^2}\right\rceil \leq \frac{16 R L \Delta}{\epsilon} + \frac{32 \sigma^2 L \Delta}{\epsilon^2}
    $$
    and
    \begin{align*}
        K \leq R \times \left\lceil\frac{K}{R}\right\rceil = R \times \left\lceil\frac{16 L \Delta}{\epsilon} + \frac{32 \sigma^2 L \Delta}{R \epsilon^2}\right\rceil.
    \end{align*}
    Using the choice of $R,$ we get
    \begin{align*}
        K \leq R \times \left\lceil\frac{48 L \Delta}{\epsilon}\right\rceil.
    \end{align*}
    In total, the algorithms will require $\left\lceil\nicefrac{48 L \Delta}{\epsilon}\right\rceil$ by $R$ consecutive updates of $x^k$ to find an $\varepsilon$--stationary point.
    Let us define $\bar{K} \eqdef \left\lceil\nicefrac{48 L \Delta}{\epsilon}\right\rceil$.
    Using \Cref{ringmaster:lem:time_R_dyn}, we know that \ringmaster requires at most 
    $$
        T^{1} \eqdef T(R,0) 
            = \min \left\{T \geq 0 : \sum_{i=1}^{n} \flr{\frac{1}{4} \int_{0}^{T} v_i(\tau) d \tau} \geq R\right\}
    $$
    seconds to finish the first $R$ consecutive updates of the iterates.
    Since the algorithms will finish the \emph{first} $R$ consecutive updates after at most $T^{1}$ seconds, they will start the iteration $R + 1$ before time $T^{1}$.
    Thus, using \Cref{ringmaster:lem:time_R_dyn} again, they will require at most 
    $$
        T^{2} \eqdef T(R,T^{1}) = \min \left\{T \geq 0 : \sum_{i=1}^{n} \flr{\frac{1}{4} \int_{T^{1}}^{T} v_i(\tau) d \tau} \geq R\right\}
    $$
    seconds to finish the first $2 \times R$ consecutive updates.
    Using the same reasoning, they will finish the first $\left\lceil \nicefrac{48 L \Delta}{\epsilon} \right\rceil \times R$ consecutive updates after at most
    $$
        T_{\bar{K}} \eqdef T(R,T_{\bar{K} - 1}) = \min \left\{T \geq 0 : \sum_{i=1}^{n} \flr{\frac{1}{4} \int_{T_{\bar{K} - 1}}^{T} v_i(\tau) d \tau} \geq R\right\}
    $$
    seconds.
\end{proof}

\section{Derivations for the Example from Section~\ref{ringmaster:sec:prel}}
\label{ringmaster:sec:deriv}

Let $\tau_i = \sqrt{i}$ for all $i \in [n],$ then
$$
    T_{\textnormal{R}} = \Theta\left(\min\limits_{m \in [n]} \left(\frac{1}{m}\sum\limits_{i=1}^m \frac{1}{\sqrt{i}}\right)^{-1} \left(\frac{L \Delta}{\varepsilon} + \frac{\sigma^2 L \Delta}{m \varepsilon^2}\right)\right).
$$
Using 
$$
    \sum_{i=1}^m \frac{1}{\sqrt{i}} = \Theta\left(\sqrt{m}\right)
$$ 
for $m \geq 1$, we simplify the term:
$$
    \left(\frac{1}{m}\sum\limits_{i=1}^m \frac{1}{\sqrt{i}}\right)^{-1} = \Theta\(\sqrt{m}\)~,
$$
$$
    T_{\textnormal{R}} = \Theta \left(\min\limits_{m \in [n]} \sqrt{m} \left(\frac{L \Delta}{\varepsilon} + \frac{\sigma^2 L \Delta}{m \varepsilon^2}\right) \right) = \Theta \left(\min\limits_{m \in [n]} \left(\frac{L \Delta \sqrt{m}}{\varepsilon} + \frac{\sigma^2 L \Delta}{\sqrt{m} \varepsilon^2}\right)\right).
$$
The minimum is achieved when the two terms are balanced, i.e., at 
$$
    m = \min\left\{\left\lceil \frac{\sigma^2}{\varepsilon} \right\rceil, n \right\}.
$$ 
Substituting this value of $m$, we obtain:
$$
    T_{\textnormal{R}} = \Theta \left(\max\left[\frac{\sigma L \Delta}{\varepsilon^{3/2}}, \frac{\sigma^2 L \Delta}{\sqrt{n} \varepsilon^2}\right]\right).
$$
We now consider $T_{\textnormal{A}}$:
$$
    T_{\textnormal{A}} = \Theta\left(\left(\frac{1}{n} \sum\limits_{i=1}^{n} \frac{1}{\tau_{i}}\right)^{-1} \left(\frac{L \Delta}{\varepsilon} + \frac{\sigma^2 L \Delta}{n \varepsilon^2}\right)\right).
$$
Using 
$$
    \sum_{i=1}^n \frac{1}{\sqrt{i}} = \Theta\left(\sqrt{n}\right)
$$
for $n \geq 1$, we simplify the term:
$$
    \left(\frac{1}{n} \sum\limits_{i=1}^n \frac{1}{\sqrt{i}}\right)^{-1} = \Theta\left(\sqrt{n}\right).
$$
Substituting this result into $T_{\textnormal{A}}$, we have:
\begin{align*}    
    T_{\textnormal{A}} 
    &= \Theta\left(\sqrt{n} \left(\frac{L \Delta}{\varepsilon} + \frac{\sigma^2 L \Delta}{n \varepsilon^2}\right)\right) = \Theta\left(\frac{L \Delta \sqrt{n}}{\varepsilon} + \frac{\sigma^2 L \Delta}{\sqrt{n} \varepsilon^2}\right) \\
    &= \Theta\left(\max\left[\frac{L \Delta \sqrt{n}}{\varepsilon}, \frac{\sigma^2 L \Delta}{\sqrt{n} \varepsilon^2}\right]\right).
\end{align*}
\section{When the Initial Point is an $\varepsilon$--Stationary Point}
\label{ringmaster:sec:l_init}
Under the assumption of $L$--smoothness (Assumption~\ref{ringmaster:ass:lipschitz_constant}), we have:
\begin{align*}
    f(y) \leq f(x) + \inp{\nabla f(x)}{y - x} + \frac{L}{2} \norm{y - x}^2
\end{align*}
for all $x, y \in \R^d$.
Taking $y = x - \frac{1}{L} \nabla f(x),$
\begin{align*}
    f\(x - \frac{1}{L} \nabla f(x)\) \leq f(x) - \frac{1}{2 L} \norm{\nabla f(x)}^2.
\end{align*}
Since $f\(x - \frac{1}{L} \nabla f(x)\) \geq f^*$ and taking $x = x^0,$ we get
\begin{align*}
    \norm{\nabla f(x^0)}^2 \leq 2 L \Delta~.
\end{align*}
Thus, if $2 L \Delta \leq \varepsilon$, then $\norm{\nabla f(x^0)}^2 \leq \varepsilon$.
            \clearpage
\newpage

\chapter{Appendix for \Cref{chapter:Ringleader}}
\label{appendix_ringleader}
\thispagestyle{empty}
\section{Arbitrarily Changing Computation Times}\label{sec:arbitrary_time}
In practice, the \emph{fixed computation model} \eqref{eq:fixed_time} is often not satisfied.
The compute power of devices can vary over time due to temporary disconnections, hardware or network delays, fluctuations in processing capacity, or other transient effects \citep{maranjyan2025mindflayer}.

In this section we extend our theory to the more general setting of arbitrarily varying computation times.

\subsection{Universal Computation Model}
\label{sec:universal_computation_model}
To formalize this setting, we adopt the \emph{universal computation model} introduced by \citet{tyurin2024tighttimecomplexitiesparallel}.

For each worker $i \in [n]$, we define a \emph{compute power} function
$$
    p_i : \R_{+} \to \R_{+}~,
$$
assumed nonnegative and continuous almost everywhere (countably many jumps allowed).
For any $T_2 \ge T_1 \ge 0$, the number of stochastic gradients \emph{completed} by worker $i$ on $[T_1, T_2]$ is
$$
    \#\text{gradients in }[T_1, T_2] \;=\; \left\lfloor \int_{T_1}^{T_2} p_i(t)\,dt \right\rfloor.
$$
Here, $p_i(t)$ models the worker's time-varying computational ability: smaller values over an interval yield fewer completed gradients, and larger values yield more.

For instance, if worker $i$ remains idle for the first $T$ seconds and then becomes active, this corresponds to $p_i(t) = 0$ for $t \leq T$ and $p_i(t) > 0$ for $t > T$.
More generally, $p_i(t)$ may follow periodic or irregular patterns, leading to bursts of activity, pauses, or chaotic changes in compute power.
The process $p_i(t)$ may even be random, and all results hold conditional on the realized sample paths of $\{p_i\}$.

The \emph{universal computation model} reduces to the \emph{fixed computation model} \eqref{eq:fixed_time} when $p_i(t) = \nicefrac{1}{\tau_i}$ for all $t \geq 0$ and $i \in [n]$.
In this case,
$$
    \#\text{gradients in }[T_1, T_2] = \left\lfloor \frac{T_2 - T_1}{\tau_i} \right\rfloor,
$$  
meaning that worker $i$ computes one stochastic gradient after $T_1 + \tau_i$ seconds, two gradients after $T_1 + 2\tau_i$ seconds, and so on.

\subsection{Toward an Optimal Method}

In the general setting of arbitrarily varying computation times, \Cref{algo:Ringleader} is not optimal.
To see why, consider the following adversarial timing pattern.

Suppose there are two workers.
During one gradient computation by the slower worker, the faster worker computes $s$ gradients.
Immediately afterwards, they switch roles: the previously fast worker slows down by a factor of $s$, while the previously slow one speeds up by the same factor.
This pattern repeats each time the slower worker finishes a gradient computation.

In this setting, if we run \Cref{algo:Ringleader}, the server waits in each Phase~1 for a single gradient from every worker.
Thus, the slower worker always contributes only one gradient, and the harmonic mean of the batch sizes satisfies 
$$
    1 \;\le\; B^k \;\le\; 2~.
$$ 
From \Cref{theorem:convergence}, the iteration complexity is
$$
    \cO\!\left(
        \frac{nL\Delta}{\varepsilon} \left( 1 + \frac{\sigma^2}{Bn\varepsilon} \right)
    \right).
$$
When $\nicefrac{\sigma^2}{n \varepsilon}$ is much larger than $B$, this dependence can be highly suboptimal.

Instead, suppose the server waits until one full round of the above process completes, collecting $s+1$ gradients from each worker.
Then the harmonic mean satisfies $B^k \ge s+1$, which can be arbitrarily larger than~2.
Since in practice both $s$ and $\nicefrac{\sigma^2}{n\varepsilon}$ can be very large, the naive strategy of waiting for only one gradient per worker (as in \Cref{algo:Ringleader}) cannot be optimal in the arbitrary-time setting.

\subsection{An Optimal Method}

The solution is simple and follows directly from the iteration complexity bound.
From
$$
    \cO\!\left(
        \frac{nL\Delta}{\varepsilon} \left( 1 + \frac{\sigma^2}{Bn\varepsilon} \right)
    \right),
$$
we see that to balance the terms it suffices to ensure
$$
    B \;\ge\; \frac{\sigma^2}{n\varepsilon}~.
$$
Accordingly, we modify the stopping condition in Phase~1 of \Cref{algo:Ringleader}.
Instead of requiring the server to receive at least one gradient from each worker, we require the stronger condition used in \malenia, namely
\begin{equation}\tag{\ref{eq:malenia_condition}}
    \left(\frac{1}{n} \sum_{i=1}^n \frac{1}{b_i} \right)^{-1} 
    \;\;\ge\;\; \max\left\{1, \frac{\sigma^2}{n\varepsilon}\right\}~,
\end{equation}
where $b_i$ is the number of gradients received from worker~$i$.

In the low-noise regime, where $\nicefrac{\sigma^2}{n\varepsilon} \le 1$, the condition reduces to requiring $b_i \ge 1$ for all $i$, so the algorithm coincides with the original \Cref{algo:Ringleader}.
In the high-noise regime, the algorithm collects more gradients in Phase~1, ensuring that $B$ is sufficiently large for optimal convergence.

With this change, Phase~1 of our algorithm matches that of \malenia.
The difference lies in Phase~2: our algorithm continues to use the ongoing gradient computations from all workers to perform $n$ updates, while \malenia discards any unfinished gradients, performs a single update, and then proceeds to the next round.

The following theorem establishes the time complexity of our algorithm under the universal computation model.
\begin{boxedtheorem}
    Under Assumptions~\ref{ass:stochastic_variance_bounded}, \ref{ass:lipschitz_constant}, and \ref{ass:lower_bound}, let the stepsize in \ringleader be
    $$
        \gamma = \frac{1}{10nL}~.
    $$
    Then, under the \emph{universal computation model}, \ringleader finds an $\varepsilon$--stationary point within at most $T^{\bar{K}}$ seconds, where
    $$
        \bar{K} \;\eqdef\; \left\lceil \frac{160 L \Delta}{\varepsilon} \right\rceil,
    $$
    and $T^{\bar{K}}$ denotes the $\bar{K}$-th element of the recursively defined sequence
    \begin{align*}
        T^k \;=\; \min \left\{ T \ge 0 :
            \left( \frac{1}{n}\sum_{i=1}^n 
            \left\lfloor \int_{T_{k-1}}^{T} p_i(t) \, dt \right\rfloor^{-1} \right)^{-1} 
            \;\;\ge\;\; \max\!\left\{1, \frac{\sigma^2}{n\varepsilon}\right\} 
        \right\}~,
    \end{align*}
    for all $k \ge 1$, with initialization $T^0 = 0$.
\end{boxedtheorem}
This result matches the lower bound derived by \citet{tyurin2024tighttimecomplexitiesparallel}, and therefore the proposed method is optimal.  
\begin{proof}
    Under the condition in \eqref{eq:malenia_condition}, each gradient-type step of the algorithm satisfies
    $$
        B^k = \left( \frac{1}{n} \sum_{i=1}^n \frac{1}{b_i^k} \right)^{-1}
        \;\;\ge\;\; \max\left\{1,\frac{\sigma^2}{n\varepsilon}\right\}.
    $$
    In \Cref{theorem:convergence}, instead of using $B$ we can substitute any valid lower bound.
    Here we choose
    $$
        B = \max\left\{1,\frac{\sigma^2}{n\varepsilon}\right\}.
    $$
    With this substitution, the iteration complexity becomes
    $$
        K = \frac{80\,nL\Delta}{\varepsilon}~.
    $$
    To derive the time complexity, consider the time required to perform $n$ iterations.
    Each block of $n$ updates occurs in Phase~2 following the Phase~1 gradient collection.
    Starting from time $T=0$, Phase~1 ends once the accumulated number of gradients satisfies condition~\eqref{eq:malenia_condition}, which occurs at time
    $$
        T_+^1 = \min \left\{T \ge 0 : 
            \left( \frac{1}{n}\sum_{i=1}^n 
            \left\lfloor \int_{0}^{T} p_i(t) \, dt \right\rfloor^{-1} \right)^{-1} 
            \;\;\ge\;\; \max\left\{1, \frac{\sigma^2}{n\varepsilon}\right\} 
        \right\}.
    $$
    After Phase~1, to complete $n$ updates in Phase~2 we must wait for the ongoing computations to finish.
    This requires at most
    $$
        T^1 = \min \left\{T \ge 0 : 
            \left( \frac{1}{n}\sum_{i=1}^n 
            \left\lfloor \int_{T_+^1}^{T} p_i(t) \, dt \right\rfloor^{-1} \right)^{-1} 
            \;\;\ge\;\; 1
        \right\}.
    $$
    Thus, the total time to complete all $K$ iterations is bounded by
    $$
        T^{\left\lceil \nicefrac{2K}{n} \right\rceil}~,
    $$
    where the sequence $\{T^k\}_{k \ge 0}$ is defined recursively as
    $$
        T^k = \min \left\{T \ge 0 : 
            \left( \frac{1}{n}\sum_{i=1}^n 
            \left\lfloor \int_{T_{k-1}}^{T} p_i(t) \, dt \right\rfloor^{-1} \right)^{-1} 
            \;\;\ge\;\; \max\left\{1, \frac{\sigma^2}{n\varepsilon}\right\} 
        \right\}~.
    $$
    and $T^0 = 0$.
\end{proof}
\section{Auxiliary Lemmas}
Here we provide proofs of lemmas omitted from the main text, along with auxiliary results that will be used later.
\subsection{Proof of \Cref{lemma:smoothness_relation}}
\label{proof:smoothness}
We begin with a lemma relating the different smoothness constants.
\begin{restate-boxedlemma}{\ref{lemma:smoothness_relation}}[Smoothness Bounds]
    Let $L_f$ denote the smoothness constant of $f$, $L_{f_i}$ the smoothness constant of $f_i$, and $L$ the constant from \Cref{ass:lipschitz_constant}.
    We have
    $$
        L_f \;\le\; L \;\le\; \sqrt{\frac{1}{n}\sum_{i=1}^n L_{f_i}^2} \;\le\; \max_{i \in [n]} L_{f_i} =: L_{\max}~.
    $$
    Moreover, if all $f_i$ are identical, i.e., $f_i = f$ for all $i \in [n]$, then $L = L_f$.
\end{restate-boxedlemma}
Recall from \Cref{ass:lipschitz_constant} that we assumed the following generalized smoothness condition:  
for some constant $L>0$ and for all $x \in \R^d$ and $y_1,\dots,y_n \in \R^d$,  
\begin{align}\label{eq:smoothness}
    \sqnorm{\nabla f(x) - \frac{1}{n}\sum_{i=1}^n \nabla f_i(y_i)}
    \;\le\; \frac{L^2}{n} \sum_{i=1}^n \sqnorm{x - y_i}~.
\end{align}
Recall that a function $\phi$ is called $L_\phi$--smooth if
$$
    \norm{\nabla \phi(x) - \nabla \phi(y)} \le L_\phi \norm{x - y}, \quad \forall x,y \in \R^d~.
$$
Here $L_\phi$ denotes the minimal such constant.
We are ready to prove the lemma.
\begin{proof}
For the first inequality, take $y_1 = \dots = y_n = y$.  
Then \eqref{eq:smoothness} reduces to
$$
    \sqnorm{\nabla f(x) - \nabla f(y)} \le L^2 \sqnorm{x-y},
$$
so $f$ is $L$--smooth.
By definition of $L_f$ as the minimal smoothness constant, this implies $L_f \le L$.

For the second inequality, by the triangle inequality, then by the smoothness of each $f_i$, and finally by Cauchy--Schwarz,
\begin{align*}
    \norm{\nabla f(x) - \frac{1}{n}\sum_{i=1}^n \nabla f_i(y_i)}
        &\le \frac{1}{n}\sum_{i=1}^n \norm{\nabla f_i(x) - \nabla f_i(y_i)}
        \le \frac{1}{n}\sum_{i=1}^n L_{f_i} \norm{x - y_i} \\
        &\le \sqrt{ \frac{1}{n}\sum_{i=1}^n L_{f_i}^2 } \ \sqrt{ \frac{1}{n}\sum_{i=1}^n \sqnorm{x-y_i} }~.
\end{align*}
Squaring both sides shows that \eqref{eq:smoothness} holds with 
$L = \sqrt{ \frac{1}{n}\sum_{i=1}^n L_{f_i}^2 }$~.

Finally, suppose all $f_i$ are identical: $f_i \equiv f$ for all $i$.
Then
\begin{align*}
    \norm{\nabla f(x) - \frac{1}{n}\sum_{i=1}^n \nabla f(y_i)}
    &\le \frac{1}{n}\sum_{i=1}^n \norm{\nabla f(x) - \nabla f(y_i)}
    \le \frac{L_f}{n}\sum_{i=1}^n \norm{x - y_i} \\
    &\le L_f \sqrt{\frac{1}{n}\sum_{i=1}^n \sqnorm{x - y_i}}~,
\end{align*}
where the last step uses Cauchy--Schwarz.
Squaring both sides yields
$$
    \sqnorm{\nabla f(x) - \frac{1}{n}\sum_{i=1}^n \nabla f(y_i)}
    \le \frac{L_f^2}{n}\sum_{i=1}^n \sqnorm{x - y_i},
$$
i.e., \eqref{eq:smoothness} holds with $L \le L_f$.
Combined with $L_f \le L$, we conclude $L = L_f$.
\end{proof}

\subsection{Variance Term}
The next lemma provides an upper bound on the variance of the gradient estimator used in \ringleader.
\begin{boxedlemma}[Variance Bound]\label{lemma:variance}
    Under \Cref{ass:stochastic_variance_bounded}, the gradient estimator used in \Cref{algo:Ringleader} satisfies the following variance-type inequality
    $$
        \E{\sqnorm{ \bar g^k - \frac{1}{n}\sum_{i=1}^n \nabla f_i\(x^{k-\delta_i^k}\) }}
        \le \frac{\sigma^2}{B^k n}~.
    $$
\end{boxedlemma}
\begin{proof}
    Recall that the gradient estimator is defined as
    $$
        \bar g^k
        = \frac{1}{n} \sum_{i=1}^n \bar g_i^k
        = \frac{1}{n} \sum_{i=1}^n g_i^{k,j}
        = \frac{1}{n} \sum_{i=1}^n \nabla f_i\(x^{k-\delta_i^k};\xi_i^{k-\delta_i^k, j}\)~.
    $$
    Let $\mathcal{F}^k$ denote the sigma-field containing all randomness up to the start of the current round, i.e., up to iteration $k-(k \bmod n)$.
    Conditioning on $\mathcal{F}^k$, the evaluation points $x^{k-\delta_i^k}$ are deterministic, and the stochastic gradients $g_i^{k,j}$ are independent across both workers $i$ and samples $j$.

    Using the law of total expectation and the independence of stochastic gradients, we have
    \begin{align*}
        \E{\sqnorm{ \bar g^k - \frac{1}{n}\sum_{i=1}^n \nabla f_i\(x^{k-\delta_i^k}\) }}
        &= \E{\ExpCond{\sqnorm{ \bar g^k - \frac{1}{n}\sum_{i=1}^n \nabla f_i\(x^{k-\delta_i^k}\)}}{\mathcal{F}^k}} \\
        &= \E{\frac{1}{n^2} \sum_{i=1}^{n} \ExpCond{\sqnorm{ \bar g_i^k - \nabla f_i\(x^{k-\delta_i^k}\)}}{\mathcal{F}^k}}.
    \end{align*}
    For each worker $i$, the conditional variance of the minibatch gradient estimator is
    \begin{align*}
        \ExpCond{\sqnorm{ \bar g_i^k - \nabla f_i\(x^{k-\delta_i^k}\)}}{\mathcal{F}^k} 
        &= \ExpCond{\sqnorm{ \frac{1}{b_i^k} \sum_{j=1}^{b_i^k} g_i^{k,j} - \nabla f_i\(x^{k-\delta_i^k}\)}}{\mathcal{F}^k}  \\
        &= \frac{1}{b_i^k} \ExpCond{\sqnorm{ g_i^{k,1} - \nabla f_i\(x^{k-\delta_i^k}\)}}{\mathcal{F}^k}  \le \frac{\sigma^2}{b_i^k}~,
    \end{align*}
    where the last inequality follows from \Cref{ass:stochastic_variance_bounded}.
    
    Combining these results, we get
    \begin{align*}
        \E{\sqnorm{ \bar g^k - \frac{1}{n}\sum_{i=1}^n \nabla f_i\(x^{k-\delta_i^k}\) }}
        &\le \frac{1}{n^2} \sum_{i=1}^{n} \frac{\sigma^2}{b_i^k}
            = \frac{\sigma^2}{n} \frac{1}{n} \sum_{i=1}^n \frac{1}{b_i^k} 
            = \frac{\sigma^2}{B^k n}~,
    \end{align*}
    where the last equality uses the definition of the harmonic mean
    $$
        B^k = \(\frac{1}{n} \sum_{i=1}^n \frac{1}{b_i^k}\)^{-1} ~.
    $$
\end{proof}
\subsection{Proof of \Cref{lemma:descent}}
\label{proof:descent}
We now prove the descent lemma.
\begin{restate-boxedlemma}{\ref{lemma:descent}}[Descent Lemma]
    Under Assumptions~\ref{ass:stochastic_variance_bounded} and~\ref{ass:lipschitz_constant}, if the stepsize in \Cref{algo:Ringleader} satisfies $\gamma \le \nicefrac{1}{4L}$, then the following inequality holds
    \begin{align*}
        \E{f\(x^{k+1}\)}
        &\le \E{f\(x^{k}\)}
            - \frac{\gamma}{2} \E{ \sqnorm{\nabla f\(x^{k}\)} }
            - \frac{\gamma}{4} \E{\sqnorm{\frac{1}{n} \sum_{i=1}^n \nabla f_i\(x^{k-\delta_i^k}\)}} \\
            &\quad + \frac{\gamma L^2}{2n} \sum_{i=1}^n \E{\sqnorm{x^{k} - x^{k-\delta_{i}^k}}}
                    + \frac{3\gamma^2 L \sigma^2}{2B} \\
            &\quad + \gamma^2 L \sum_{\ell = k-(k \bmod n)}^{k-1} \E{\sqnorm{\frac{1}{n}\sum_{i=1}^n \nabla f_i\(x^{\ell-\delta_i^\ell}\)}}.
    \end{align*}
\end{restate-boxedlemma}
\begin{proof}
    Some proof techniques are adapted from the works of \citet{maranjyan2025ringmaster} and \citet{wang2025incremental}.
    
    From \Cref{ass:lipschitz_constant} and \Cref{lemma:smoothness_relation}, we know that $f$ is $L$--smooth. 
    Therefore, the following standard inequality holds \citep{nesterov2018lectures}
    \begin{equation}\label{eq:l-smoothness}
        \E{f(x^{k+1})}
        \le \E{ f(x^{k}) 
            - \gamma \inp{\nabla f(x^{k})}{\bar g^k}
            + \frac{L\gamma^2}{2} \sqnorm{\bar g^k} }.
    \end{equation}
    Recall that the gradient estimator is defined as
    $$
        \bar g^k
        = \frac{1}{n} \sum_{i=1}^n \bar g_i^k
        = \frac{1}{n} \sum_{i=1}^n \frac{1}{b_i^k} \sum_{j=1}^{b_i^k} g_i^{k,j} ~.
    $$
    Let $\mathcal{F}^k$ denote the sigma-field containing all randomness up to the start of the current Phase~2, i.e., up to iteration $k - (k \bmod n)$.
    A key observation is that all gradients in the current gradient table were computed and received during the current round.
    Since these gradients were computed at points from previous iterations within the current round, we have $k - \delta_i^k \le k - (k \bmod n)$ for all $i\in[n]$.
    Conditioning on $\mathcal{F}^k$, the points $x^{k-\delta_i^k}$ are deterministic.
    Therefore, we can compute the conditional expectation of the gradient estimator:
    $$
        \ExpCond{\bar g^k}{\mathcal{F}^k}
        = \frac{1}{n} \sum_{i=1}^n \frac{1}{b_i^k} \sum_{j=1}^{b_i^k} \ExpCond{g_i^{k,j}}{\mathcal{F}^k}
        = \frac{1}{n} \sum_{i=1}^n \nabla f_i\(x^{k-\delta_i^k}\).
    $$
    The last equality follows from the unbiasedness of the stochastic gradient estimator (\Cref{ass:stochastic_variance_bounded}).
    
    Using this conditional expectation and the law of total expectation, we can now simplify the inner product term in \eqref{eq:l-smoothness}:
    \begin{align*}
        &\E{\inp{\nabla f\(x^{k}\)}{\bar g^k}}
        = \E{\inp{\nabla f\(x^{k}\)}{\frac{1}{n} \sum_{i=1}^n \nabla f_i\(x^{k-\delta_i^k}\)}} \\
            &\quad\qquad\qquad + \E{\inp{\nabla f\(x^{k}\)}
            {\bar g^k - \frac{1}{n} \sum_{i=1}^n  \nabla f_i\(x^{k-\delta_i^k}\) }} \\
        &\qquad\qquad= 
        \E{\inp{\nabla f\(x^{k}\)}{\frac{1}{n} \sum_{i=1}^n \nabla f_i\(x^{k-\delta_i^k}\)}} \\
            &\quad\qquad\qquad + \E{\inp{\nabla f\(x^{k}\) - \nabla f\(x^{k-(k \bmod n)}\)}
            {\bar g^k - \frac{1}{n} \sum_{i=1}^n  \nabla f_i\(x^{k-\delta_i^k}\) }} \\
                &\qquad\qquad\qquad + \E{\inp{\nabla f\(x^{k-(k \bmod n)}\)}
                {\bar g^k - \frac{1}{n} \sum_{i=1}^n \nabla f_i\(x^{k-\delta_i^k}\) }} \\
        &\qquad\qquad= 
        \underbrace{ \E{\inp{\nabla f\(x^{k}\)}{\frac{1}{n} \sum_{i=1}^n \nabla f_i\(x^{k-\delta_i^k}\)}} }_{T_1} \\
            &\qquad\qquad\qquad + \underbrace{ \E{\inp{\nabla f\(x^{k}\) - \nabla f\(x^{k-(k \bmod n)}\)}
            {\bar g^k - \frac{1}{n} \sum_{i=1}^n  \nabla f_i\(x^{k-\delta_i^k}\) }} }_{T_2} .
    \end{align*}
    Next, using \Cref{ass:lipschitz_constant}, we have
    \begin{align*}
        2T_1
        &= \E{2 \inp{\nabla f\(x^{k}\)}{ \frac{1}{n} \sum_{i=1}^n \nabla f_i\(x^{k-\delta_i^k}\) }} \\
        &= \E{ \sqnorm{\nabla f\(x^{k}\)} 
            + \sqnorm{ \frac{1}{n} \sum_{i=1}^n \nabla f_i\(x^{k-\delta_i^k}\) } } \\
            &\quad - \E{\sqnorm{\nabla f\(x^{k}\) - \frac{1}{n} \sum_{i=1}^n \nabla f_i\(x^{k-\delta_i^k}\) }} \\
        &\ge \E{ \sqnorm{\nabla f\(x^{k}\)} }
            + \E{ \sqnorm{ \frac{1}{n} \sum_{i=1}^n \nabla f_i\(x^{k-\delta_i^k}\) } } \\
            &\quad - \frac{L^2 }{n} \sum_{i=1}^n \E{\sqnorm{x^{k} - x^{k-\delta_i^k}}}.
    \end{align*}
    Next, we analyze $T_2$
    \begin{align*}
        T_2
        &= \E{\inp{\nabla f\(x^{k}\) - \nabla f\(x^{k-(k \bmod n)}\)}
        {\bar g^k - \frac{1}{n} \sum_{i=1}^n  \nabla f_i\(x^{k-\delta_i^k}\) }} \\
        &\ge - \E{\norm{\nabla f\(x^{k}\) - \nabla f\(x^{k-(k \bmod n)}\)} 
            \norm{\bar g^k - \frac{1}{n} \sum_{i=1}^n  \nabla f_i\(x^{k-\delta_i^k}\)}} \\
        &\ge - L \E{\norm{x^{k} - x^{k-(k \bmod n)} }
            \norm{\bar g^k - \frac{1}{n} \sum_{i=1}^n  \nabla f_i\(x^{k-\delta_i^k}\)}} \\
        &= - L \E{\norm{\gamma \sum_{\ell = k-(k \bmod n)}^{k-1} \bar g^\ell }
            \norm{\bar g^k - \frac{1}{n} \sum_{i=1}^n  \nabla f_i\(x^{k-\delta_i^k}\)}} \\
        &\ge - L\gamma \sum_{\ell = k-(k \bmod n)}^{k-1} \E{\norm{ \bar g^\ell }
            \norm{\bar g^k - \frac{1}{n} \sum_{i=1}^n  \nabla f_i\(x^{k-\delta_i^k}\)}} \\
        &\ge - L\gamma \sum_{\ell = k-(k \bmod n)}^{k-1} \frac{1}{2}\(\E{\sqnorm{ \bar g^\ell }} 
            + \E{\sqnorm{\bar g^k - \frac{1}{n} \sum_{i=1}^n  \nabla f_i\(x^{k-\delta_i^k}\)}} \) \\
        &\ge - \frac{L\gamma}{2} \sum_{\ell = k-(k \bmod n)}^{k-1} \E{\sqnorm{ \bar g^\ell } }
            - (k \bmod n) \frac{L\gamma \sigma^2}{2B^k n} \\
        &\ge - \frac{L\gamma}{2} \sum_{\ell = k-(k \bmod n)}^{k-1} \E{\sqnorm{ \bar g^\ell } }
            - \frac{L\gamma \sigma^2}{2B^k}~.
    \end{align*}
    The inequalities follow from the Cauchy-Schwarz inequality, $L$--smoothness of $f$, the triangle inequality, Young's inequality, \Cref{lemma:variance}, and finally $(k \bmod n) \le n-1 < n$.

    It remains to bound the term $\E{\sqnorm{ \bar g^k } }$.
    Using Young's inequality, we have
    \begin{align*}
        \E{\sqnorm{\bar g^k}} 
        &\le 2\E{\sqnorm{\frac{1}{n}\sum_{i=1}^n \nabla f_i\(x^{k-\delta_i^k}\)}} 
            + 2\E{\sqnorm{ \bar g^k - \frac{1}{n}\sum_{i=1}^n \nabla f_i\(x^{k-\delta_i^k}\)}} \\
        &\le 2\E{\sqnorm{\frac{1}{n}\sum_{i=1}^n \nabla f_i\(x^{k-\delta_i^k}\)}} 
            + \frac{2 \sigma^2}{B^k n}~,
    \end{align*}
    where in the last step we used \Cref{lemma:variance}.

    Now, by combining all terms in \eqref{eq:l-smoothness}, we obtain
    \begin{align*}
        \E{f\(x^{k+1}\)}
        &\le \E{f\(x^{k}\)}
            - \frac{\gamma}{2} \E{ \sqnorm{\nabla f\(x^{k}\)} }
            - \frac{\gamma}{2} \E{\sqnorm{\frac{1}{n} \sum_{i=1}^n \nabla f_i\(x^{k-\delta_i^k}\)}} \\
            &\quad + \frac{\gamma L^2}{2n} \sum_{i=1}^n \E{\sqnorm{x^{k} - x^{k-\delta_{i}^k}}} \\
            &\quad + \frac{\gamma^2 L}{2} \sum_{\ell = k-(k \bmod n)}^{k-1} \E{\sqnorm{ \bar g^\ell } }
            + \frac{\gamma^2 L \sigma^2}{2 B^k} + \frac{\gamma^2L}{2} \E{\sqnorm{\bar g^k}} \\
        &\le \E{f\(x^{k}\)}
            - \frac{\gamma}{2} \E{ \sqnorm{\nabla f\(x^{k}\)} }
            - \frac{\gamma}{2} \E{\sqnorm{\frac{1}{n} \sum_{i=1}^n \nabla f_i\(x^{k-\delta_i^k}\)}} \\
            &\quad + \frac{\gamma L^2}{2n} \sum_{i=1}^n \E{\sqnorm{x^{k} - x^{k-\delta_{i}^k}}} \\
            &\quad + \gamma^2 L \sum_{\ell = k-(k \bmod n)}^{k-1} \E{\sqnorm{\frac{1}{n}\sum_{i=1}^n \nabla f_i\(x^{\ell-\delta_i^\ell}\)}} \\
            &\quad + \gamma^2L \E{\sqnorm{\frac{1}{n}\sum_{i=1}^n \nabla f_i\(x^{k-\delta_i^k}\)}} \\
            &\quad + \gamma^2 L \sum_{\ell = k-(k \bmod n)}^{k-1}\frac{\sigma^2}{B^\ell n} 
                + \frac{\gamma^2 L \sigma^2 }{2 B^k} 
                + \frac{\gamma^2L \sigma^2}{B^k n} ~.
    \end{align*}
    This completes the proof under the stepsize condition $\gamma \le \nicefrac{1}{4L}$ and $B \eqdef \inf_{k \ge 0} B^k$.
\end{proof}
\subsection{Proof of \Cref{lemma:residual}}
\label{proof:residual}
The following lemma provides an upper bound on the residual error due to delays.
\begin{restate-boxedlemma}{\ref{lemma:residual}}[Residual Estimation]
    Under \Cref{ass:stochastic_variance_bounded}, the iterates of \ringleader (\Cref{algo:Ringleader}) with stepsize $\gamma \le \nicefrac{1}{4nL}$ satisfy the following bound
    \begin{equation*}
    \frac{1}{K} \sum_{k=0}^{K-1} \frac{1}{n} \sum_{i=1}^n \E{\sqnorm{x^{k} - x^{k-\delta_i^k}}}
        \le \frac{2\gamma n}{LK}\sum_{k=0}^{K-1} \E{\sqnorm{ \frac{1}{n}\sum_{j=1}^n \nabla f_j\(x^{k-\delta_j^k}\) }}
            + \frac{2\gamma \sigma^2}{LB}~.
    \end{equation*}
\end{restate-boxedlemma}
\begin{proof}
    By Young's inequality, we have
    \begin{align*}
    \E{\sqnorm{x^{k} - x^{k-\delta_i^k}}}
        &= \E{\sqnorm{\gamma \sum_{\ell = k-\delta_i^k}^{k-1} \bar g^\ell}} \\
        &\le 2\gamma^2 \E{\sqnorm{\sum_{\ell = k-\delta_i^k}^{k-1} \frac{1}{n}\sum_{j=1}^n \nabla f_j\(x^{\ell-\delta_j^\ell}\) }} \\
        &\quad + 2\gamma^2 \E{\sqnorm{\sum_{\ell = k-\delta_i^k}^{k-1} \( \bar g^\ell - \frac{1}{n}\sum_{j=1}^n \nabla f_j\(x^{\ell-\delta_j^\ell}\) \)}} \\
        &\le 2\gamma^2 \underbrace{\delta_i^k \sum_{\ell = k-\delta_i^k}^{k-1} \E{\sqnorm{ \frac{1}{n}\sum_{j=1}^n \nabla f_j\(x^{\ell-\delta_j^\ell}\) }}}_{T_{ik}} + 2\gamma^2 (\delta_i^k)^2 \frac{\sigma^2}{Bn}~.
    \end{align*}
    In the last inequality, we used Jensen's inequality and \Cref{lemma:variance}.
    
    Next, we estimate the sum of $T_{ik}$
    \begin{align*}
    \sum_{k=0}^{K-1} \frac{1}{n} \sum_{i=1}^n T_{ik}
        &= \sum_{k=0}^{K-1} \frac{1}{n} \sum_{i=1}^n \delta_i^k \sum_{\ell = k-\delta_i^k}^{k-1} \E{\sqnorm{ \frac{1}{n}\sum_{j=1}^n \nabla f_j\(x^{\ell-\delta_j^\ell}\) }} \\
        &= \frac{1}{n} \sum_{i=1}^n \sum_{k=0}^{K-1} \delta_i^k \sum_{\ell = k-\delta_i^k}^{k-1} \E{\sqnorm{ \frac{1}{n}\sum_{j=1}^n \nabla f_j\(x^{\ell-\delta_j^\ell}\) }} \\
        &\le 2 \sum_{i=1}^n \sum_{k=0}^{K-1} \sum_{\ell = k-\delta_i^k}^{k-1} \E{\sqnorm{ \frac{1}{n}\sum_{j=1}^n \nabla f_j\(x^{\ell-\delta_j^\ell}\) }} \\
        &\le 2 \sum_{i=1}^n \delta_i^{\max} \sum_{k=0}^{K-1} \E{\sqnorm{ \frac{1}{n}\sum_{j=1}^n \nabla f_j\(x^{k-\delta_j^k}\) }} \\
        &\le 4n^2 \sum_{k=0}^{K-1} \E{\sqnorm{ \frac{1}{n}\sum_{j=1}^n \nabla f_j\(x^{k-\delta_j^k}\) }} .
    \end{align*}
    In the first and last inequality, we used the bound $\delta_i^{\max} \le 2n$ from \Cref{lemma:delay}.
    Finally, applying the stepsize condition $\gamma \le \nicefrac{1}{4nL}$ yields the result.
\end{proof}

\section{Time Complexity of \iasgdtitle}\label{sec:iasgd_time_complexity}
The iteration complexity of \iasgd \citep{wang2025incremental} is
$$
    K = \cO\!\left(\frac{\delta^{\max}L\Delta}{\varepsilon} \left(1 + \frac{\sigma^2}{n\varepsilon}\right)\right).
$$
We now analyze the corresponding wall-clock time under the \textit{fixed computation model} \eqref{eq:fixed_time}.
Since the algorithm performs an update whenever a single worker finishes a computation, we seek the minimal time $T$ such that
$$
    \sum_{i=1}^n\left\lfloor\frac{T}{\tau_i}\right\rfloor \;\ge\; K~.
$$
Observe that
$$
    \sum_{i=1}^n \frac{T}{\tau_i}
        \;\ge\; \sum_{i=1}^n\left\lfloor\frac{T}{\tau_i}\right\rfloor.
$$
Hence, if we define $T^\prime$ by
$$
    \sum_{i=1}^n \frac{T^\prime}{\tau_i} \;=\; K~,
$$
then 
$$
    T^\prime = \left(\sum_{i=1}^n \frac{1}{\tau_i}\right)^{-1} K~.
$$
It follows that the minimal time $T$ is necessarily larger than $T^\prime$.

It remains to bound $\delta^{\max}$.
At initialization, all workers start computing their first gradients simultaneously.
By the time the slowest worker completes its first gradient (at time $\tau_n$), the other workers may each have completed multiple gradients.
In particular,
$$
    \delta^{\max} \;\ge\; \sum_{i=1}^n \left\lfloor \frac{\tau_n}{\tau_i}\right\rfloor.
$$
Combining this with the iteration complexity bound, we obtain that the total runtime satisfies
$$
    T \;\ge\; c \times \frac{\tau_n L\Delta}{\varepsilon} \left(1 + \frac{\sigma^2}{n\varepsilon}\right),
$$
for some universal constant $c > 0$.

Note that the expression above should not be viewed as an exact upper bound on the runtime.
It is better understood as a simplified estimate of~$T$, which is sufficient for our purposes
and provides a cleaner basis for comparison.

\section{Improved \maleniatitle}\label{sec:malenia_param_free}
\malenia has the following iteration complexity \citep{tyurin2024optimal}
$$
    K \;\ge\; \frac{12\Delta L_f}{\varepsilon} \;+\; \frac{12 \Delta L_f \sigma^2}{\varepsilon^2 n S}~,
$$
where $S$ is a lower bound on the harmonic mean of the batch sizes, i.e.,
$$
    \left(\frac{1}{n}\sum_{i=1}^n \frac{1}{b_i^k}\right)^{-1} \;\ge\; S~,
$$
for all iterations $k$.  
In the original \malenia analysis \citep{tyurin2024optimal}, this bound follows from the condition in \eqref{eq:malenia_condition}, which fixes the same value of $S$ across all iterations.

In the fixed-time regime \eqref{eq:fixed_time}, however, this condition is no longer necessary.  
By adopting the same strategy as \ringleader (\Cref{algo:Ringleader})\textemdash namely, waiting for at least one gradient from each worker\textemdash we effectively replace $S$ with $B$ in the rate.  
This yields the following time complexity
$$
    \tau_n K
    \;=\; \frac{12\tau_n\Delta L_f}{\varepsilon}\left(1 + \frac{\sigma^2}{\varepsilon n B}\right).
$$
Substituting the expression for $B$ from \Cref{lemma:time_for_n_iter} and proceeding as in the proof of \Cref{thm:time_complexity}, we obtain the same overall time complexity as before\textemdash this time \emph{without} requiring condition \eqref{eq:malenia_condition}, which depends on knowing $\sigma$ and fixing $\varepsilon$ in advance.

Finally, note that this improvement is only valid in the fixed-time regime.  
In the setting with arbitrarily varying computation times, the same optimization cannot be applied, for the same reasons discussed for \ringleader in \Cref{sec:arbitrary_time}.
            \newpage

\chapter{Appendix for \Cref{chapter:ATA}}
\label{appendix_ata}
\thispagestyle{empty}

\section{Additional Experiments}
\label{ata:section:additional_experiments}

The objective function is a convex quadratic function $f \,:\, \R^d \to \R$ defined as  
$$
    f(x) = \frac{1}{2} x^\top A x - b^\top x~,
$$
where
\begin{align*}
    A = \frac{1}{4}
    \begin{bmatrix}
    2 & -1 &  & 0 \\
    -1 & \ddots & \ddots &  \\
    & \ddots & \ddots & -1 \\
    0 & & -1 & 2 \\
    \end{bmatrix}
    \in \R^{d \times d} ~,
    \quad \text{and} \quad 
    b = \frac{1}{4}
    \begin{bmatrix}
    -1 \\
    0 \\
    \vdots \\
    0 \\
    \end{bmatrix}
    \in \R^d~.
\end{align*}
We denote $f^*$ as the minimum value of the function $f$.
Each of the $n$ workers is able to calculate unbiased stochastic gradients $g(x)$ that satisfy
$$
    \E{\sqnorm{g(x) - \nabla f(x)}} \le 0.01^2~.
$$
This is achieved by adding Gaussian noise to the gradients of $f$.

The experiments were implemented in Python.
The distributed environment was emulated on machines with Intel(R) Xeon(R) Gold 6248 CPU @ 2.50GHz.
\subsection{Linear Noise}
\label{ata:sec:linear_noise}
In this section we model the computation time for worker $i$ by the distribution
$$
    \nu_i = 29i + \mathrm{Exp}(29i), \quad \text{for all} \quad i \in [n]~.
$$
The expected value of this distribution is $\mu_i = 2 \cdot 29 i$~.
Furthermore, the Orlicz norm satisfies the bound $\alpha_i \le 2\mu_i$.

We again set $B = 23$ and run simulations similar to those in \Cref{ata:section:experiments}.
The results are shown in \Cref{ata:fig:linear}.
\begin{figure*}[thb]
    \centering
	{ %
    \setlength{\tabcolsep}{-4pt} 
    \renewcommand{\arraystretch}{0} 
    \begin{tabular}{@{}cccc@{}}
        \includegraphics[width=0.27\textwidth]{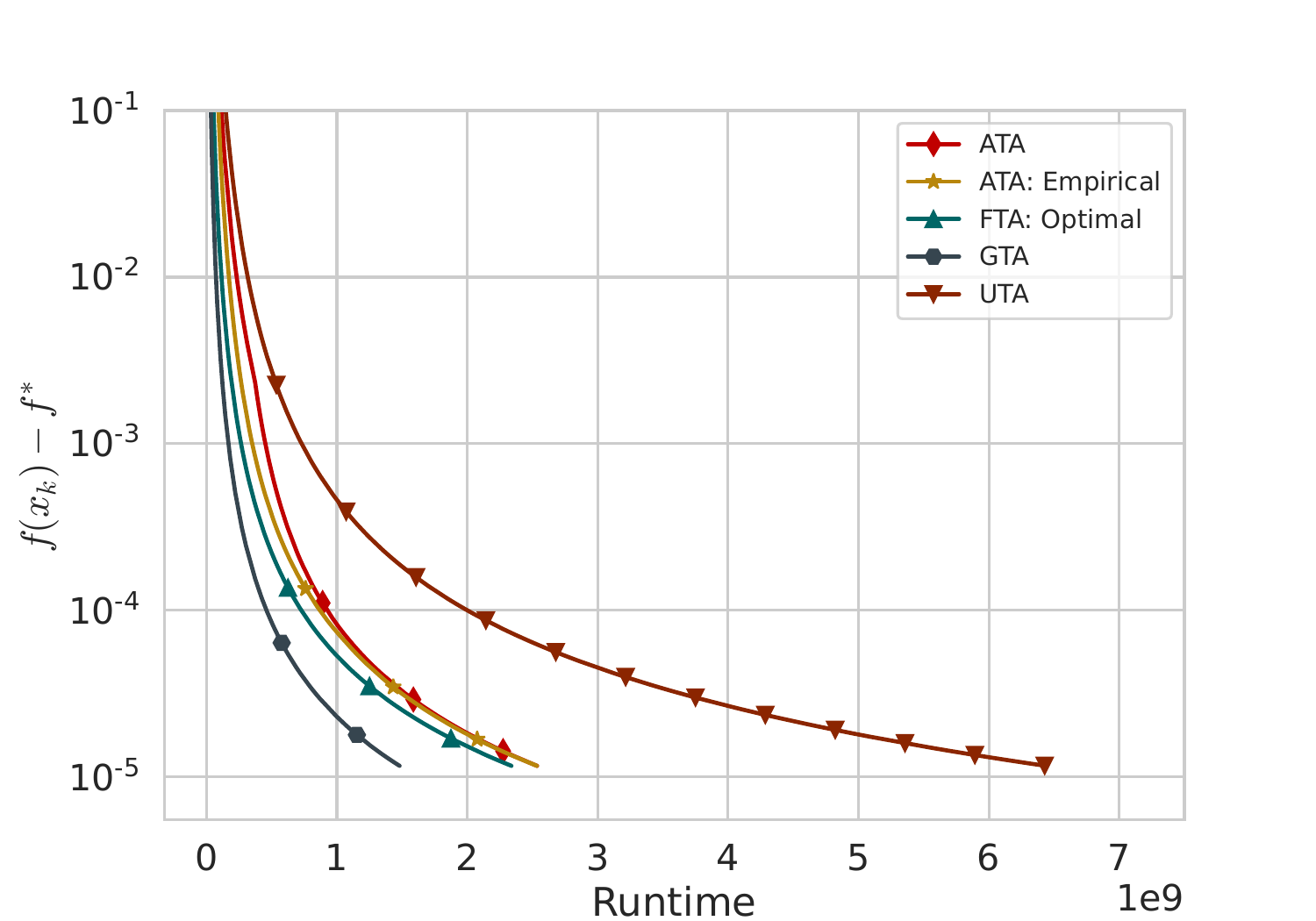} &
        \includegraphics[width=0.27\textwidth]{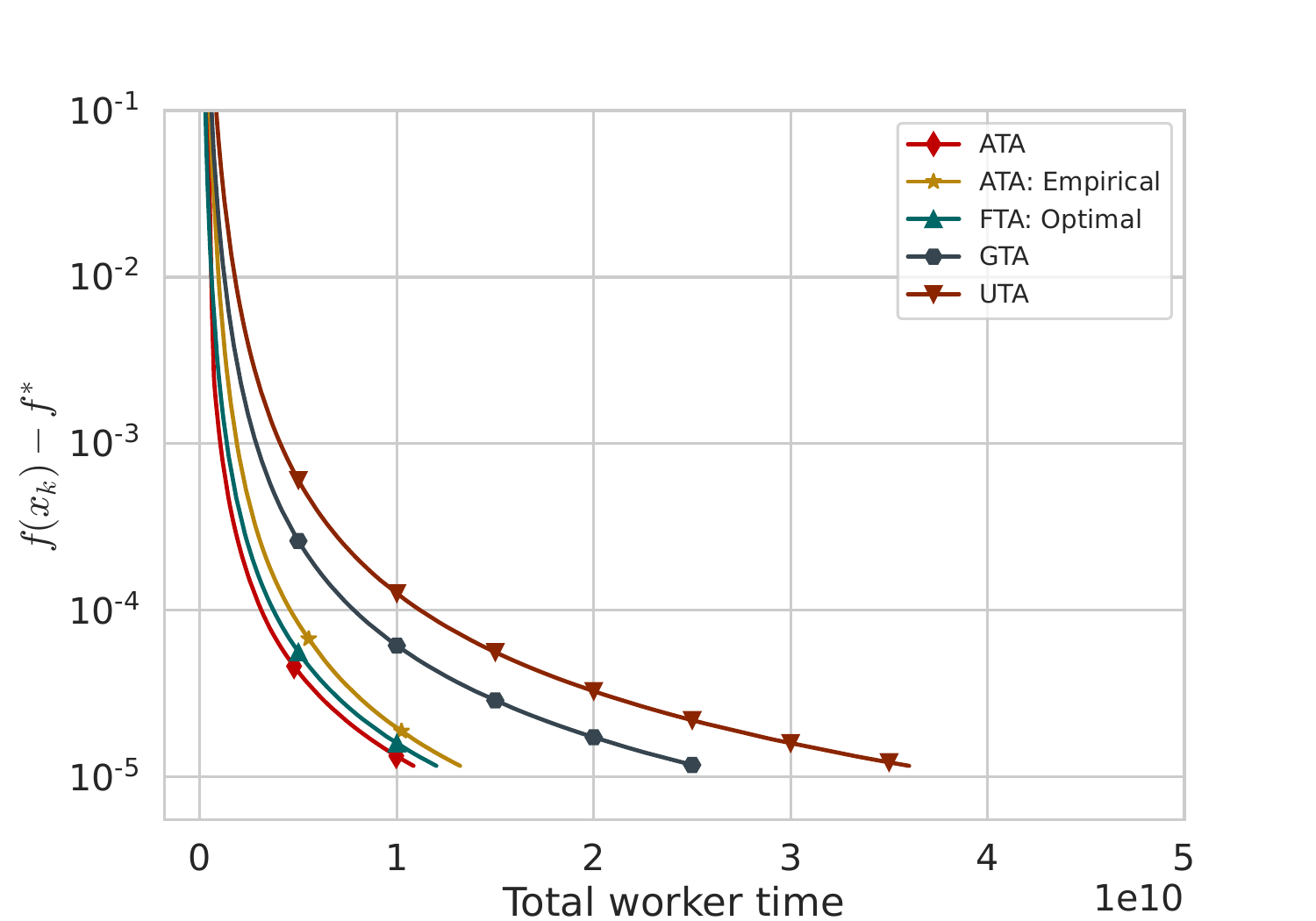} &
        \includegraphics[width=0.27\textwidth]{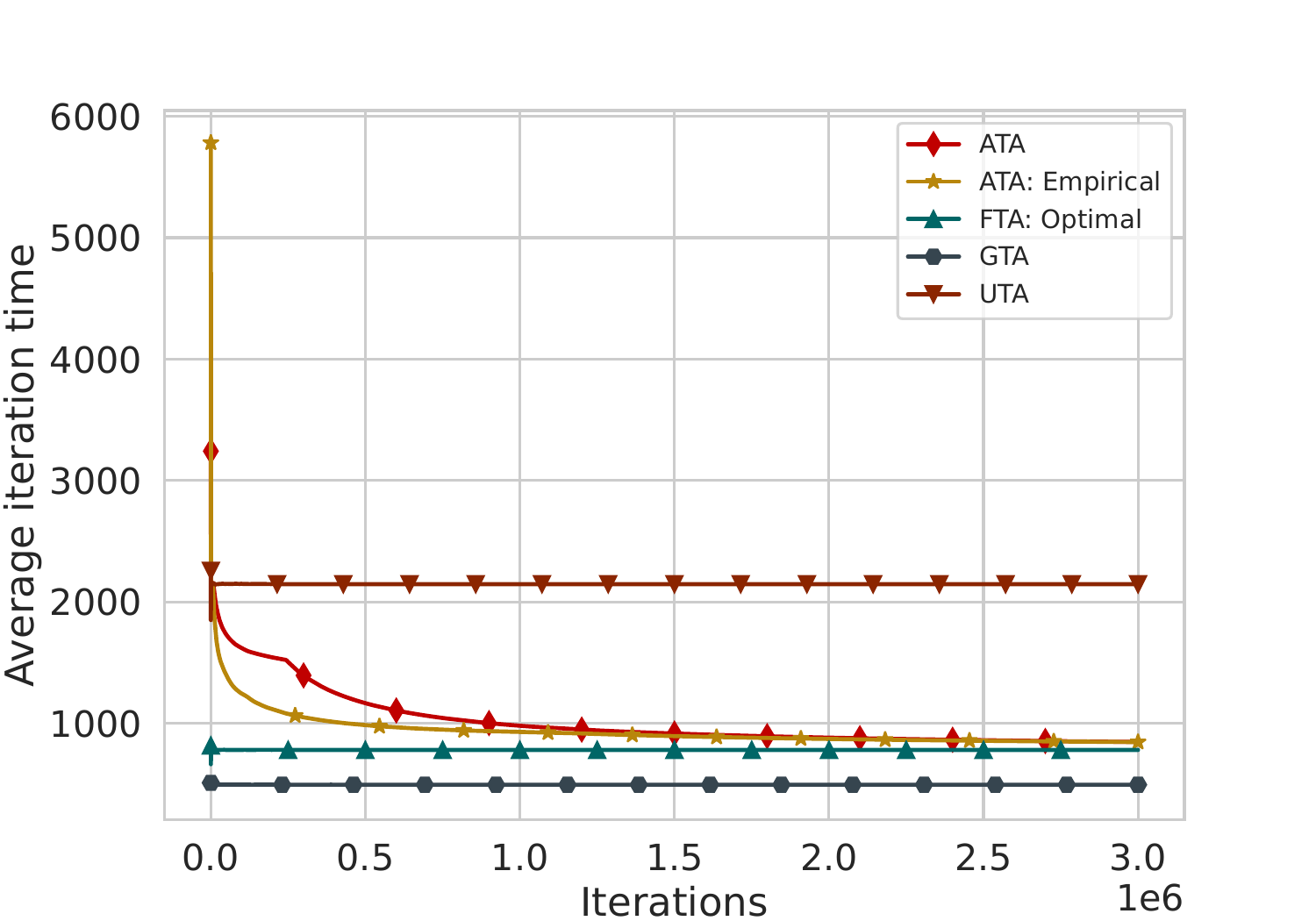} &
        \includegraphics[width=0.27\textwidth]{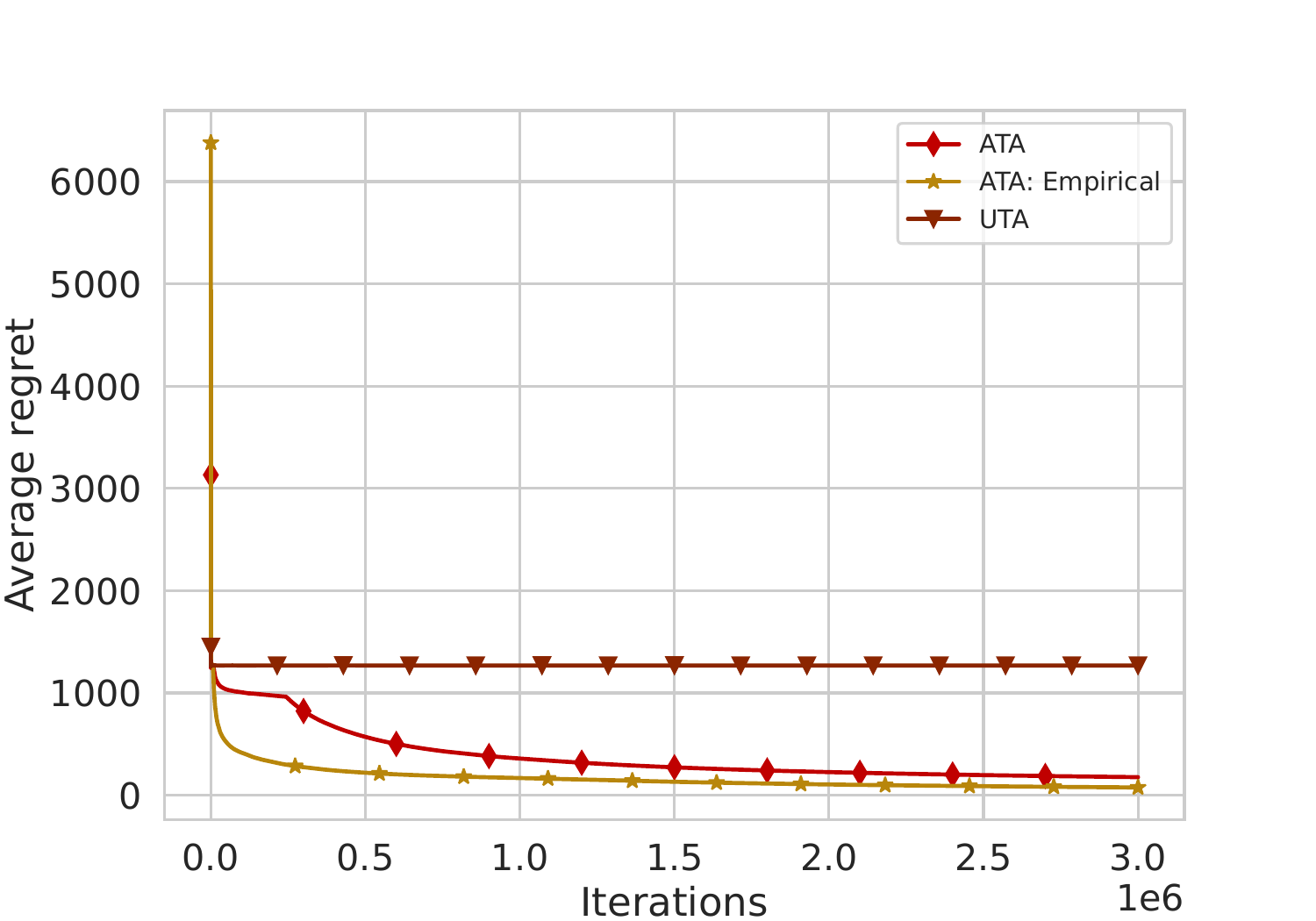} \\
        \includegraphics[width=0.27\textwidth]{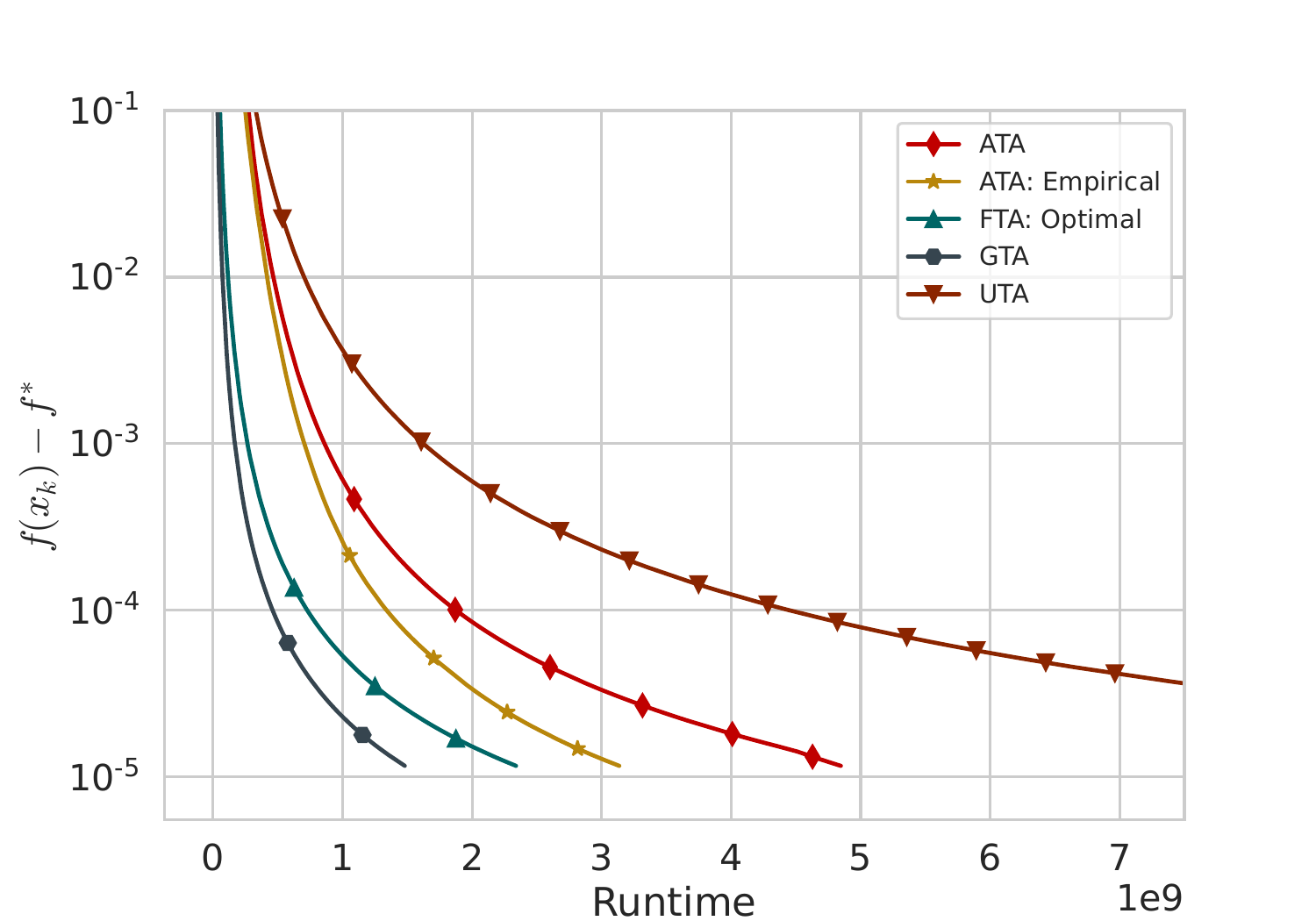} &
        \includegraphics[width=0.27\textwidth]{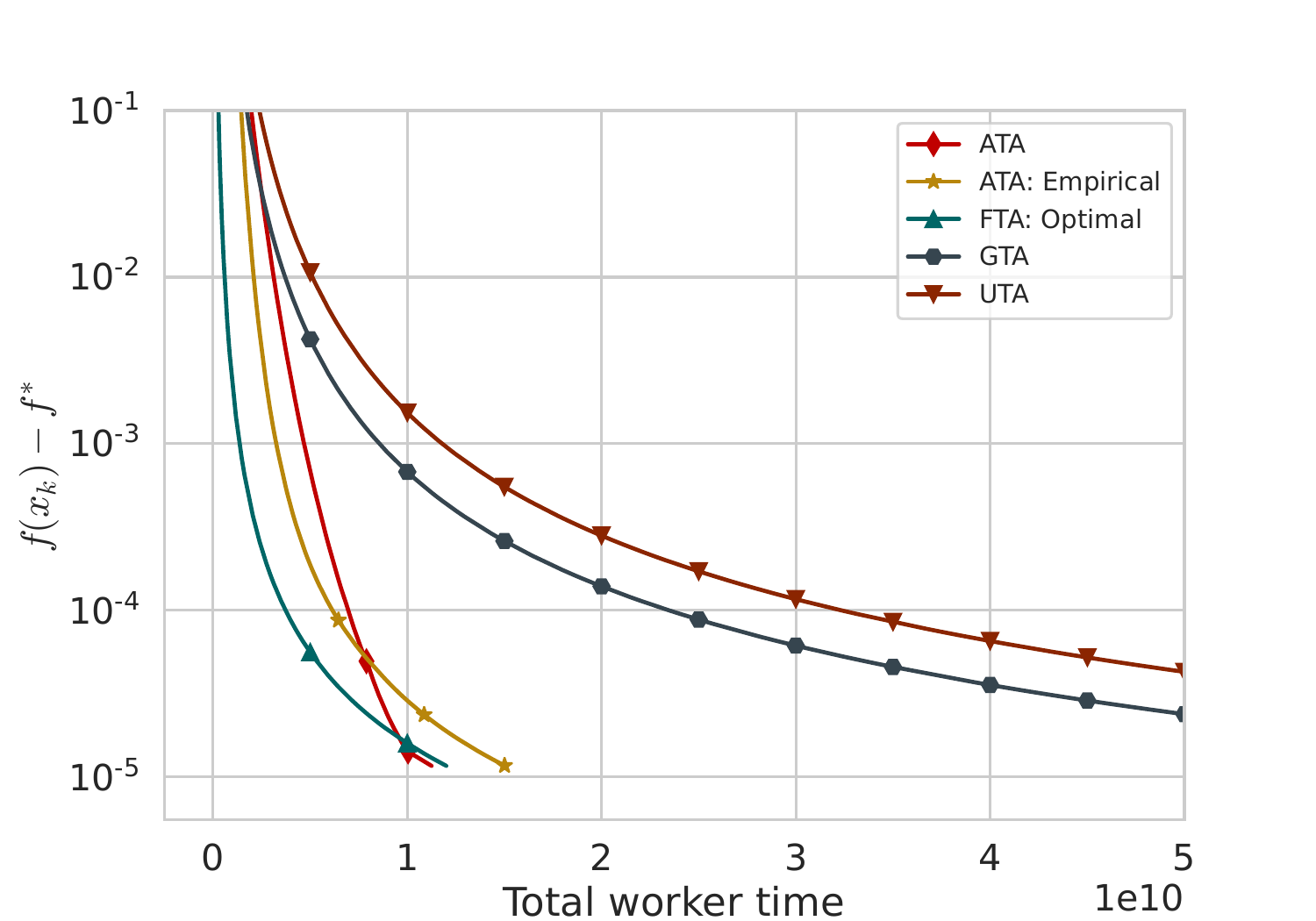} &
        \includegraphics[width=0.27\textwidth]{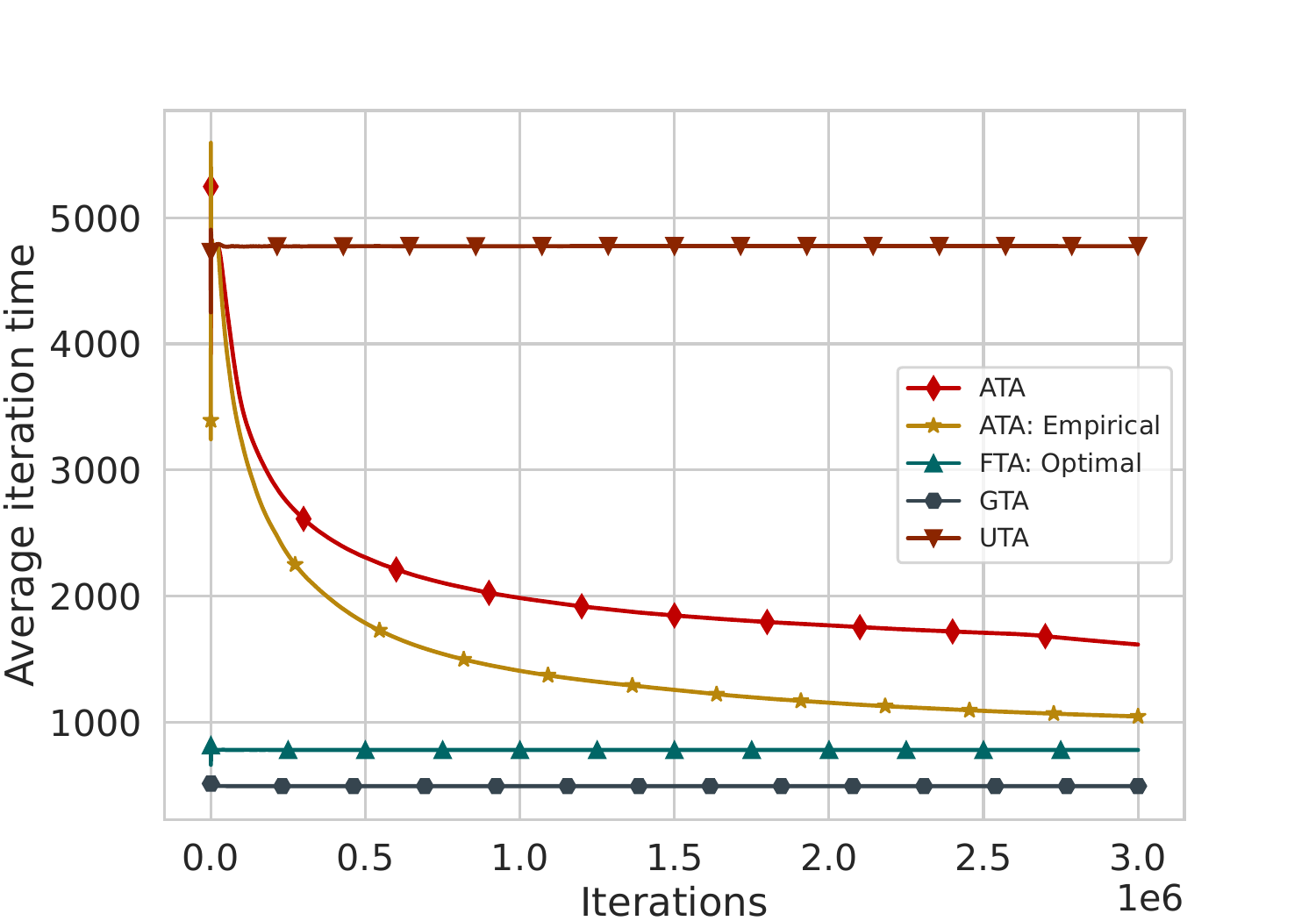} &
        \includegraphics[width=0.27\textwidth]{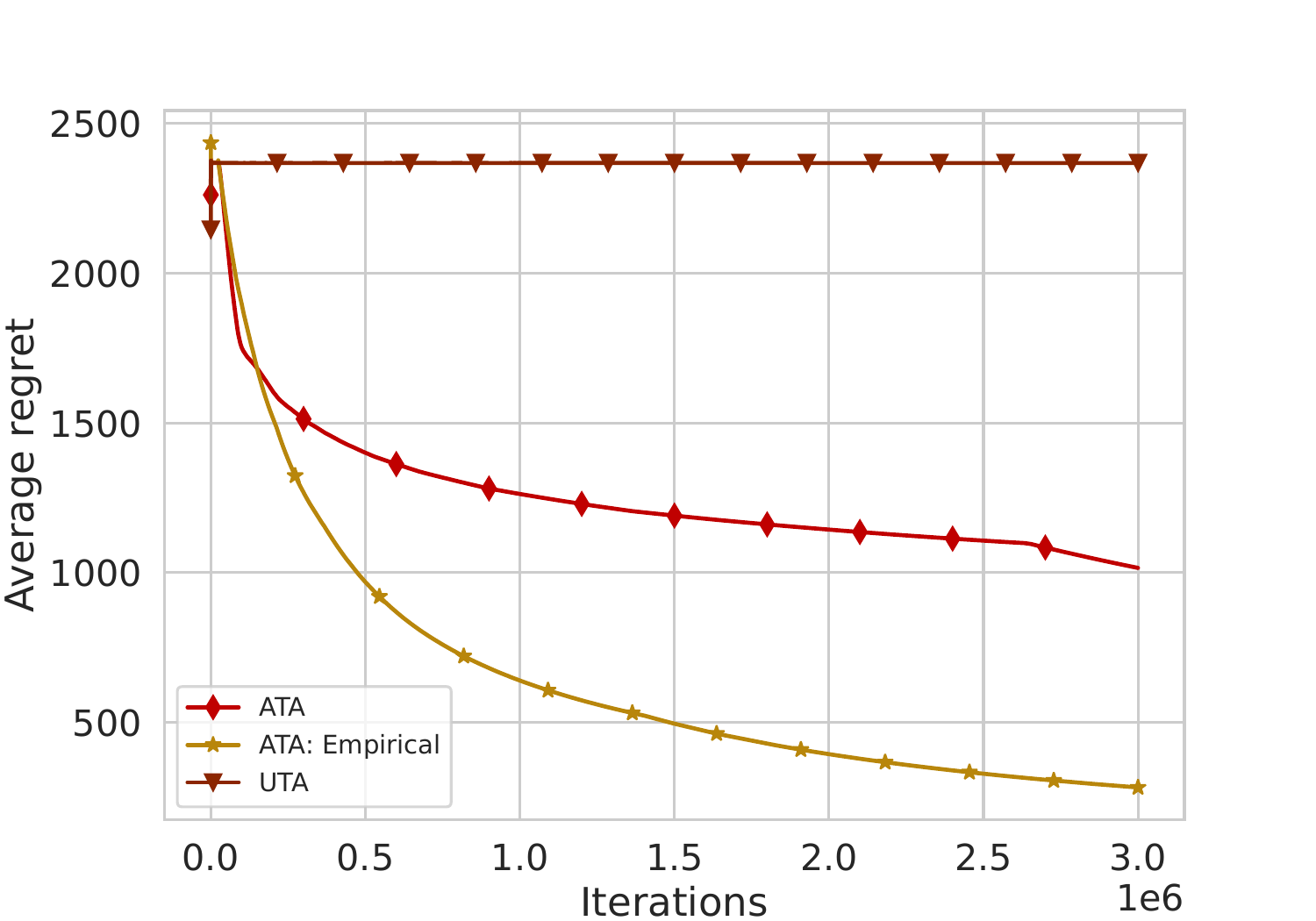} \\
        \includegraphics[width=0.27\textwidth]{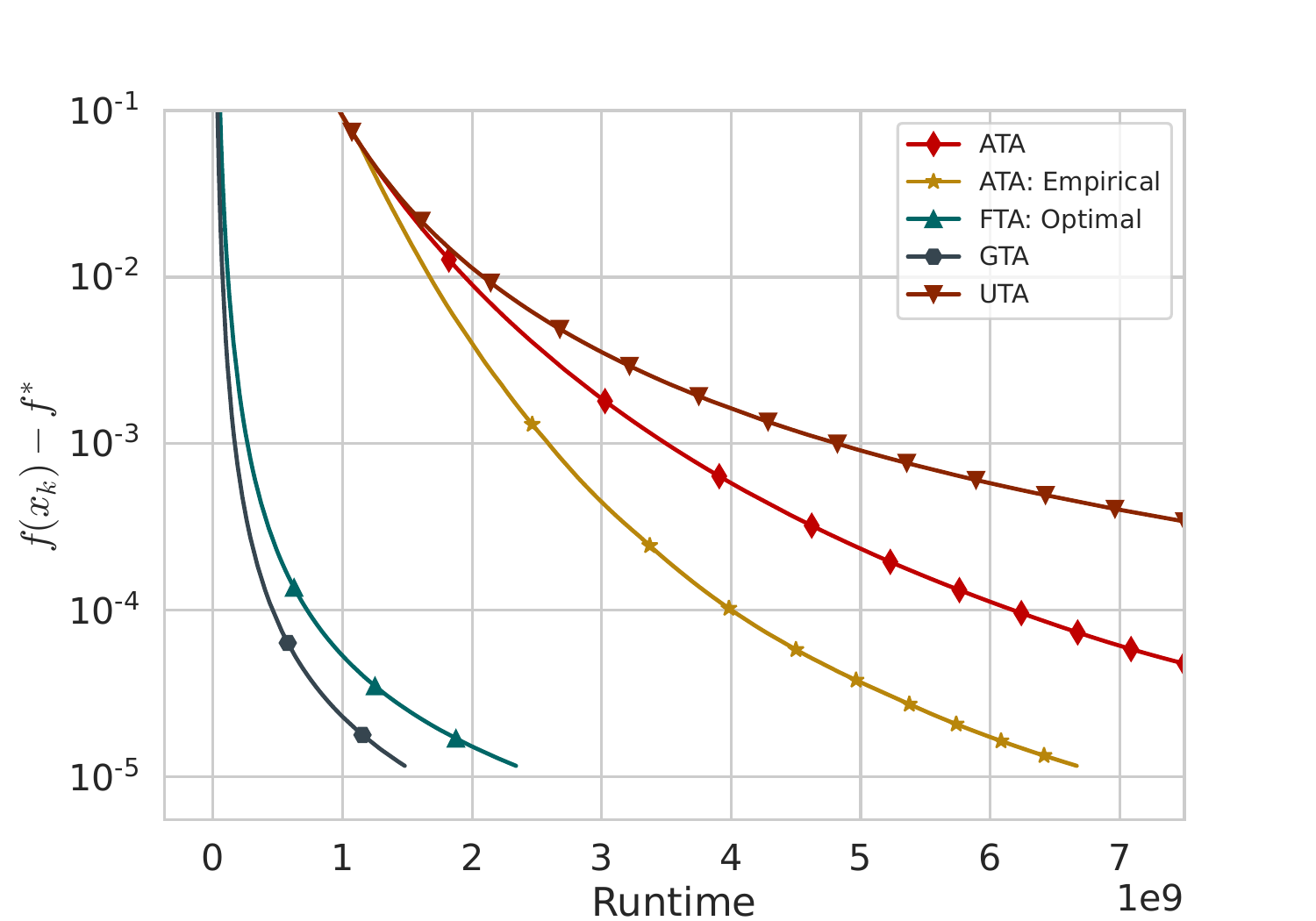} &
        \includegraphics[width=0.27\textwidth]{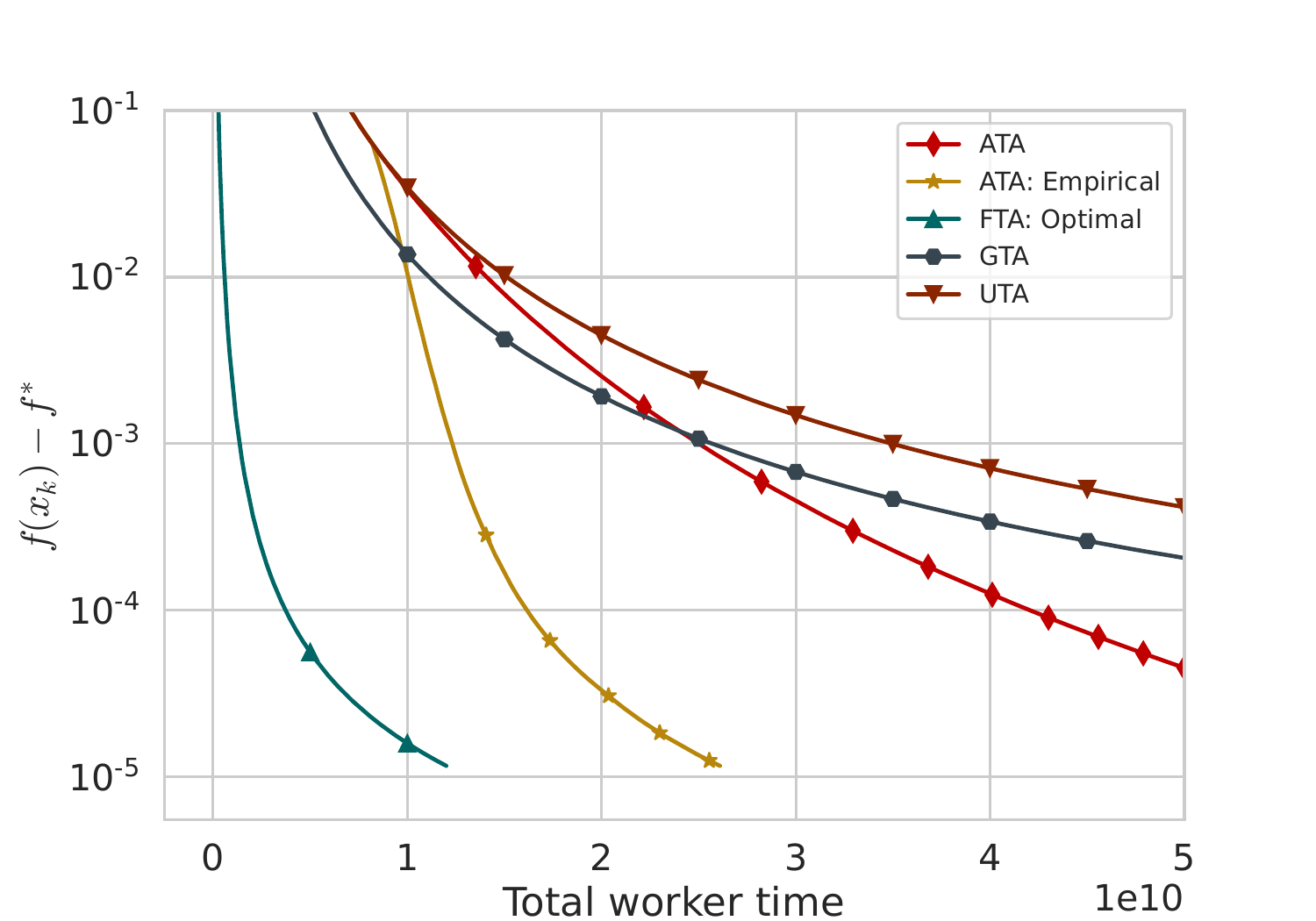} &
        \includegraphics[width=0.27\textwidth]{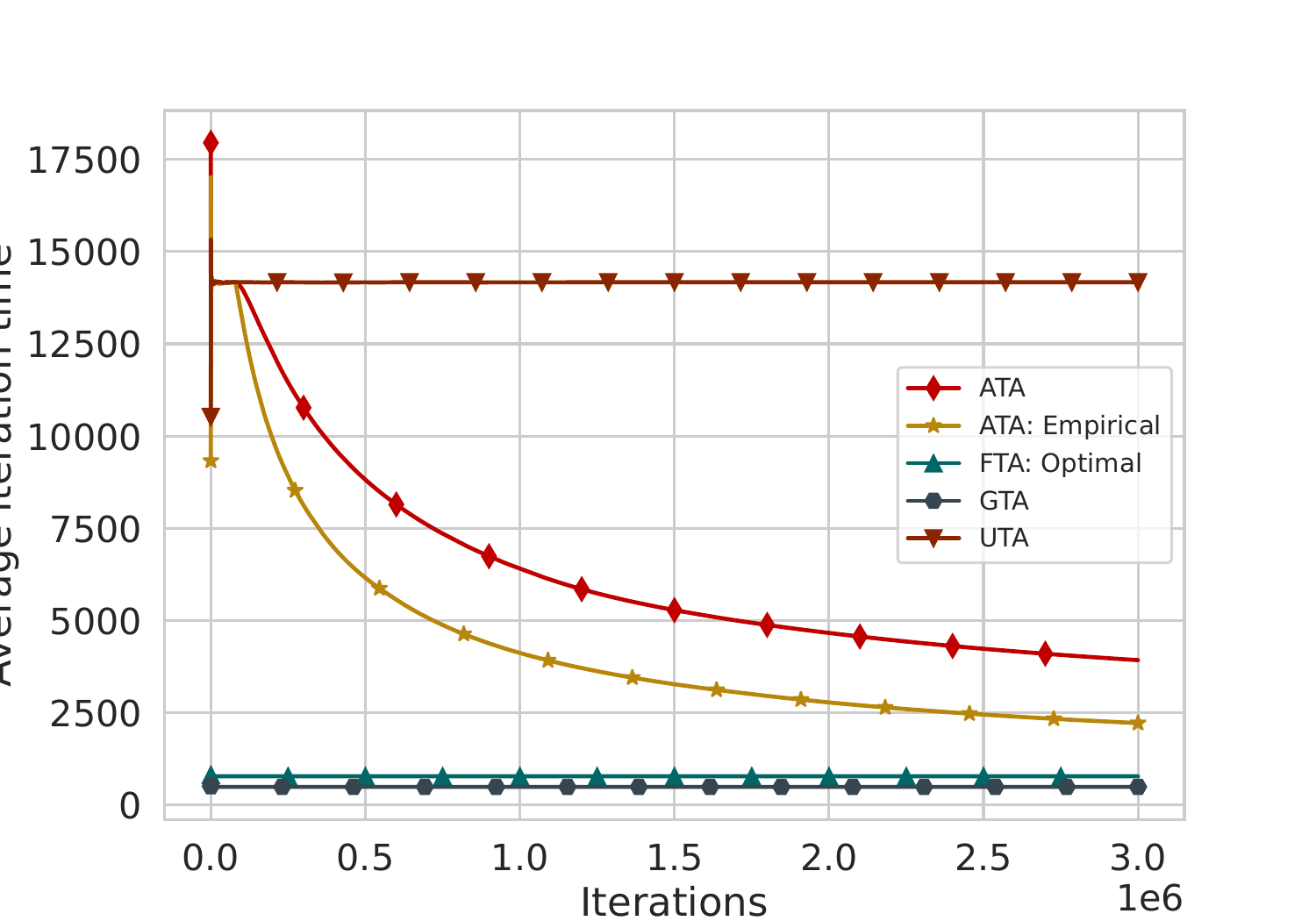} &
        \includegraphics[width=0.27\textwidth]{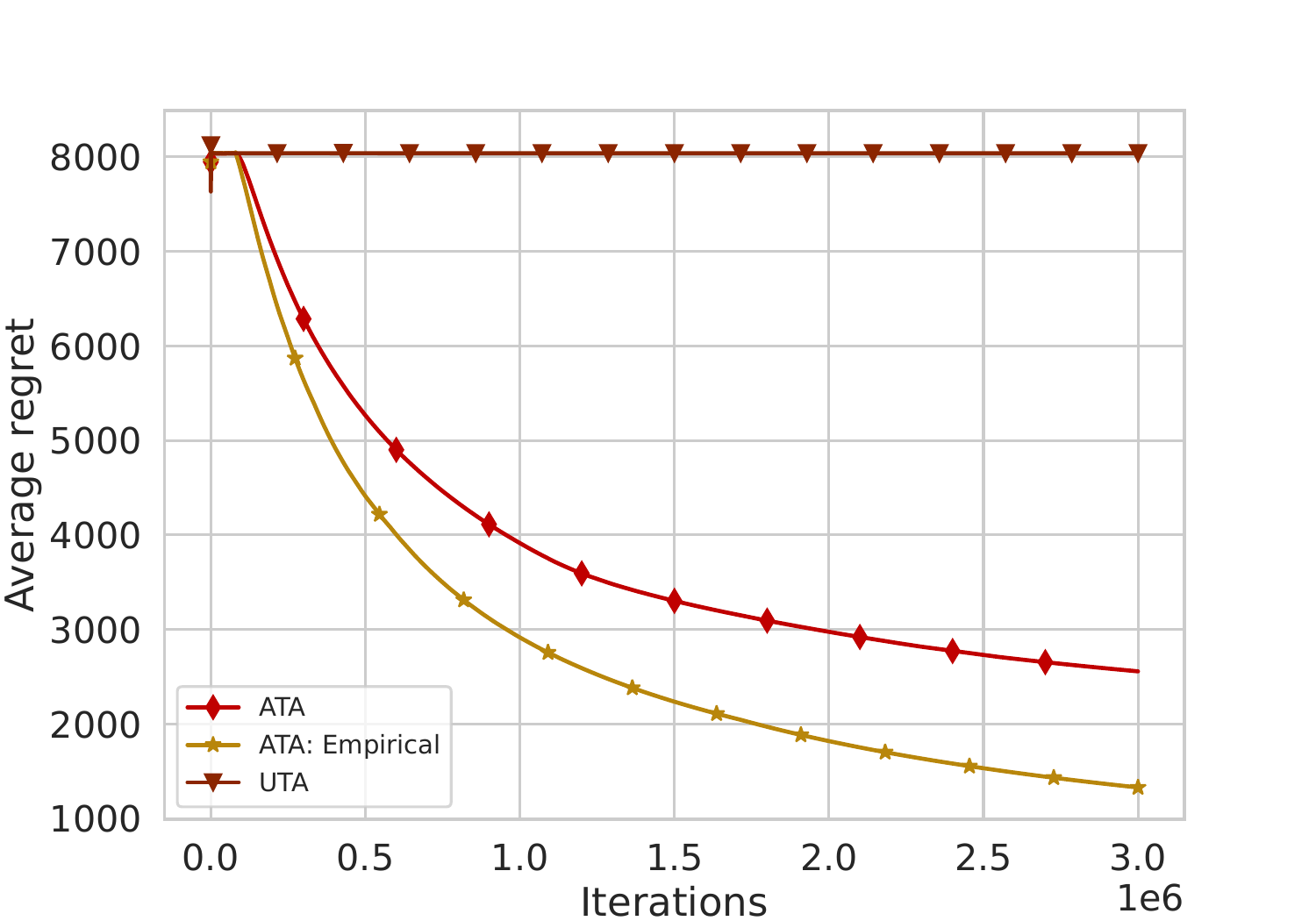}
    \end{tabular}
	}
    \caption{
		Each row increases the number of workers by a factor of 3, starting from $17$, that is, $n = 17, 51, 153, 459$ from top to bottom.
        The first column shows runtime vs. suboptimality.  
        The second column also plots suboptimality, but against total worker time, i.e., $\sum_{i=1}^n T_i^k$ in \Cref{ata:alg:ata}.
        The third column presents the average iteration time, given by $\nicefrac{C_k}{k}$ over all iterations $k$.  
        The last column displays the averaged cumulative regret, as defined in \eqref{ata:eq:proxy_loss}.
    }
    \label{ata:fig:linear}
\end{figure*}

The important difference to the previous \Cref{ata:fig:sqrt} is that here \algname{ATA-Empirical} performs better than \algname{ATA}.
This is because the Orlicz norm $\alpha = 4\cdot29n$ is much larger.

Similarly, we provide a numerical comparison in \Cref{ata:table:linear}.
\begin{table}[thb]
	\caption{
		Ratios of total worker times and runtimes required to achieve $f(x) - f^{*} < 10^{-5}$.
        For total worker time, we divide the total worker time of \algname{GTA} by the corresponding total worker times of the other algorithms listed.
        For runtime, we do the opposite, dividing the runtime of the other algorithms by the runtime of \algname{GTA}, since \algname{GTA} is the fastest.
        }
	\label{ata:table:linear}

	\begin{center}
	\begin{small}
	\begin{sc}
	\begin{tabular}{l|ccc|ccc}
	\toprule
	\multirow{2}{*}{$n$} & \multicolumn{3}{c|}{Total worker time ratio} & \multicolumn{3}{c}{Runtime ratio} \\
	\cmidrule(lr){2-4} \cmidrule(lr){5-7}
	& \algname{ATA} & \algname{ATA-Empirical} & \algname{OFTA} & \algname{ATA} & \algname{ATA-Empirical} & \algname{OFTA} \\
	\midrule
	$17$ & $2.32$ & $1.91$ & $2.1$ & $1.71$ & $1.71$ & $1.58$ \\
	$51$ & $6.71$ & $5.02$ & $6.29$ & $3.27$ & $2.12$ & $1.58$ \\
	$153$ & $3.41$ & $8.68$ & $18.87$ & $7.96$ & $4.5$ & $1.58$ \\
	\bottomrule
	\end{tabular}
	\end{sc}
	\end{small}
	\end{center}
\end{table}

\subsection{Heterogeneous Time Distributions}
\label{ata:sec:heterogeneous_distributions}

So far, we have only considered cases where clients follow the same distributions but with different means.
In this section, we extend our experiments to cases where the distributions themselves differ. 
We consider five distributions: Exponential, Uniform, Half-Normal, lognormal, and Gamma.
We group five workers with these five distributions so that each group has the same mean, then vary the mean across different groups.
More concretely, we use:
\begin{itemize}
    \item $\mathrm{Exp}(c(5g+1))$,
    \item $\mathrm{Uniform}\left(\frac{c(5g+1)}{2}, 3\frac{c(5g+1)}{2}\right)$,
    \item $\left|\cN\left(0, c(5g+1) \sqrt{\frac{\pi}{2}}\right)\right|$,
    \item $\mathrm{Lognormal}\left(\log(c(5g+1))/2, \sqrt{\log(c(5g+1))}\right)$,
    \item $\mathrm{Gamma}\left((c(5g+1))^2, \frac{1}{c(5g+1)}\right)$ with shape and scale parameters.
\end{itemize}
Next, we add a constant $c(5g+1)$ to all the distributions, where $c = 29$, and $g$ represents the group number, starting from 0. 
The clients are divided into $\nicefrac{n}{5}$ groups.

The results of the experiments are shown in \Cref{ata:fig:hetero}.
The plots demonstrate that the algorithms are robust across different distributions.

\begin{figure*}[thb]
    \centering
    { %
    \setlength{\tabcolsep}{-4pt} 
    \renewcommand{\arraystretch}{0} 
    \begin{tabular}{@{}cccc@{}}
        \includegraphics[width=0.27\textwidth]{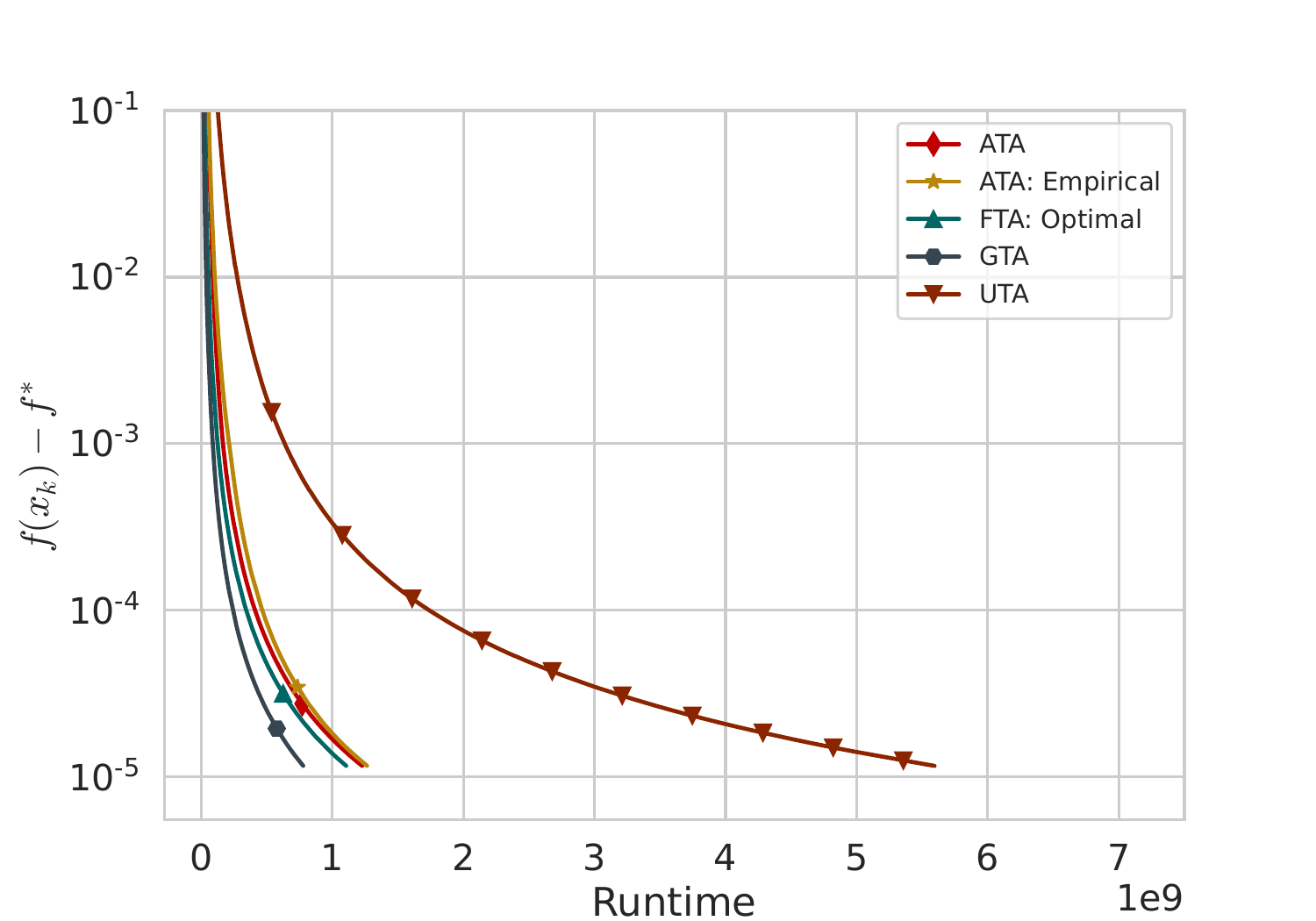} &
        \includegraphics[width=0.27\textwidth]{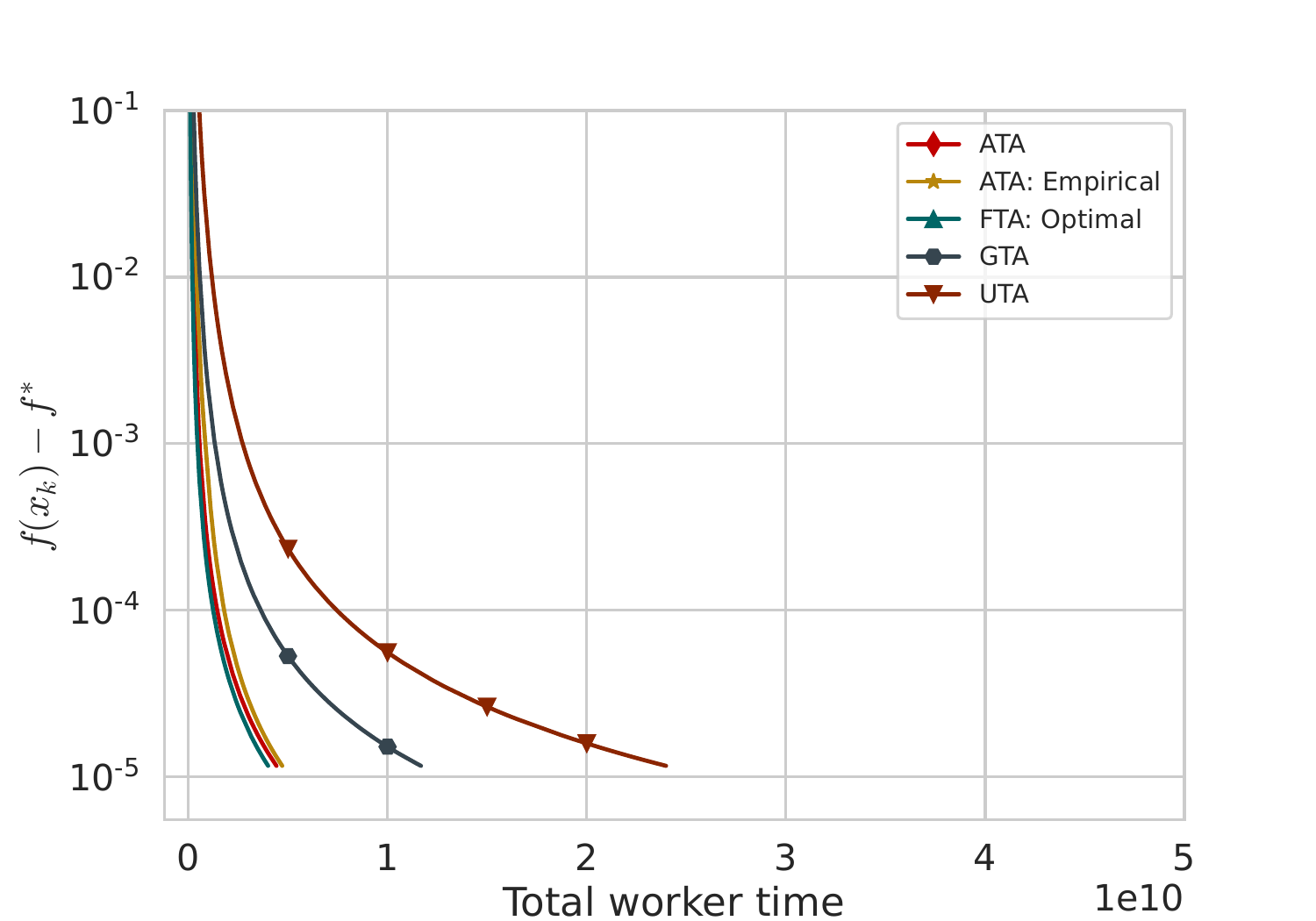} &
        \includegraphics[width=0.27\textwidth]{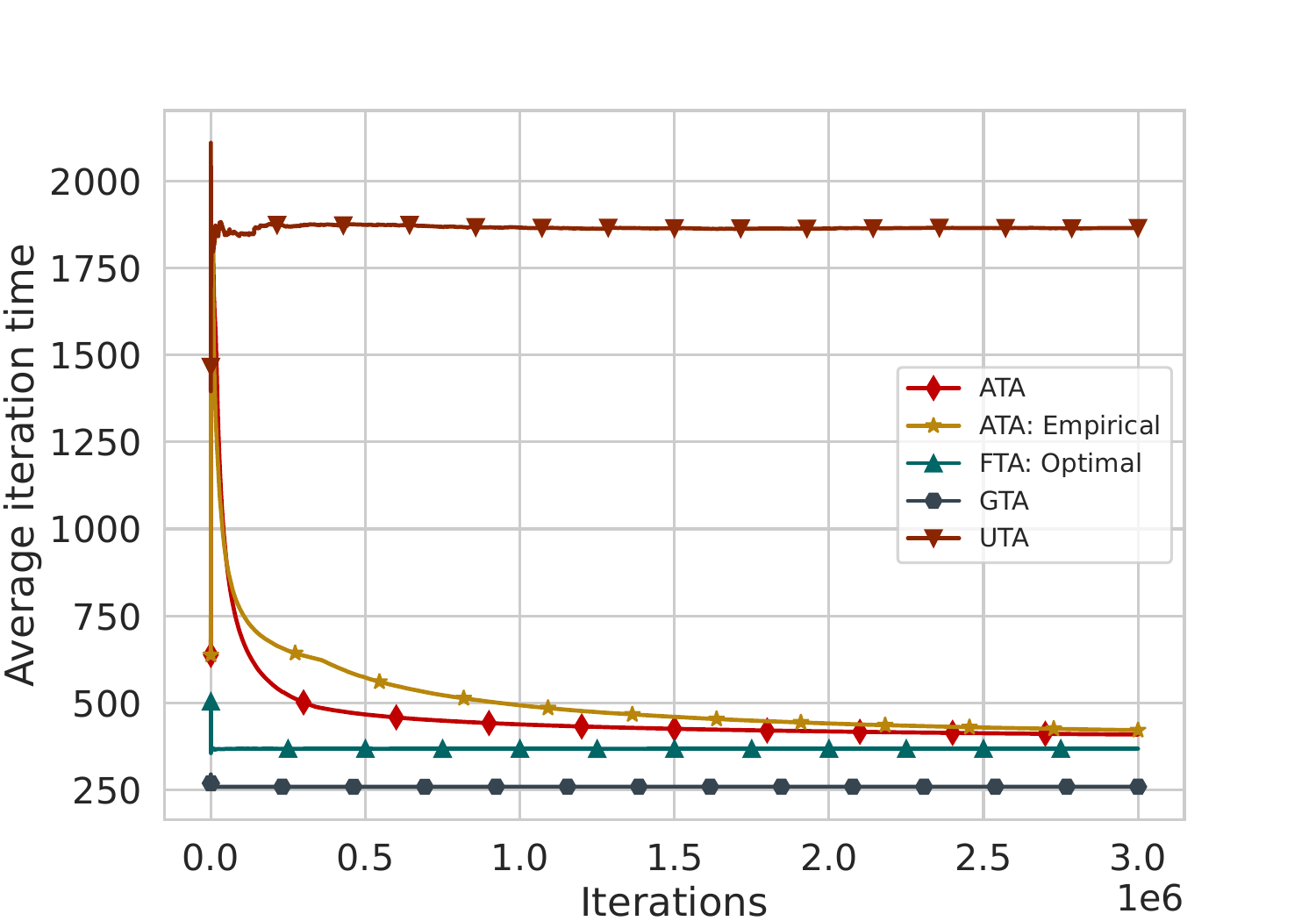} &
        \includegraphics[width=0.27\textwidth]{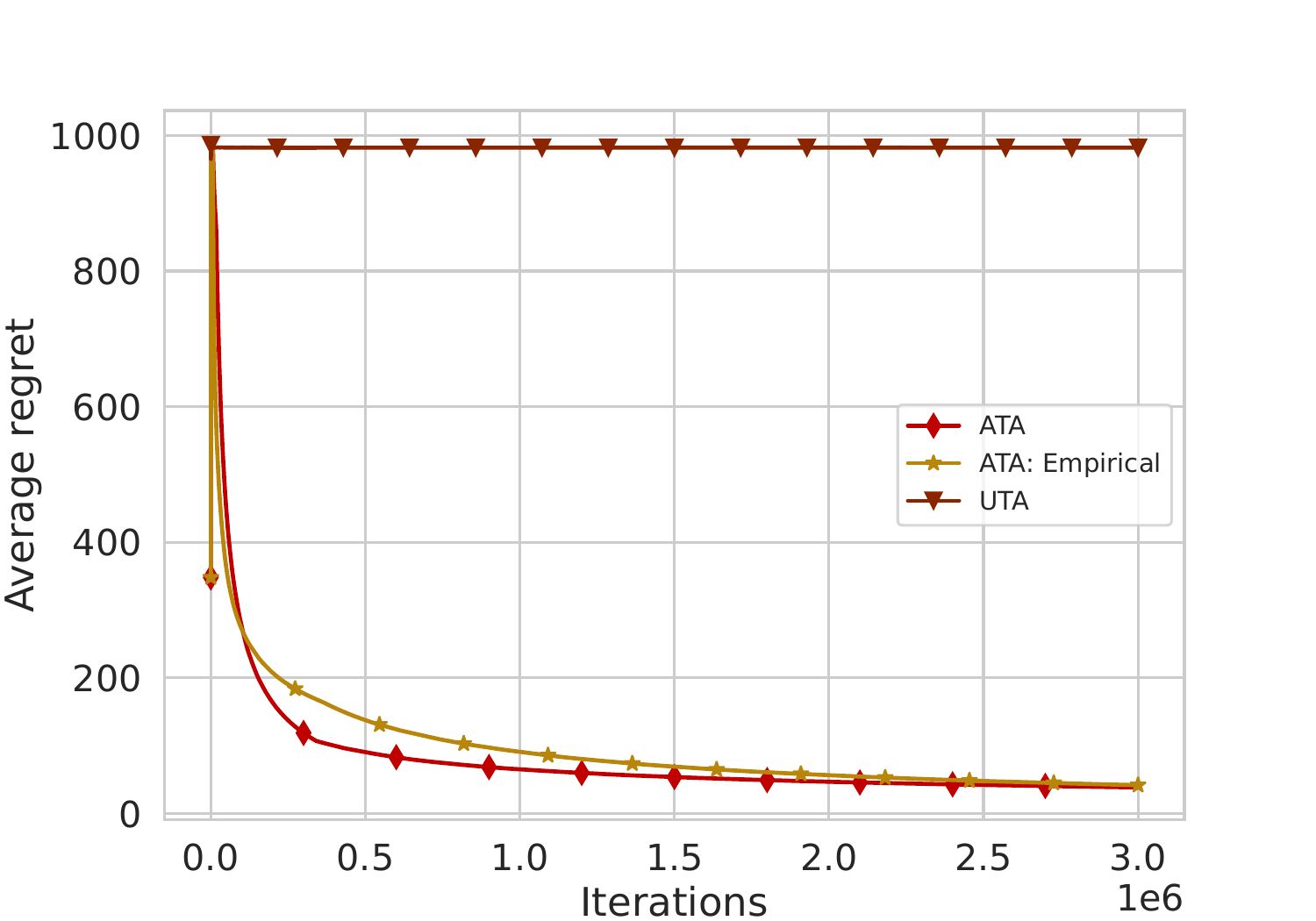} \\
        \includegraphics[width=0.27\textwidth]{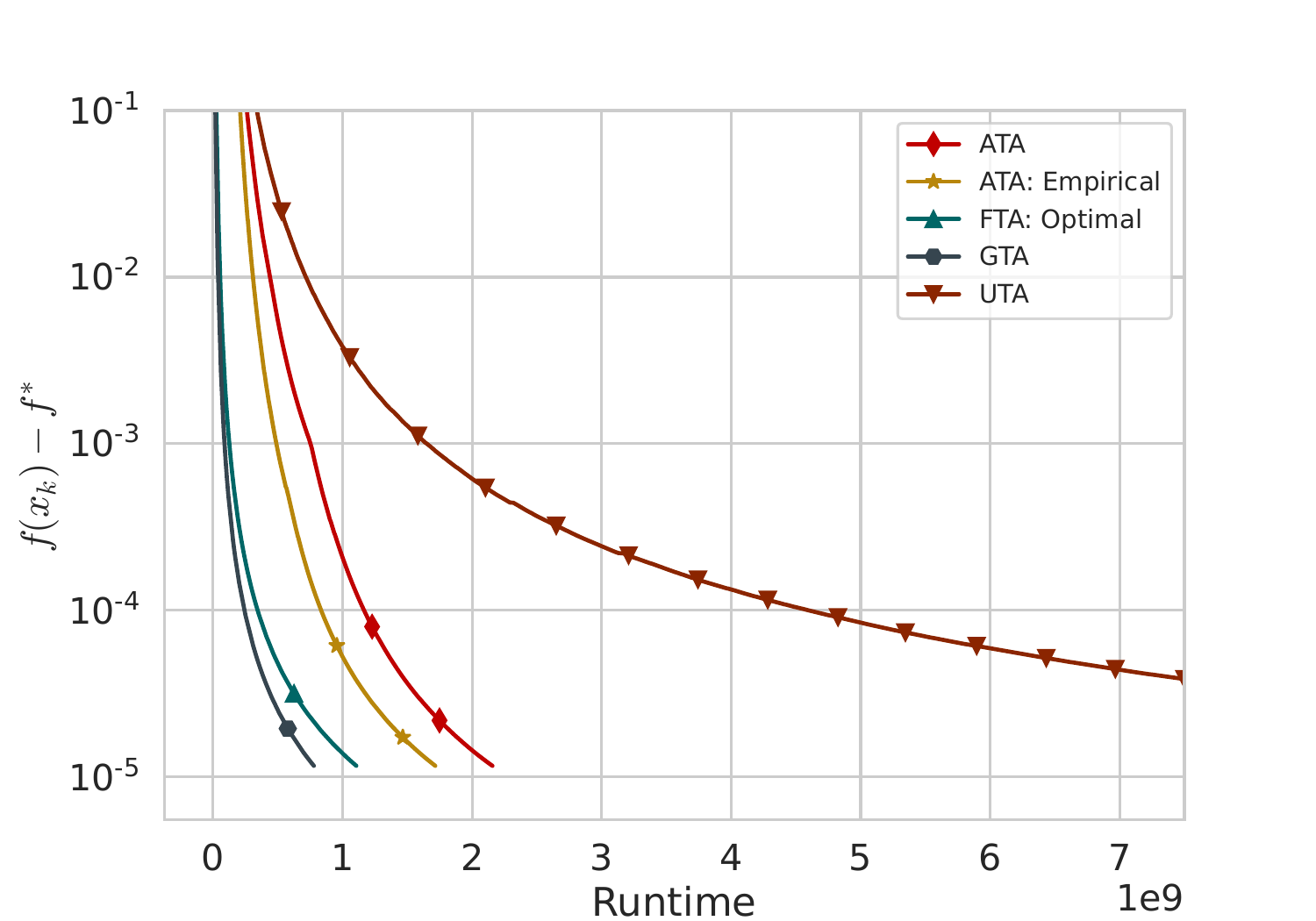} &
        \includegraphics[width=0.27\textwidth]{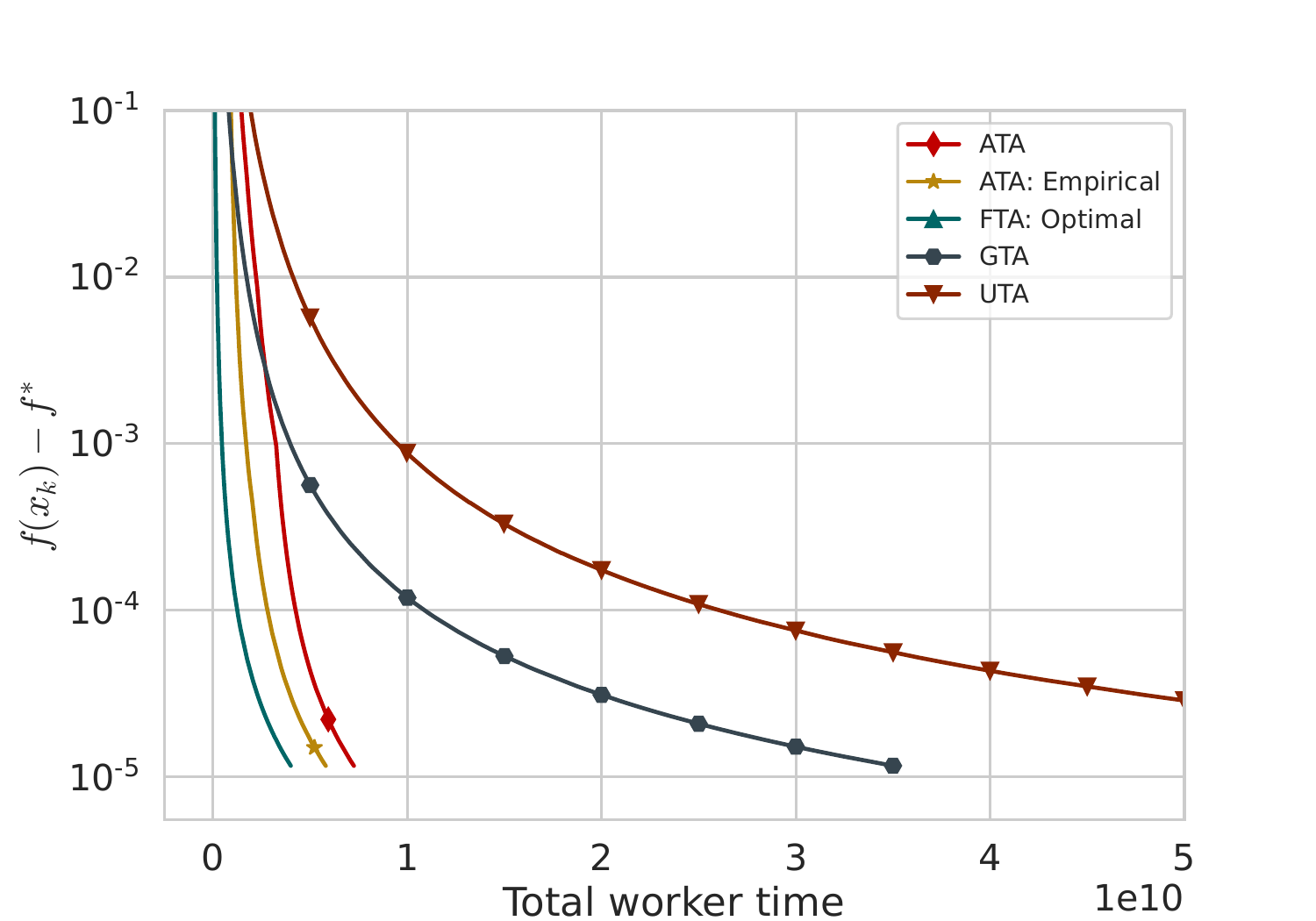} &
        \includegraphics[width=0.27\textwidth]{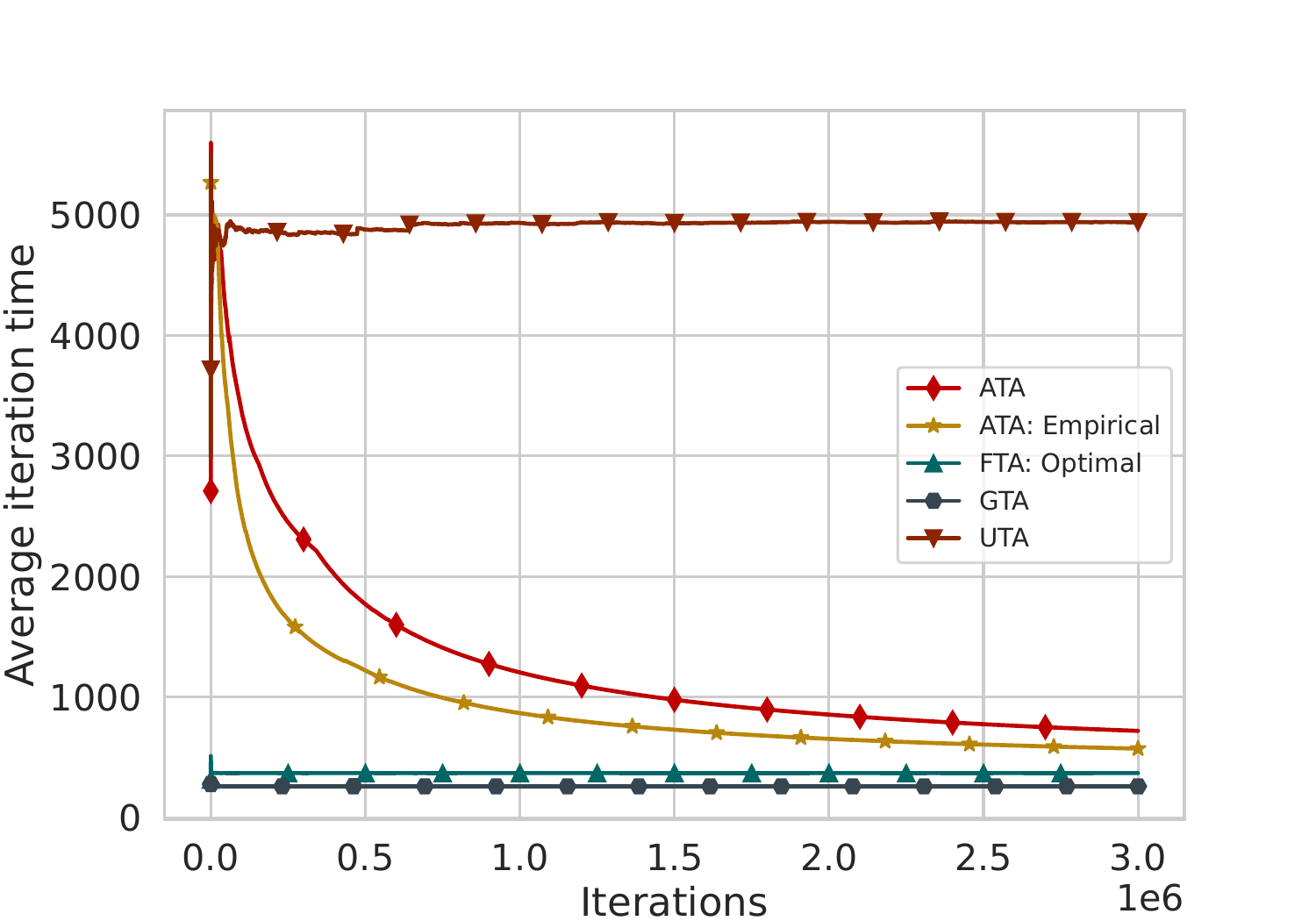} &
        \includegraphics[width=0.27\textwidth]{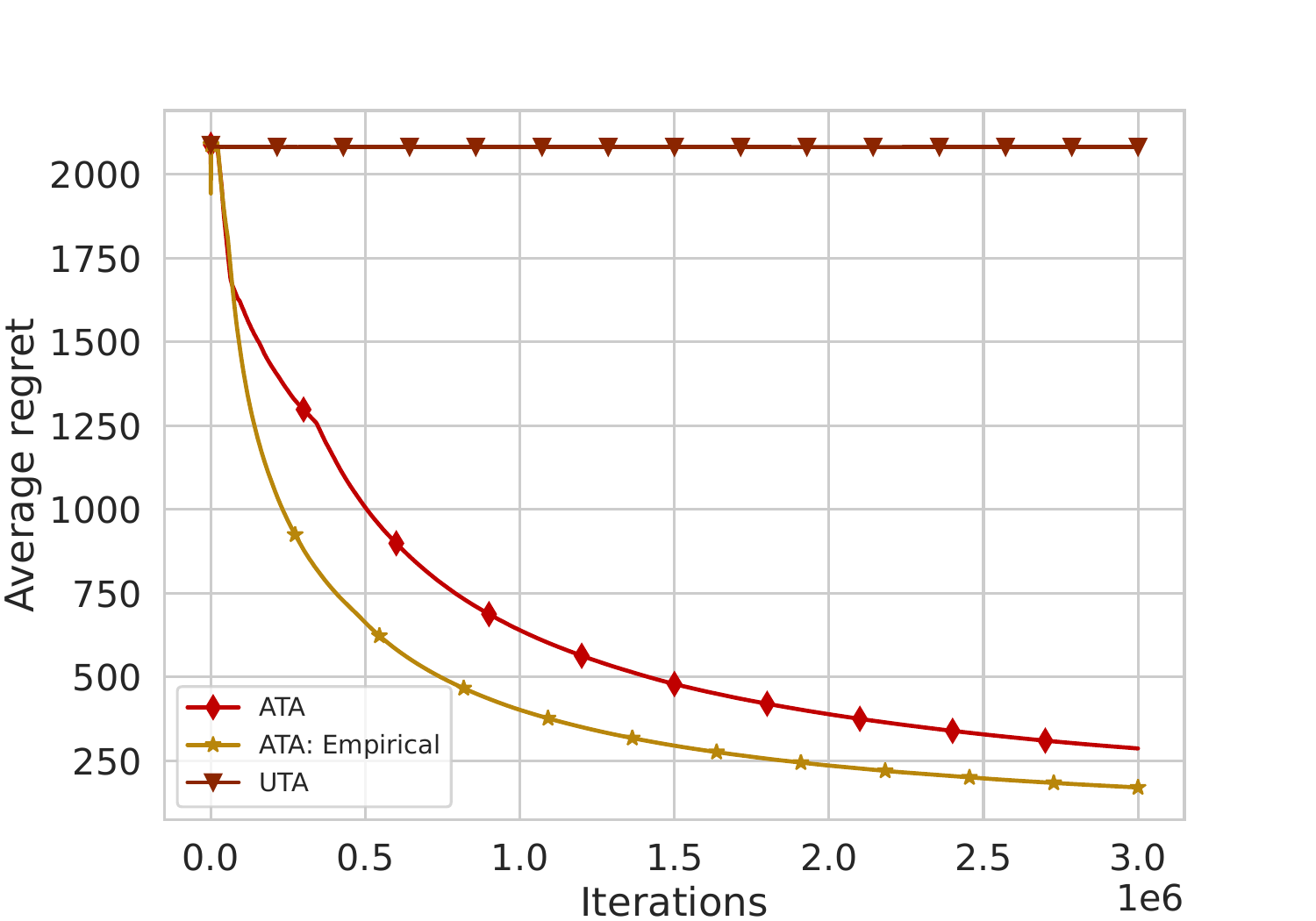} \\
        \includegraphics[width=0.27\textwidth]{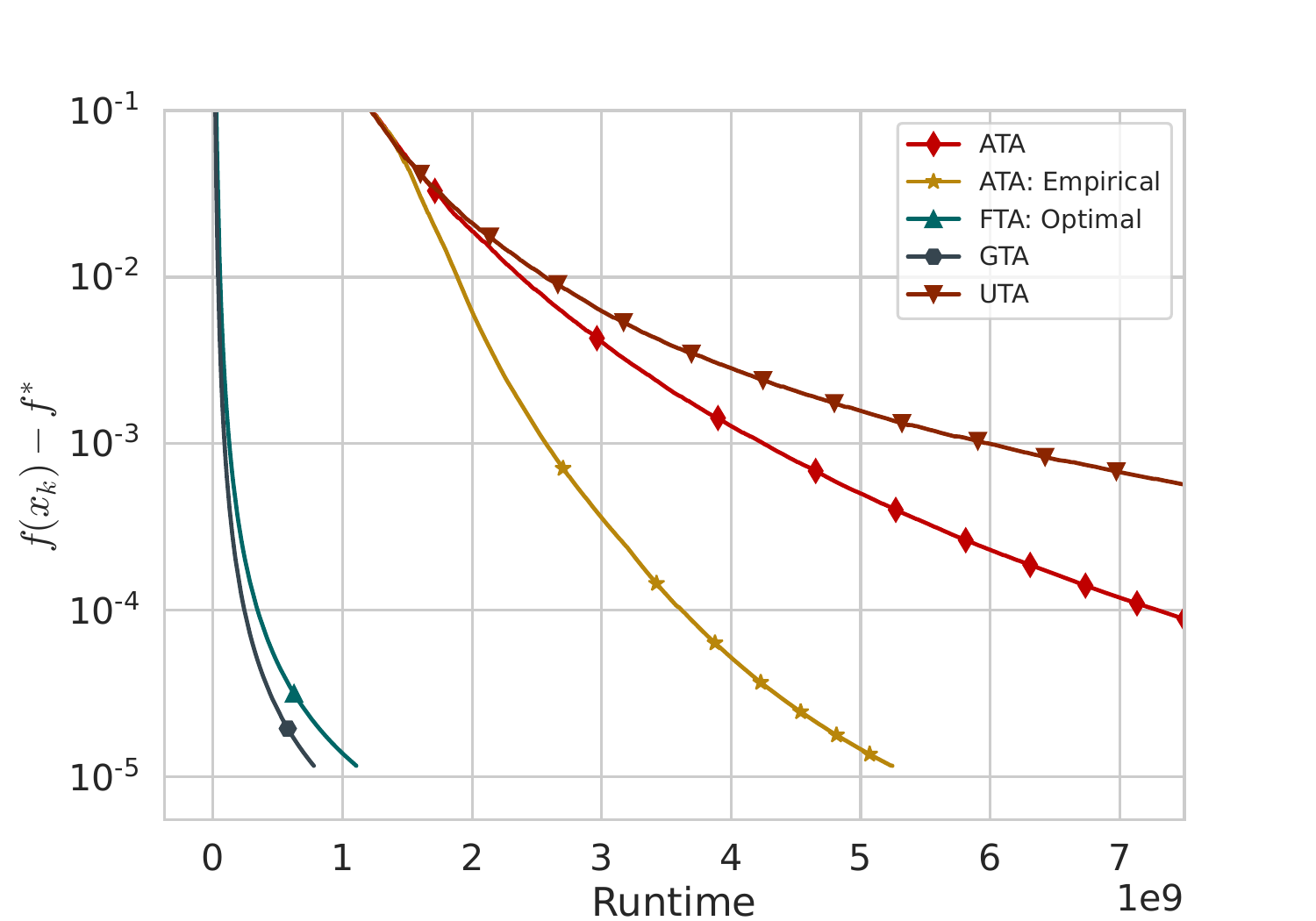} &
        \includegraphics[width=0.27\textwidth]{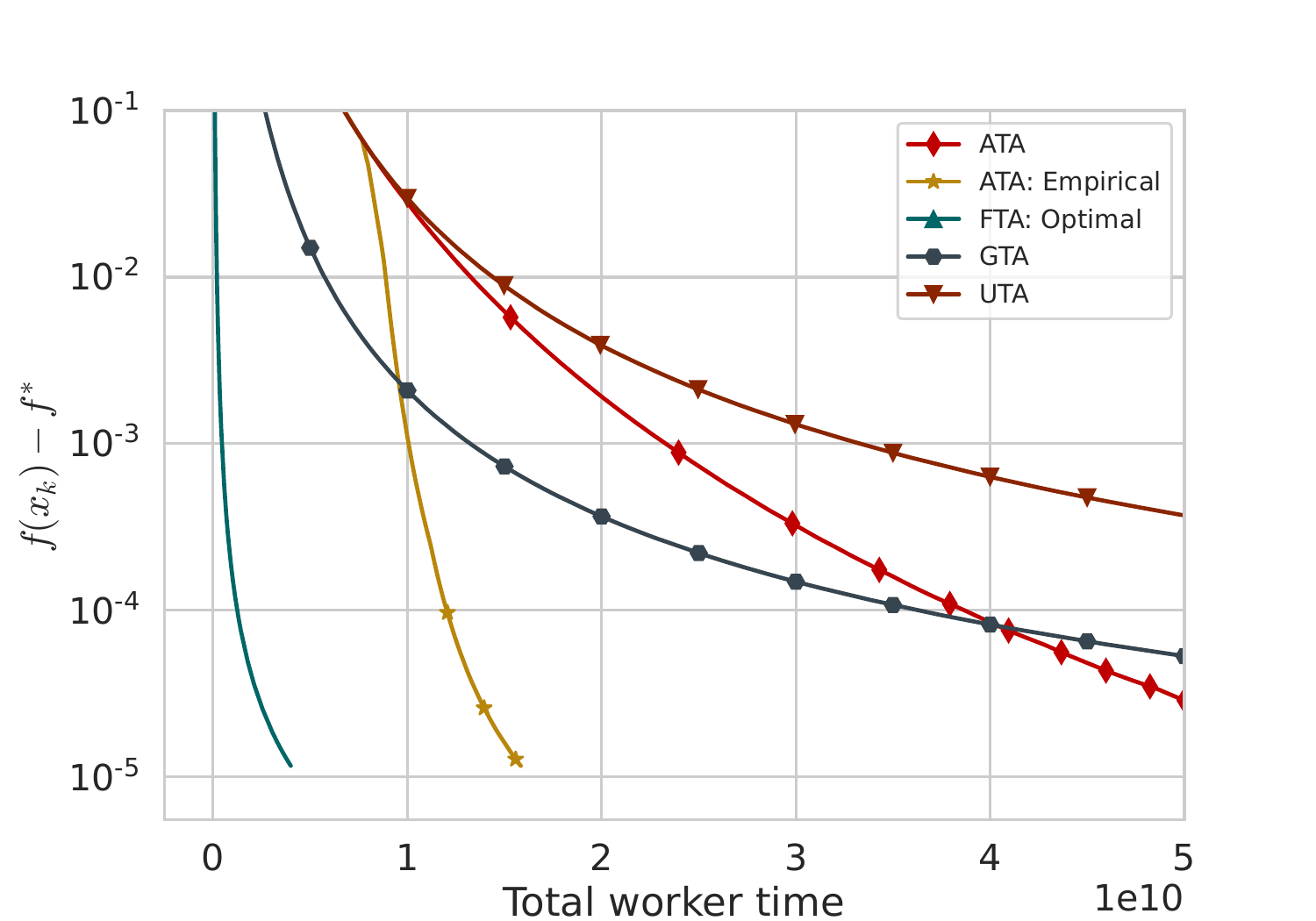} &
        \includegraphics[width=0.27\textwidth]{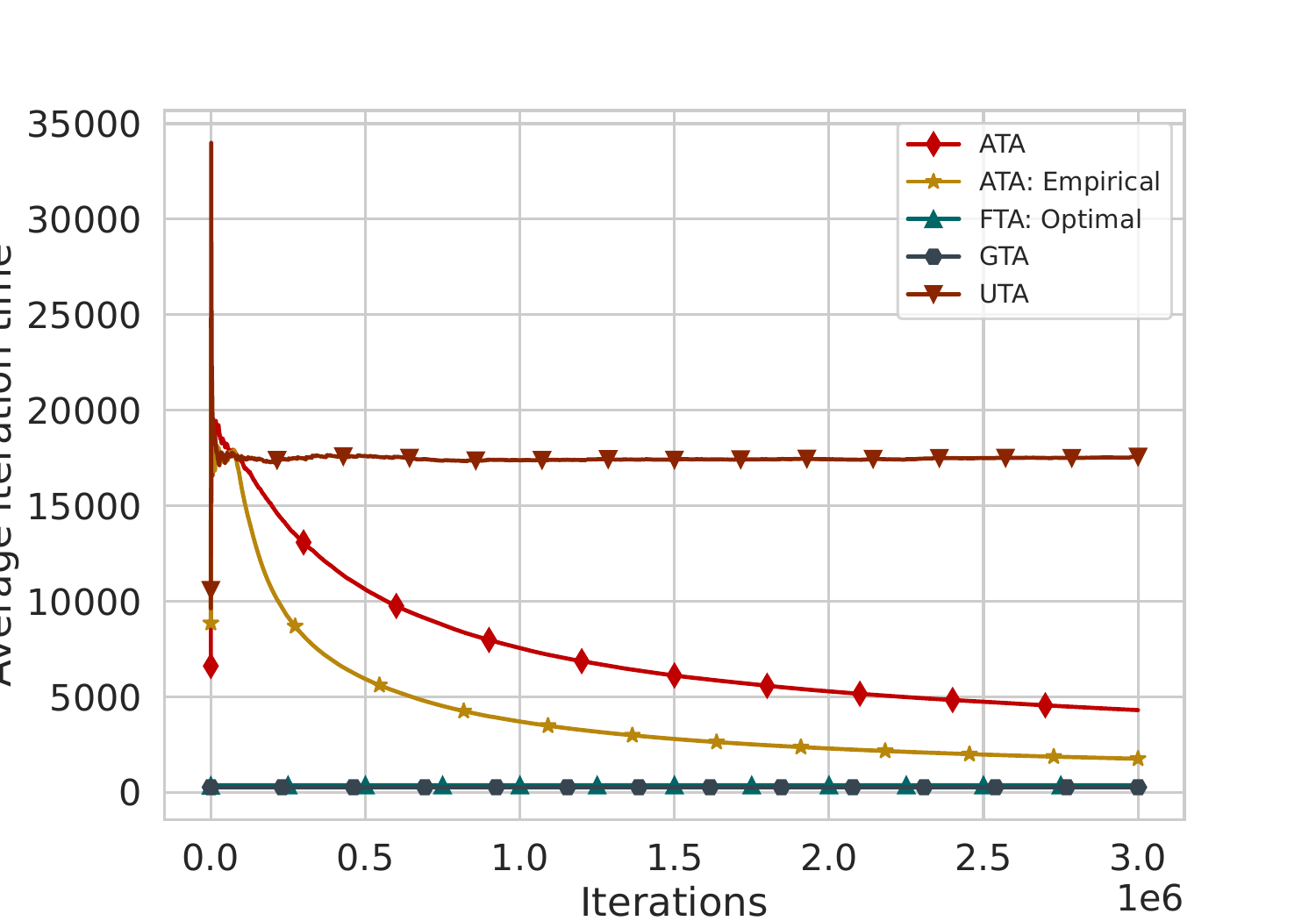} &
        \includegraphics[width=0.27\textwidth]{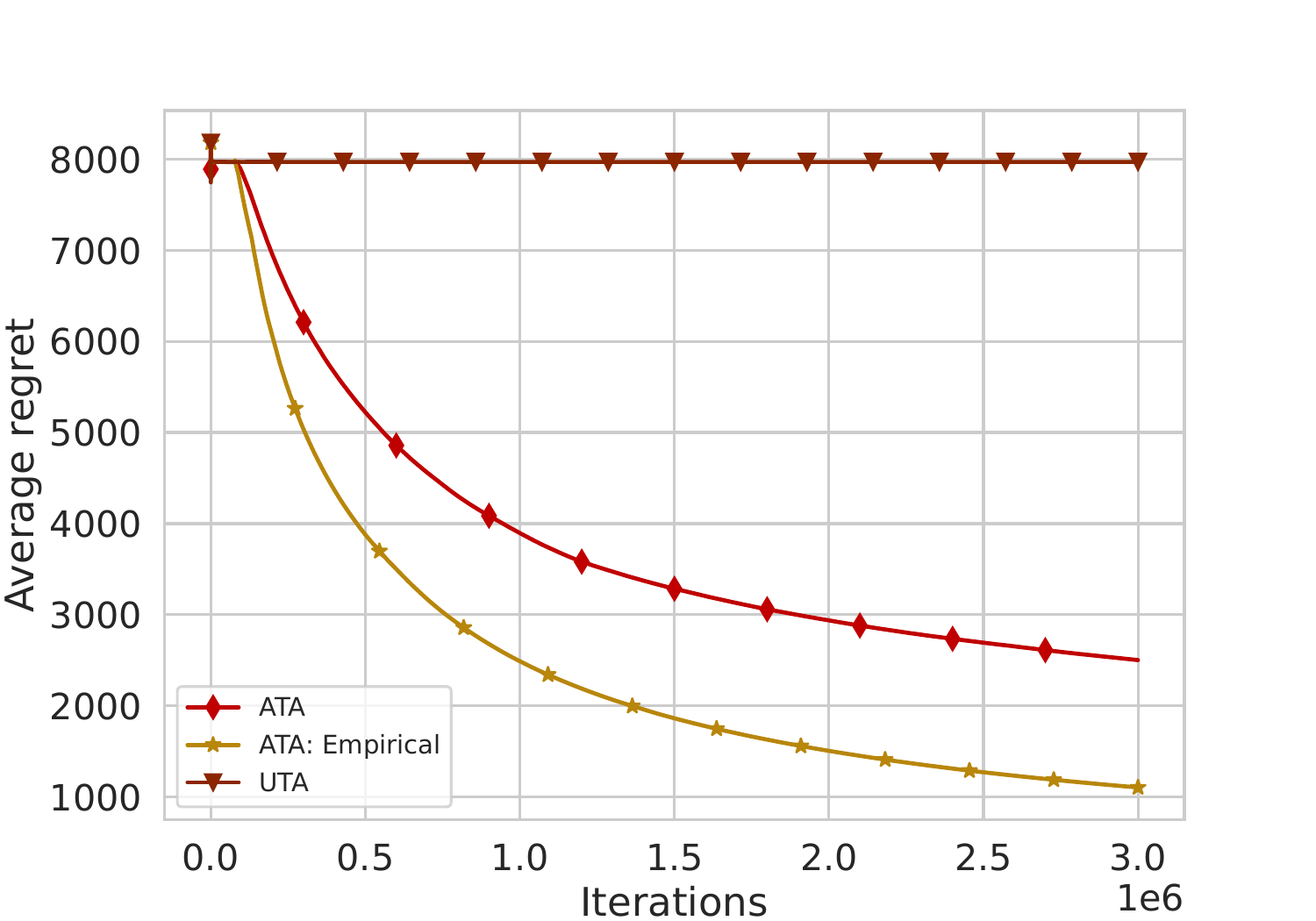}
    \end{tabular}
	}
    \caption{
        Each row corresponds to an increasing number of workers, with $n = 15, 45, 150$ from top to bottom.
        We consider five distributions---Exponential, Uniform, Half-Normal, Lognormal, and Gamma---grouping them to have the same mean and then varying the mean across different groups.
        The results demonstrate that the algorithms remain robust across different distributions.
        The columns represent the same as in \Cref{ata:fig:linear}.
    }
    \label{ata:fig:hetero}
\end{figure*}

\subsection{Regret}
\label{ata:sec:regret}

In this section, we verify \Cref{ata:thm:main} and \ref{ata:thm:main2} on regret through simulations.
\begin{figure*}[h]
    \centering
    \includegraphics[width=0.5\textwidth]{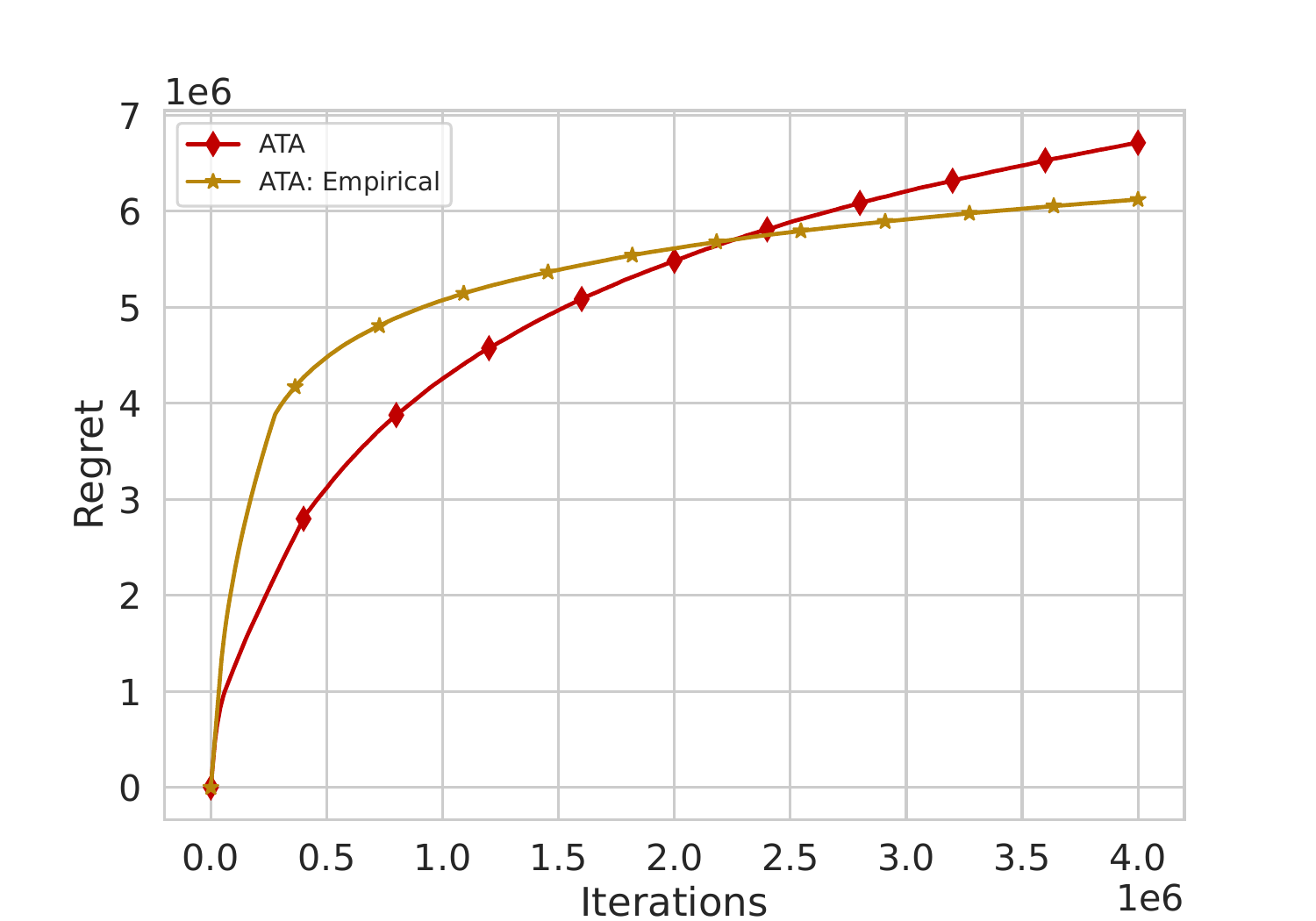}
    \caption{Regret growth over iterations.}
    \label{ata:fig:regret}
\end{figure*}
We set $n = 20$ and $B = 5$, with the computation time for worker $i$ following the distribution
$$
	\nu_i = \mathrm{Exp}(2i), \quad \text{for all} \quad i \in [n]~.
$$
We ran the simulation five times, and the plots include standard deviation bars, although they are not visible.
The results are presented in \Cref{ata:fig:regret}.

As expected, the regret grows logarithmically.

\subsection{Real Dataset}
\label{ata:sec:real_dataset}

In this section, we present an experiment where we train a convolutional neural network (CNN) on the CIFAR-100 dataset \citep{krizhevsky2009learning}.
The network consists of three convolutional layers and two fully connected layers, with a total of 160k parameters.

We use the Adam optimizer \citep{kingma2014adam} with a constant step size of $8 \cdot 10^{-5}$.
The computation time of the workers follows the same setup as in \Cref{ata:fig:linear}.
The results are shown in \Cref{ata:fig:real}.

\begin{figure*}[thb]
    \centering
	{ %
    \renewcommand{\arraystretch}{0} 
    \begin{tabular}{@{}cc@{}}
        \includegraphics[width=0.40\textwidth]{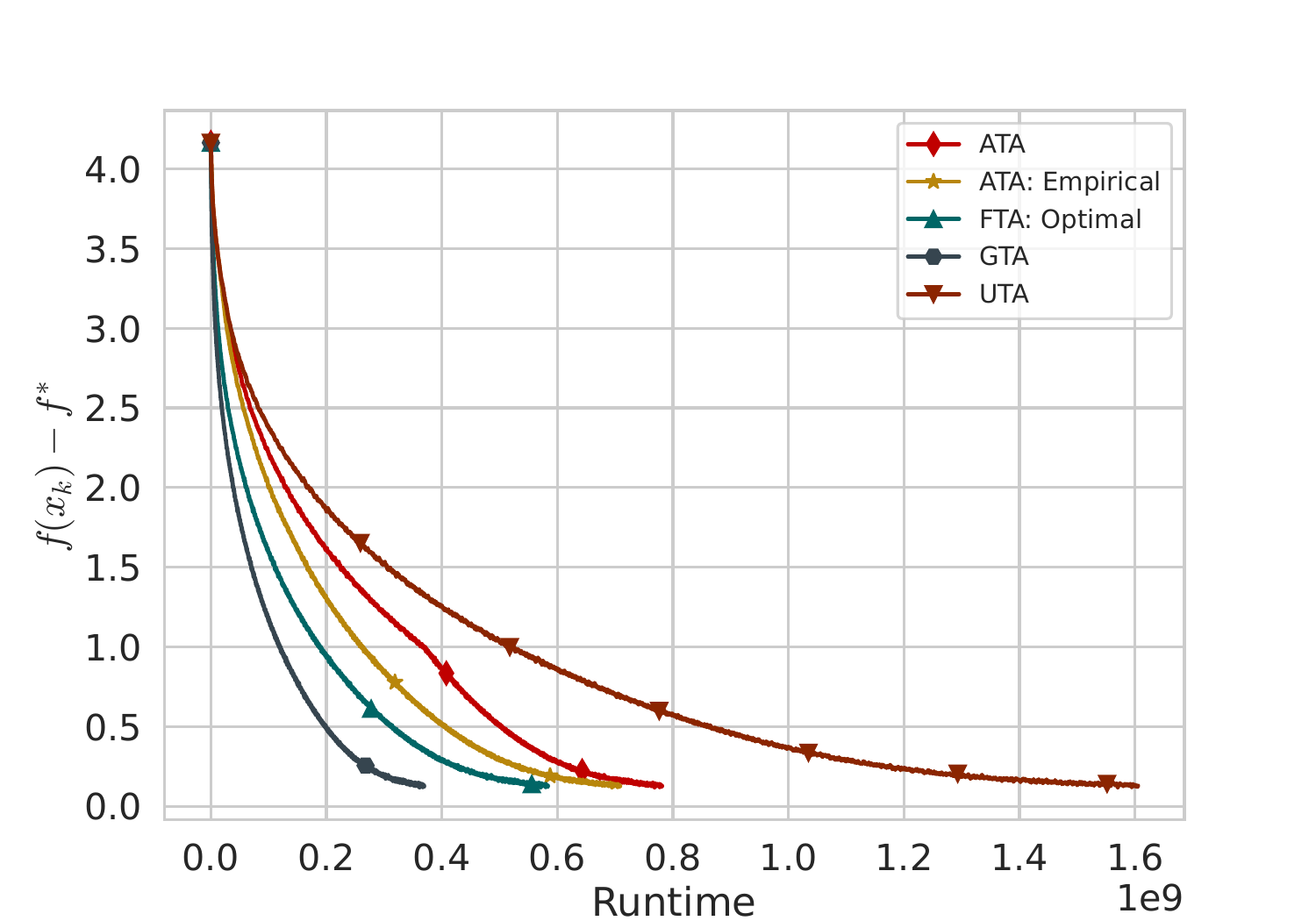} &
        \includegraphics[width=0.40\textwidth]{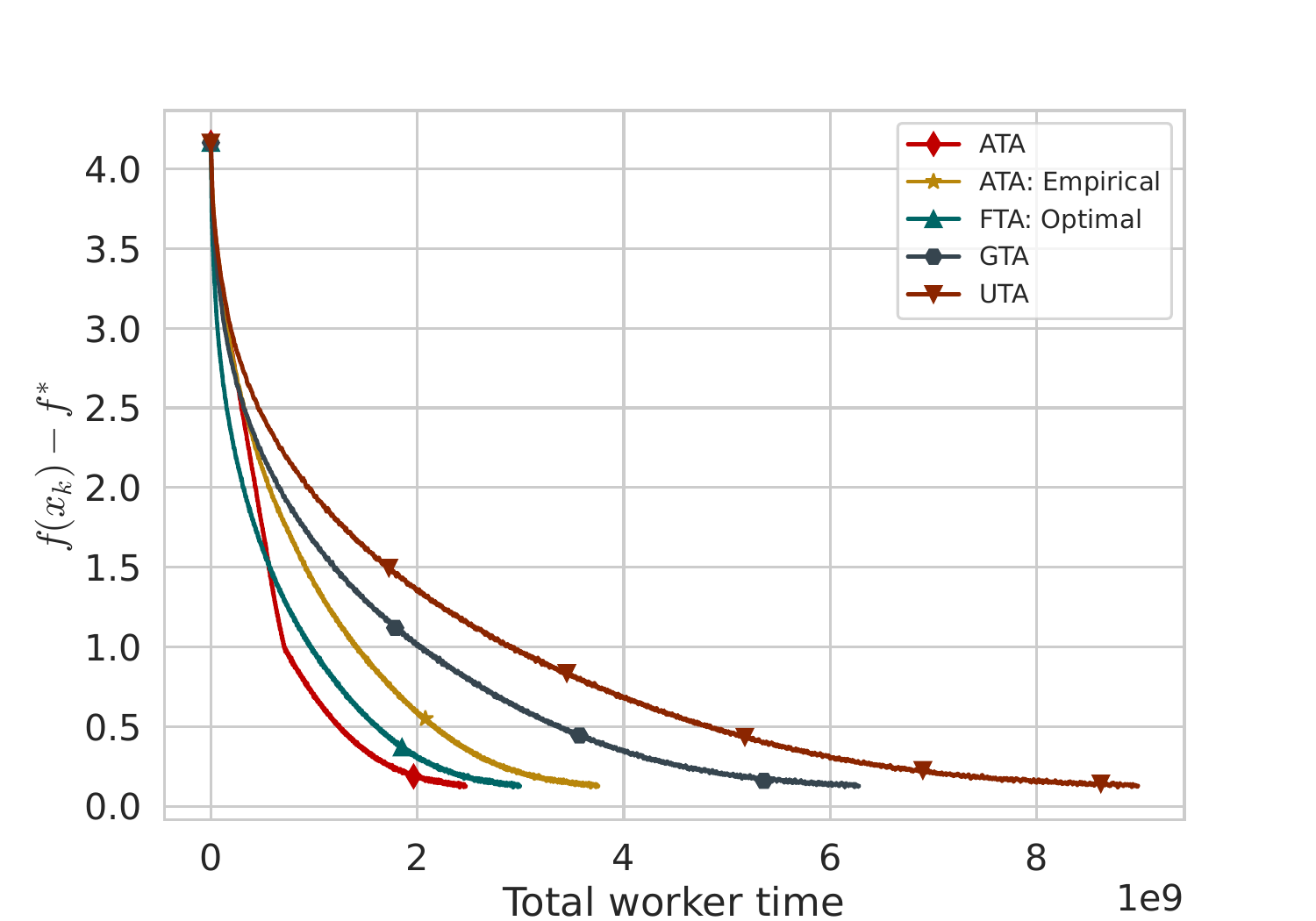} \\
        \includegraphics[width=0.40\textwidth]{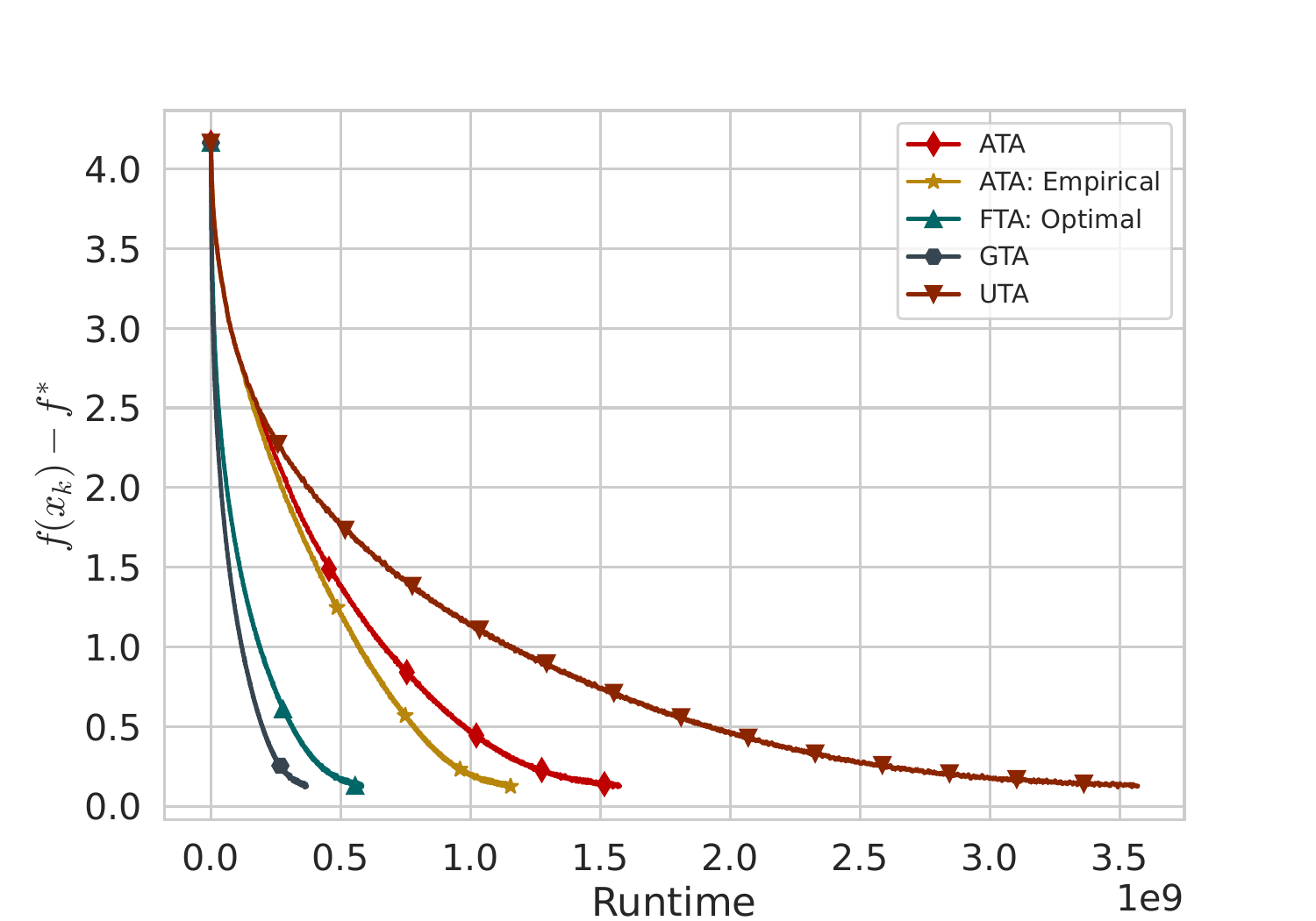} &
        \includegraphics[width=0.40\textwidth]{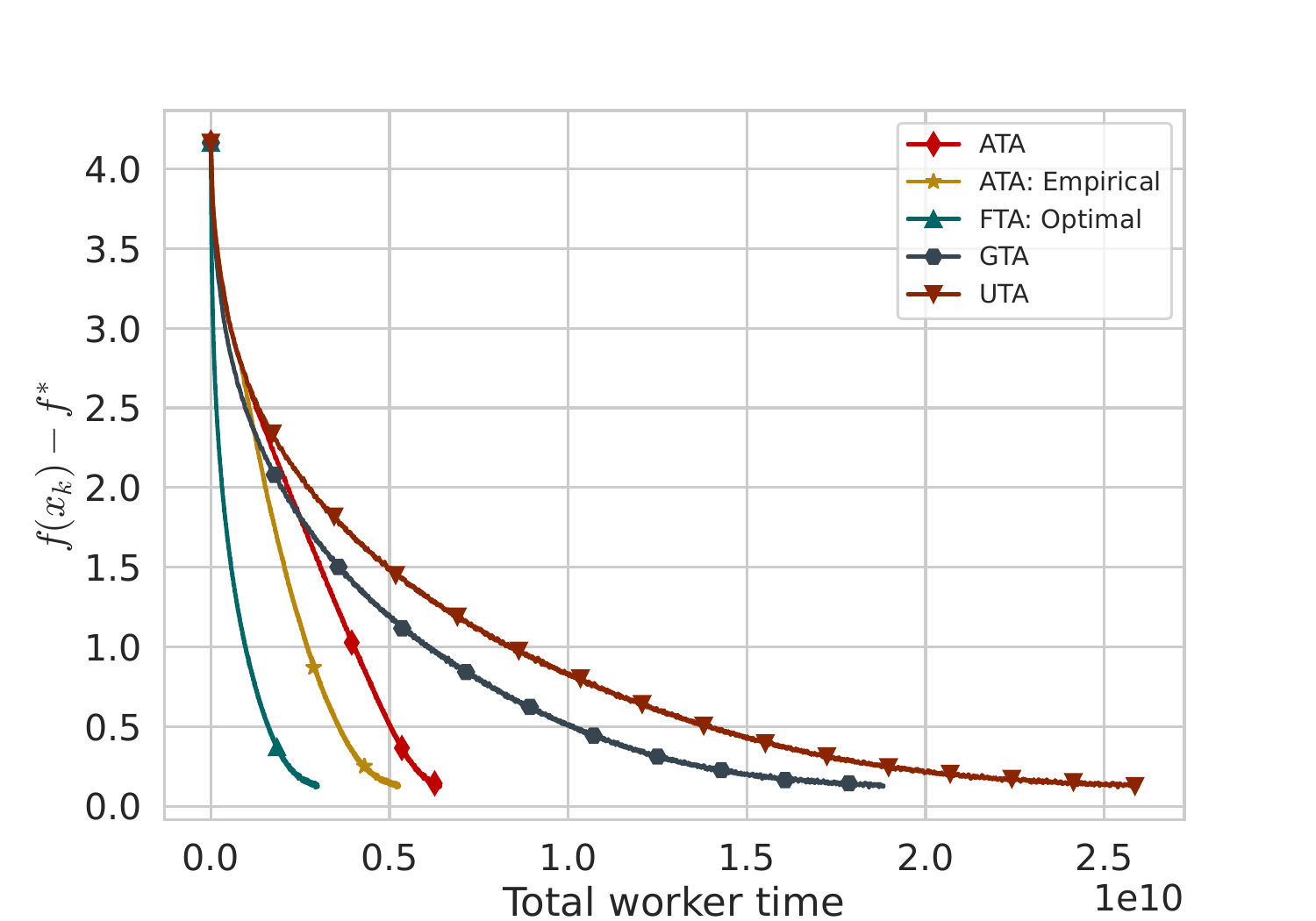} \\
        \includegraphics[width=0.40\textwidth]{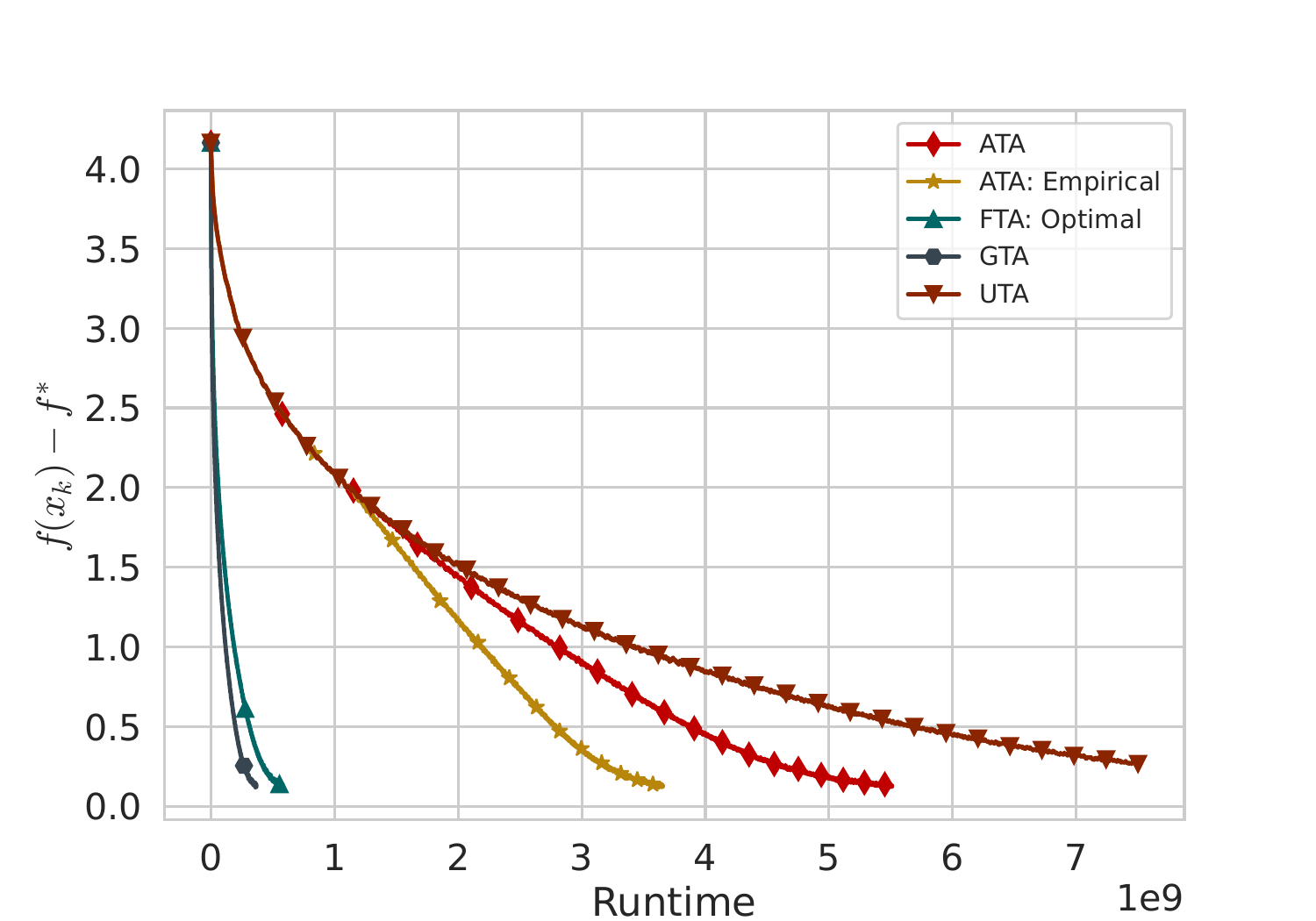} &
        \includegraphics[width=0.40\textwidth]{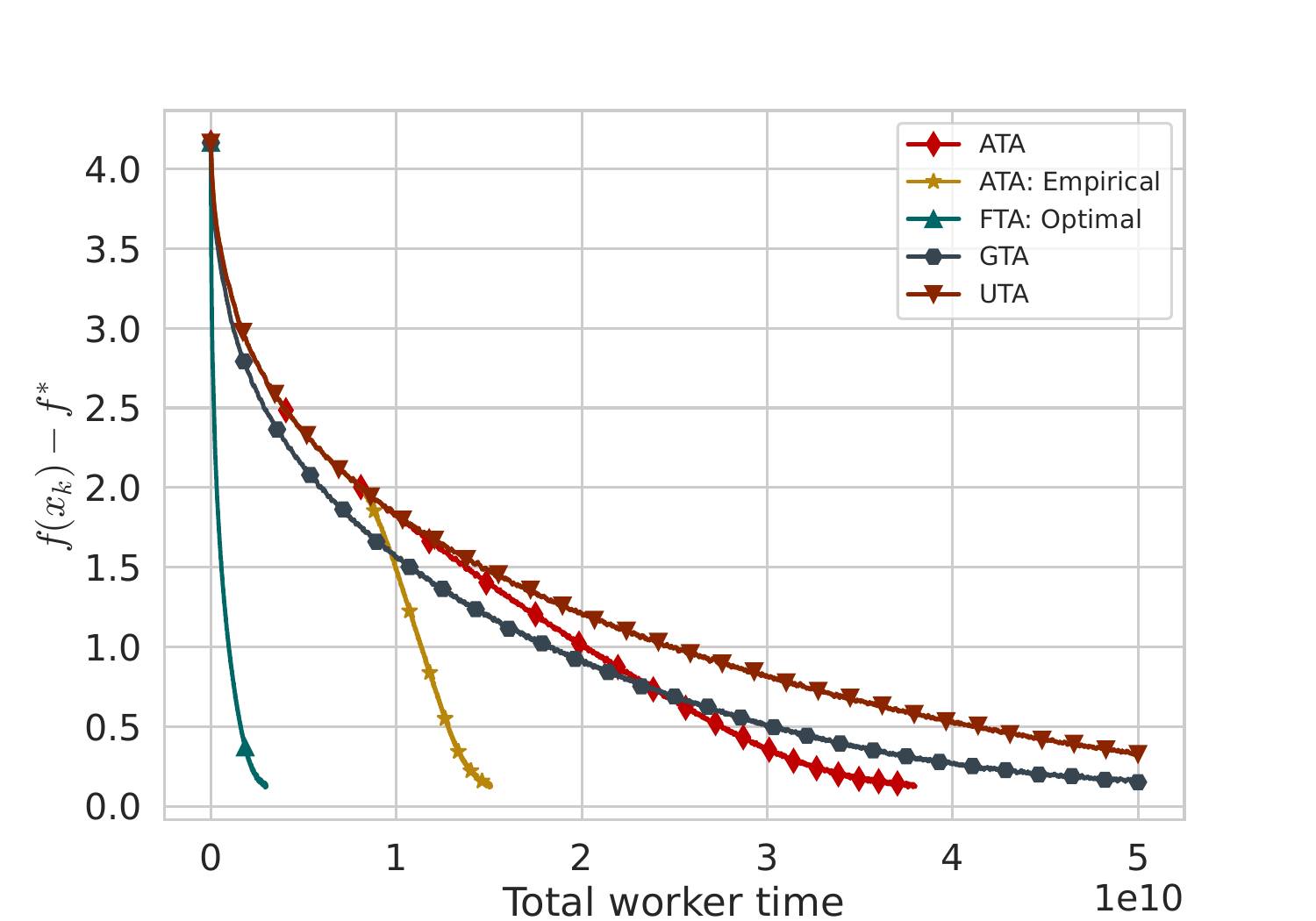}
    \end{tabular}
	}
    \caption{
        We use the CIFAR-100 dataset \citep{krizhevsky2009learning}.
        The model is a CNN with three convolutional layers and two fully connected layers, totaling 160k parameters.
        The Adam optimizer \citep{kingma2014adam} is used with a constant step size of $8 \cdot 10^{-5}$.
        The computation time of the workers follows the same setup as in \Cref{ata:fig:linear}, where the mean time increases linearly.
        The batch size remains the same at $B=23$.
        Each row corresponds to a different number of workers, with $n = 17, 51, 153$ from top to bottom.
    }
    \label{ata:fig:real}
\end{figure*}

\subsection{Impact of Prior Knowledge on Time Distributions}
\label{ata:sec:prior_knowledge}

In real-world systems where multiple machine learning models are trained, estimates of computation times from previous runs may be available. 
With this prior knowledge, \algname{ATA} and \algname{ATA-Empirical} can be much faster, as they spend less time on exploration and quickly approach the performance of \algname{OFTA}.

To illustrate this, we vary the number of prior runs, $P$.
Since our algorithms operate independently of the underlying optimization process, we first focus solely on the bandit component, updating the confidence scores of machines over several iterations.
We then apply the loss curves to different segments of the bandit phase and compare the results as $P$ increases.
A larger $P$ yields more accurate estimates.

The optimization setup remains the same as in \Cref{ata:fig:real}, with $B=23$ and $n=51$.
The results are presented in \Cref{ata:fig:prior}.

\begin{figure*}[thb]
    \centering
    { %
    \renewcommand{\arraystretch}{0} 
    \begin{tabular}{@{}cc@{}}
        \includegraphics[width=0.40\textwidth]{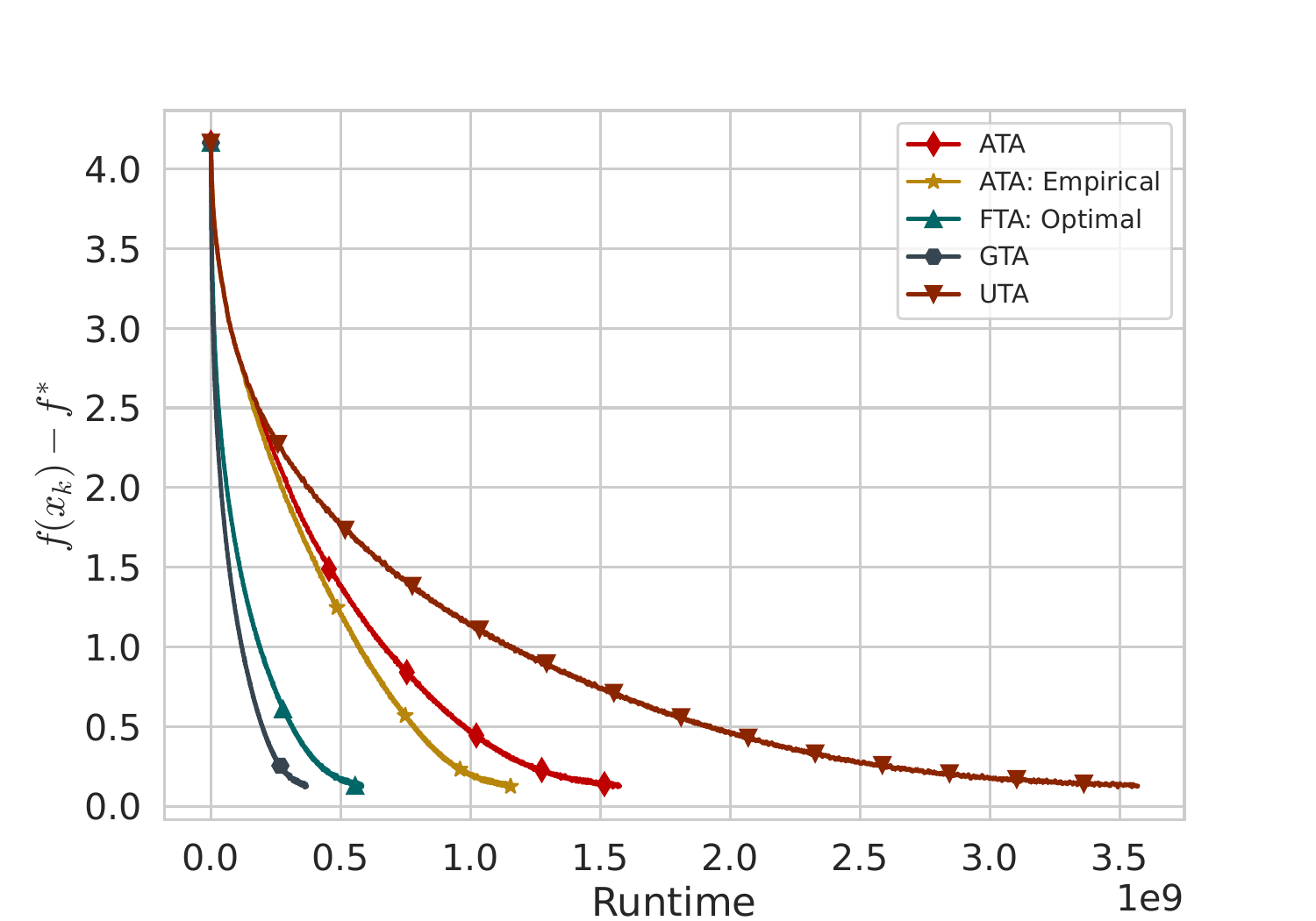} &
        \includegraphics[width=0.40\textwidth]{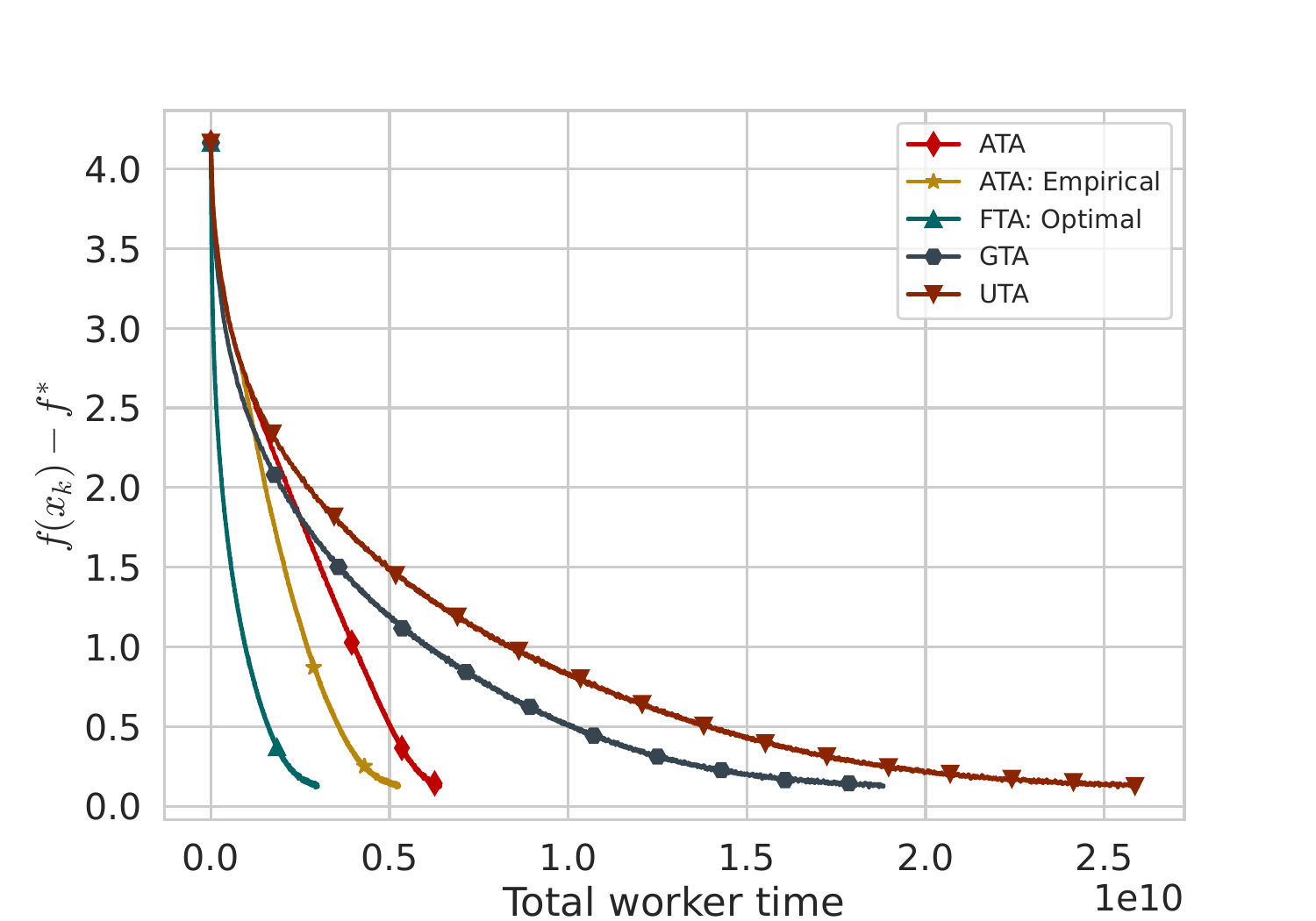} \\
        \includegraphics[width=0.40\textwidth]{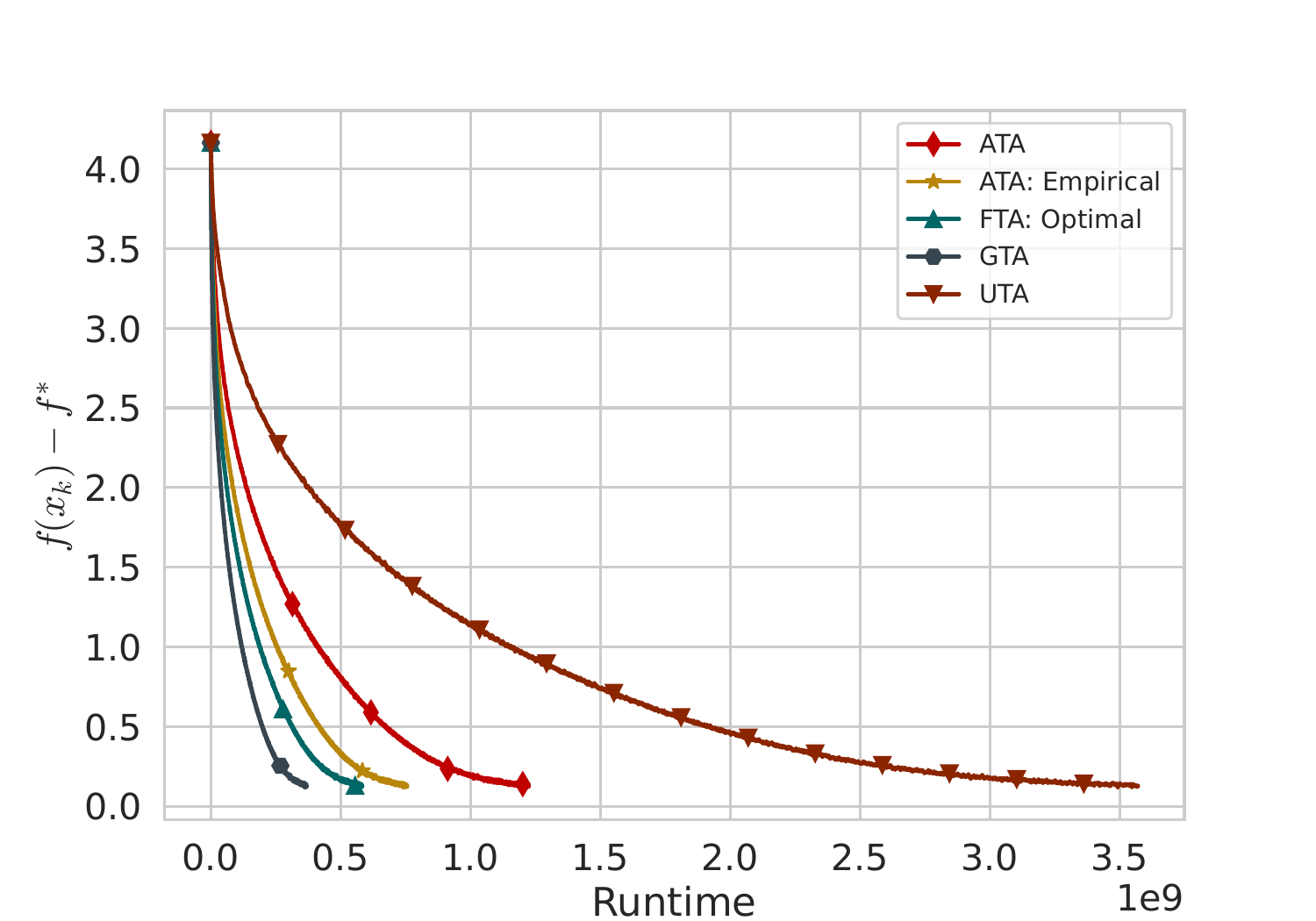} &
        \includegraphics[width=0.40\textwidth]{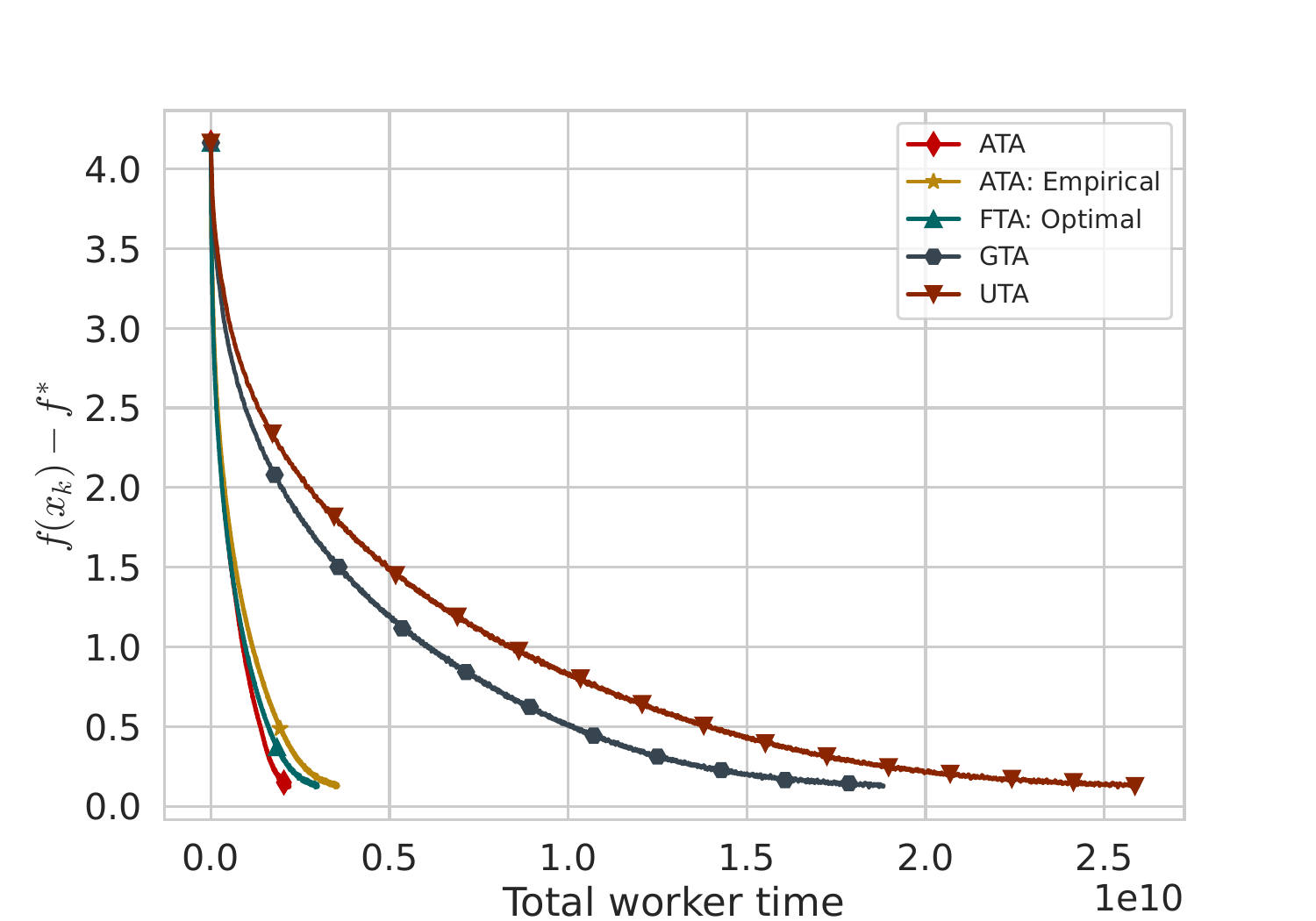} \\
        \includegraphics[width=0.40\textwidth]{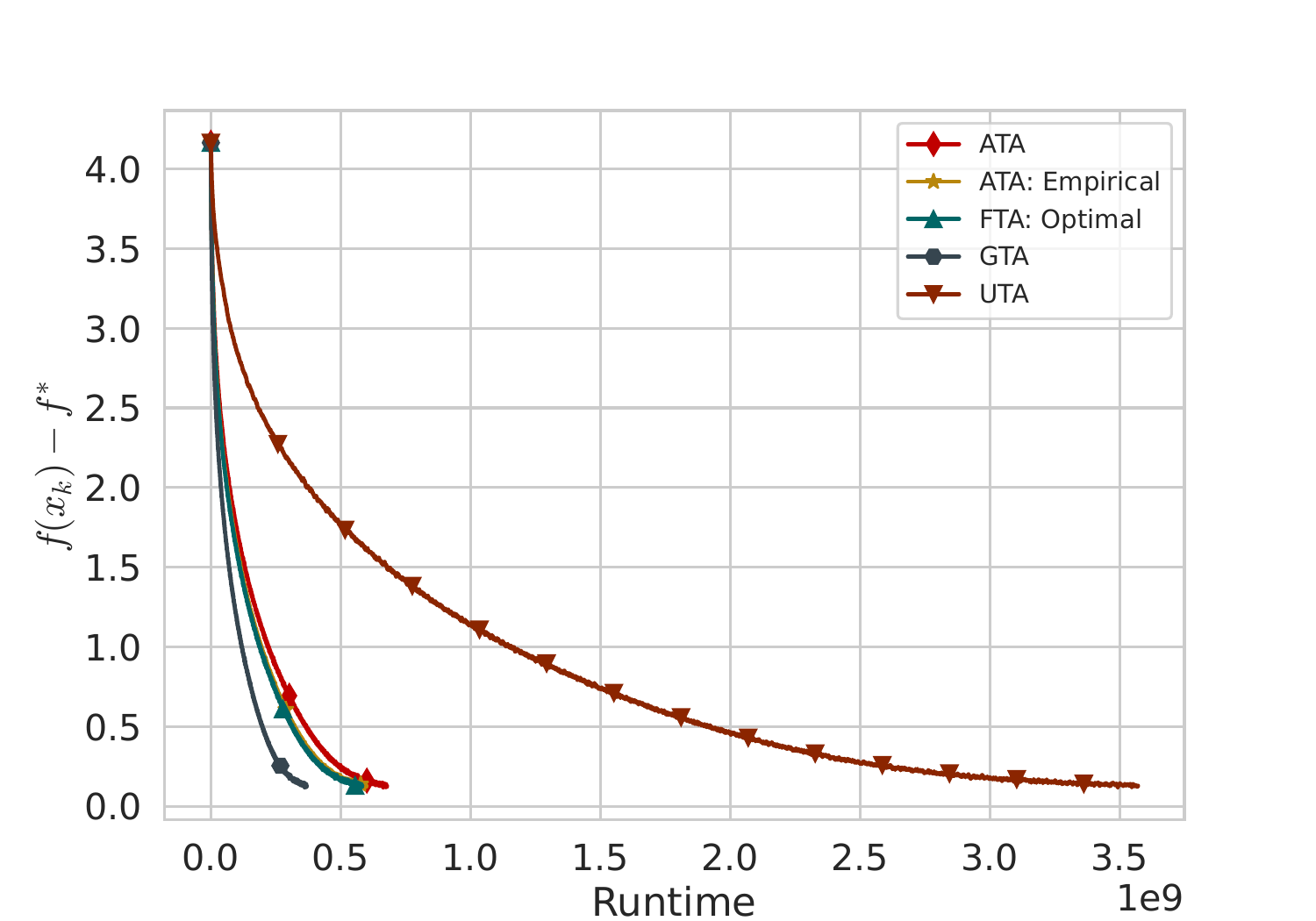} &
        \includegraphics[width=0.40\textwidth]{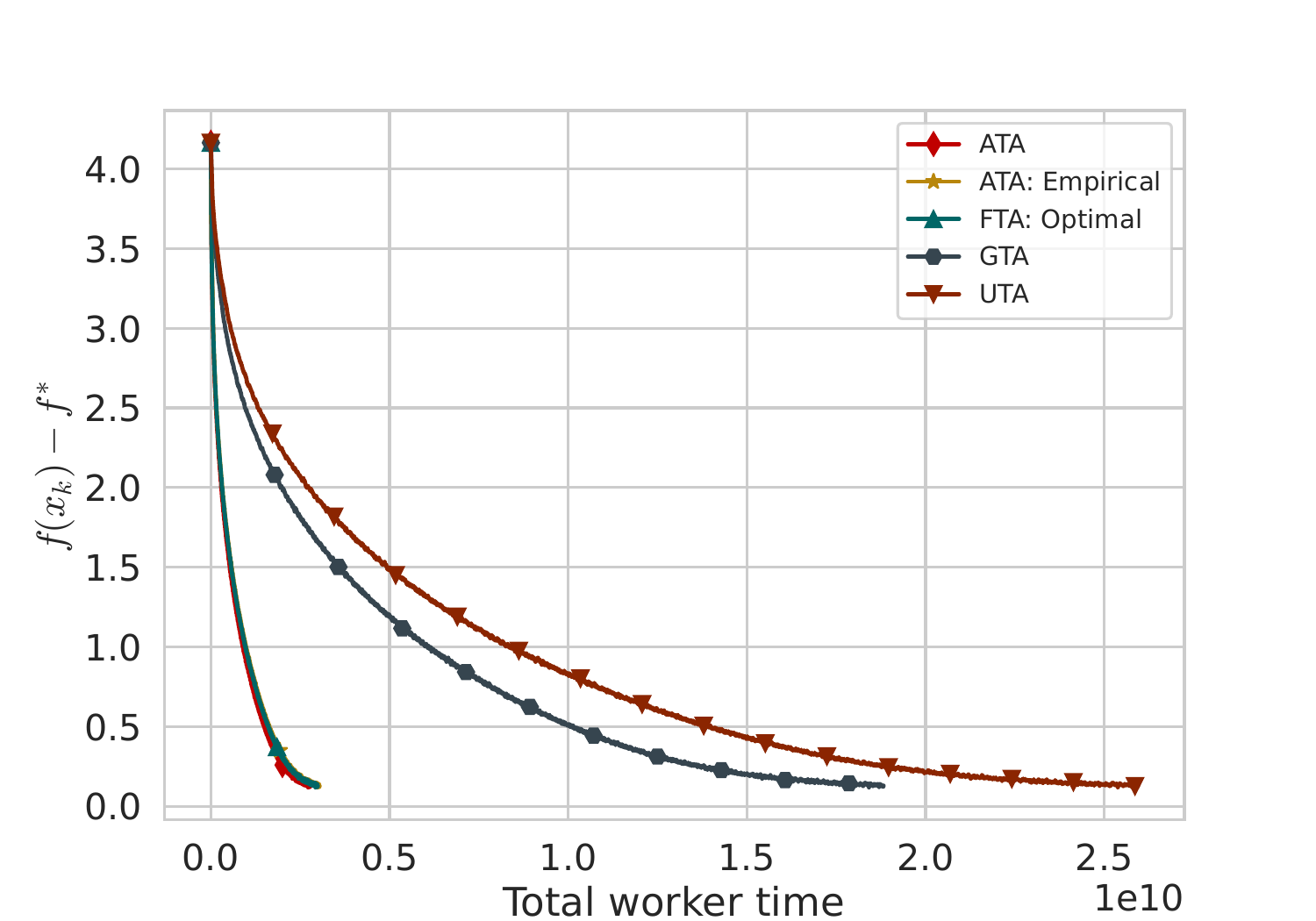}
    \end{tabular}
	}
    \caption{
        We use the same optimization setup as in \Cref{ata:fig:real}, with $B=23$ and $n=51$.
        The number of prior iterations, $P$, varies across rows, starting from the top with $P = 0, 5 \cdot 10^5, 5 \cdot 10^6$.
        As the number of prior iterations increases, we observe that the training of \algname{ATA} and \algname{ATA-Empirical} accelerates, bringing their performance closer to the optimal performance of \algname{OFTA}.
    }
    \label{ata:fig:prior}
\end{figure*}

\section{Concrete Optimization Methods}
\label{ata:section:other_methods}

In this section, we provide concrete examples of optimization algorithms using the \algname{ATA} and \algname{GTA} allocation strategies.

For optimization problems, we focus on \algname{SGD} and \algname{Asynchronous SGD}.
Other methods, such as stochastic proximal point methods and higher-order methods, can be developed in a similar fashion.

\subsection{Stochastic Gradient Descent}
Consider the problem of finding an approximate stationary point of the optimization problem
\begin{equation}
    \label{ata:eq:homo_problem}
    \minimize\limits_{x \in \R^d} \ \left\{f(x) \eqdef \ExpSub{\xi \sim {\cal D}}{f(x;\xi)}\right\}.
\end{equation}
We assume that each worker is able to compute stochastic gradient $f(x;\xi)$ satisfying $\mathbb{E}_{\xi \sim {\cal D}}\[ \|f(x;\xi) - \nabla f(x)\|^2\] \leq \sigma^2$ for all $x\in \R^d$.

In this case, \algname{SGD} with allocation budget $B$ becomes \algname{Minibatch SGD} with batch size $B$.
The next step is determining how the batch is collected.
For \algname{ATA}, we refer to this method as \algname{SGD-ATA}, as described in \Cref{ata:alg:sgd-ata}.

\begin{algorithm}[H]
	\caption{\algname{SGD-ATA}}
    \label{ata:alg:sgd-ata}
	\begin{algorithmic}[1]
		\STATE \textbf{Optimization inputs}: initial point $x^1 \in \R^d$, stepsize $\gamma > 0$
        \STATE \textbf{Allocation inputs}: allocation budget  $B$
        \STATE \textbf{Initialize}: empirical means $\hmu_i^1 = 0$, usage counts $K_i^1 = 0$, and usage times $T_i^1 = 0$, for all $i \in [n]$
		\FOR{$k = 1,\ldots, K$}
			\STATE Compute LCBs $(s_i^k)$ for all $i \in [n]$
			\STATE Find allocation:
				$$
					a^k \in  \argmin_{a\in \mathcal{A}} \ell(a, s^k)
				$$
			\STATE Allocate $a_i^k$ tasks to each worker $i \in [n]$
			\STATE Update $x$:
				$$
					x^{k+1} = x^k - \frac{\gamma}{B} \sum_{i=1}^n \sum_{j=1}^{a_i^k} \nabla f\(x^k;\xi_i^j\)
				$$
			\STATE For all $i$ such that $a_i^k \neq 0$, update:
			\begin{align*}
				K_i^{k+1} &= K_i^k + a_i^k \\
				T_i^{k+1} &= T_i^k + \sum_{j=1}^{a_i^k} X_i^{k,j} \\
				\hmu_i^{k+1} &= \frac{T_i^{k+1}}{K_i^{k+1}}
			\end{align*}
		\ENDFOR
	\end{algorithmic}
\end{algorithm}

In this case, each task consists in calculating the gradient using the device's local data, which is assumed to have the same distribution as the data on all other devices.
Because of this, it does not matter which device performs the task.
The method then averages these gradients to obtain an unbiased gradient estimator and performs a gradient descent step.

Now, let us give the version of \algname{Minibatch SGD} using greedy allocation \Cref{ata:alg:sgd-gta}.

\begin{algorithm}[H]
	\caption{\algname{SGD-GTA}}
    \label{ata:alg:sgd-gta}
	\begin{algorithmic}[1]
		\STATE \textbf{Input}: initial point $1 \in \R^d$, stepsize $\gamma > 0$, allocation budget $B$
		\FOR{$k = 1,\ldots, K$}
			\STATE $b=0$
			\STATE Query single gradient from each worker $i \in [n]$ 
			\WHILE{$b<B$}
				\STATE Gradient $\nabla f\big(x^k; \xi^{k,b}\big)$ arrives from worker $i^{k,b}$
				\STATE $g^k = g^k + \nabla f\big(x^k; \xi^{k,b}\big)$
				\STATE $b = b+1$
				\STATE Query gradient at $x^k$ from worker $i^{k,b}$
			\ENDWHILE
			\STATE Update the model: 
				$$
					x^{k+1} = x^k - \gamma \frac{g^k}{B}
				$$
		\ENDFOR
	\end{algorithmic}
\end{algorithm}

\Cref{ata:alg:sgd-gta} corresponds exactly to the \algname{Rennala SGD} method proposed by \citet{tyurin2024optimal}, which has optimal time complexity when the objective function is non-convex and smooth.

If the computation times are deterministic, then \algname{GTA} makes the same allocation in each iteration.
In that case, \algname{SGD-ATA} will converge to this fixed allocation.
If the times are random, the allocation found by \algname{GTA} may vary in each iteration.
In this case, \algname{SGD-ATA} will approach the best allocation for the expected times.

\subsection{Asynchronous SGD}

Here, we focus on the homogeneous problem given in \Cref{ata:eq:homo_problem}.
The greedy variant, \algname{Ringmaster ASGD}, was proposed by \citet{maranjyan2025ringmaster} and, like \algname{Rennala SGD}, achieves the best runtime.

We now present its version with \algname{ATA}, given in \Cref{ata:alg:asgd-ata}.

\begin{algorithm}[H]
	\caption{\algname{ASGD-ATA}}
    \label{ata:alg:asgd-ata}
	\begin{algorithmic}[1]
		\STATE \textbf{Optimization inputs}: initial point $x^1 \in \R^d$, stepsize $\gamma > 0$
        \STATE \textbf{Allocation inputs}: allocation budget  $B$
        \STATE \textbf{Initialize}: empirical means $\hmu_i^1 = 0$, usage counts $K_i^1 = 0$,\\
		and usage times $T_i^1 = 0$, for all $i \in [n]$
		\FOR{$k = 1,\ldots, K$}
        \STATE Compute LCBs $(s_i^k)$ for all $i \in [n]$
        \STATE Find allocation:
        $
        	a^k = \algname{RAS} (s^k; B)
        $
        \STATE Update $x^k$ using \Cref{ata:alg:asgd} with allocation $a^k$
        \STATE For all $i$ such that $a_i^k \neq 0$, update:
        \begin{align*}
            K_i^{k+1} &= K_i^k + a_i^k \\
            T_i^{k+1} &= T_i^k + \sum_{j=1}^{a_i^k} X_i^{k,j} \\
            \hmu_i^{k+1} &= \frac{T_i^{k+1}}{K_i^{k+1}}
        \end{align*}
		\ENDFOR
	\end{algorithmic}
\end{algorithm}

\begin{algorithm}[H]
	\caption{\algname{ASGD} updates}
	\label{ata:alg:asgd}
	\begin{algorithmic}[1]
		\STATE \textbf{Input:} Initial point $x^1 \in \R^d$, stepsize $\gamma > 0$,\\
		allocation vector $a$ with $\norm{a}_1 = B$
		\STATE Workers with $a_i > 0$ start computing stochastic gradients at $x^1$
		\FOR{$s = 1, 2, \ldots, B$}
			\STATE Receive gradient $\nabla f\big(x^{s+\delta^s}; \xi_{i^s}^{s+\delta^s}\big)$ from worker $i^s$
			\STATE Update: $x^{s+1} = x^{s} - \gamma \nabla f\big(x^{s+\delta^s}; \xi_{i^s}^{s+\delta^s}\big)$
			\IF{$a_{i^s} > 0$}
				\STATE Worker $i^s$ begins computing $\nabla f\big(x^{s+1}; \xi_{i^s}^{s+1}\big)$
				\STATE Decrease remaining allocation for worker $i^s$ by one: $a_{i^s} = a_{i^s} - 1$
			\ENDIF
		\ENDFOR
		\STATE \textbf{return:} $x^{B+1}$
	\end{algorithmic}
	\vspace{0.2cm}
	The sequence $\{\delta^s\}$ represents delays, where $\delta^s \geq 0$ is the difference between the iteration when worker $i^s$ started computing the gradient and iteration $s$, when it was applied.
\end{algorithm}

Here, the task remains gradient computation, but each worker's subsequent tasks use different points for computing the gradient.
These points depend on the actual computation times and the asynchronous nature of the method, hence the name \asyncsgd.

\section{Recursive Allocation Selection Algorithm}
\label{ata:sec:RAS}

In this section, we introduce an efficient method for finding the best allocation.
Given LCBs $s^k$ and allocation budget $B$, each iteration of \algname{ATA} (\Cref{ata:alg:ata}) determines the allocation by solving
$$
    a^k \in \argmin_{a\in \mathcal{A}} \ \ell(a, s^k)~,
$$
where 
$$
    \ell(a,\mu) \eqdef \max_{i \in [n]} \  a_{i} \mu_i  = \infnorm{a\odot \mu}~,
$$
with $\odot$ denoting the element-wise product.
When clear from context, we write $\ell(a)$ instead of $\ell(a, \mu)$.

In the early iterations, when some $s_i$ values are $0$, \algname{ATA} allocates uniformly across these arms until all $s_i$ values become positive.
After that, the allocation is determined using the recursive routine in \Cref{ata:alg:RAS}.

\begin{algorithm}
    \caption{Recursive Allocation Selection (\algname{RAS})}
    \label{ata:alg:RAS}
    \begin{algorithmic}[1]
        \STATE \textbf{Input:} Scores $s_1, \dots, s_n$, allocation budget $B$
        \STATE Assume without loss of generality that $s_1 \leq s_2 \leq \dots \leq s_n$ (i.e., sort the scores)
        \IF{$B = 1$}
            \STATE \textbf{return:} $(1, 0, \dots, 0)$
        \ENDIF
        \STATE Find the previous best allocation:
        $$
            a = (a_1, \dots, a_n) = \algname{RAS}\(s_1, \dots, s_n; B-1\)
        $$
        \STATE Determine the first zero allocation:
        \begin{equation}
            \label{ata:eq:r}
            r = 
            \begin{cases}
                \min\{i \mid a_i = 0\}, & \text{if }\ a_n = 0, \\
                n, & \text{otherwise}
            \end{cases}
        \end{equation}
        \STATE Find the best next query allocation set:
        $$
            M = \argmin_{i \in [r]} \  \infnorm{(a + e_i) \odot s},
        $$
        where $e_i$ is the unit vector in direction $i$.
        \STATE Select $j \in M$ such that the cardinality of 
        $$
            \argmax_{i \in [r]} \ (a_i + e_{j,i}) s_i
        $$ 
        is minimized
        \STATE \textbf{return:} $a + e_j$
    \end{algorithmic}
\end{algorithm}
\begin{remark}
    The iteration complexity of \algname{RAS} is 
	$$
		\cO \(n \ln\(\min\{B, n\}\) + \min\{B, n\}^2\) ~.
	$$
    In fact, the first $n \ln(\min\{B, n\})$ term arises from identifying the smallest $B$ scores.
    For the second term, note that in \eqref{ata:eq:r}, we have $r \leq \min\{B, n\}$.
\end{remark}
\subsection{Optimality}
We now prove that \algname{RAS} finds the optimal allocation, as stated in the following lemma.
\begin{boxedlemma}
    \label{ata:thm:RAS_optimality}
    For positive scores $0<s_1 \le s_2 \le \ldots \le s_n$, \algname{RAS} (\Cref{ata:alg:RAS}) finds an optimal allocation $h\in \cA$, satisfying
    $$
    h\in \argmin_{a\in \cA} \  \infnorm{a\odot s} ~.
    $$
\end{boxedlemma}
\begin{proof}
    We prove the claim by induction on the allocation budget $B$.
    
    \textbf{Base Case ($B = 1$):}  
    When $B = 1$, \algname{RAS} (\Cref{ata:alg:RAS}) allocates the task to worker with the smallest score (line 9).
    Thus, the base case holds.

    \textbf{Inductive Step:}  
    Assume that \algname{RAS} finds an optimal allocation for budget $B - 1$, denoted by
    $$
        \bar{h} = \algname{RAS}(s_1, \ldots, s_n; B-1)~.
    $$
    We need to prove that the solution returned for budget $B$, denoted by $h= \bar{h} + e_r$, is also optimal.

    Assume, for contradiction, that there exists $a\in \cA$ such that $a\neq h$ and $\ell(a) < \ell(h)$. 
    Write $a= \bar{a} + e_q$ for some $q \in [n]$. Observe that $\|\bar{a}\|_1=B-1$ because $a\in \mathcal{A}$.

    We consider two cases based on the value of $\ell\(\bar{h} + e_r\)$:

    \begin{itemize}
        \item $\ell\(\bar{h} + e_r\) = h_k s_k$ for some $k \neq r$.  
        In this case, adding one unit to index $r$ does not change the maximum value, i.e., $\ell\(\bar{h}\) = \ell\(\bar{h} + e_r\)$. 
        By the inductive hypothesis, $\bar{h}$ minimizes $\ell(\bx)$ for budget $B - 1$. 
        Therefore, we have
        $$
        	\ell(a) \geq \ell\(\bar{a}\) \geq \ell\(\bar{h}\) = \ell\(\bar{h} + e_r\) =\ell(h)~,
        $$
        which contradicts the assumption that $\ell(a) < \ell(h)$.

        \item $\ell\(\bar{h} + e_r\) = \(\bar{h}_r + 1\)s_r$.  
        By the algorithm's logic, $\(\bar{h}_r + 1\)s_r \leq \(\bar{h}_i + 1\)s_i$ for all $i \neq r$.
        Since $\ell(\bar{h}+e_r)\leq \ell(\bar{h}+e_q)$ and we assumed $\ell(\bar{a}+e_q)=\ell(a)<\ell(h)=\ell(\bar{h}+e_r)$, then $\bar{a} \neq \bar{h}$ otherwise $\ell(\bar{a}+e_q)<\ell(\bar{a}+e_r)$.
        Given that $\|\bar{h}\|_1=\|\bar{a}\|_1$, this implies that there exists some $u \in [n]$ such that $0\le\bar{a}_u \leq \bar{h}_u - 1$ and another index $v \in [n]$ where $\bar{a}_v \geq \bar{h}_v + 1$.

        In addition, note that $r$ is chosen such that $\ell\(\bar{h} + e_r\)$ is minimum. Using the fact that $\ell\(\bar{h} + e_r\) = \(\bar{h}_r + 1\)s_r$, we have that for any index $q$, we also necessarily have $\ell\(\bar{h} + e_q\) = \(\bar{h}_q + 1\)s_q$.
        Using this, we deduce
        $$
			\ell(h)
			=\ell\(\bar{h} + e_r\)
			\leq \ell\(\bar{h} + e_v\)
			= \(\bar{h}_v + 1\)s_v
			\leq \max_i \  \bar{a}_i s_i
			= \ell\(\bar{a}\) \leq \ell(a)~,
        $$
        where in the second inequality we used the fact that $\bar{a}_v \geq \bar{h}_v + 1$ and in the last inequality we used the fact that the loss is not decreasing for we add one element to the vector.
        This chain of inequalities again contradicts the assumption that $\ell(a) < \ell(h)$.
    \end{itemize}
    Since both cases lead to contradictions, we conclude that no $a\in \cA$ exists with $\ell(a) < \ell(h)$. 
    Thus, \algname{RAS} produces an optimal allocation for budget $B$.
\end{proof}

\subsection{Minimal Cardinality}

Among all possible allocations \algname{RAS} choose one that always minimizes the cardinality of the set:
$$
    \argmax_{i \in [n]} \  a_i s_i~.
$$
The reason for this choice is just technical as it allows the \Cref{ata:lem:1} to be true. 
\begin{boxedlemma}
    \label{ata:thm:minimal_cardinality}
    The output of $\algname{RAS}$ ensures the smallest cardinality of the set:
    $$
        \argmax_{i \in [n]} \ a_i s_i
    $$
    among all the optimal allocations $a$.
\end{boxedlemma}
\begin{proof}
    This proof uses similar reasoning as the one before.

    Let $h= \algname{RAS}(s;B)$, and denote the cardinality of the set $\argmax_{i \in [n]} \ a_i s_i $ for allocation $a$ by
    $$
        C_B(a) = \left| \argmax_{i \in [n]} \ a_i s_i  \right| \geq 1~.
    $$  
    We prove the claim by induction on $B$.

    \textbf{Base Case ($B=1$):}  
    For $B=1$, there is a single coordinate allocation, thus $C_1(h) = 1$, which is the smallest possible cardinality.

    \textbf{Inductive Step:}  
    Assume that $\algname{RAS}$ finds an optimal allocation for budget $B-1$ with the smallest cardinality, denote its output by
    $$
        \bar{h} = \algname{RAS}(s_1, \ldots, s_n; B-1)~.
    $$  
    We need to prove that $h= \bar{h} + e_r$ minimizes $C_B(a)$ among all optimal allocations for budget $B$.

    Assume, for contradiction, that there exists $a\in \cA$ such that $a\neq h$, $\ell(a) = \ell(h)$, and $C_B(a) < C_B(h)$. 
    Write $a= \bar{a} + e_q$ for some $q \in [n]$.
    We consider three cases:

    \begin{itemize}
        \item 
        $C_B(h) = 1$.  
        Since the minimum cardinality is exactly 1, we must have $C_B(a) \ge 1 = C_B(h)$, that contradicts our assumption.
        \item 
        $C_B(h) = C_{B-1}\(\bar{h}\)>1$.
        This occurs when $\ell(h) = \ell\(\bar{h}\) \ne \(\bar{h}_r + 1\)s_r$. 
        By the optimality of $h$, we have $\ell\(\bar{h}\) \le \ell\(\bar{a}\) \le \ell(a) = \ell(h)=\ell(\bar{h})$, which implies $\ell\(\bar{a}\) = \ell(a)$.
        Therefore, $C_{B-1}\(\bar{a}\) \le C_B(a)$. 
        Since the induction hypothesis holds for $B-1$, we have $C_{B-1}\(\bar{h}\) \le C_{B-1}\(\bar{a}\)$. 
        Thus,
        $$
        	C_B(h) = C_{B-1}\(\bar{h}\) \le C_{B-1}\(\bar{a}\) \le C_B(a)~,
        $$
        which leads to a contradiction.
        \item 
        $C_B(h) = C_{B-1}\(\bar{h}\) + 1$.  
        This occurs when $\ell(h) = \ell\(\bar{h}\) = \(\bar{h}_r + 1\)s_r$. 
        Proceeding as in the previous case, we have $\ell\(\bar{a}\) = \ell(a)$, and hence $C_{B-1}\(\bar{a}\) \le C_B(a)$.
        Since the induction hypothesis holds for $B-1$, we know $C_{B-1}\(\bar{h}\) \le C_{B-1}\(\bar{a}\)$.

        We now have additional cases:
        \begin{itemize}
			\item If $C_{B-1}\(\bar{a}\) = C_{B-1}\(\bar{h}\) + 1$, then
			$$
				C_B(h) = C_{B-1}\(\bar{h}\) + 1 = C_{B-1}\(\bar{a}\) \le C_B(a)~,
			$$
			which leads to a contradiction.
			\item Now assume $C_{B-1}\(\bar{a}\) = C_{B-1}\(\bar{h}\)$.
			We will show that in this case, $C_B(a) = C_{B-1}\(\bar{a}\) + 1$. 
			By contradiction, suppose $C_B(a) = C_{B-1}\(\bar{a}\)$, which implies $(\bar{a}_q + 1)s_q < \ell(a)$.
			Let $k$ be an index such that $\bar{a}_k s_k = \ell(a)$.
			Construct a new allocation $a^{\prime} = \bar{a} + e_q - e_k$. 
			Then,
			$$
				C_{B-1}\(a^{\prime}\) = C_{B-1}\(\bar{a}\) - 1 < C_{B-1}\(\bar{h}\)~,
			$$
			which contradicts the induction hypothesis. 
			Thus, 
			$$
				C_B(a) = C_{B-1}\(\bar{a}\) + 1 ~.
			$$
			Using this, we have
			$$
				C_B(h) = C_{B-1}\(\bar{h}\) + 1 = C_{B-1}\(\bar{a}\) + 1 = C_B(a)~,
			$$
			which again contradicts $C_B(a) < C_B(h)$.
        \end{itemize}
    \end{itemize}
    This concludes the proof.
\end{proof}
\section{Proofs of \Cref{ata:thm:main}, \Cref{ata:thm:main2}, and \Cref{ata:cor:main}}
\label{ata:sec:proof_1}
We start by recalling the notation.
For $i\in [n]$ and $k \in [K]$, $(X^{(u)}_{i,k})_{u \in [B]}$ denote $B$ independent samples at round $k$ from distribution $\nu_i$.
When using an allocation vector $a^k \in \mathcal{A}$, the total computation time of worker $i$ at round $k$ is $\sum_{u=1}^{a_i^k} X^{(u)}_{i,k}$, when $a_i^k >0$. $\mu = (\mu_1, \dots, \mu_K)$ is the vector of means.
For each $k \in [K]$, when using the allocation vector $a^k$, we recall the definition of the proxy loss $\ell: \mathcal{A}\times \R_+^n \to \R_+$ by
$$
	\ell(a^k, \lambda) = \max_{i\in [n]} \ a_i^k \lambda_i~,
$$
where $\lambda = (\lambda_1, \dots, \lambda_n)$ is a vector of non-negative components. When $\lambda = \mu$, we drop the dependence on the second input of $\ell$.
For each $\lambda$, let $\bar{a}_{\lambda}\in \mathcal{A}$, be the action minimizing this loss
$$
	\bar{a}_{\lambda} \in \argmin_{a \in \mathcal{A}} \ \ell(a, \lambda)~.
$$
We drop the dependency on $\mu$ from $\bar{a}_{\mu}$ to ease notation.
The actual (random) computation time at round $k$ is denoted by $C: \mathcal{A} \to \R_+$:
\begin{equation}\tag{\ref{ata:eq:def_c}}
	C(a^k) := \max_{i\in [n]} \ \sum_{u=1}^{a_i^k} X_i^{k,u}~.
\end{equation}
Let $a^*$ be the action minimizing the expected time
$$
    a^* \in \argmin_{a \in \mathcal{A}} \ \E{C(a)}~.
$$
The expected regret after $K$ rounds is defined as follows
$$
    \cR_K := \sum_{t=1}^{K} \E{\ell(a^k)-\ell(\bar{a})}~.
$$
For the remainder of this analysis we consider $\bar{a} \in \argmin_{a \in \mathcal{A}} \ \ell(a)$ found using the \algname{RAS} procedure.
For each $i\in [n]$, recall that $k_i$ is the smallest integer such that
\begin{equation}
    \label{ata:eq:def_n}
	(\bar{a}_i+k_i)\mu_i > \ell(\bar{a})~.
\end{equation}
Below we present a technical lemma used in the proofs of Theorems~\ref{ata:thm:main}~and~\ref{ata:thm:main2}.
\begin{boxedlemma}
	\label{ata:lem:1}
	Let $x=(x_1, \dots, x_n) \in \R_{\ge 0}^n$.
    Let $a$ be the output of $\algname{RAS}(x; B)$.
    For each $i, j \in [n]$, we have
	$$ 
	a_{ j} x_j \le \left(a_{ i}+1 \right) x_i~.
	$$
\end{boxedlemma}
\begin{proof}
	Fix $x \in \R_+^n$, and let $a = \algname{RAS}(x;B)$.
	The result is straightforward when $\min\limits_{i\in [n]}{x_i} = 0$.
	
	Suppose that $x_i >0$ for all $i \in [n]$.
	Let $s\ge 1$ denote the cardinality
	$$
	    s:= \abs{\argmax_{i \in [n]} \  a_{i} x_i }~.
	$$
	Fix $i,j \in [n]$, let $k \in \argmax_{i \in [n]} \  a_{i} x_i $. We need to show that
	$$
	    a_{k} x^k \le (a_{i}+1)x_i~.
	$$
	We use a proof by contradiction.
	Suppose that we have $ a_{k} x^k > (a_{i}+1)x_i$ consider the allocation vector $a^{\prime}\in \cA$ given by $a^{\prime}_k = a_{k}-1$, $a^{\prime}_i = a_i+1$ and $a^{\prime}_u = \bar{a}_u$ when $u \notin \{i,k\}$.
    Let $R := \max_{u \neq i,k} \{a_u x_u \}$.
    We have
	\begin{align*}
		\ell(a^{\prime},x)
		= \max_{u \in [n]} \ a^{\prime}_u x_u
		= \max \{(a_i+1) x_i, (a_k-1) x^k, R \}~.
	\end{align*}
	We consider two cases:
	\begin{itemize}
		\item Suppose that $s=1$ (i.e., the only element in $[n]$ such that $a_ux_u=\ell(a, x)$ is $k$), then we have necessarily $R< a_k x^k $. Moreover, by the contradiction hypothesis, $(a_i+1)x_i < a_k x^k$.
		Therefore,
		\begin{align*}
			\ell(a^{\prime}, x) = \max\{ (a_i+1)x_i, (a_k-1)x^k, R\}
			    < a_k x^k = \ell(a, x)~,
		\end{align*}
		which contradicts the definition of $a$.
		\item Suppose that $s\ge 2$, since by hypothesis $ a_k x^k > (a_i+1)x_i$, we clearly have $a_ix_i < \ell(a, x)$ therefore among the set $[n]\setminus \{k,i\}$ there are exactly $s-1$ elements such that $a_u x_u = \ell(a, x)$. In particular, this gives
		\begin{align*}
			\ell(a^{\prime},x)
			= \max_{u \in [n]} \ \{(a_i+1) x_i, (a_k-1) x^k, R \}
			= R = \ell(a, x)~.
		\end{align*}
		Therefore, $a^{\prime} \in \argmin_{a \in \mathcal{A}} \ \ell(a, x)$ and the number of elements such that $a^{\prime}_i x_i = \ell(a^{\prime}, x)=\ell(a, x)$ is at most $s-1$, which contradicts the fact that $s$ is minimal given the \algname{RAS} choice and \Cref{ata:thm:minimal_cardinality}.
	\end{itemize}
	As a conclusion we have $a_k x^k\le (a_i+1)x_i$.
\end{proof}
\begin{remark}
	Recall that \Cref{ata:lem:1} guarantees that $k_i$ defined in \eqref{ata:eq:def_n} satisfy: $k_i \in \{1, 2\}$ for each $i \in [n]$.
\end{remark}

\subsection{Proof of \Cref{ata:thm:main}}
\label{ata:proof:thm:main}
Below we restate the theorem.
\begin{restate-boxedtheorem}{\ref{ata:thm:main}}
	Suppose that \Cref{ata:a:sube} holds.
    Let $\bar{a} \in \argmin_{a \in \mathcal{A}} \ell(a)$, in case of multiple optimal actions, we consider the one output by \algname{RAS} when fed with $\mu$.
	Then, the expected regret of \algname{ATA} with inputs $(B, \alpha)$  satisfies
	$$
        \cR_K
        \le 2n\max_{i \in [n]} \{B\mu_i -\ell(\bar{a})\}+c \cdot\sum_{i=1}^{n} \frac{\alpha^2(\bar{a}_i+k_i)(B \mu_i - \ell(\bar{a})) }{\left((\bar{a}_i+k_i)\mu_i - \ell(\bar{a})\right)^2}\cdot \ln K~,
	$$
	where $\alpha = \max_{i \in [n]} \|X_i-\mu_i\|_{\psi_1}$, and $c$ is a constant.
\end{restate-boxedtheorem}
\begin{proof}
	Let $K_i^k$ be the number of rounds where arm $i$ was queried prior to round $k$ (we take $K_i^1=0$).
    Recall that we chose the following confidence bound: if $K_i^k \ge 1$, then
	$$
	    \text{conf}(i,k) = 2\alpha\sqrt{\frac{ \ln(2k^2)}{K_i^k}}+2\alpha\frac{ \ln(2k^2)}{K_i^k}~,
	$$
	and $\text{conf}(i,k) = \infty$ otherwise.
    Recall that $\hat{\mu}_i^k$ denotes the empirical mean of samples from $\nu_i$ observed prior to $k$ if $K_i^k\ge 0$ and $\hat{\mu}_i^k=0$ if $K_i^k=0$.
    Let $s_i^k$ denote the lower confidence bound used in the algorithm:
	$$
	    s_i^k = \left(\hat{\mu}_i^k -\text{conf}(i,k)\right)_{+}~.
	$$
	We introduce the events $\cE_i^k$ for $i \in [n]$ and $k \in [K]$ defined by
	$$
	    \cE_i^k := \left\lbrace \abs{\hat{\mu}_i^k-\mu_i} > \text{conf}(i,k)\right\rbrace.
	$$
	Let
	$$
	    \cE_k = \bigcup\limits_{i \in [n]} \cE_i^k~.
	$$
	Let us prove that for each $k \in [K]$ and $i \in [n]$: $\mathbb{P}\left(\cE_i^k\right) \le \nicefrac{1}{k^2}$, which gives using a union bound $\mathbb{P}(\cE_k) \le \nicefrac{n}{k^2}$.
	Let $i \in [n]$, using \Cref{ata:prop:concentration} and taking $\delta = \nicefrac{1}{k^2}$, we have
	\begin{align*}
		\mathbb{P}(\cE_i^k)
			= \mathbb{P}\left\{\abs{\hat{\mu}_i^k-\mu} > \text{conf}(i,k)\right\}
			\le \frac{1}{k^2} ~.
	\end{align*}
	We call a ``bad round", a round $k$ where we have $\ell(a^k) > \ell(\bar{a})$. Let us upper bound the number of bad rounds. 
	
	\noindent Observe that in a bad round there is necessarily an arm $i \in [K]$ such that $a_i^k \mu_i > \ell(\bar{a})$.
    For each $i\in [n]$, let $N_i(k)$ denote the number of rounds $q\in \{1,\dots, k\}$ where $a_i^q \mu_i > \ell(\bar{a})$ and $i \in \argmax_{j \in [n]} \ a_j^q \mu_j$ (this corresponds to a bad round triggered by worker $i$)
	$$
	    N_i(k) := \abs{\left\lbrace q \in \{1, \dots, k\}: a_i^q\mu_i > \ell(\bar{a}) \text{ and } a_i^q\mu_i = \ell(a_q) \right\rbrace}~.
	$$
	We show that in the case of $\ell(a^k) > \ell(\bar{a})$, the following event will hold: there exists $i \in [n]$ such that 
	$$
		E_{i,k} := \cE_k \text{ or }\left\lbrace N_i(k-1) \le  \frac{24\alpha^2 (\bar{a}_i+k_i) \ln(2K^2)}{\left((\bar{a}_i+k_i)\mu_i - \ell(\bar{a})\right)^2}  \right\rbrace~.
	$$
	To prove this we use a contradiction argument.
    Suppose that for each $i \in [n]$, $\neg E_{i,k}$ holds and that $\ell(a^k) > \ell(\bar{a})$.
    This means that $k$ is a bad round, let $i$ be an arm that triggered this bad round (i.e., $i \in \argmax_{j \in [n]} \ a_j^k \mu_j$).
    Event $\neg E_{i,k}$ gives in particular
	\begin{equation}
        \label{ata:eq:ni}
		N_i(k-1) > \frac{24\alpha^2 (\bar{a}_i+k_i) \ln(2K^2)}{\left((\bar{a}_i+k_i)\mu_i - \ell(\bar{a})\right)^2}~.
	\end{equation}
	Observe that in each round where $N_i(\cdot)$ is incremented, the number of samples received from the distribution $\nu_i$ increases by at least $\bar{a}_{i}+k_i$. 
	Therefore, we have \eqref{ata:eq:ni} implies
	\begin{align*}
		K_i^k
		> \frac{24\alpha^2(\bar{a}_i+k_i)^2 \ln(2K^2)}{\left((\bar{a}_i+k_i)\mu_i - \ell(\bar{a})\right)^2}
		= \frac{24\alpha^2 \ln(2K^2)}{\left(\mu_i - \frac{\ell(\bar{a})}{\bar{a}_i+k_i}\right)^2}~.
	\end{align*}
	
	\noindent Then we have, using the expressions of $\text{conf}(\cdot)$ and the bound above
	\begin{align}
		2\text{conf}(i,k) &=  4\alpha\sqrt{\frac{ \ln(2k^2)}{K_i^k}}+4\alpha\frac{ \ln(2k^2)}{K_i^k}\nonumber\\
		&\le \mu_i - \frac{\ell(\bar{a})}{\bar{a}_i+k_i}~.\label{ata:eq:conf}
	\end{align}
	The contradiction hypothesis gives that $a_i^k \mu_i > \ell(\bar{a})$, then we have, using the definition of $k_i$ that $a_i^k \ge \bar{a}_i + k_i$.
    Therefore, \eqref{ata:eq:conf} gives
	\begin{equation}
        \label{ata:eq:conf2}
		2\text{conf}(i,k) < \mu_i - \frac{\ell(\bar{a})}{a_i^k}~.
	\end{equation}
	Observe that in each round $\|a^k\|_1 = B$, therefore if we have $a_i^k \ge \bar{a}_i+k_i > \bar{a}_i$ for some $i$, we necessarily have that there exists $j \in [n]\setminus \{i\}$ such that $a_j^k \le \bar{a}_j-1$.
    Using the fact that $\ell(\bar{a}) \ge \bar{a}_j \mu_j$ with \eqref{ata:eq:conf2}, we get
	\begin{equation}
        \label{ata:eq:e1}
		a_i^k(\mu_i-2\text{conf}(i,k)) > \bar{a}_j \mu_j~.
	\end{equation}
	Since both $\neg \cE_i^k$ and $\neg \cE_{j,k}$ hold (because $\neg E_{i,k}$ implies $\neg \cE_k$), we have that
	\begin{align}\label{ata:eq:mu2}
		\mu_i - 2\text{conf}(i,k) &\le \hat{\mu}_{i, k} - \text{conf}(i,k) \le s_i^k~,
	\end{align} 
	and $\mu_j \ge \hat{\mu}_{j,k} - \text{conf}(j,k)$.
	Recall that $\mu_j \ge 0$, therefore
	\begin{align}
		\mu_j 
		\ge \left(\hat{\mu}_{j,k} - \text{conf}(j,k)\right)_{+}
		= s_{j,k}~.\label{ata:eq:mu3}
	\end{align}
	Using the bounds \eqref{ata:eq:mu2} and \eqref{ata:eq:mu3} in \eqref{ata:eq:e1}, we have
	$$
	    a_i^k s_i^k > \bar{a}_j s_{j,k} \ge (a_j^k+1) s_{j,k}~,
	$$
	where we used the definition of $j$ in the second inequality.
	This contradicts the statement of \Cref{ata:lem:1}, which concludes the contradiction argument.
    Therefore, the event that $k$ is a bad round implies that $E_{i,k}$ holds for at least one $i\in [n]$.
	We say that a bad round was triggered by arm $i$, a round where $N_i(\cdot)$ was incremented. 
	Observe that if $k \in [K]$ is not a bad round then $\E{\ell(a^k)}-\ell(\bar{a})=0$, otherwise if $k$ is a bad round triggered by $i \in [n]$ then 
	$$
		\E{\ell(a^k)}-\ell(\bar{a}) \le B\mu_i-\ell(\bar{a})~.
	$$
	To ease notation we introduce for $i\in [n]$
	$$
	    H_i := \frac{24\alpha^2 (\bar{a}_i+k_i) \ln(2K^2)}{\left((\bar{a}_i+k_i)\mu_i - \ell(\bar{a})\right)^2}~.
	$$
	We denote by $\mathds{1}(\cdot)$ the indicator function, which takes the value $1$ if the condition inside the braces is true and $0$ otherwise.

	The expected regret satisfies
	\begin{align*}
		\cR_K &= \sum_{i=1}^{K} \mathbb{E}\left[\ell(a^k)-\ell(\bar{a})\right]\\
		&\le \sum_{i=1}^{n} (B\mu_i-\ell(\bar{a}))\mathbb{E}[N_i(K)]\\ 
		&= \sum_{i=1}^{n}\sum_{k=1}^{K} (B\mu_i-\ell(\bar{a}))\mathbb{E}\left[\mathds{1}(k \text{ is a bad round triggered by }i)\right]\\
		&\le \max_{i\in [n]}\{(B\mu_i-\ell(\bar{a}))\}\cdot\sum_{t=1}^{K} \mathbb{P}(\cE_k) \\
            &\quad + \sum_{i=1}^{n}(B\mu_i-\ell(\bar{a}))\sum_{k=1}^{K} \mathbb{E}\left[\mathds{1}(k \text{ is a bad round triggered by }i) \mid \neg \cE_k\right]\\
		&\le \max_{i\in [n]}\{(B\mu_i-\ell(\bar{a}))\}\cdot\sum_{t=1}^{K} \mathbb{P}(\cE_k) \\
            &\quad+ \sum_{i=1}^{n}(B\mu_i-\ell(\bar{a}))\sum_{k=1}^{K} \mathbb{E}\left[\mathds{1}(N_i(k)=1+N_i(k-1) \text{ and } N_i \le H_i ) \mid \neg \cE_k\right]\\
		&\le \max_{i\in [n]}\{(B\mu_i-\ell(\bar{a}))\}\cdot\sum_{k=1}^{K} \mathbb{P}(\cE_k)+ \sum_{i=1}^{n} (B\mu_i-\ell(\bar{a}))H_i\\
		&\le 2n \max_{i\in [n]}\{(B\mu_i-\ell(\bar{a}))\}+  \sum_{i=1}^{n} \frac{24\alpha^2 (\bar{a}_i+k_i)(B\mu_i-\ell(\bar{a})) \ln(2K^2)}{\left((\bar{a}_i+k_i)\mu_i - \ell(\bar{a})\right)^2}~.
	\end{align*}
\end{proof}

\subsection{Proof of \Cref{ata:thm:main2}}
\label{ata:proof:thm:main2}
We restate below the theorem that provides the regret bound for \algname{ATA-Empirical} and then present its proof.
\begin{restate-boxedtheorem}{\ref{ata:thm:main2}}
	Suppose that \Cref{ata:a:sube} holds.
    Let $\bar{a} \in \argmin_{a \in \mathcal{A}} \ell(a)$, in case of multiple optimal actions, we consider the one output by \algname{RAS} when fed with $\mu$.
    Then, the expected regret of \algname{ATA-Empirical} with the empirical confidence bounds using the inputs $(B, \eta)$  satisfies
    \begin{align*}
        \cR_K &\le 2n\max_{i \in [n]} \{B\mu_i -\ell(\bar{a})\}
        +c \cdot\sum_{i=1}^{n} \frac{\eta^2 \mu_i^2(\bar{a}_i+k_i)(B \mu_i - \ell(\bar{a})) }{\left((\bar{a}_i+k_i)\mu_i - \ell(\bar{a})\right)^2}\cdot \ln K ~,
    \end{align*}
    where $\eta = \max_{i \in [n]} \nicefrac{\alpha_i}{\mu_i}$, and $c$ is a constant.
\end{restate-boxedtheorem}
\begin{proof}
	We build on the techniques used in the proof of \Cref{ata:thm:main}. Recall the expression of $\eta$:
	$$
	    \eta = \max_{i \in [n]}~ \frac{\alpha_i}{\mu_i}~.
	$$
	Define the quantities $C_{i,k}$ by
	$$
	    C_{i,k} = 2 \sqrt{\frac{\ln(2k^2)}{K_i^k}}+2 \frac{\ln(2k^2)}{K_i^k}~.
	$$
	Recall that the lower confidence bounds used here are defined as
	$$
	    \hat{s}_i^k = \hat{\mu}_i^k \left(1-\eta C_{i,k}\right)_{+}~.
	$$
	We additionally define the following quantities
	\begin{equation*}
		\hat{u}_{i,k} := \hat{\mu}_i^k \left(1+\frac{4}{3}\eta C_{i,k}\right)~.
	\end{equation*}
	\noindent We introduce the events $\cE_i^k$ for $i \in [n]$ and $k \in [K]$ defined by
	$$
		\cE_i^k := \left\lbrace \mu_i < \hat{s}_i^k \right\rbrace \quad \text{ or } \quad \left\lbrace \eta C_{i,k} \le \frac{1}{4} \quad \text{and} \quad \mu_i > \hat{u}_{i,k} \right\rbrace ~.
	$$
	Let 
	$$
		\cE_k = \bigcup_{i \in [n]} \cE_i^k ~.
	$$
	We have, using Lemma~\ref{ata:lem:conc2}, for each $k \in [K]$ and $i \in [n]$: $\mathbb{P}\left(\cE_i^k\right) \le \nicefrac{1}{k^2}$, which gives using a union bound $\mathbb{P}(\cE_k) \le \nicefrac{n}{k^2}$. 
	
	Following similar steps as in the proof of \Cref{ata:thm:main}, we call a ``bad round", a round $k$ where we have $\ell(a^k) > \ell(\bar{a})$.
    Let us upper bound the number of bad rounds. 
	
	Observe that in a bad round there is necessarily an arm $i \in [K]$ such that $a_i^k \mu_i > \ell(\bar{a})$.
    For each $i\in [n]$, let $N_i(k)$ denote the number of rounds $q\in \{1,\dots, k\}$ where $a_i^q \mu_i > \ell(\bar{a})$ and $i \in \argmax_{j \in [n]} \{ a_j^q \mu_j\}$ (this corresponds to a bad round triggered by worker $q$):
	$$
	    N_i(k) := \abs{\left\lbrace q \in \{1, \dots, k\}: a_i^q\mu_i > \ell(\bar{a}) \text{ and } a_i^q\mu_i = \ell(a_q) \right\rbrace}~.
	$$
	We show that in the case of $\ell(a^k) > \ell(\bar{a})$, the following event will hold: there exists $i \in [n]$ such that 
	$$
	    E_{i,k} := \cE_k \text{ or }\left\lbrace N_i(k-1) \le  \frac{185\eta^2\mu_i^2 (\bar{a}_i+k_i) \ln(2K^2)}{\left((\bar{a}_i+k_i)\mu_i - \ell(\bar{a})\right)^2}  \right\rbrace~.
	$$
	To prove this, suppose for a contradiction argument that we have for some $i \in [n]$: $\neg E_{i,k}$ and that $k$ is a bad round triggered by arm $i$ (i.e., $\ell(a^k) > \ell(\bar{a})$ and $i \in \argmax_{j \in [n]} {a_j^k\mu_j}$).
	
	This gives in particular
	\begin{equation}\label{ata:eq:ni2}
		N_i(k-1) >  \frac{185\eta^2 \mu_i^2 (\bar{a}_i+k_i) \ln(2K^2)}{\left((\bar{a}_i+k_i)\mu_i - \ell(\bar{a})\right)^2}~.
	\end{equation}
	Observe that in each round where $N_i(\cdot)$ is incremented, the number of samples received from the distribution $\nu_i$ increases by at least $\bar{a}_{i}+k_i$. 
	Therefore, we have \eqref{ata:eq:ni2} implies
	\begin{align*}
		K_i^k
		>  \frac{185\eta^2 \mu_i^2 (\bar{a}_i+k_i)^2 \ln(2K^2)}{\left((\bar{a}_i+k_i)\mu_i - \ell(\bar{a})\right)^2}
		=  \frac{185\eta^2 \mu_i^2  \ln(2K^2)}{\left(\mu_i - \frac{\ell(\bar{a})}{\bar{a}_i+k_i}\right)^2}~.
	\end{align*}
	Therefore, we have, using the expression of $C_{i,k}$ and the bound above
	\begin{align}
		C_{i,k} &= 2 \sqrt{\frac{\ln(2k^2)}{K_i^k}}+2 \frac{\ln(2k^2)}{K_i^k}\nonumber\\
		&\le \frac{3}{19\eta \mu_i} \left(\mu_i - \frac{\ell(\bar{a})}{\bar{a}_i+k_i}\right)~.\label{ata:eq:Ci}
	\end{align}
	The last bound implies in particular that $\eta C_{i,k} \le \nicefrac{3}{19}$, hence $(1-\eta C_{i,k})_{+} = 1-\eta C_{i,k}$.
	
	We have
	\begin{align}
		\hat{\mu}_i^k -\hat{s}_i^k &= \hat{\mu}_i^k \left(1-\left(1-\eta C_{i,k}\right)_{+}\right)\nonumber\\
		&\le \eta C_{i,k}\hat{\mu}_i^k\nonumber\\
		&\le \eta C_{i,k}\frac{\mu_i}{1-\eta C_{i,k}}\nonumber\\
		&\le \frac{1}{5} \left(\mu_i - \frac{\ell(\bar{a})}{\bar{a}_i+k_i}\right)~,\label{ata:eq:b10}
	\end{align}
	where we used the event $\neg \cE_i^k$ in the penultimate inequality (in particular $\hat{s}_i^k = \hat{\mu}_i^k(1-\eta C_{i,k})_{+}\le \mu_i$), and the bound \eqref{ata:eq:Ci} in the last inequality.

    Given the hypothesis that $\ell(a^k)>\ell(\bar{a})$ and, $\ell(a^k) = a_i^k \mu_i$, we necessarily have $a_i^k \ge \bar{a}_i+k_i$.
    Therefore, bound \eqref{ata:eq:b10}
	$$
	    5\hat{\mu}_i^k - 5\hat{s}_i^k < \mu_i - \frac{\ell(\bar{a})}{a_i^k}~.
	$$
	Observe that in each round $\norm{a^k}_1 = B$, therefore if we have $a_i^k \ge \bar{a}_i+k_i > \bar{a}_i$ for some $i$, we necessarily have that there exists $j \in [n]\setminus \{i\}$ such that $a_j^k \le \bar{a}_j-1$.
	Therefore, using the fact that $\ell(\bar{a}) \ge \bar{a}_j \mu_j$, we obtain
	\begin{equation}\label{ata:eq:e12}
		a_i^k(\mu_i+ 5\hat{s}_i^k-5\hat{\mu}_i^k) > \ell(\bar{a}) \ge \bar{a}_j \mu_j~.
	\end{equation}
	Since both $\neg \cE_i^k$ and $\neg \cE_{j,k}$ hold (because $\neg E_{i,k}$ implies $\neg \cE_k$), we have that
	\begin{align*}
		\mu_i +5\hat{s}_i^k-5\hat{\mu}_i^k &= \hat{s}_i^k + \mu_i - \hat{\mu}_i^k + 4 \left(\hat{s}_i^k - \hat{\mu}_i^k \right)\\
		&= \hat{s}_i^k + \mu_i - \hat{\mu}_i^k + 4\hat{\mu}_i^k \left((1-\eta C_{i,k})_{+}-1 \right)\\
		&\le \hat{s}_i^k + \mu_i - \hat{\mu}_i^k - 4\hat{\mu}_i^k \eta C_{i,k}\\
		&\le \hat{s}_i^k + \mu_i - \hat{u}_{i,k}~,
	\end{align*} 
	To conclude, observe that given $\eta C_{i,k} \le \nicefrac{3}{19}$, event $\neg \cE_i^k$ implies that $\mu_i \le \hat{u}_{i,k}$, therefore 
	$$
		\mu_i +5\hat{s}_i^k-5\hat{\mu}_i^k \le \hat{s}_i^k~.
	$$
	Since $\neg \cE_{j,k}$ holds, we also have
	\begin{equation*}
		\mu_j \ge \hat{s}_{j,k}~.
	\end{equation*}
	Using the two last bounds in \eqref{ata:eq:e12}, we have
	$$
	    a_i^k \hat{s}_i^k > \bar{a}_j \hat{s}_{j,k} \ge (a_j^k+1) \hat{s}_{j,k}~,
	$$
	where we used the definition of $j$, as an arm satisfying $\bar{a}_j \ge 1+a_j^k$, in the second inequality.
	This contradicts the statement of \Cref{ata:lem:1}, which concludes the contradiction argument.
    Therefore, the event that $k$ is a bad round implies that $E_{i,k}$ holds for at least one $i\in [n]$.
	We say that a bad round was triggered by arm $i$, a round where $N_i(\cdot)$ was incremented. 
	Observe that if $k \in [K]$ is not a bad round then $\E{\ell(a^k)}-\ell(\bar{a})=0$, otherwise if $k$ is a bad round triggered by $i \in [n]$ then $\E{\ell(a^k)}-\ell(\bar{a}) \le B\mu_i-\ell(\bar{a})$.
	To ease notation we introduce for $i\in [n]$
	$$
	    H_i := \frac{185 \eta^2\mu_i^2 (\bar{a}_i+k_i) \ln(2K^2)}{\left((\bar{a}_i+k_i)\mu_i - \ell(\bar{a})\right)^2}~.
	$$
	The expected regret satisfies $\mathds{s}$
	\begin{align*}
		\cR_K &= \sum_{i=1}^{K} \mathbb{E}\left[\ell(a^k)-\ell(\bar{a})\right]\\
		&\le \sum_{i=1}^{n} (B\mu_i-\ell(\bar{a}))\mathbb{E}[N_i(K)]\\ 
		&= \sum_{i=1}^{n}\sum_{k=1}^{K} (B\mu_i-\ell(\bar{a}))\mathbb{E}\left[\mathds{1}(k \text{ is a bad round triggered by }i)\right]\\
		&\le \max_{i\in [n]}\{(B\mu_i-\ell(\bar{a}))\}\cdot\sum_{t=1}^{K} \mathbb{P}(\cE_k) \\ 
            &\quad+ \sum_{i=1}^{n}(B\mu_i-\ell(\bar{a}))\sum_{k=1}^{K} \mathbb{E}\left[\mathds{1}(k \text{ is a bad round triggered by }i) \mid \neg \cE_k\right]\\
		&\le \max_{i\in [n]}\{(B\mu_i-\ell(\bar{a}))\}\cdot\sum_{t=1}^{K} \mathbb{P}(\cE_k) \\
            &\quad+ \sum_{i=1}^{n}(B\mu_i-\ell(\bar{a}))\sum_{k=1}^{K} \mathbb{E}\left[\mathds{1}(N_i(k)=1+N_i(k-1) \text{ and } N_i \le H_i ) \mid \neg \cE_k\right]\\
		&\le \max_{i\in [n]}\{(B\mu_i-\ell(\bar{a}))\}\cdot\sum_{k=1}^{K} \mathbb{P}(\cE_k)+ \sum_{i=1}^{n} (B\mu_i-\ell(\bar{a}))H_i\\
		&\le 2n \max_{i\in [n]}\{(B\mu_i-\ell(\bar{a}))\}+  \sum_{i=1}^{n} \frac{185 \eta^2\mu_i^2 (\bar{a}_i+k_i)(B \mu_i - \ell(\bar{a})) \ln(2K^2)}{\left((\bar{a}_i+k_i)\mu_i - \ell(\bar{a})\right)^2}~. 
	\end{align*}
\end{proof}

\subsection{Proof of \Cref{ata:cor:main}}
\label{ata:sec:proof_2}

Let us first restate the theorem.
\begin{restate-boxedtheorem}{\ref{ata:cor:main}}
	Suppose Assumption~\ref{ata:a:sube} holds and let $\eta := \max_{i \in [n]} \nicefrac{\alpha_i}{\mu_i}$, where $\alpha_i = \|X_i-\mu_i\|_{\psi_1}$.
	Then, the total expected computation time after $K$ rounds, using the allocation prescribed by \algname{ATA} with inputs $(B, \alpha)$ satisfies
	$$
	\cC_K \le \left(1+4\eta \ln(B)\right)\cC_K^* + \mathcal{O}(\ln K)~.
	$$
\end{restate-boxedtheorem}
\begin{proof}
	Let $\mathbb{E}_k$ be the expectation with respect to the variables observed up to and including $k$ and $\mathcal{F}_k$ the corresponding filtration. Using the tower rule, we have
	$$
		\sum_{k=1}^{K}\E{C(a^k)} = \E{ \sum_{k=1}^{K} \ExpSub{k-1}{C(a^k)} }.
	$$
	Consider round $k \in [K]$, let us upper bound $\mathbb{E}_{k-1}[C(a_t)]$ using $\mathbb{E}_{k-1}[\ell(a^k)]$. Recall that $a^k \in \mathcal{F}_{k-1}$, let $Y_i = \sum_{u=1}^{a_i^k} X^{(u)}_{i,k}$, since $Y_i$ is the sum of $a_i^k$ i.i.d samples we have that $\mathbb{E}_{k-1}[Y_i] = a_i^k \mu_i$ and $\|Y_i - a_i^k \mu_i\|_{\psi_1} \le a_i^k \|X_i-\mu_i\|_{\psi_1}$.
	Thus, using \Cref{ata:lem:exp_max}, we get
	\begin{align*}
		\mathbb{E}_{k-1}\left[C(a^k) \right] &= \mathbb{E}_{k-1}\left[ \max_{i \in \text{supp}(a^k)} \sum_{u=1}^{a_i^k} X^{(u)}_{i,k} \right]\\
		&\le \max_{i \in \text{supp}(a^k)} a_i^k \mu_i  + 4\max_{i \in \text{supp}(a^k)} a_i^k \alpha_i \cdot \ln(B)\\
		&\le \max_{i \in \text{supp}(a^k)} a_i^k \mu_i  + 4\max_{i \in \text{supp}(a^k)} a_i^k\, \eta\mu_i \cdot \ln(B)\\
		&= \left(1+4\eta \ln(B)\right)\max_{i \in \text{supp}(a^k)} a_i^k \mu_i .
	\end{align*}
	Moreover, using Jensen's inequality, we have
	\begin{align*}
		\max_{i \in [n]} a^*_i \mu_i
		\le \E{ \max_{i \in [n]} \sum_{u=1}^{a_{k,i}} X^{(u)}_{i,k} }
		= \E{C(a^*)} ~.
	\end{align*}
	Using the last two bounds with the result of \Cref{ata:thm:main}, we get the result.
\end{proof}
\section{Technical Lemmas}
\label{ata:sec:technical}
The lemma below gives a concentration bound on sub-exponential variables.
Note that this result can be inferred from Proposition 7 in the work of \citet{maurer2021concentration}, although applying the last result directly requires assuming the variables are positive, this is not needed in their proof in the one dimensional case.
For completeness, we present the full proof below.
\begin{boxedlemma}\label{ata:prop:concentration}
	Let $Y_1, \dots, Y_n$ be iid random variables with $\mathbb{E}[Y_1] = 0$ and $\alpha = \|Y_1\|_{\psi_1} < +\infty$.
    Then for any $\delta \in (0,1)$, we have with probability at least $1-\delta$
	$$
		\abs{\frac{1}{n}\sum_{i=1}^n Y_i} \le 2\alpha \left(\sqrt{\frac{\ln(\nicefrac{2}{\delta})}{n}}+ \frac{\ln(\nicefrac{2}{\delta})}{2n}\right)~.
	$$
\end{boxedlemma}
\begin{proof}
	Let $v := 2n \alpha^2$.
    We have using Lemma~\ref{ata:lem:mgf} that
	$$
		\sum_{i=1}^n \E{Y_i^2} \le 2n \alpha^2 \quad \text{and} \quad \sum_{i=1}^{n} \E{(Y_i)_{+}^q} \le \frac{q!}{2} v \alpha^{q-2}~.
	$$
	Therefore, using Bernstein concentration inequality (Proposition~\ref{ata:prop:bernstein}) we obtain that
	$$
		\mathbb{P}\left(\abs{\sum_{i=1}^n Y_i} \ge 2\alpha \sqrt{nt}+ \alpha t\right) \le 2\exp(-t)~.
	$$
	Choosing $t=\ln(\nicefrac{2}{\delta})$, we obtain the result.
\end{proof}

\begin{proposition}[Theorem 2.10 in the work of \citet{boucheron2005moment}]
    \label{ata:prop:bernstein}
	Let $X_{1},\dots,X_{n}$ be independent real-valued random variables.
	Assume there exist positive numbers $v$ and $c$ such that
	$$
		\sum_{i=1}^{n}\mathbb{E}\bigl[X_{i}^{2}\bigr] \le v \quad\text{and}\quad \sum_{i=1}^{n}\mathbb{E}\bigl[(X_{i})_{+}^{\,q}\bigr] \le \frac{q!}{2}\,v\,c^{\,q-2}~, \qquad\text{for all integers } q\ge 3~,
	$$
	where $x_{+}:=\max\{x,0\}$.
	Define the centered sum
	$$
		S :=\sum_{i=1}^{n}\bigl(X_{i}-\mathbb{E}X_{i}\bigr).
	$$
	Then, for every $t>0$,
	$$
		\mathbb{P}\left( S \ge \sqrt{2vt}+ct\right) \le e^{-t}~.
	$$
\end{proposition}
\begin{boxedlemma}\label{ata:lem:conc2}
	Let $X_1, \dots, X_n$ be i.i.d positive random variables with mean $\mu$ and $\|X_1\|_{\psi_1}<+\infty$. Denote $\alpha := \|X_1-\mu\|_{\psi_1}$.
	Denote $\hat{X}_n = \frac{1}{n}\sum_{i=1}^n X_i$ and let $\delta \in (0,1)$. Define $\eta$, $C_{\cdot, \cdot}$ by:
	\begin{equation*}
		\eta :=  \frac{\alpha}{\mu} \qquad \text{and} \qquad 
		C_{n,\delta} := 2  \sqrt{\frac{\ln(\nicefrac{2}{\delta})}{n}}+2 \frac{\ln(\nicefrac{2}{\delta})}{n}~.
	\end{equation*}
	Then with probability at least $1-\delta$ we have
	$$
		\mu \ge \hat{X}_n (1-\eta \cdot C_{n, \delta})_{+}~,
	$$
	where we use the notation $(a)_+ = \max\{0,a\}$. Moreover, if $\eta~C_{n,\delta} \le \nicefrac{1}{4}$, then with probability least $1-\delta$
	$$
		\hat{X}_n (1-\eta \cdot C_{n,\delta})_{+} \le \mu \le \hat{X}_n \left(1+\frac{4}{3}\eta\cdot  C_{n,\delta}\right)~.
	$$
\end{boxedlemma}
\begin{proof}
	Fix $n, \delta$.
	We work on the event 
	$$
		\cE_{n,\delta} = \left\lbrace \abs{\hat{X}_n - \mu} \le \alpha \cdot  C_{n, \delta} \right\rbrace
	$$
	that holds with probability at least $1-\delta$ if we apply \Cref{ata:prop:concentration} to $X_i-\mu$.
	
	\noindent \textbf{Proof of $\mu \ge \hat{X}_n (1-\eta \cdot C_{n,\delta})_{+}$:}
	
	If $\eta C_{n,\delta} \ge 1$, we have $(1-\eta \cdot C_{n,\delta})_{+} = 0$ and the result follows from the fact that $X$ is non-negative which gives $\mu \ge 0$. 
	
	\noindent Suppose now that $\eta C_{n,\delta} < 1$. Recall that event $\cE_{n,\delta}$ implies that 
	$$
	    \hat{X}_n \le \mu (1+\eta C_{n,\delta})~.
	$$
	Therefore, we have
	$$
		\frac{\hat{X}_n}{1+\eta C_{n,\delta}} \le \mu ~.
	$$
	Using $1-\eta C_{n,\delta} \le \frac{1}{1+\eta C_{n,\delta}}$ with the bound above, we obtain
	$$
		\hat{X}_n (1-\eta C_{n,\delta})_{+} = \hat{X}_n (1-\eta C_{n,\delta}) \le \frac{\hat{X}_n}{1+\eta \cdot C_{n,\delta}} \le \mu~.
	$$
	\noindent \textbf{Proof of $ \mu \le \left(1+\frac{4}{3} \eta \cdot C_{n,\delta}\right).$} 
	Recall that event $\cE_{n,\delta}$ gives
	$$
		\hat{X}_n \ge \mu - \alpha C_{n,\delta} = \mu (1-\eta C_{n,\delta})~. 
	$$
	Suppose that $\eta C_{n,\delta} \le \nicefrac{1}{4}$.
    We therefore have
	$$
		\mu \le \frac{\hat{X}_n}{1-\eta C_{n,\delta}}~.
	$$ 
	Next, we use the fact that for any $x \in [0, \nicefrac{1}{4}]$, we have
	$$
		\frac{1}{1-x} \le 1+\frac{4}{3}x~,
	$$
	which gives
	$$
		\mu \le \hat{X}_n \left(1+\frac{4}{3}\eta \cdot C_{n,\delta}\right)~,
	$$
	when $\eta C_{n,\delta} \le \nicefrac{1}{4}$.
\end{proof}
The following lemma provides a useful bound on the expectation of the maximum of independent sub-exponential random variables.
\begin{boxedlemma}\label{ata:lem:exp_max}
	Let $X_1, \dots, X_n$ be a sequence of independent random variables with finite Orlicz norm $\|X_i\|_{\psi_1} <+\infty$ and let $\E{X_i} = \mu_i$.
    Then we have
	$$
		\E{\max_{i \in [n]} X_i } \le \max_{i \in [n]} \mu_i + 4\alpha \ln(n)~,
	$$
where $\alpha = \max_{i \in [n]} \|X_i - \mu_i\|_{\psi_1}$. 
\end{boxedlemma}
\begin{proof}
	If $n=1$ the bound is straightforward, suppose that $n\ge 2$.
	Let $Y_i := X_i- \mu_i$, then $\alpha = \max_{i \in [n]} \|Y_i\|_{\psi_1}$.
    Let $\lambda \in (0, 1/\alpha)$, we have
	\begin{align*}
		\max_{i\in [n]} Y_i &= \frac{1}{\lambda} \ln\left(\exp\left(\lambda \max_{i\in [n]} Y_i\right)\right)\\
		&\le \frac{1}{\lambda} \ln\left(\sum_{i=1}^n \exp\left(\lambda Y_i\right)\right)~.
	\end{align*}
	Taking the expectation and using \Cref{ata:lem:mgf}, we have
	\begin{align*}
		\E{\max_{i\in [n]} Y_i} &\le \frac{1}{\lambda} \ln\left(\sum_{i=1}^n \E{\exp\left(\lambda Y_i\right)}\right)\\
		&\le \frac{1}{\lambda} \ln\left(\frac{n}{1-\lambda \alpha}\right)~.
	\end{align*}
	We choose 
	$$
		\lambda = \frac{1 - 1 /n}{\alpha} ~,
	$$
	which gives
	\begin{align}
		\E{\max_{i\in [n]} Y_i} &\le \frac{\alpha}{1-\frac{1}{n}} \ln(n^2)\nonumber\\
		&= 2\alpha \frac{n}{n-1} \ln(n)\nonumber\\
		&\le 4\alpha \ln(n)~.\label{ata:eq:y}
	\end{align}
	Let $i^* \in \argmax_{i \in [n]} X_i$, we have 
	\begin{align*}
		\max_{i \in [n]} X_i - \max_{i\in [n]} \mu_i &= X_{i^*} - \max_{i\in [n]} \mu_i\\
		&\le X_{i^*} - \mu_{i^*} \le \max_{i \in [n]} \(X_i - \mu_i\) = \max_{i\in [n]} Y_i~.
	\end{align*}
	Combining the last bound with \eqref{ata:eq:y} we obtain
	$$
		\E{\max_{i \in [n]} X_i } \le \max_{i \in [n]} \mu_i + 4\alpha \ln(n)~.
	$$
\end{proof}
Lemma below is based on a standard argument we give here for completeness.
\begin{boxedlemma}\label{ata:lem:mgf}
	Let $Y$ be a variable such that $\alpha = \|Y\|_{\psi_1} <+\infty$.
    Then we have for any $\lambda \in \left(-\nicefrac{1}{\alpha}, \nicefrac{1}{\alpha}\right)$
	$$
		\E{\exp\left(\lambda Y\right)} \le \frac{1}{1-\abs{\lambda}\alpha}~.
	$$
	Moreover, we have for any $q\ge 3$
	$$
		\E{Y^2} \le 2\alpha^2 \quad \text{and} \quad \E{(Y)_{+}^q} \le \frac{q!}{2} \cdot (2\alpha^2) \cdot \alpha^{q-2}~.
	$$
\end{boxedlemma}
\begin{proof}
	Let $Z = \nicefrac{\abs{Y}}{\alpha}$.
    First observe that we have
	$$
		\sum_{k \ge 0} \frac{\E{Z^k}}{k!}
		= \E{\exp(Z)}
		= \E{\exp\(\frac{\abs{Y}}{\alpha}\)} \le 2 ~,
	$$
	so,
	$$
		\sum_{k \ge 1} \frac{\E{Z^k}}{k!} \le 1 ~.
	$$
	This implies 
	$$
		\E{\abs{Y}^k} \le k! \alpha^k
	$$
	for all $k \ge 1$.
	Using this bound, we estimate the moment generating function
	\begin{align*}
		\E{\exp(\lambda Y)}
		&\le \E{\exp\left(\abs{\lambda}\abs{Y}\right)}\\
		&= \sum_{k \ge 0} \frac{\abs{\lambda}^k\E{\abs{Y}^k}}{k!}\\
		&\le 1 + \sum_{k \ge 1} \abs{\lambda}^k\alpha^k\\
		&= \frac{1}{1-\abs{\lambda}\alpha}~.
	\end{align*}
	The remaining bounds follow from 
	$$
		\E{\abs{Y}^k} \le k! \alpha^k ~.
	$$
\end{proof}

\vspace{8em}
\begin{center}
    {\Huge \textit{Khalas}}
\end{center}
		\endgroup


\end{document}